%% file: gener.tex
\documentstyle[a4,txmac,%
amssymb,%
case,%
myrot,%
mypic,%
epsf,%
twoside,%
nocaphead,%
amssymb,times,mathptm]{article}

\makeatletter

\advance\oddsidemargin by -1.7cm
\advance\evensidemargin by -1.7cm
\advance\textwidth by 3.4cm

\def\mynewtheo#1#2{%
\newtheorem{@#1}{#2}[section]%
\newenvironment{#1}{\begin{@#1}\rm}{\end{@#1}}}
 
\mynewtheo{lemma}{Lemma}
\mynewtheo{exer}{Exercise}
\mynewtheo{theo}{Theorem}
\mynewtheo{rem}{Remark}
\mynewtheo{defi}{Definition}
\mynewtheo{conj}{Conjecture}
\mynewtheo{corr}{Corollary}
\mynewtheo{prop}{Proposition}
\mynewtheo{question}{Question}
\mynewtheo{exam}{Example}

\def\inx{\mathop {\operator@font ind}\mathord{}}
\def\conv{\mathop {\operator@font conv}\mathord{}}
\def\mpb{\mathop {\operator@font mpb}\mathord{}}
\def\mwf{\mathop {\operator@font mwf}\mathord{}}
\def\spn{\mathop {\operator@font span}\mathord{}}
\def\len{\mathop {\operator@font len}\mathord{ }}
\def\str{\mathop {\operator@font star}\mathord{ }}

\def\ffrac#1#2{\mbox{\small$\ds\frac{#1}{#2}$}}

\newenvironment{theorem}{\begin{theo}}{\end{theo}}

\def\eqref{\@ifnextchar[{\@eqref}{\@eqref[]}}
\def\@eqref[#1]#2{\mbox{(\protect\reference{#2}#1})}
\def\proof{\@ifnextchar[{\@proof}{\@proof[\unskip]}}
\def\@proof[#1]{\noindent{\bf Proof #1.}\enspace}

\begin{document}

{\def\thefootnote{\fnsymbol{footnote}}
\author{A. Stoimenow 
\\[2mm]
\small Department of Mathematics, Keimyung University,\\
\small Dalseo-Gu, Sindang-Dong 1000, Daegu 704-701, Korea \\[1mm] 
\small e-mail: {\tt stoimeno@stoimenov.net},\\[0mm]
\small WWW: {\hbox{\web|http://stoimenov.net/stoimeno/homepage/|}}
}

}

\title{\large\bf
\uppercase{Diagram genus, generators and applications}
\\[4mm] \phantom{\small\it This is preprint. I would be grateful
for any comments and corrections.}}

\date{}

\maketitle

\makeatletter
\def\chrd#1#2{\picline{1 #1 polar}{1 #2 polar}}
\def\arrow#1#2{\picvecline{1 #1 polar}{1 #2 polar}}

\def\labch#1#2#3{\chrd{#1}{#2}\picputtext{1.3 #2 polar}{$#3$}}
\def\labar#1#2#3{\arrow{#1}{#2}\picputtext{1.3 #2 polar}{$#3$}}
\def\labbr#1#2#3{\arrow{#1}{#2}\picputtext{1.3 #1 polar}{$#3$}}

\def\rottab#1#2{%
\expandafter\advance\csname c@table\endcsname by -1\relax
\centerline{%
\rbox{\centerline{\vbox{\setbox1=\hbox{#1}%
\centerline{\mbox{\hbox to \wd1{\hfill\mbox{\vbox{{%
\caption{#2}}}}\hfill}}}%
\vskip9mm
\centerline{
\mbox{\copy1}}}}%
}%
}%
}

\def\GD{\szCD{6mm}}
\def\szCD#1#2{{\let\@nomath\@gobble\small\diag{#1}{2.4}{2.4}{
  \picveclength{0.27}\picvecwidth{0.1}
  \pictranslate{1.2 1.2}{
    \piccircle{0 0}{1}{}
    #2
}}}}
\def\CD{\szCD{4mm}}

\def\point#1{{\picfillgraycol{0}\picfilledcircle{#1}{0.08}{}}}
\def\labpt#1#2#3{\pictranslate{#1}{\point{0 0}\picputtext{#2}{$#3$}}}

\def\vrt#1{{\picfillgraycol{0}\picfilledcircle{#1}{0.09}{}}}
\def\svrt#1{{\picfillgraycol{0}\picfilledcircle{#1}{0.06}{}}}
\def\lvrt#1#2#3{
  {
     \pictranslate{#1}{
       \picfillgraycol{0}
       \picfilledcircle{0 0}{0.09}{}
       \picputtext{#2}{#3}
     }
  }
}

\def\@dcont{}
\def\svCD#1{\ea\glet\csname #1\endcsname\@dcont}
\def\rsCD#1{\ea\glet\ea\@dcont\csname #1\endcsname\ea\glet
\csname #1\endcsname\relax}

\def\addCD#1{\ea\gdef\ea\@dcont\ea{\@dcont #1}}
\def\drawCD#1{\szCD{#1}{\@dcont}}

\def\noloop{{\diag{0.5cm}{0.5}{1}{\picline{0.25 0}{0.25 1}}}}

\def\vrt#1{{\picfillgraycol{0}\picfilledcircle{#1}{0.09}{}}}

\def\ReidI#1#2{
  \diag{0.5cm}{0.9}{1}{
    \pictranslate{0.4 0.5}{
      \picscale{#11 #21}{
        \picmultigraphics[S]{2}{1 -1}{
	  \picmulticurve{-6 1 -1 0}{0.5 -0.5}{0.5 0}{0.1 0.3}{-0.2 0.3}
	} 
	\piccirclearc{-0.2 0}{0.3}{90 270}
      }
    }
  }
}
\def\Pos#1#2{{\diag{#1}{1}{1}{#2
\picmultiline{-5 1 -1 0}{0 1}{1 0}
\picmultiline{-5 1 -1 0}{0 0}{1 1}
}}}
\def\Neg#1#2{{\diag{#1}{1}{1}{#2
\picmultiline{-5 1 -1 0}{0 0}{1 1}
\picmultiline{-5 1 -1 0}{0 1}{1 0}
}}}
\def\Nul#1#2{{\diag{#1}{1}{1}{#2
\piccirclearc{0.5 1.4}{0.7}{-135 -45}
\piccirclearc{0.5 -0.4}{0.7}{45 135}
}}}
\def\Inf#1#2{{\diag{#1}{1}{1}{#2
\piccirclearc{0.5 1.4 x}{0.7}{135 -135}
\piccirclearc{0.5 -0.4 x}{0.7}{-45 45}
}}}
\def\pos{\Pos{0.5em}{\piclinewidth{10}}}
\def\neg{\Neg{0.5em}{\piclinewidth{10}}}
\def\nul{\Nul{0.5em}{\piclinewidth{10}}}
\def\iinf{\Inf{0.5em}{\piclinewidth{10}}}

\def\pt#1{{\picfillgraycol{0}\picfilledcircle{#1}{0.1}{}}}
\def\ppt#1{{\picfillgraycol{0}\picfilledcircle{#1}{0.06}{}}}

\def\@curvepath#1#2#3{%
  \@ifempty{#2}{\piccurveto{#1 }{@stc}{@std}#3}%
    {\piccurveto{#1 }{#2 }{#2  #3  0.5 conv}
    \@curvepath{#3}}%
}
\def\curvepath#1#2#3{%
  \piccurve{#1 }{#2 }{#2 }{#2  #3  0.5 conv}%
  \picPSgraphics{/@stc [ #1  #2  -1 conv ] $ D /@std [ #1  ] $ D }%
  \@curvepath{#3}%
}

\def\@opencurvepath#1#2#3{%
  \@ifempty{#3}{\piccurveto{#1 }{#1 }{#2 }}%
    {\piccurveto{#1 }{#2 }{#2  #3  0.5 conv}\@opencurvepath{#3}}%
}
\def\opencurvepath#1#2#3{%
  \piccurve{#1 }{#2 }{#2 }{#2  #3  0.5 conv}%
  \@opencurvepath{#3}%
}

\def\@br#1#2#3#4#5#6{
  \pictranslate{#1}{
     \picrotate{#2}{
        #6{
          \picmultigraphics[S]{2}{-1 1}{
            \picmulticurve{-7 1 -1 0}
	      {#3 -.5 * #3 .5 *}{#3 -.2 * #3 .5 *}
              {#3 .2 * #3 -.5 *}{#3 .5 * #3 -.5 *}
          }
        }
        \picmultigraphics[S]{2}{-1 1}{
	   \@ifempty{#5}{
             \picmultigraphics[S]{2}{1 -1}{
               \picline{#3 .5 * d}{#4 .5 * #3 .5 *}
	     }
	   }{
	     \picscale{1 -1}{\picline{#3 .5 * d}{#4 .5 * #3 .5 *}}
	     \picvecline{#3 .5 * d}{#4 #3 d + + 6 : #3 .5 *}
	     \piclineto{#4 .5 * #3 .5 *}
	   }
        }
     }
  }
}

\def\lbr#1#2#3#4{\@br{#1}{#2}{#3}{#4}{}{\picscale{-1 1}}}
\def\lvbr#1#2#3#4{\@br{#1}{#2}{#3}{#4}{a}{\picscale{-1 1}}}
\def\rbr#1#2#3#4{\@br{#1}{#2}{#3}{#4}{}{}}
\def\rvbr#1#2#3#4{\@br{#1}{#2}{#3}{#4}{a}{}}

\long\def\@makecaption#1#2{%
   \vskip 10pt
   {\let\label\@gobble
   \let\ignorespaces\@empty
   \xdef\@tempt{#2}%
   }%
   \ea\@ifempty\ea{\@tempt}{%
   \setbox\@tempboxa\hbox{%
      \fignr#1#2}%
      }{%
   \setbox\@tempboxa\hbox{%
      {\fignr#1:}\capt\ #2}%
      }%
   \ifdim \wd\@tempboxa >\captionwidth {%
      \rightskip=\@captionmargin\leftskip=\@captionmargin
      \unhbox\@tempboxa\par}%
   \else
      \hbox to\captionwidth{\hfil\box\@tempboxa\hfil}%
   \fi}%
\def\fignr{\small\sffamily\bfseries}%
\def\capt{\small\sffamily}%

\newdimen\@captionmargin\@captionmargin2\parindent\relax
\newdimen\captionwidth\captionwidth\hsize\relax

\def\nin{\not\in}
\def\bQ{{\Bbb Q}}
\def\bC{{\Bbb C}}
\def\bR{{\Bbb R}}
\def\bN{{\Bbb N}}
\def\bZ{{\Bbb Z}}
\def\cE{{\cal E}}
\def\cL{{\cal L}in}
\def\cK{{\cal K}}
\def\cI{{\cal I}}
\def\cR{{\cal R}}
\def\cV{{\cal V}}
\def\cO{{\cal O}}
\def\hD{{\hat D}}
\def\cD{{\cal D}}
\def\cX{{\cal X}}
\def\hK{{\hat K}}
\def\bm{\bar t'_2}
\let\bt\bm 
\def\R{\Re e\,}
\def\I{\Im m\,}
\def\pr{\text{\rm pr}\,}
\def\ncap{\not\mathrel{\cap}}
\def\|{\mathrel{\kern1.5pt\Vert\kern1.5pt}}
\def\So{\Rightarrow}
\def\lra{\longrightarrow}
\let\ds\displaystyle
\let\reference\ref
\def\so{\Rightarrow}
\let\x\exists
\let\fa\forall
\let\es\enspace

\let\x\exists
\let\sg\sigma
\let\Sg\Sigma
\let\tl\tilde
\let\ap\alpha
\let\dl\delta
\let\Dl\Delta
\let\be\beta
\let\gm\gamma
\let\Gm\Gamma
\let\nb\nabla
\let\Lm\Lambda
\let\lm\lambda
\let\sm\setminus
\let\vn\varnothing
\let\dt\det
\let\pa\partial 
\let\bd\pa
\def\tG{\tl G}
\def\tD{\tl D}
\def\tS{\tl S}

\let\eps\varepsilon
\let\ul\underline
\let\ol\overline
\def\md{\min\deg}
\def\Md{\max\deg}
\def\Mcf{\max{\operator@font cf}}
\let\Mc\Mcf
\def\sgn{{\operator@font sgn}}
\def\_#1{/_{\!\!e}\,}
\def\^#1{/\!\!/_{\!\!e}\,}

\def\ssim{\stackrel{\ds \sim}{%
  \vbox{\vskip-0.2em\hbox{$\scriptstyle *$}}}}
\def\TM{$^\text{\raisebox{-0.2em}{${}^\text{TM}$}}$}

\def\bysame{\same[\kern2cm]\,}
\def\qed{\hfill\@mt{\Box}}
\def\@mt#1{\ifmmode#1\else$#1$\fi}

\def\proof{\@ifnextchar[{\@proof}{\@proof[\unskip]}}
\def\@proof[#1]{\noindent{\bf Proof #1.}\enspace}

\def\myfrac#1#2{\raisebox{0.2em}{%
  \small$#1$}\!/\!\raisebox{-0.2em}{\small$#2$}}
\def\ffrac#1#2{\mbox{\small$\ds\frac{#1}{#2}$}}

\newcommand{\mybr}[2]{\text{$\Bigl\lfloor\mbox{%
\small$\displaystyle\frac{#1}{#2}$}\Bigr\rfloor$}}
\def\mybrtwo#1{\mbox{\mybr{#1}{2}}}

\def\br#1{\left\langle#1\right\rangle}
\def\BR#1{\left\lceil#1\right\rceil}
\def\Br#1{\left\lfloor#1\right\rfloor}

\def\epsfs#1#2{{\ifautoepsf\unitxsize#1\relax\else
\epsfxsize#1\relax\fi\epsffile{#2.eps}}}
\def\epsfsv#1#2{{\vcbox{\epsfs{#1}{#2}}}}
\def\epsfsb#1#2{{\bcbox{\epsfs{#1}{#2}}}}
\def\vcbox#1{\setbox\@tempboxa=\hbox{#1}\parbox{\wd\@tempboxa}{\box
  \@tempboxa}}
\def\bcbox#1{\setbox\@tempboxa=\hbox{#1}\parbox[b]{\wd\@tempboxa}{\box
  \@tempboxa}}

\def\@test#1#2#3#4{%
  \let\@tempa\go@
  \@tempdima#1\relax\@tempdimb#3
    \@tempdima\relax\@tempdima#4\unitxsize\relax
  \ifdim \@tempdimb>\z@\relax
    \ifdim \@tempdimb<#2%
      \def\@tempa{\@test{#1}{#2}}%
    \fi
  \fi
  \@tempa
}

\def\go@#1\@end{}
\newdimen\unitxsize
\newif\ifautoepsf\autoepsftrue

\unitxsize4cm\relax
\def\epsfsize#1#2{\epsfxsize\relax\ifautoepsf
  {\@test{#1}{#2}{0.1 }{4   }
		{0.2 }{3   }
		{0.3 }{2   }
		{0.4 }{1.7 }
		{0.5 }{1.5 }
		{0.6 }{1.4 }
		{0.7 }{1.3 }
		{0.8 }{1.2 }
		{0.9 }{1.1 }
		{1.1 }{1.  }
		{1.2 }{0.9 }
		{1.4 }{0.8 }
		{1.6 }{0.75}
		{2.  }{0.7 }
		{2.25}{0.6 }
		{3   }{0.55}
		{5   }{0.5 }
		{10  }{0.33}
		{-1  }{0.25}\@end
		\ea}\ea\epsfxsize\the\@tempdima\relax
		\fi
		}

\input{myeqn.tex}

\let\diagram\diag

\def\boxed#1{\diagram{1em}{1}{1}{\picbox{0.5 0.5}{1.0 1.0}{#1}}}

\def\rato#1{\hbox to #1{\rightarrowfill}}
\def\arrowname#1{{\enspace
\setbox7=\hbox{F}\setbox6=\hbox{%
\setbox0=\hbox{\footnotesize $#1$}\setbox1=\hbox{$\to$}%
\dimen@\wd0\advance\dimen@ by 0.66\wd1\relax
$\stackrel{\rato{\dimen@}}{\copy0}$}%
\ifdim\ht6>\ht7\dimen@\ht7\advance\dimen@ by -\ht6\else
\dimen@\z@\fi\raise\dimen@\box6\enspace}}

\def\contr{\diagram{1em}{0.6}{1}{\piclinewidth{35}%
\picstroke{\picline{0.5 1}{0.2 0.4}%
\piclineto{0.6 0.6}\picveclineto{0.3 0}}}}

\newcounter{pp}%
\newenvironment{mylist}[1]{%
\begin{list}{#1{pp})}%
{\usecounter{pp}\setlength{\labelwidth}{4mm}%
\setlength{\leftmargin}{0.6cm}\setlength{\itemsep}{1mm}%
\setlength{\topsep}{1mm}}%
\gdef\myitem{\item\xdef\@currentlabel{#1{pp})}}%
}{\end{list}}

\def\abstractname{}

\parskip5pt plus 1pt minus 2pt
\parindent\z@

{\let\@noitemerr\relax
\vskip-2.7em\kern0pt\begin{abstract}
\noindent{\bf Abstract.}\enspace
We continue the study of the genus of knot diagrams, deriving
a new description of generators using Hirasawa's algorithm.
This description leads to good estimates on the maximal number
of crossings of generators and allows us to complete their
classification for knots of genus 4.

As applications of the genus 4 classification, we establish
non-triviality of the skein polynomial on $k$-almost positive
knots for $k\le 4$, and of the Jones polynomial for $k\le 3$.
For $k\le 4$, we classify the occurring achiral knots, and
prove a trivializability result for $k$-almost positive unknot
diagrams. This yields also estimates on the number of unknotting
Reidemeister moves. We describe the positive knots of signature
(up to) 4.

Using a study of the skein polynomial, we
prove the exactness of the Morton-Williams-Franks braid
index inequality and the existence of a minimal string 
Bennequin surface for alternating knots up to genus $4$.
We also prove for such knots conjectures of Hoste and Fox
about the roots and coefficients of the Alexander polynomial.
\\[1mm]
\end{abstract}
}
\vspace{1mm}

{\parskip0.2mm\tableofcontents}
\vspace{7mm}

\section{Introduction}

Introducing the \em{genus} $g(K)$ of a knot $K$, Seifert
\cite{Seifert} gave a construction of compact oriented surfaces
in 3-space bounding the knot (\em{Seifert surface}) by an
algorithm starting with some diagram of the knot (see
\cite[\S 4.3]{Adams} or \cite{Rolfsen}). The surface given
by this algorithm is called \em{canonical}. A natural problem
is when the diagram is \em{genus-minimizing}, or \em{of
minimal genus}, that is, its canonical Seifert surface
has minimal genus among all Seifert surfaces of the knot.
This problem has been studied over a long period. First, the
minimal genus property was shown for alternating diagrams,
independently by Crowell \cite{Crowell} and Murasugi \cite{Murasugi2}.
Their proof is algebraical, using the Alexander polynomial
$\Dl$ \cite{Alexander} and the inequality $\Md\Dl\le g$
(which thus they prove to be exact for alternating knots).
Later Gabai \cite{Gabai} developed a geometric method using
foliations and showed that this method is likewise successful
for alternating diagrams. Murasugi \cite{Murasugi3} introduced
the operation *-product. It was shown to behave naturally
w.r.t.\ the Alexander polynomial by himself, and later
by Gabai \cite{Gabai2,Gabai3} in the geometrical context.
These results imply the extension of the minimal genus property
of alternating diagrams to the homogeneous diagrams of \cite{Cromwell}.
An important other subclass of the class of homogeneous diagrams
and links are the positive diagrams and links. Such links
have been considered (in general or in special cases) independently
before. (See e.g. \cite{CochranGompf,MorCro,WilFr,Ozawa,Rudolph,%
Traczyk,Yokota,Zulli}.) The minimal genus property for
such diagrams follows from yet a different source, the work
of Bennequin \cite{Bennequin} on contact structures. His inequality
(theorem 3 in that paper; stated as theorem \reference{Biq} below)
in fact allows to estimate the difference 
between the genus of the diagram and the genus of the knot in
terms of the number of positive or negative crossings.

The prospect of applications led to the treatment of canonical
surfaces, and the establishment of the \em{canonical genus}, the
minimal genus of all such surfaces for a given knot or link, in its
own right. In this paper, we will continue the study of canonical
Seifert surfaces from the combinatorial point of view. This study
was initiated my myself in \cite{gen1}, and independently by
Brittenham \cite{Brittenham}, and set forth in \cite{gen2}, and
later in \cite{STV,SV}. (Some explanation of this work is given
also in section 5.3 of Cromwell's recent book \cite{Cromwell2}.) 
Using the theoretical insight gained there we had methods
efficient enough to complete the calssification up to genus 3,
in terms of the list of ``generators''. In \cite{SV} the relation
to certain algebraic objects named Wicks forms was discussed,
and applied to the enumeration of alternating knots by genus.
We gave an extensive list of other applications in \cite{gen2},
and later for example in \cite{gbi2,gsigex}.

The previous method of \cite{SV} does not apply well for
several components, and thus we develop an alternative
approach, which bases on the special diagram algorithm
first found by Hirasawa \cite{Hirasawa2} (and rediscovered
a little later independently in \cite[\S 7]{bseq}). We will
work out inequalities for the crossing number and number of
$\sim$-equivalence classes of $\bm$-irreducible link diagrams
(generators). So far we carried this out only for knots. The
present approach allows to improve what we know in the knot
case and extend it to links.

Then we treat the description, and applications, of
the compilation of knot generators of genus $4$. As has already
become apparent in the preceding generator compilations for genus
$2$ and $3$, and then also from the result of \cite{SV}, the growth
of the number of generators is enormous. In order to push the
task back within the limits of (reasonable) computability, we need
important new theoretical knowledge. It relies heavily on the
special diagram algorithm of Hirasawa and the previous careful
analysis of this procedure we carry out. 

The first application addresses the identification of the unknot.
This a basic problem in knot theory. There are general methods
using braid foliations \cite{BirMen2}, and Haken theory \cite{HL}.
However, these methods are difficult to use in practice. So one is
interested in applicable criteria, at least for special classes
of diagrams. 

It is known, for example from Crowell \cite{Crowell} and Murasugi
\cite{Murasugi2}, that in alternating unknot diagrams all crossings
are nugatory. For positive, as for alternating, diagrams, one can
observe the same phenomenon using the Alexander polynomial (as in
\cite{Cromwell}). It is also an application of the Bennequin
inequality (theorem \reference{Biq} below), in its extended form
based on the Vogel algorithm \cite{Vogel} (see \cite{pos}). For almost
alternating diagrams a trivializibility result using flypes and
``tongue moves'' (see figure \reference{fig1}) was proved recently
by Tsukamoto \cite{Ts}, confirming a previous conjecture of Adams. 
Here we prove several related results, among others the following:

\begin{theorem}\label{th1.6}
For $k\le 4$, all $k$-almost positive unknot diagrams are trivializable
by crossing number reducing wave moves and factor slides.
\end{theorem}

A wave move is shown in figure \ref{fwv}, and a factor slide
in figure \reference{figtan} (b).

For $k$-almost positive diagrams the description in the case $k=1$
(unknotted twist knot diagrams and possible nugatory crossings)
was written down in \cite{apos}. It was, though, previously known
to Przytycki, who observed it as a consequence of Taniyama's
result \cite{Taniyama}. Their joint work had been announced long
ago, but the full paper was not finished until very recently
\cite{PrzTan}. This description includes and concretifies
Hirasawa's result \cite{Hirasawa} on almost special alternating
diagrams. For 2-almost positive diagrams (in particular 2-almost
special alternating diagrams) the result is given in \cite{gen2}.
Now we settle the cases $k\le 4$ in the stated way.

This theorem is also quite analogous to the result\footnote{I was
informed, though, of a possible gap in its proof.} in \cite{NO} for
3-bridge diagrams, and for arborescent (or Conway-algebraic) diagrams,
which is a consequence of Bonahon-Siebenmann's classification of
arborescent links. A partial writeup of their (still unpublished)
work, which covers the treatment of the unknot, is given in \cite{FG}.

A tongue move is a special type of wave move, preserving the property
`almost alternating'. From this point of view, we can extend
Tsukamoto's (and Hirasawa's) result for $k=1$ in special diagrams
to $k\le 5$ (see corollary \ref{cr!} and theorem \ref{t5p}). In
fact, this extension was one of the motivations for our attention
to this type of problem. We can apply our work also to give
polynomial estimates on the number of Reidemeister moves needed for
unknotting, which has been recently studied by various authors
\cite{HL,HN,Hayashi} (see proposition \reference{pRei}).

{}From a different perspective, a factor slide is a (very special)
type of preserving wave move. Then theorem \ref{th1.6} confirms in
our case, in a stronger form, a conjecture (conjecture \ref{cu}),
stating that preserving and reducing wave moves in combination
suffice to trivialize unknot diagrams. Example \ref{xzq} shows that
our stronger statement is not true in general, and the assumption
$k\le 4$ in the theorem cannot be improved. 

The work on theorem \ref{th1.6} allows us to address two related
problems which have been given some attention in the literature.
The following result relates to the non-triviality problem of link
polynomials.


\begin{theorem}\label{th1.7}
Any non-trivial $k$-almost positive knot has non-trivial skein
polynomial $P$ for $k\le 4$ and non-trivial Jones polynomial
$V$ for $k\le 3$.
\end{theorem}

So far, only the cases $k=0,1$ were known (see \cite{restr}). For
$P$, the result follows from the proof of theorem \ref{th1.6} almost
directly. For $V$, we combine the skein and Kauffman bracket properties
of the Jones polynomial in \cite{Bmo,restr,adeq} with some arguments
about the signature and Gabai's geometric work \cite{Gabai2,Gabai3}.
(They allow us to reduce the necessary computations to a minimum.)

The proof of theorem \reference{th1.6} can be applied in some form
not only to the trivial, but more generally to an amphicheiral
(achiral) knot (see corollary \reference{Cz}). We can then obtain
the classification of these knots.

\begin{theorem}\label{th1.8}
The $k$-almost positive achiral knots for $k\le 4$ are those
of $2k$ crossings, that is, $4_1$ for $k=2$; $6_3$ and $3_1\#!3_1$
for $k=3$; and $8_3$, $8_9$, $8_{12}$, $8_{17}$ $8_{18}$, and
$4_1\#4_1$ for $k=4$. 
\end{theorem}

This property was known for $k\le 2$ by \cite{PrzTan}; prior to
latter's completion, written account was given in \cite{restr} for
$k\le 1$ and \cite{gen2} for $k=2$. Theorem \ref{th1.8} adds the
cases $k=3,4$. It was obtained in \cite{restr} for alternating prime
knots (independently provable using the work in \cite{Thistle}),
and checked there also for prime knots of $\le 16$ crossings.
The proof of theorem \ref{th1.8} will consist in extending the
arguments for theorem \ref{th1.6} and simplifying diagrams of such
knots to these checked low crossing cases (see corollary \ref{Cz}).

A final application concerning positivity is the determination
of the positive knots of signature 4 (theorem \ref{tsg4}). Again
the result for signature 2 (namely, that these are precisely the
positive knots of genus 1) is a consequence of Taniyama's work
\cite{Taniyama}. Our theorem is the first non-trivial explicit step
beyond Taniyama toward the general case of the conjecture that
positive knots of only finitely many genera have given signature;
see \cite{gsig,gsigex}.

The later sections of the paper extend the applications of generators
and regularization to alternating knots. In \S\reference{S6}, we
use a regularization of the skein polynomial and the work of
Murasugi-Przytycki \cite{MP} to show

\begin{theorem}\label{tBi}
The Morton-Williams-Franks braid index inequality MWF
is exact on alternating knots of genus up to 4.
\end{theorem}

We will also confirm a conjecture of Murasugi-Przytycki
(conjecture \reference{C7}) for such knots. 
Theorem \ref{tBi} complements the
results in \cite{Murasugi} for fibered and 2-bridge knots;
it comes close to the examples of strict MWF inequality of a
4-component genus 3 (alternating) link and a genus 6 knot
found in \cite{MP}. Evidence from some computations during
its proof led us to conjecture an improvement of Ohyama's
\cite{Ohyama} braid index inequality for special alternating
links (\S\ref{OhS}). We confirm this stronger inequality
(apart from knots of genus up to 4) also for arborescently
alternating links.

Our study of Murasugi-Przytycki's index requires also to clarify
and correct an inexactness in \cite{MP}, which affects the proof
of their main upper braid index estimate (see \S\ref{S61}).
We revealed a slight discrepancy between their definition
of graph index and the diagrammatic move they introduce to
reduce the number of Seifert circles. We will be able to
justify their estimates, but still there is some cost, in that
the correspondence between
a diagram and (its Seifert) graph becomes (in general) lost
during the recursive calculation of the index (using their old
definition). This oversight seems to propagate to other papers
and may cause a problem at some point. We will explain how
to define, on the level of Seifert graphs, the ``right'' index
w.r.t. their move.

In \S\reference{SMS} we develop a method of constructing a
Bennequin surface of links for a minimal genus canonical surface.
The number of strings of the resulting braid is generally
low, and can be calculated by modifying our corrected version
of graph index. We apply this construction to alternating
knots of genus up to 4, and show that one can span for these
knots a Bennequin surface on a minimal string braid. Bennequin
\cite{Bennequin} proved such a result for all 3-braid links,
and Hirasawa (unpublished) for the 2-bridge links. In contrast,
knots lacking such surfaces are known for braid index 4 and
genus 3. We will see how to apply the Murasugi-Przytycki move
keeping a minimal genus surface, finally converting it into a
minimal string Bennequin surface.

The final section \S\reference{SHo_} deals with the Alexander
polynomial of an alternating knot. We consider first a conjecture
of Hoste about the roots of the Alexander polynomial. We develop
several tests based on the generator description, and apply them to
confirm this conjecture on knots up to genus 4 (theorem \ref{thHo}).

Another conjecture we treat is the $\log$-concavity conjecture
of \cite{spec}. It states that the square of each coefficient of
the Alexander polynomial of an alternating link is not less than
the product of its two neighbors. This implies Fox's ``trapezoidal''
conjecture \ref{Fxc}. We will verify both conjectures for knots
up to genus 4 (theorem \ref{tcc}). In this case it was possible
to extend the result to links (theorem \ref{lcL}). 

As another qualitative extension of the Fox conjecture, we will
conclude that the convex hull of alternating knot polynomials
of given genus is a polytope, and we can determine this polytope
(i.e.\ the complete set of linear inequalities satisfied by
the coefficients of Alexander polynomials of alternating knots)
up to genus 4 in \S\ref{rcs}. Apart from an improvement of
the trapezoidal inequalities, we will compare our result also
to the inequalities obtained by Ozsv\'ath-Szab\'o using knot
Floer homology \cite{OS}. (Their inequalities give an alternative
proof of the trapezoidal conjecture for genus 2.) In-Dae
Jong \cite{Jong} used the generator description (theorem \ref
{thgen}) to prove the $\log$-concavity conjecture and then,
following the discussion in \S\ref{rcs}, to give the complete
set of linear inequalities for genus 2.

\section{Preliminaries}

\subsection{Graphs}

A graph $G$ will have for us possibly multiple edges (edges
connecting the same two vertices), but usually no
loop edges (edges connecting one and the same vertex).
By $V(G)$ we will denote the set of vertices of $G$, and by $E(G)$ the
set of edges of $G$ (each multiple edge counting as a set of
single edges); $v(G)$ and $e(G)$ will be the number of
vertices and edges of $G$ (thus counted), respectively.

For a graph, let the operation
\[
\diag{1cm}{1}{0.5}{
  \picfillgraycol{0}
  \picfilledcircle{0 0.25}{0.08}{}
  \picfilledcircle{1 0.25}{0.08}{}
  \picline{0 0.25}{1 0.25}
}
\qquad\lra\qquad
\diag{1cm}{1.5}{0.2}{
  \picfillgraycol{0}
  \picfilledcircle{0 0.1}{0.08}{}
  \picfilledcircle{0.75 0.1}{0.08}{}
  \picfilledcircle{1.5 0.1}{0.08}{}
  \picline{0 0.1}{1.5 0.1}
}
\]
(adding a vertex of valence $2$) be called \em{bisecting} and
its inverse (removing such a vertex) \em{unbisecting} (of an
edge). We call a graph $G'$ a \em{bisection} of a graph $G$
with no valence-2-vertices,
if $G'$ is obtained from $G$ by a sequence of edge bisections.
We call a bisection $G'$ \em{reduced}, if it has no adjacent
vertices of valence $2$ (that is, each edge of $G$ is bisected
at most once). Contrarily, if $G'$ is a graph, its \em{unbisected
graph} $G$ is the graph with no valence-2-vertices, of which
$G'$ is a bisection.

Similarly, a \em{contraction} is  the operation
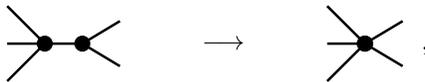
\[
\diag{1cm}{3.5}{2}{
  \vrt{1.5 1}
  \picline{1.5 1}{1 1.5}
  \picline{1.5 1}{1 1}
  \picline{1.5 1}{1 0.5}
  \picline{1.5 1}{2 1}
  \pictranslate{-1.5 0}{
    \vrt{3.5 1}
    \picline{3.5 1}{4 1.3}
    \picline{3.5 1}{4 0.7}
  }
}\lra
\diag{1cm}{2.2}{2}{
  \vrt{1.5 1}
  \picline{1.5 1}{1 1.5}
  \picline{1.5 1}{1 1}
  \picline{1.5 1}{1 0.5}
  \pictranslate{-2 0}{
    \vrt{3.5 1}
    \picline{3.5 1}{4 1.3}
    \picline{3.5 1}{4 0.7}
  }
}\,,
\]
and a \em{decontraction} its inverse.

The \em{doubling} of an edge consists in adding a new edge
connecting the same two vertices.

A graph is \em{$n$-connected}, if $n$ is the minimal number of
edges needed to remove from it to disconnect it. (Thus connected
means $1$-connected.) Such a collection of edges is called an
\em{$n$-cut}.

Hereby, when we delete an edge, we understand that a vertex
it is incident to is \em{not} to be deleted too. In that sense, the
set $I_v$ of edges incident to a given vertex $v$ always forms a
cut~-- if we delete these edges, $v$ gets an isolated component,
so the graph is not connected anymore.
A \em{cut vertex} is a vertex which disconnects a graph,
when removed \em{together} with all its incident edges.

A graph $G$ is \em{planar} if it is embeddable in the plane \em{and
equipped} with a fixed such planar embedding.
Observe that there is a natural bijection of edges between a planar
graph $G$ and its dual graph $G^*$; in that sense we can talk of
the dual $e^*\in E(G^*)$ of an edge $e\in E(G)$. 
The operations doubling and bisection become dual to each other.

\subsection{Knots and diagrams\label{kkd}}

A crossing $p$ in a knot diagram $D$ is called \em{reducible} (or
nugatory) if it looks like on the left of figure \reference{figred}.
$D$ is called reducible if it has a reducible crossing, else it is
called \em{reduced}. The reducing of the reducible crossing $p$
is the move depicted on figure \reference{figred}. Each diagram 
$D$ can be (made) reduced by a finite number of these moves.

\begin{figure}[htb]
\begin{eqn}\label{eqred}
\diag{6mm}{4}{2}{
  \picrotate{-90}{\rbraid{-1 2}{1 1.4}}
  \picputtext[u]{2 0.7}{$p$}
  \picscale{1 -1}{
    \picfilledcircle{0.7 -1}{0.8}{$P$}
  }
  \picfilledcircle{3.3 1}{0.8}{$Q$}
} \qquad\lra\qquad
\diag{6mm}{4}{2}{
  \piccirclearc{2 0.5}{1.3}{45 135}
  \piccirclearc{2 1.5}{1.3}{-135 -45}
  \picfilledcircle{1 1}{0.8}{$P$}
  \picfilledcircle{3 1}{0.8}{$Q$}
} 
\end{eqn}
\caption{\label{figred}}
\end{figure}

We assume in the following all diagrams reduced, unless otherwise
stated.

\begin{figure}[htb]
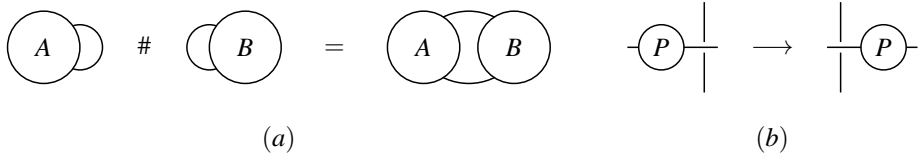

\[
\begin{array}{c@{\qquad\es}c}
\diag{6mm}{3}{2}{
  \piccirclearc{1.8 1}{0.5}{-120 120}
  \picfilledcircle{1 1}{0.8}{$A$}
}\, \#\,
\diag{6mm}{3}{2}{
  \piccirclearc{1.2 1}{0.5}{60 300}
  \picfilledcircle{2 1}{0.8}{$B$}
}\quad =\quad
\diag{6mm}{4}{2}{
  \piccirclearc{2 0.5}{1.3}{45 135}
  \piccirclearc{2 1.5}{1.3}{-135 -45}
  \picfilledcircle{1 1}{0.8}{$A$}
  \picfilledcircle{3 1}{0.8}{$B$}
} 
 &
\diag{6mm}{2}{2}{
  \picline{1.7 0}{1.7 2}
  \picmultiline{-7 1 -1 0}{0 1}{2 1}
  \picfilledcircle{0.75 1}{0.5}{$P$}
} \quad\lra\quad
\diag{6mm}{2}{2}{
  \picline{0.3 0}{0.3 2}
  \picmultiline{-7 1 -1 0}{0 1}{2 1}
  \picfilledcircle{1.25 1}{0.5}{$P$}
} \\[8mm]
(a) & (b) 
\end{array}
\]
\caption{Diagram connected sum and factor slide move\label{figtan}}
\end{figure}

Part (a) of figure \reference{figtan} displays the
\em{connected sum} $D=A\# B$ of the diagrams $A$ and $B$.
Latter are called \em{factors} of $D$. The connected sum on
diagrams is uniquely determined only up to the move shown in
part (b) of the figure (and the mirror image of that move).
We will call it a \em{factor slide}.
If a diagram $D$ can be represented as the connected sum of 
diagrams $A$ and $B$, such that both $A$ and $B$ have at least one
crossing, then $D$ is called \em{composite}, otherwise it is
called \em{prime}. A knot of link $K$ is \em{prime} if whenever
$D=A\# B$ is a composite diagram of $K$, one of $A$ and $B$
represent an unknotted arc (but not both; the unknot is
not prime per convention).

A diagram $D$ is \em{connected} if its curve is a connected set in
$\bR^2$, that is, there is no closed curve $\gm$ disjoint from $D$
such that both the interior and exterior of $\gm$ have non-trivial
intersection with $D$. Otherwise $D$ is \em{disconnected} or
\em{split}. A link is \em{split} if it has a split diagram, and
otherwise \em{non-split}.

\begin{theorem} (\cite{Menasco})\label{ThM}
If $D$ has a prime alternating non-trivial diagram of $K$, then
$K$ is prime.
\end{theorem}

The \em{(Seifert) genus} $g(K)$ resp.\ \em{Euler characteristic}
$\chi(K)$ of a knot or link $K$ is said to be the minimal genus
resp.\ maximal Euler characteristic of Seifert surface of $K$.
For a diagram $D$ of $K$, $g(D)$ is defined to be the genus
of the Seifert surface obtained by Seifert's algorithm on $D$,
and $\chi(D)$ its Euler characteristic. Let $c(D)$ denote the
\em{number of crossings} of $D$ and $n(D)=n(K)$ the \em{number
of components} of $D$ or $K$ (so $n(K)=1$ if $K$ is a knot).
Write $s(D)$ for the \em{number of Seifert circles} of $D$.
Then $\chi(D)=s(D)-c(D)$ and $2g(D)=2-n(D)-\chi(D)$.

\begin{theorem}
(see \cite{Murasugi2,Crowell,Gabai,Cromwell})
If $D$ has an alternating or positive diagram of $K$, then
$g(K)=g(D)$ and $\chi(D)=\chi(K)$.
\end{theorem}

The \em{crossing number} $c(K)$ is the minimal crossing number
of all diagrams $D$ of $K$.
The \em{canonical genus} $\tl g(K)$ resp.\ \em{canonical
Euler characteristic} $\tl\chi(K)$ is defined as the minimal genus
resp.\ maximal Euler characteristic 
of all diagrams of $K$. In general we can have $g(K)<\tl g(K)$,
that is, no diagrams of $K$ of minimal genus (see \cite{Morton}).

The \em{writhe}, or \em{(skein) sign}, is a number $\pm1$,
assigned to any crossing in a link
diagram. A crossing as in figure \ref{figwr}(a) has writhe 1 and
is called \em{positive}. A crossing as in figure \ref{figwr}(b) has
writhe $-1$ and is called \em{negative}. The \em{writhe} of a link
diagram is the sum of writhes of all its crossings.

Let $c_{\pm}(D)$ be the number of positive, respectively negative
crossings of a diagram $D$, so that $c(D)=c_+(D)+c_-(D)$ and
$w(D)=c_+(D)-c_-(D)$. Let $c_{\pm}(K)$ for a knot $K$ denote
the minimal number of positive resp. negative 
crossings of a diagram of $K$.

\begin{figure}[htb]
\[
\begin{array}{c@{\qquad}c}
\diag{6mm}{1}{1}{
\picmultivecline{0.18 1 -1.0 0}{1 0}{0 1}
\picmultivecline{0.18 1 -1.0 0}{0 0}{1 1}
} &
\diag{6mm}{1}{1}{
\picmultivecline{0.18 1 -1 0}{0 0}{1 1}
\picmultivecline{0.18 1 -1 0}{1 0}{0 1}
}
\\[2mm]
(a) & (b)
\end{array}
\]
\caption{\label{figwr}}
\end{figure}

A diagram is \em{positive} if all its crossings are positive.
A diagram is \em{almost positive} if all its crossings are positive
except exactly one. A knot is \em{positive} if it has a positive
diagram. (See e.g. \cite{MorCro,Ozawa,Yokota,Zulli}.) It is \em{almost
positive} if it is not positive but has an almost positive diagram.
More generally a diagram $D$ is \em{$k$-almost positive} if it has
exactly $k$ negative crossings, i.e.\ $c_-(D)=k$, and a knot $K$
is $k$-almost positive if it has a $k$-almost positive, but
no $k-1$-almost positive diagram, i.e.\ $c_-(K)=k$. We will
sometimes call the switch of crossings of a diagram $D$ so that
all become positive the \em{positification} of $D$.

\em{Bennequin's inequality} (theorem 3 in \cite{Bennequin}) can
be stated, using (as explained in \cite{pos}) the work by Vogel
\cite{Vogel}, thus:

\begin{theorem}\label{Biq}
If $D$ is a diagram of a knot $K$, then $g(K)\ge g(D)-c_-(D)$.
\end{theorem}

In particular, if $D$ is a diagram of the unknot, then
$c_-(D)\ge g(D)$.

We call a crossing $p$ connected to a Seifert circle $s$ also \em{%
adjacent} or \em{attached} to $s$. The \em{valency} of a Seifert
circle $s$ is the number of crossings attached to $s$. We call
a Seifert circle \em{negative} if only negative crossings are
attached to it. Let $s_-(D)$ be the number of negative Seifert circles
of $D$.

Bennequin's inequality was improved by Rudolph \cite{Rudolph}.

\begin{theorem}\label{Biq2}
If $D$ is a diagram of a knot $K$, 
then 
\begin{eqn}\label{RBi}
g(K)\,\ge\, g_s(K)\,\ge\,g(D)-c_-(D)+s_-(D)\,.
\end{eqn}
Here $g_s(K)$ is the \em{smooth slice genus} of $K$.
\end{theorem}

We will refer to \eqref{RBi} as the \em{Rudolph-Bennequin inequality}.

A diagram is \em{special} if no Seifert circle contains other
Seifert circles in both regions it separates the plane into.
Such Seifert circles are called \em{separating}. It is an
easy observation that for connected diagrams two of the 
properties alternating, positive/negative and special imply
the third. A diagram with these properties is called \em{special
alternating}. A knot is special alternating if it has a special 
alternating diagram. Such knots were introduced and studied
by Murasugi \cite{Murasugi} and have a series of special features. 
Contrarily, all knots have a special (not necessarily alternating) 
diagram. Hirasawa \cite{Hirasawa2} shows how to a modify any knot 
diagram $D$ into a special diagram $D'$ so that $g(D)=g(D')$.
(Actually, the canonical surfaces of $D$ and $D'$ are isotopic.)

A diagram is \em{almost alternating} \cite{Adams,Adams2,Adams3,GHY}
if it can be turned by one crossing change into an alternating one.
A knot is \em{almost alternating}, if it has an almost alternating
diagram, but is not alternating.

The \em{index}
$\inx(s)$ of a separating Seifert circle $s$ is defined as follows:
denote for an inner crossing (i.e., attached from the inside) of $s$
a letter `$i$', and for an outer crossing a letter `$o$' cyclically
along $s$. Then $\inx(s)$ is by definition the minimal number of
disjoint subwords of the from $i^n$ or $o^n$ ($n>0$) of this cyclic
word. For example, the Seifert circle $s$ in figure \reference{_fig1}
on page \pageref{_fig1} has index $4$.

\subsection{Diagrammatic moves}

\begin{defi}
A \em{flype} is a move on a diagram shown in figure \reference{fig_}.
We say that a crossing \em{admits} a flype if it can be represented
as the distinguished crossing in the diagrams in the figure, and both
tangles have at least one crossing.

\begin{figure}[htb]
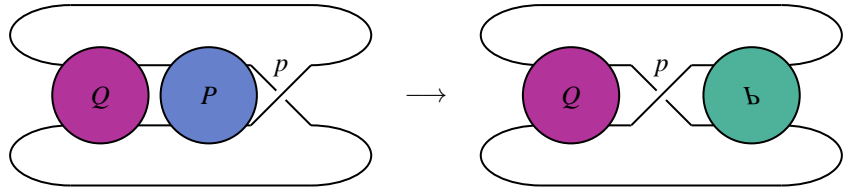

\[ 
\diag{8mm}{6}{3}{
    \pictranslate{3 1.5}{
       \picputtext[d]{1.5 0.4}{$p$}
       \picmultigraphics[S]{2}{1 -1}{
           \picmultiline{0.22 1 -1.0 0}{2 -0.5}{1 0.5}
           \picmultigraphics[S]{2}{-1 1}{
                \picellipsearc{-2 -1.0}{1 0.5}{90 270}
           }
           \picline{-2 -1.5}{2 -1.5}
           \picline{-2 -0.5}{1 -0.5}
      }
      \picfillgraycol{0.4 0.5 0.8}
      \picfilledcircle{0.3 0}{0.8}{$P$}
      \picfillgraycol{0.7 0.2 0.6}
      \picfilledcircle{-1.5 0}{0.8}{$Q$}
   }
}\quad\lra\quad
\diag{8mm}{6}{3}{
    \pictranslate{3 1.5}{
       \picputtext[d]{0 0.4}{$p$}
       \picmultigraphics[S]{2}{1 -1}{
           \picmultiline{0.22 1 -1.0 0}{0.5 -0.5}{-0.5 0.5}
           \picmultigraphics[S]{2}{-1 1}{
                \picellipsearc{-2 -1.0}{1 0.5}{90 270}
                \picline{-2 -0.5}{-0.5 -0.5}
           }
           \picline{-2 -1.5}{2 -1.5}
      }
      \picscale{1 -1}{
           \picfillgraycol{0.3 0.7 0.6}
           \picfilledcircle{1.5 0}{0.8}{$P$}
      }
      \picfillgraycol{0.7 0.2 0.6}
      \picfilledcircle{-1.5 0}{0.8}{$Q$}
   }
}
\] 
\caption{A flype near the crossing $p$\label{fig_}}
\end{figure}
\end{defi}

By the fundamental work of Menasco-Thistlethwaite, we have
a proof of the Tait flyping conjecture.

\begin{theorem}(\cite{MenThis})\label{thfl}
For two alternating diagrams of the same prime alternating link,
there is a sequence of flypes taking the one diagram into the other.
\end{theorem}

\begin{figure}[htb]
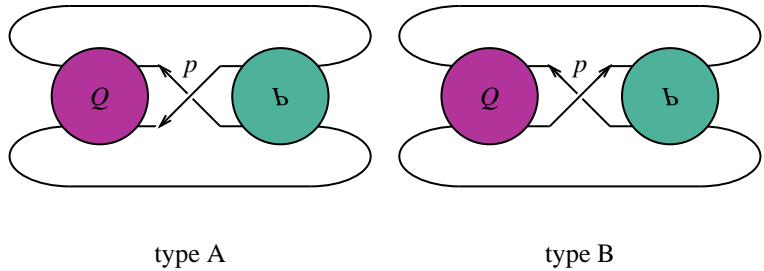

\[
\begin{array}{cc}
\diag{8mm}{6}{3}{
    \pictranslate{3 1.5}{
       \picputtext[d]{0 0.4}{$p$}
       \picmultigraphics[S]{2}{1 -1}{
           \picmultivecline{0.12 1 -1.0 0}{0.5 -0.5}{-0.5 0.5}
           \picmultigraphics[S]{2}{-1 1}{
                \picellipsearc{-2 -1.0}{1 0.5}{90 270}
                \picline{-2 -0.5}{-0.5 -0.5}
           }
           \picline{-2 -1.5}{2 -1.5}
      }
      \picscale{1 -1}{
           \picfillgraycol{0.3 0.7 0.6}
           \picfilledcircle{1.5 0}{0.8}{$P$}
      }
      \picfillgraycol{0.7 0.2 0.6}
      \picfilledcircle{-1.5 0}{0.8}{$Q$}
   }
}
&
\diag{8mm}{6}{3}{
    \pictranslate{3 1.5}{
       \picputtext[d]{0 0.4}{$p$}
       \picmultigraphics[S]{2}{1 -1}{
           \picmultiline{0.12 1 -1.0 0}{0.5 -0.5}{-0.5 0.5}
           \picmultigraphics[S]{2}{-1 1}{
                \picellipsearc{-2 -1.0}{1 0.5}{90 270}
                \picline{-2 -0.5}{-0.5 -0.5}
           }
           \picline{-2 -1.5}{2 -1.5}
       }
       \picvecline{-0.3 0.3}{-0.5 0.5}
       \picvecline{0.3 0.3}{0.5 0.5}
      \picscale{1 -1}{
           \picfillgraycol{0.3 0.7 0.6}
           \picfilledcircle{1.5 0}{0.8}{$P$}
      }
      \picfillgraycol{0.7 0.2 0.6}
      \picfilledcircle{-1.5 0}{0.8}{$Q$}
   }
}\\[17mm]
\mbox{type A} & \mbox{type B} 
\end{array}
\]
\caption{A flype of type A and B\label{ffl2}}
\end{figure}

We introduced (see \cite{SV}) a distinction of flypes
according to the orientation near the crossing $p$ at which
the flype is performed. See figure \reference{ffl2}.
An important observation is that each crossing admits at most
one of the types A and B of flypes, and this remains so after
applying any sequence on flypes on the diagram.

A \em{bridge} is a piece of a strand of a knot diagram
that passes only crossings from above.
The \em{length} of the bridge is the number of crossings passed 
by it (\em{excluding} the initiating and terminating underpass).
A \em{tunnel} is the mirror image of a bridge.

\begin{figure}[htb]
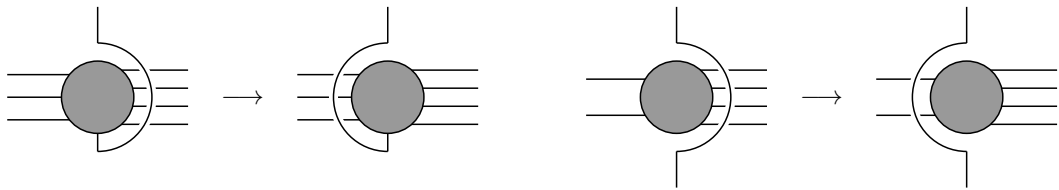

\[
\diag{6mm}{4}{4}{
  \picfillgraycol{0.6}
  \picmultigraphics{3}{0 -0.5}{
    \picline{0 2.5}{2 2.5}
  }
  \picmultigraphics{4}{0 -0.4}{
    \picline{4 2.6}{2 2.6}
  }
  \picfilledcircle{2 2}{0.8}{}
  \picmulticirclearc{-6 1 -1 0}{2 2}{1.2}{-90 90}
  \picline{2 4}{2 3.2}
  \picline{2 1.2}{2 0.8}
}\quad\lra\quad
\diag{6mm}{4}{4}{
  \picfillgraycol{0.6}
  \picmultigraphics{3}{0 -0.5}{
    \picline{0 2.5}{2 2.5}
  }
  \picmultigraphics{4}{0 -0.4}{
    \picline{4 2.6}{2 2.6}
  }
  \picfilledcircle{2 2}{0.8}{}
  \picmulticirclearc{-6 1 -1 0}{2 2}{1.2}{90 -90}
  \picline{2 4}{2 3.2}
  \picline{2 1.2}{2 0.8}
}\quad\quad\quad\quad
\diag{6mm}{4}{4}{
  \picfillgraycol{0.6}
  \picmultigraphics{2}{0 -0.8}{
    \picline{0 2.4}{2 2.4}
  }
  \picmultigraphics{4}{0 -0.4}{
    \picline{4 2.6}{2 2.6}
  }
  \picfilledcircle{2 2}{0.8}{}
  \picmulticirclearc{-6 1 -1 0}{2 2}{1.2}{-90 90}
  \picline{2 4}{2 3.2}
  \picline{2 0}{2 0.8}
}\quad\lra\quad
\diag{6mm}{4}{4}{
  \picfillgraycol{0.6}
  \picmultigraphics{2}{0 -0.8}{
    \picline{0 2.4}{2 2.4}
  }
  \picmultigraphics{4}{0 -0.4}{
    \picline{4 2.6}{2 2.6}
  }
  \picfilledcircle{2 2}{0.8}{}
  \picmulticirclearc{-6 1 -1 0}{2 2}{1.2}{90 -90}
  \picline{2 4}{2 3.2}
  \picline{2 0}{2 0.8}
}
\]
\caption{Wave-moves. The number of strands on left and right
of the shaded cicrle may vary. It is only important that the
parities are equal resp. different, and that the left-outgoing
strands are fewer that the right-outgoing ones.\label{fwv}}
\end{figure}

A \em{wave move} is a replacement of a
bridge $a$ of length $l_1$ by another one $b$ of length $l_2$.
We will assume throughout that $l_2\le l_1$.
The move is (crossing number) \em{reducing}, if $l_2<l_1$ 
(see figure \ref{fwv} or \cite{bridge} for example), and
(crossing number) \em{preserving}, if $l_2=l_1$.
Elsewhere a wave move is also called a \em{$(l_1,l_2)$-pass}.
Note that a $(1,0)$-pass is the removal of a nugatory crossing
(much like in figure \ref{figred}, except that $P$ is flipped),
and a $(1,1)$-pass is exactly the factor slide of figure
\reference{figtan} (b). We will usually consider only reducing
wave moves, except for the case $l_1=l_2=1$.
Obviously, such a move works for tunnels instead of bridges
and we \em{will not distinguish between both}, unless clearly
stated.

\begin{figure}[htb]
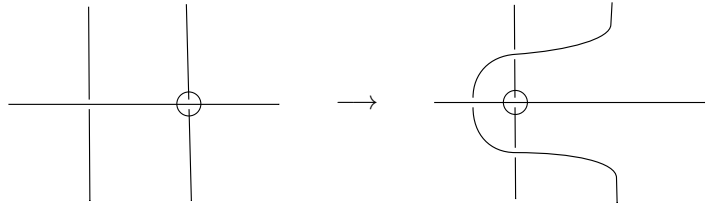

\[
\diag{1cm}{3.5}{2.5}{
  \picputtext{1.75 1.25}{ \epsfs{3cm}{t1-genmin1} }
  \piclinewidth{70}
  \piccircle{2.35 1.23}{0.16}{}
}
\qquad\lra\qquad
\diag{1cm}{3.5}{2.5}{
  \picputtext{1.75 1.25}{ \epsfs{3cm}{t1-genmin2} }
  \piclinewidth{70}
  \piccircle{1.03 1.25}{0.16}{}
}
\]
\caption{\label{fig1}A tongue move between two almost
alternating diagrams. The encircled crossing is to be switched
to obtain an alternating diagram.}
\end{figure}

Adams introduced a \em{tongue move} allowing to build more complicated
almost alternating diagrams from a given one. This move is
shown in figure \reference{fig1}. Herein the crossing needed
to be switched to obtain an alternating diagram is encircled; we
call this crossing \em{dealternator}.

Adams formulated a conjecture, stating how to recognize the unknot in
almost alternating diagrams, which was proved recently by Tsukamoto.

\begin{theorem}(Tsukamoto \cite{Ts})\label{TTs}
An almost alternating unknot diagram is trivializable by tongue moves,
flypes, and crossing number reducing Reidemeister I and II moves.
\end{theorem}

A tongue move is a special type of wave move, preserving the property
almost alternating. In that sense theorems \ref{th1.6} and \ref{t5p}
can be thought of as extending Tsukamoto's result in special diagrams.

For $\le 10$ crossings we use the numbering of prime knots of
\cite{Rolfsen}. For knots from 11 to 16 crossings our numbering
is that of KnotScape \cite{KnotScape}. Latter is reorganized so that
non-alternating knots are appended after alternating ones (of the
same crossing number), instead of using `a' and `n' superscripts.

We write $!D$ for the \em{mirror image} of $D$, and $!K$ denotes
the mirror image of $K$. Clearly $g(!D)=g(D)$ (and therefore $\tl
g(K)=\tl g(!K)$), and $g(!K)=g(K)$.

\subsection{Link polynomials\label{SLP}}

Let $X\in\bZ[t,t^{-1}]$. The \em{minimal} or \em{maximal degree} 
$\md X$ or $\Md X$ is the minimal resp.\ maximal exponent of $t$ 
with non-zero coefficient in $X$. Let $\spn_tX=\Md_tX-\md_tX$,
the \em{span} of \em{breadth} of $X$.
The coefficient in degree $d$ of $t$ in $X$ is denoted $[X]_{t^d}$
or $[X]_{d}$. The \em{leading coefficient} $\Mcf\,X$ of $X$ is its
coefficient in degree $\Md X$. If $X\in\bZ[x_1^{\pm 1},x_2^{\pm 1}]$,
then $\Md_{x_1}X$ denotes the maximal degree in $x_1$. Minimal
degree and coefficients are defined similarly, and of course
$[X]_{x_1^k}$ is regarded as a polynomial in $x_2^{\pm 1}$.

The \em{skein polynomial} $P$ \cite{HOMFLY,LickMil} is a Laurent
polynomial in two variables $l$ and $m$ of oriented knots and links
and can be defined by being $1$ on the unknot and the (skein) relation
\begin{equation}\label{1}
l^{-1}\,P\Bigl(\,
\diag{5mm}{1}{1}{
\picvecwidth{0.09}
\picmultivecline{-7.5 1 -1.0 0}{1 0}{0 1}
\picmultivecline{-7.5 1 -1.0 0}{0 0}{1 1}
}\,
\Bigr)\,+\,
l \,P\Bigl(\,
\diag{5mm}{1}{1}{
\picvecwidth{0.09}
\picmultivecline{-7.5 1 -1 0}{0 0}{1 1}
\picmultivecline{-7.5 1 -1 0}{1 0}{0 1}
}\,
\Bigr)\,=\,
-m\,P\Bigl(\,
\diag{5mm}{1}{1}{
\picvecwidth{0.09}
\piccirclearc{1.35 0.5}{0.7}{-230 -130}\picvecrlineto{0.01 0.011}
\piccirclearc{-0.35 0.5}{0.7}{310 50}\picvecrlineto{-.01 0.011}
}\,
\Bigr)\,.
\end{equation}
(The convention differs from \cite{LickMil} by the interchange of 
$l$ and $l^{-1}$.) We will denote in each triple as in \eqref{1} 
the diagrams (from left to right) by $D_+$, $D_-$ and $D_0$.
For a diagram $D$ of a link $L$, we will use all of the notations
$P(D)=P_D=P_D(l,m)=P(L)$ etc.\ for its skein polynomial, with
the self-suggestive meaning of indices and arguments.

The \em{Jones polynomial} \cite{Jones} $V$, and (one variable)
\em{Alexander polynomial} \cite{Alexander} $\Dl$ are obtained
from $P$ by the substitutions
\begin{eqnarray}
V(t) & = & P(-it,i(t^{-1/2}-t^{1/2}))\,, \label{PtoV}\\
\Delta(t) & = & P(i,i(t^{1/2}-t^{-1/2})) \label{PtoDl}\,,
\end{eqnarray}
hence these polynomials also satisfy corresponding skein relations.
(In algebraic topology, the Alexander polynomial is usually
defined only up to units in $\bZ[t,t^{-1}]$; the present
normalization is so that $\Dl(t)=\Dl(1/t)$ and $\Dl(1)=1$.)

We will use sometimes instead of $\Dl$ also the \em{Conway polynomial}
\cite{Conway} $\nb(z)$ with $\Dl(t)=\nb(t^{1/2}-t^{-1/2})$.
$\nb$ satisfies the skein relation $\nb(D_+)-\nb(D_-)=z\nb(D_0)$.
Note that $\Md\nb=2\Md\Dl$, which we use in particular to implicitly
restate some of the results of \cite{Cromwell} in the sequel.


The \em{Kauffman polynomial} \cite{Kauffman} $F$ is usually defined
via a regular isotopy invariant $\Lm(a,z)$ of unoriented links.
We use here a slightly different convention for the variables
in $F$, differing from \cite{Kauffman,Thistle} by the interchange
of $a$ and $a^{-1}$. Thus in particular we have for a link diagram $D$
the relation $F(D)(a,z)=a^{w(D)}\Lm(D)(a,z)$, where $\Lm(D)$ is the
writhe-unnormalized
version of the polynomial, given in our convention by the properties
\[
\begin{array}{c}
\Lm\bigl(\Pos{0.5cm}{}\bigr)\ +\ \Lm\bigl(\Neg{0.5cm}{}\bigr)\ =\ z\ \bigl(\ 
\Lm\bigl(\Nul{0.5cm}{}\bigr)\ +\ \Lm\bigl(\Inf{0.5cm}{}\bigr)\ \bigr)\,,\\[2mm]
\Lm\bigl(\ \ReidI{-}{-}\bigr) = a^{-1}\ \Lm\bigl(\noloop\bigr);\quad
\Lm\bigl(\ \ReidI{-}{ }\bigr) = a\ \Lm\bigl(\noloop\bigr)\,,\\[2mm]
\Lm\bigl(\,\mbox{\Large $\bigcirc$}\,\bigr) = 1\,.
\end{array}
\]

The Jones polynomial $V$ is obtained from $F$ (in
our convention) by the substitution (see \cite[\S III]{Kauffman2})
\[
V(t) \quad = \quad F(-t^{3/4},t^{1/4}+t^{-1/4})\,.
\]

\begin{figure}[htb]
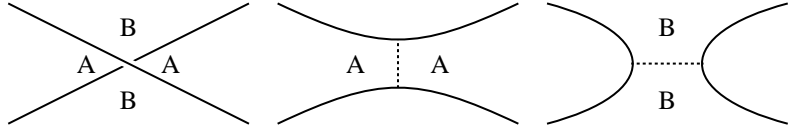

\[
\diag{8mm}{4}{2}{
   \picline{0 0}{4 2}
   \picmultiline{-5.0 1 -1.0 0}{0 2}{4 0}
   \picputtext{2.7 1}{A}
   \picputtext{1.3 1}{A}
   \picputtext{2 1.6}{B}
   \picputtext{2 0.4}{B}
} \quad
\diag{8mm}{4}{2}{
   \pictranslate{2 1}{
       \picmultigraphics[S]{2}{1 -1}{
           \piccurve{-2 1}{-0.3 0.2}{0.3 0.2}{2 1}
       }
       {\piclinedash{0.05 0.05}{0.01}
        \picline{0 -0.4}{0 0.4}
       }
   }
   \picputtext{2.7 1}{A}
   \picputtext{1.3 1}{A}
} \quad
\diag{8mm}{4}{2}{
   \pictranslate{2 1}{
       \picmultigraphics[S]{2}{-1 1}{
           \piccurve{2 -1}{0.1 -0.5}{0.1 0.5}{2 1}
       }
       {\piclinedash{0.05 0.05}{0.01}
        \picline{0 -0.6 x}{0 0.6 x}
       }
   }
   \picputtext{2 1.6}{B}
   \picputtext{2 0.4}{B}
}
\]
\caption{\label{figsplit}The A- and B-corners of a
crossing, and its both splittings. The corner A (resp. B)
is the one passed by the overcrossing strand when rotated 
counterclockwise (resp. clockwise) towards the undercrossing 
strand. A type A (resp.\ B) splitting is obtained by connecting 
the A (resp.\ B) corners of the crossing. It is useful to
put a ``trace'' of each splitted crossing as an arc connecting
the loops at the splitted spot.}
\end{figure}

An alternative description of $V$ is given by the Kauffman bracket in
\cite{Kauffman2}, which we recall next. The Kauffman bracket $[D]$
of a diagram $D$ is a Laurent polynomial in a variable
$A$, obtained by summing over all states the terms
\begin{eqn}\label{eq_12}
A^{\#A-\#B}\,\left(-A^2-A^{-2}\right)^{|S|-1}\,.
\end{eqn}
A \em{state} is a choice of splittings of type $A$ or 
$B$ for any single crossing (see figure \ref{figsplit}), 
$\#A$ and $\#B$ denote the number of
type A (resp. type B) splittings and $|S|$ the
number of (disjoint) circles obtained after all
splittings in a state.

The Jones polynomial of a link $L$ is related to the
Kauffman bracket of some diagram $D$ of $L$ by
\begin{eqn}\label{conv}
V_L(t)\,=\,\left(-t^{-3/4}\right)^{-w(D)}\,[D]
\raisebox{-0.6em}{$\Big |_{A=t^{-1/4}}$}\,.
\end{eqn}

Let the \em{$A$-state} of $D$ be the state where all
crossings are $A$-spliced; similarly define the $B$-state.
We call a diagram \em{$A$-(semi)adequate} if in the $A$-state
no crossing trace (one of the dotted lines in figure
\reference{figsplit}) connects a loop with itself.
We call such a trace a \em{self-trace}.
Similarly we define $B$-(semi)adequate. A diagram is
\em{semiadequate} if it is $A$- or $B$-semiadequate, and
\em{adequate} if it is simultaneously $A$- and $B$-semiadequate.
A link is adequate/semiadequate if it has an adequate/semiadequate
diagram.

\begin{theorem}(\cite{ntriv})\label{Vn1}
For a semiadequate knot (or link) diagram $D$, we have $V(D)\ne 1$
(or $V(D)\ne (-t^{1/2}-t^{-1/2})^{n(D)-1}$).
\end{theorem}

After the conversion \eqref{conv}, we see that
the minimal degree of $t$ to which the sum \eqref{eq_12}
contributes occurs in the $A$-state $A(D)$. For our
subsequent arguments, this degree can be more conveniently
written for a knot diagram $D$ as
\begin{eqn}\label{XA}
m(D):=g(D)-\frac{3}{2}c_-(D)+\frac{s(D)-|A(D)|}{2}\,.
\end{eqn}

If a diagram is $A$-adequate, only the $A$-state contributes in
degree $m(D)$, and the coefficient is $\pm 1$ (see \cite
{LickThis}). In Bae-Morton \cite{Bmo}, the contribution of the sum
of \eqref{eq_12} in degree $m(D)$ for general diagrams $D$ was
studied. The following easy property will be useful.

Call a self-trace \em{isolated}, if it does not pair up with
another self-trace like
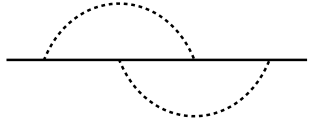
\[
\diag{1cm}{4}{2}{
  \picline{1 0 x}{1 4 x}
  \piclinedash{0.05 0.05}{0.01}
  \piccurve{0.5 1}{0.9 2}{2.1 2}{2.5 1}
  \piccurve{1.5 1}{1.9 0}{3.1 0}{3.5 1}
}
\]

\begin{lemma}\label{lBM}(Bae-Morton \cite{Bmo})
If the $A$-state of $D$ has an isolated self-trace, then the
contribution of the sum of \eqref{eq_12} in degree $m(D)$
vanishes, that is, $\md V(D)>m(D)$.
\end{lemma}

Let $F(K)\in \bZ[a^{\pm 1},z]$ be the Kauffman polynomial of a knot
$K$ and $Q_K(z)=F_K(1,z)$ be the BLMH polynomial \cite{BLM,Ho}. Let
further
\begin{eqn}\label{a+-}
a_-(K)\,=\,\max\,\{\,m-l\,:\,[F(K)]_{a^lz^m}\ne 0\,\}\quad\mbox{and}\quad
a_+(K)\,=\,\max\,\{\,m+l\,:\,[F(K)]_{a^lz^m}\ne 0\,\}\,.
\end{eqn}
It easily follows that $\Md_zF(K)\le a_+(K)+a_-(K)$. Moreover, if the r.h.s.
is positive, strict inequality holds, because the substitution $F(i,z)=1$
shows that $\Mc_zF(K)$ cannot be a single monomial (in $a$).

In \cite{Thistle}, Thistlethwaite proves that 
\begin{eqn}\label{capm}
c_{\pm}(K)\ge a_{\pm}(K)\,,
\end{eqn}
with equality if and only if $K$ is $A$ resp. $B$-semiadequate.
Let us say in the following that a (non-strict) inequality of
the form `$A\ge B$' is \em{sharp} or \em{exact} if $A=B$, 
and \em{strict} otherwise (i.e.\ if $A>B$).

Adequate (in particular alternating) diagrams make both
inequalities in \eqref{capm} simultaneously sharp. It follows 
from these inequalities that for every non-trivial knot $K$,
\begin{eqn}\label{QFc}
\Md Q(K)\,\le\,\Md_z F(K)\,<\,c_+(K)+c_-(K)\,.
\end{eqn}
(If $c_-(K)=c_+(K)=0$, then $K$ is simultaneously positive and negative,
and so is trivial; see for example \cite{restr}, or the below remarks on
the signature.) 

Note that for $P$ and $F$ there are several other
variable conventions, differing from each other by possible inversion
and/or multiplication of some variable by some fourth root of unity.

\subsection{The signature\label{Sz}}

The \em{signature} $\sg$ is a $\bZ$-valued invariant of knots and links.
Originally it was defined terms of Seifert matrices \cite{Rolfsen}.
We have that $\sg(L)$ has the opposite parity to the number of
components of a link $L$, whenever the \em{determinant}
$\dt(L)=\big|\,\Dl_L(-1)\,\big|\ne 0$. 
This in particular always happens for
$L$ being a knot (since $\Dl_L(-1)$ is always odd in this
case), so that $\sg$ takes only even values on knots.
Most of the early work on the signature was done by
Murasugi \cite{Murasugi5}, who showed several properties
of this invariant.

Then, for links $L_{\pm,0}$ with diagrams as in \eqref{1},
we have
\begin{eqnarray}
\sg(L_+)-\sg(L_-) & \in & \{0,1,2\} \label{2a} \\[2mm]
\sg(L_\pm)-\sg(L_0) & \in & \{-1,0,1\}\,. \label{2b}
\end{eqnarray}
(Note: In the first property one can also have $\{0,-1,-2\}$ istead
of $\{0,1,2\}$, since other authors, like Murasugi, take $\sg$ to
be with opposite sign. Thus \eqref{2a} not only defines a property,
but also specifies our sign convention for $\sg$.)

Further, Murasugi found the following important relation between
$\sg(K)$ and $\dt(K)$ for a knot $K$.
\begin{eqn}\label{3}
\begin{array}{*3c}
\sg(K)\equiv 0\,(4) & \iff & \dt(K)\equiv 1\,(4) \\[2mm]
\sg(K)\equiv 2\,(4) & \iff & \dt(K)\equiv 3\,(4)
\end{array}
\end{eqn}
If $\dt=1$, then even $8\mid \sg$, because of the property of
signatures of unimodular quadratic forms.

These conditions, together with the initial value $\sg(\bigcirc)=0$
for the unknot, and the additivity of $\sg$ under split union (denoted
by `$\sqcup$') and connected sum (denoted by `$\#$')
\[
\sg(L_1\# L_2)\,=\,\sg(L_1\sqcup L_2)\,=\,\sg(L_1)+\sg(L_2)\,,
\]
allow one to calculate $\sg$ for most links (incl. all
knots). The following further property is very useful:
$\sg(!L)=-\sg(L)$, where $!L$ is the mirror image of $L$.

If $K$ is a positive knot, then $\sg(K)>0$ \cite{CochranGompf}. 
Przytycki's concern was to improve and extend this result.
In particular, $\sg(K)\ge 4$ if $g(K)\ge 2$; also $\sg(K)>0$
if $K$ is almost positive or 2-almost positive, except a twist
knot. See \cite{apos,2apos,gen2} or \S\reference{Ss} below
for written (though not identical) proofs.

\subsection{Braid index and skein polynomial\label{SBR}}

The \em{braid group} $B_n$ on $n$ strands (or strings) is considered
to be generated by
the Artin \em{standard generators} $\sg_i$ for $i=1,\dots,n-1$. These
are subject to relations of the type $[\sg_i,\sg_j]=1$ for $|i-j|>1$,
which we call \em{commutativity relations} (the bracket denotes the
commutator) and $\sg_{i+1}\sg_i\sg_{i+1}= \sg_i\sg_{i+1}\sg_i$,
which we call \em{Yang-Baxter} (or shortly YB) \rm{relations}.

The \em{braid index} $b(L)$ of a link $L$ is the smallest
number of strands of a braid $\be$ whose closure $\hat\be$ is
$L$. See \cite{Morton,WilFr,Murasugi4}. (Alexander's theorem
asserts that $L=\hat\be$ for some braid $\be$.)

In \cite{Morton,WilFr} it was proved that
\begin{eqn}\label{MWF}
\mwf(L):=\,\frac 12\spn_lP(L)+1\,\le b(L)\,,
\end{eqn}
the \em{Morton-Williams-Franks} (MWF) \em{inequality}. Since
we will need the left-hand side of this inequality later,
let us write for it $\mwf(L)$ and call it the
\em{Morton-Williams-Franks bound} for $b(L)$. The inequality
\eqref{MWF} results from two other inequalities, due to Morton,
namely that for a diagram $D$, we have
\begin{eqn}\label{mq}
1-s(D)+w(D)\,\le\,\md_lP(D)\,\le\,\Md_lP(D)\,\le\,s(D)-1+w(D)\,.
\end{eqn}
Williams-Franks showed these inequalities for the case of
braid representations. Later it was observed from the
algorithm of Yamada \cite{Yamada} and Vogel \cite{Vogel}
that the braid version is actually equivalent to, and not just
a special case of, the diagram version. (These algorithms
allow to turn any diagram $D$ into a braid diagram without
altering $s(D)$ and $w(D)$.) Nonetheless we will
refer below to \eqref{mq} as `Morton's inequalities'. 

These inequalities were later improved in \cite{MP} in a way 
that allows to settle the braid index problem for many links 
(see \S\ref{S6} or also \cite{Ohyama}). For this purpose, 
Murasugi-Przytycki developed the concept of \em{index of a graph}.
We recall some main points of Murasugi-Przytycki's work, referring
to \cite{MP} for further details, and cautioning to the
correction explained in \S\ref{mhs}~--\ref{mht} below.

\begin{defi}
Let $G$ be a signed graph (each edge carries a sign $+$ or 
$-$). For a vertex $v$ in $G$ let $G_v=G/v$ be the graph obtained 
by contracting the \em{star} $\str v$ of $v$, that is, the set of
edges of which $v$ is one of the endpoints. Let $G^v$ be the graph
obtained from $G$ by deleting all these edges and $v$.

We say that $v$ is a \em{cut vertex} if $G^v$ is disconnected.
A cut vertex decomposes $G$ into a \em{block sum} or \em{join}
$G_1* G_2$ of two graphs $G_{1,2}$. If $G=G_1*\dots*G_n$ and
all $G_i$ have no cut vertex, the $G_i$ will be
called \em{block components} or \em{join factors} of $G$.
\end{defi}

\begin{defi}\label{d2.3}
We define (recursively) a sequence of edges $\mu=(e_1,\dots,e_n)$ 
to be \em{independent} in a graph $G$, if the following conditions 
are satisfied. 
\begin{enumerate}
\item The empty (edge) sequence is independent per definition. 
\item 
  Let $e_1$ connect vertices 
  $v_{1,2}$. Then we demand that $e_1$ is \em{simple}, i.e.\ there 
is no other edge connecting $v_{1,2}$, and that $e_2,\dots,e_n$ is 
independent in (one of) $G_{v_1}$ or $G_{v_2}$ (i.p. $e_2,\dots,e_n$
are disjoint from the set of edges incident to $v_1$ or $v_2$ resp.). 
\end{enumerate}
An \em{independent set} is a set of edges admitting an ordering
as an independent sequence.

The \em{index} $\inx(G)$, resp. \em{positive index} $\inx_+(G)$ and 
\em{negative index} $\inx_-(G)$ of $G$ are defined as the maximal
length of an independent edge set (or sequence), resp. independent
positive or negative edge set/sequence in $G$. A sequence is
\em{maximal independent} if it realizes the index of $G$.
\end{defi}

(This index is to be kept strictly apart from the index of a
Seifert circle of \S\reference{kkd}. Both are not unrelated,
but we do not investigate about their relation, and use them
in quite separate contexts.)

Now to each link diagram $D$ we associate its \em{Seifert graph}
$G=\Gm(D)$, which is a plane bipartite signed graph. It consists
of a vertex for each Seifert circle in $D$ and an edge for each
crossing, connecting two Seifert circles. Each edge is signed by the
skein sign of the crossing it represents. We will for convenience
sometimes identify crossings/Seifert circles of $D$ with edges/vertices
of $G$. Let also $\inx_{(\pm)\, }(D)=\inx_{(\pm)\, }(\Gm(D))$. 
 
\begin{prop}(see \cite[(8.4) and (8.8)]{MP})\label{p71}
If $D$ is a diagram of an oriented link $L$, then
\begin{eqnarray}
\label{q1}\Md_l P(L) & \le & w(D)+s(D)-1-2\inx_+(D) \\
\label{q2}\md_l P(L) & \ge & w(D)-s(D)+1+2\inx_-(D) \\
\label{q3}    b(L) & \le & \mpb(D):=s(D)-\inx(D)\,.
\end{eqnarray}
\end{prop}

For any diagram $D$, we have
\begin{eqn}\label{q4}
\inx_+(D)+\inx_-(D)\ge \inx(D)\,.
\end{eqn}
For alternating (and more generally homogeneous \cite{Cromwell})
diagrams $D$ equality holds, because each join factor of $\Gm(D)$
contains only edges of the same sign. This implies that if in such
diagrams \eqref{q1}, \eqref{q2} are sharp, then \eqref{MWF} and
\eqref{q3} become sharp, too.

It is clear that one can reconstruct $D$ from $\Gm(D)$
(when latter is given in a planar embedding). If $\Gm(D)$
has block components $G_i$, then we call the diagrams
$D_i$ with $\Gm(D_i)=G_i$ the \em{Murasugi atoms} of $D$.
An alternative way to specify $D_i$ is to say that they
are the prime factors of the Murasugi summands of $D$
(see definition \reference{Mz} below, and \cite{QW}).

The following theorem of Murasugi-Przytycki is important:

\begin{theorem}(theorem 2.5 in \cite{MP})\label{T25}
The index is additive under block sum of bipartite graphs,
i.e.\ $\inx(G_1* G_2)=\inx(G_1)+\inx(G_2)$ if $G_{1,2}$ are bipartite.
\end{theorem}

We will treat below, simultaneously to (the sharpness of) MWF, 
the following conjecture of theirs.

\begin{conj}\label{C7}(Murasugi-Przytycki)
If $D$ is an alternating diagram of a link $L$, then $b(L)=\mpb(D)$.
\end{conj}

\subsection{Genus generators\label{gg}}

Now let us recall, from \cite{gen1,gen2}, some basic facts concerning
knot generators of given genus. An explanation is given also in
section 5.3 in \cite{Cromwell2}. (There are several equivalent
forms of these definitions, and we choose here one that closely leans
on the terminology of Gau\ss{} diagrams; see for example \cite{pos}.)
We will also set up some notations and conventions used below.
The situation for links is discussed only briefly here, and in much
more detail in \S\reference{S33}.

\begin{defi}
We call two crossings $p$, $q$ in a knot diagram \em{linked}, and
write $p\cap q$, if passing their crossingpoints along the orientation
of $D$ we have the cyclic order $pqpq$, and not $ppqq$.
\end{defi}

\begin{defi}\label{dqa}
Let $D$ be a knot diagram, and $p$ and $q$ be crossings.
\def\labelenumi{(\roman{enumi})}
\begin{enumerate}
\item We call $p$ and $q$ \em{(twist) equivalent},
$q\simeq p$, if for all $r\ne p,q$ we have $r\cap p\iff r\cap q$.
\item We call $p$ and $q$ \em{$\sim$-equivalent}
and write $p\sim q$ if $p$ and $q$ are equivalent and $p\ncap q$.
\item Similarly $p$ and $q$ are called \em{$\ssim$-equivalent},
$p\ssim q$, if $p$ and $q$ are equivalent and $p\cap q$.
\item Finally, call two crossings $p$ and $q$ \em{Seifert equivalent},
if they connect the same two Seifert circles.
\end{enumerate}
\end{defi}

Let us record the following easy but useful observation.
 
\begin{lemma}\label{l45}
$\ssim$-equivalent crossings are Seifert equivalent. The converse
is true in special diagrams. \qed
\end{lemma}

\begin{defi}
A $\sim$-equivalence class consisting of one crossing is called
\em{trivial}, a class of more than one crossing \em{non-trivial}.
A $\sim$-equivalence class is \em{reduced} if it has at most
two crossings; otherwise it is \em{non-reduced}.
\end{defi}

\begin{defi}
Let $t(D)$ be the number of $\sim$-equivalence classes of $D$.
For $i=1,\dots,t(D)$ let $t_i(D)$ be the number of crossings
in the $i$-th $\sim$-equivalence class. Then for
$i=1,\dots,t(D)$ and $j=1,\dots,t_i(D)$ let $p(D,i,j)$
be the $j$-th crossing in the $i$-th $\sim$-equivalence class
of $D$.
\end{defi}


\begin{defi}
A \em{$\bt$ move} or \em{twist} at a crossing $x$ in a diagram $D$ is
a move, which creates a pair of $\sim$-equivalent crossings to $x$. 
(This is well-defined up to flypes.) We will call it positive
or negative, depending on whether it acts on a crossing of the
according sign.
\end{defi}

\begin{defi}\label{dbr}
An alternating diagram $D$ is called \em{$\bt$ 
irreducible} or \em{generating diagram}, if all $\sim$-equivalence
classes are reduced, that is, $t_i(D)\le 2$ for $i=1,\dots,
t(D)$. An alternating knot $K$ is called \em{generator} 
if some of its alternating diagrams is generating. The diagrams 
obtained from $D$ by $\bt$ moves and crossing changes form the 
\em{sequence} or \em{series} $\br{D}$ of $D$.
\end{defi}

(In many cases it will be useful to restrict oneself in $\br{D}$
to alternating or positive diagrams.)

Observe that, since $\sim$-equivalence is invariant under flypes,
theorem \ref{thfl} implies that some alternating diagram of $K$ is
generating if and only if all its alternating diagrams are so.

A flype in figure \reference{fig_} is called \em{trivial} if one
of the tangles $P,Q$ contains only crossings equivalent to
the crossing admitting the flype. A flype is of type A if the
strand orientation is so that strands on the left/right side 
of each tangle are directed equally w.r.t. the tangle (i.e.\ both
enter or both exit). Otherwise it is a flype of type B. Compare
\cite{SV}. So the property a crossing to admit a type A resp. type
$B$ flype is invariant of the $\ssim$ resp. $\sim$-equivalence class.

\begin{theorem}(\cite{gen1})
There exist only finitely many generators of given genus.
All diagrams of that genus can be obtained from diagrams
of these generators, under $\bt$ twists, flypes, and
crossing changes.
\end{theorem}

The finiteness of generators, together with the Flyping
theorem \cite{MenThis}, shows

\begin{theorem}\label{thang}(see \cite{gen1})\ \ 
Let $a_{n,g}$ be the number of prime alternating knots $K$ of
genus $g(K)=g$ and crossing number $c(K)=n$. Then for $g\ge 1$
\[
\sum_na_{n,g}x^n\,=\,\frac{R_g(x)}{(x^{p_g}-1)^{d_g}}\,,
\]
for some polynomial $R_g\in\bZ[x]$, and $p_g,d_g\in\bN$. Alternatively,
this statement can be written also in the following form:
there are numbers $p_g$ (period), $n_g$ (initial number of
exceptions) and polynomials $P_{g,1},\dots,P_{g,p_g}\in\bQ[n]$
with $a_{n,g}=P_{g,n\bmod p_g}(n)$ for $n\ge n_g$.
\end{theorem}

This was explained roughly in \cite{gen1}, and then in more
detail in \cite{SV}, where we made effort to characterize
the leading coefficient of these polynomials $P_{g,i}$. Even
if the polynomials vary with a very large period $p_g$, the
leading coefficients depend only on the parity of $n$.
The degrees of all $P_{g,i}$ are also the same, and equal to
$1$ less than the maximal number of $\sim$-equivalence classes
of diagrams of canonical genus $g$, with the exception
$g=1$, in which case this degree is $1$ or $2$ dependingly
on whether $i$ is even or odd.

In practice (in particular as we will see below) it is important
to obtain the list of generators for small genus. Genus one is
easy, and also observed independently.

\begin{theorem}(\cite{gen1}; see also \cite{Rudolph})\label{thge}
There are two generators of genus one, the trefoil and figure-8-knot.
\end{theorem}

Genus two and three require much more work. For suggestive
reasons, it is sufficient to find prime diagrams, and by theorem 
\reference{ThM}, prime generators.

\begin{theorem}(\cite{gen2})\label{thgen}
There are 24 prime generators of genus two,
$5_{1}$, $6_{2}$, $6_{3}$, $7_{5}$, $7_{6}$, $7_{7}$, $8_{12}$, 
$8_{14}$, $8_{15}$, $9_{23}$, $9_{25}$, $9_{38}$, $9_{39}$, $9_{41}$, 
$10_{58}$, $10_{97}$, $10_{101}$, $10_{120}$, $11_{123}$, $11_{148}$, 
$11_{329}$, $12_{1097}$, $12_{1202}$, and $13_{4233}$, and
4017 prime generators of genus 3.
\end{theorem}

A classification (by means of obtaining the list of
prime generators) for genus $g=4$ is also possible, and
will be explained below, as well as a general statement
in theorem \reference{th1gen}.

It follows from theorem \reference{thfl} that the series of
different (alternating) diagrams of the same generating knot
are equivalent \em{up to mutations}. When using tests involving
the polynomial invariants, which are invariant under mutations,
it is legitimate that \em{a priori} we fix 
\em{a single specific diagram $D$ for each generator $K$}, and work
only with (the series of) this diagram. Then we write $t(K)=t(D)$,
$t_i(K)=t_i(D)$ and $p(K,i,j)=p(D,i,j)$, and speak of the sequence
or series of $K$. With some simple observations we will be able
to comfort us with the assumption in \S\reference{S6},
while in \S\reference{S7} more care is needed.

Let us parametrize diagrams in the series of a generator $K$
as $K(x_1,\dots,x_l)$, where $l=t(K)$, the $l$ $\sim$-equivalence
classes of $K$ are ordered in some fixed way, and $x_i$ are
defined as follows.

For a trivial $\sim$-equivalence (i.e.\ $t_i=1$), $x_i\ge 1$ means
a $\sim$-equivalence class of $2x_i-1$ positive crossings, or
alternatively, the result of applying $(x_i-1)$ $\bt$-moves to
a single positive crossing in the generator. If $x_i\le 0$,
then we have a $\sim$-equivalence class of $1-2x_i$ negative
crossings, that is, a negative crossing with $(-x_i)$ $\bt$-moves
applied.

For a $\sim$-equivalence class of $t_i=2$ crossings, $x_i> 0$ 
means a $\sim$-equivalence class of $2x_i$ positive crossings, 
or alternatively, the result of applying ($x_i-1$) $\bt$-moves
to one of the two positive crossings in the generator.
$x_i=0$ means that the two crossings have opposite sign,
that is, form a trivial clasp after flypes. $x_i<0$
means $-2x_i$ negative crossings in the $\sim$-equivalence class,
that is, the result of ($-1-x_i$) $\bt$-moves on one of the
negative two crossings in the generator.

This convention will remain valid for the rest of the paper.
Note that it implies that we discard diagrams with crossings
of different sign within the same $\sim$-equivalence class
(unless $x_i=0$ and $t_i=2$).
Such diagrams have a trivial clasp after flypes, and are
of little interest. We will call a $\sim$-equivalence class
\em{positive} or \em{negative} depending on the sign of its crossings.

Turning to \em{link} diagrams, the $\bar t'_2$ twist (or
$\bar t'_2$ move) is given up to mirroring by
\begin{eqn}\label{t2}
\diag{7mm}{1}{1}{
    \picvecwidth{0.09}
    \picmultivecline{-7.5 1 -1.0 0}{1 0}{0 1}
    \picmultivecline{-7.5 1 -1.0 0}{0 0}{1 1}
}\quad\lra\quad
\diag{7mm}{3}{2}{
  \picPSgraphics{0 setlinecap}
  \picvecwidth{0.09}
  \pictranslate{1 1}{
    \picrotate{-90}{
      \rbraid{0 -0.5}{1 1}
      \rbraid{0 0.5}{1 1}
      \rbraid{0 1.5}{1 1}
      \pictranslate{-0.5 0}{
      \picvecline{0.021 1.95}{0 2}
      \picvecline{0.021 -.95}{0 -1}
    }
    }
    }
}\es.
\end{eqn}
We call diagrams that cannot be reduced by flypes and inverses of
the move \eqref{t2} \em{generating} or \em{$\bm$-irreducible}.
We know (in the case of knots, but we will extend this result
below to links) that reduced knot diagrams of given genus
$g$ (with $1-\chi=2g$) decompose into finitely many
equivalence classes  under $\bar t'_2$ twists and their
inverses. We call these collections of diagrams \em{series}.

Recall that in \cite{gen1} we called two crossings
\em{$\sim$-equivalent}, if after a sequence of flypes
they can be made to form a \em{reverse clasp}
$ \diag{7mm}{2}{1.2}{
  \picPSgraphics{0 setlinecap}
  \picvecwidth{0.09}
  \pictranslate{1 0.6}{
    \picrotate{-90}{
      \rbraid{0 -0.5}{1 1}
      \rbraid{0 0.5}{1 1}
      \pictranslate{-0.5 0}{
      \picvecline{0.021 -.95}{0 -1}
      \picvecline{0.979 .95}{1 1}
    } } } }\es; $
it is an exercise to check that this is an equivalence relation.
Similarly we call (see \cite{gen2}) two
crossings \em{$\ssim$-equivalent}, if after a sequence of flypes
they can be made to form a \em{parallel clasp}
$ \diag{7mm}{2}{1.2}{
  \picPSgraphics{0 setlinecap}
  \picvecwidth{0.09}
  \pictranslate{1 0.6}{
    \picrotate{-90}{
      \rbraid{0 -0.5}{1 1}
      \rbraid{0 0.5}{1 1}
      \pictranslate{-0.5 0}{
      \picvecline{0.021 .95}{0 1}
      \picvecline{0.979 .95}{1 1}
    } } } }\es. $
We observed (see \cite{gen1,gen2,STV}) that $\sim$-\ and $\ssim$-%
equivalent crossings of a knot diagram are linked with the same set
of other crossings in the diagram, so we have an equivalence to the
previous definition in the case of knots.


\subsection{Knots vs. links}

To conclude our setup, let us make more precise the separation
between knots and links in the following work.

The methods developed in this paper put into prospect to adapt
results about `$k$-almost positive knots' and `alternating knots
of genus $k$' to statements that apply to `$\tl k$-almost positive
$l$-component links' and `alternating $l$-component links of
genus $\tl k$' for $\tl k+l \le k+1$. We have largely waived on
treating links for a technical reason: the way we designed our
computation (see similarly remark \ref{rDn}).

The results where links are covered are those in \S\ref{S33} (here
the inclusion of links is essential and will be used elsewhere)
and \S\ref{OhS}. They rely on purely theoretical arguments, and
computation is not needed. To some extent, we managed to adapt
computations to the link case, as is shown in theorem \ref{lcL}.
However, for other parts more effort will be needed. We also expect
that some statements for links would look less pleasant than for
knots.

\section{The maximal number of generator crossings and
$\sim$-equivalence classes\label{S33}}

{
\def\@curvepath#1#2#3{%
  \@ifempty{#2}{\piccurveto{#1 polar}{@stc}{@std}#3}%
    {\piccurveto{#1 polar}{#2 polar}{#2 polar #3 polar 0.5 conv}
    \@curvepath{#3}}%
}
\def\curvepath#1#2#3{%
  \piccurve{#1 polar}{#2 polar}{#2 polar}{#2 polar #3 polar 0.5 conv}%
  \picPSgraphics{/@stc [ #1 polar #2 polar -1 conv ] $ D /@std [ #1 polar ] $ D }%
  \@curvepath{#3}%
}

\def\@opencurvepath#1#2#3{%
  \@ifempty{#3}{\piccurveto{#1 polar}{#1 polar}{#2 polar}}%
    {\piccurveto{#1 polar}{#2 polar}{#2 polar #3 polar 0.5 conv}\@opencurvepath{#3}}%
}
\def\opencurvepath#1#2#3{%
  \piccurve{#1 polar}{#2 polar}{#2 polar}{#2 polar #3 polar 0.5 conv}%
  \@opencurvepath{#3}%
}
\subsection{Generator crossing number inequalities} 

As explained in the introduction, we must first spend some
effort in estimates for the crossing number of generators
of given canonical Euler characteristic. These estimates
turn out to be rather sharp, and are a consequence of a
detailed study of the special diagram algorithm of Hirasawa
\cite{Hirasawa2} and myself \cite[\S 7]{bseq}.

\begin{theorem}\label{th1gen}
In a connected link diagram $D$ of
canonical Euler characteristic $\chi(D)\le 0$ there are at most 
\begin{eqn}\label{*+*}
t(D)\,\le\,
\left \{ \begin{array}{c@{\quad\mbox{if}\es}c} -3\chi(D) & \chi(D)<0 \\
1 & \chi(D)=0 \end{array} \right.
\end{eqn}
$\sim$-equivalence classes of crossings.
If $D$ is $\bm$-irreducible and has $n(D)$ link components, then
\begin{eqn}\label{sDi}
c(D)\,\le\,\left \{
\begin{array}{c@{\quad}l}
4 & \mbox{if\es $\chi(D)=-1$ and $n(D)=1$\,,} \\
2 & \mbox{if\es $\chi(D)=0$\,,} \\
-6\chi(D) & \mbox{if\es $\chi(D)<0$ and $n(D)=2-\chi(D)$\,,} \\
-5\chi(D)+n(D)-3 & \mbox{else}\,.
\end{array} \right.
\end{eqn}
\end{theorem}

Thus we settle the problem to determine the maximal crossing number of
a generator for knots.

\begin{corr}\label{cor4.1}
The maximal crossing number of a knot generator of genus
$g\ge 2$ is $10g-7$.
\end{corr}

\proof We know from \cite{SV} that for any $g\ge 2$, there are
examples of generator diagrams with $10g-7$ crossings. \qed

The inequality \eqref{*+*} is also optimal, and we call the
generating diagrams $D$ that make \eqref{*+*} exact \em{maximal
generating diagrams} and their knots/links \em{maximal generators}.
For knots these generators were studied in detail in \cite{SV}.

The main idea behind the description of diagrams
of given canonical Euler characteristic was to show that
they decompose into finitely many equivalence classes
under $\bar t'_2$ twists and their inverses. In \cite{gen1}
we showed this only for knots, but the case of links
can be easily recurred to it. For a diagram $D$ of a link
$L$ with $n(L)$ components, one applies $n-1=n(L)-1$ moves,
which replace a positive/negative crossing of two different
components by a parallel positive/negative clasp:
\begin{eqn}\label{y}
\diag{7mm}{1}{1}{
  \picvecwidth{0.09}
  \picmultivecline{0.18 1 -1.0 0}{0 1}{1 0}
  \picmultivecline{0.18 1 -1.0 0}{0 0}{1 1}
}\quad\lra\quad
\diag{7mm}{2}{2}{
  \picvecwidth{0.09}
  \picPSgraphics{0 setlinecap}
  \pictranslate{1 1}{
    \picrotate{-90}{
      \rbraid{0 -0.5}{1 1}
      \rbraid{0 0.5}{1 1}
      \pictranslate{-0.5 0}{
        \picvecline{0.021 .95}{0 1}
        \picvecline{0.979 .95}{1 1}
      }
    }
  }
}\es,\rx{1.5cm}
\diag{7mm}{1}{1}{
  \picvecwidth{0.09}
  \picmultivecline{0.18 1 -1.0 0}{0 0}{1 1}
  \picmultivecline{0.18 1 -1.0 0}{0 1}{1 0}
}\quad\lra\quad
\diag{7mm}{2}{2}{
  \picPSgraphics{0 setlinecap}
  \picvecwidth{0.09}
  \pictranslate{1 1}{
    \picrotate{-90}{
      \lbraid{0 -0.5}{1 1}
      \lbraid{0 0.5}{1 1}
      \pictranslate{-0.5 0}{
        \picvecline{0.021 .95}{0 1}
        \picvecline{0.979 .95}{1 1}
      }
    }
  }
}
\es.
\end{eqn}
We call this prodecure below a \em{clasping}. Thus $D$ can be
transformed into a knot diagram $D'$, with $1-\chi(D')\,=\,n-\chi(D)$. 

Although this simple argument establishes the picture qualitatively,
it is not useful for an optimal estimate, and neither was our
original approach in \cite{gen1} for knots. Then, in
subsequent work \cite{STV,SV} we established a relation
between canonical Seifert surfaces for knot diagrams
and 1-vertex triangulations of surfaces, and obtained
some partial information (certainly much better than
in \cite{gen1}) on the maximal number of crossings and
$\sim$-equivalence classes of $\bm$-irreducible knot
diagrams of given canonical genus. Unfortunately, this
method does not extent pleasantly to links. We will
thus introduce now a different approach, which is
entirely knot theoretic and circumvents this problem.
(Contrarily, there are insights of the old approach,
which will also be used later, but which cannot be recovered.)
This approach originates from an algorithm, first found by Hirasawa
\cite{Hirasawa}, to make any link diagram into a special one
without altering the canonical Euler characteristic
(in fact, even its isotopy type).

We will apply this algorithm to prove the above result theorem
\ref{th1gen}, which is likely in its optimal and most general form.
Before the proof of theorem \reference{th1gen}, first we
discuss Hirasawa's algorithm.

\subsection{An algorithm for special diagrams}

\begin{defi}\label{Mz}
Seifert circles in an arbitrary
diagram are called \em{non-separating}, if they have
empty interior or exterior; the others are called \em{separating}.
A diagram is called \em{special} if all its Seifert circles
are non-separating. Any link diagram decomposes along its
separating Seifert circles as the \em{Murasugi
sum} ($*$-product) of special diagrams (see \cite[\S 1]{Cromwell}).
\end{defi}

In \cite{BurZie} it was proved that each link has a special
diagram by a procedure how to turn any given diagram of the link into
a special one. However, the procedure in this proof alters drastically
the initial diagram and offers no reasonable control on the
complexity (canonical genus and crossing number)
of the resulting special diagram. A much more
economical procedure was given by Hirasawa \cite{Hirasawa}
and rediscovered a little later independently in \cite[\S 7]{bseq}.
Hirasawa's move consists of laying a part of a separating Seifert
circle $s$ along itself in opposite direction (we call this move
\em{rerouting}; it is the opposite of the wave move in
figure \reference{fwv}), while changing
the side of $s$ dependingly on whether interior or exterior adjacent
crossings to $s$ are passed. See figure \reference{_fig1}.

\begin{figure}[htb]
\[
\diag{1cm}{8}{8}{
  \picmultigraphics{3}{1 0}{}
  \pictranslate{4 4}{
    \picfillgraycol{0.8}
    \picputtext{3.15 120 polar}{$s$}
    \picputtext{2.85 155 polar}{$p$}
    \picputtext{2.5 175 polar}{$a$}
    \picputtext{3.3 175 polar}{$L$}
    \picputtext{2.8 193 polar}{$b$}
    \picputtext{2.6 220 polar}{$M$}
    \picrotate{20}{\lbraid{2.5 0}{1 0.6}}
    \picrotate{60}{\lbraid{3.5 0}{1 0.6}}
    \picrotate{80}{\lbraid{2.5 0}{1 0.6}}
    \picrotate{100}{\lbraid{3.5 0}{1 0.6}}
    \picrotate{170}{\lbraid{2.5 0}{1 0.6}}
    \picrotate{230}{\lbraid{3.5 0}{1 0.6}}
    \picrotate{210}{\lbraid{3.5 0}{1 0.6}}
    \picrotate{260}{\lbraid{2.5 0}{1 0.6}}
    \picrotate{280}{\lbraid{2.5 0}{1 0.6}}
    \picrotate{300}{\lbraid{3.5 0}{1 0.6}}
    \picrotate{320}{\lbraid{3.5 0}{1 0.6}}
    \picrotate{350}{\lbraid{3.5 0}{1 0.6}}
    \picfilledcircle{0 0}{2}{$T$}
    \piccircle{0 0}{3}{}
    {\piclinedash{0.1}{0.05}
    \picstroke{
      \opencurvepath{3 130}{3 135}{2.6 140}{2.3 98}{3.7 90}{3.7 75}%
      {2.3 70}{2.3 55}{3.7 40}{3.7 30}{3.7 20}%
      {3.7 10}{2.3 0}{2.3 350}{2.3 325}{2.3 300}%
      {3.7 290}{3.7 270}{3.7 250}{2.2 230}{2.2 220}{2.2 210}%
      {3.7 190}{3.7 180}{3.7 170}{3.7 160}{2.3 140}{2.8 140}{3 145}{3 150}
      {}
    }}
    {\piclinewidth{15}
     \picstroke{
       \curvepath{2.8 170}{2.83 167}{3.1 166}{3.1 164}{3.1 163}
       {3.1 160}{3.4 160}{3.5 167}{3.6 174}{3.6 180}
       {3.3 190}{3.1 195}{3.1 175}{3.0 174}{2.83 173}{}
     }
     \picstroke{
       \curvepath{3.3 210}{3.27 211}{2.9 214}{2.9 220}
       {2.9 226}{3.27 229}{3.3 230}{3.27 231}{2.9 235}{2.9 240}
       {2.7 240}{2.5 232}
       {2.3 225}{2.3 220}{2.3 215}{2.4 210}{2.6 205}{2.7 203}
       {2.9 203}{2.9 206}{3.27 209}{}
     }
    }
  }
}
\]
\caption{\label{_fig1}}
\end{figure}

The move of \cite{bseq} is similar, only that this type
of rerouting is applied to the Seifert circle connected
to $s$ by an crossing $c$ exterior to $s$. This move lowers
the canonical Euler characteristic by two, but by properly choosing
to reroute the strand above or below the rest of the diagram
(that is, such that it passes all newly created crossings
as over- or undercrossings), one obtains a trivial parallel clasp
involving $c$, whose deletion raises the canonical Euler
characteristic back by two. Then we obtain an instance of
Hirasawa's move.

Hereby, unlike in Hirasawa's original version of his algorithm, we
take the freedom to alter the signs of the new crossings, as far as
the isotopy type of the link, but \em{not} necessarily
of the canonical Seifert surface is preserved. It is of importance to
us only that the canonical Euler characteristic of the diagram is preserved. 
We assume that this freedom is given throughout the rest of this
section.

Hirasawa's algorithm is very economical~-- the number of
new crossings added is linearly bounded in the crossing number,
and even in the canonical Euler characteristic of the
diagram started with. (Note, for example, that the braid algorithms
of Yamada \cite{Yamada} and Vogel \cite{Vogel} have quadratic growth.)

Using Hirasawa's algorithm we prove now 

\begin{lemma}\label{pq}
Any connected link diagram $D$ of maximal number of 
$\sim$-equivalence classes of crossings for
given canonical Euler characteristic is special.
Moreover, this is also true for any $\bm$-irreducible
link diagram of maximal crossing number, except
if $\chi(D)=-1$ and $D$ is the 4-crossing (figure-8-knot)
diagram.
\end{lemma}

The following notions will be of particular importance in the proof:

\begin{defi}
A \em{region} of a link diagram is a connected component of the
complement of the (plane curve of) the diagram.
An \em{edge} of $D$ is the part of the plane curve of $D$
between two crossings (clearly each edge bounds two regions).
At each crossing $p$, exactly two of the four adjacent
regions contain a part of the Seifert circles near $p$.
We call these the \em{Seifert circle regions} of $p$. 
The other two regions are called the
\em{non-Seifert circle regions} of $p$. If the diagram is
special, each Seifert circle coincides with (the boundary of)
some region. We call the regions accordingly Seifert circle regions
or non-Seifert circle regions (without regard to a particular crossing).
\end{defi}

Note that two crossings are $\sim$-equivalent iff
they share the same pair of non-Seifert circle regions.

\begin{defi}
We call two crossings $a$ and $b$ in a diagram $D$ {\em neighbored}, if
they belong to a reversely oriented primitive Conway tangle in $D$, that
is, there are crossings $c_1,\dots,c_n$ with $a=c_1$ and $b=c_n$,
such that $c_i$ and $c_{i+1}$ form a reverse clasp in $D$. 

This is a similar definition to $\sim$-equivalence, but with no flypes
allowed. Thus the number of $\sim$-equivalence classes of a diagram is
not more than the number of neighbored equivalence classes of
the same diagram, or of any flyped version of it.
\end{defi}

\proof [of lemma \reference{pq}]
We consider first the second statement, that is, the one for
the crossing number of $\bm$-irreducible link diagrams.

To show this, take a non-special reduced diagram $D$
(not necessarily a generator). Apply flypes so that
all $\sim$-equivalent crossings are neighbored equivalent.
(This is possible, because flyping in a flyping circuit
is independent from the other ones. For the definition of
flyping circuits and related discussion see \cite[\S 3]{SunThis}.
Alternatively one can consult Lackenby's paper \cite{Lackenby},
where such diagrams are called `twist reduced'.)

We will show now that the application of Hirasawa moves
and appropriate flypes augments the crossing number of the generator,
in whose series the diagram lies, and preserves the condition
that $\sim$-equivalent crossings are neighbored equivalent.
The flypes we will apply always reduce the number of separating
Seifert circles, and so does any Hirasawa move. Thus we will
be done by induction on the number of separating Seifert circles.

Clearly, one can assume that $D$ is prime. The composite
case follows easily from the prime one.

Consider again the picture of the Hirasawa move,
figure \reference{fig4}. Here we avoid the creation of
the nugatory crossing $p$ in figure \reference{_fig1}.
Let $s$ be the separating
Seifert circle of $D$ on which the move is performed.

We distinguish three cases according to $\inx(s)$.

\begin{caselist}
\case If $\inx(s)=1$, $D$ is not prime. 

\case Now assume that $\inx(s)\ge 3$. Then
$s$ has $\inx(s)$ adjacent regions from inside and outside.
In figure \reference{fig4}, the inner regions are called
$A$, $A_1$, $A_2$ and $A_3$, and the outer regions are called
$B$, $B_1$, $B_2$ and $B_3$. The Hirasawa move splits up from
every such region $R$ a small part $R'$, containing a new
Seifert circle. (In figure \reference{fig4} four of the
8 such $R'$ are displayed.) Then it joins an inner and outer
region (here $A$ and $B$) to a new region we call $AB$. It adds
$2\inx(s)-1$ crossings. We call these crossings new, the others,
existing already before the move, old. Call the new diagram $D'$.

\begin{figure}[htb]
\[
\diag{1.6cm}{8}{8}{
  \picmultigraphics{3}{1 0}{}
  \pictranslate{4 4}{
    \picfillgraycol{0.8}
    \picputtext{3.15 120 polar}{$s$}
    \picputtext{2.5 175 polar}{$q_2$}
    \picputtext{3.5 105 polar}{$q_1$}
    \picputtext{2.8 90 polar}{$p_1$}
    \picputtext{3.5 130 polar}{$B$}
    \picputtext{3.8 80 polar}{$B_1$}
    \picputtext{3.9 30 polar}{$B_2$}
    \picputtext{3.3 24 polar}{$B_2'$}
    \picputtext{3.7 -80 polar}{$B_3$}
    \picputtext{2.7 120 polar}{$A'$}
    \picputtext{3.3 175 polar}{$B'$}
    \picputtext{2.8 193 polar}{$p_2$}
    \picputtext{2.6 220 polar}{$A_3'$}
    \picputtext{2.3 190 polar}{$A_3$}
    \picputtext{2.6 330 polar}{$A_2$}
    \picputtext{2.6 40 polar}{$A_1$}
    \picputtext{2.6 150 polar}{$A$}
    \picrotate{20}{\lbraid{2.5 0}{1 0.6}}
    \picrotate{60}{\lbraid{3.5 0}{1 0.6}}
    \picrotate{80}{\lbraid{2.5 0}{1 0.6}}
    \picrotate{100}{\lbraid{3.5 0}{1 0.6}}
    \picrotate{170}{\lbraid{2.5 0}{1 0.6}}
    \picrotate{230}{\lbraid{3.5 0}{1 0.6}}
    \picrotate{210}{\lbraid{3.5 0}{1 0.6}}
    \picrotate{260}{\lbraid{2.5 0}{1 0.6}}
    \picrotate{280}{\lbraid{2.5 0}{1 0.6}}
    \picrotate{300}{\lbraid{3.5 0}{1 0.6}}
    \picrotate{320}{\lbraid{3.5 0}{1 0.6}}
    \picrotate{350}{\lbraid{3.5 0}{1 0.6}}
    \picfilledcircle{0 0}{2}{$T$}
    {\piclinedash{0.1}{0.05}
     \piccircle{0 0}{3}{}
    }
    \picstroke{
      \opencurvepath{3 130}{3 135}{2.6 140}{2.3 98}{3.7 90}{3.7 75}%
      {2.3 70}{2.3 55}{3.7 40}{3.7 30}{3.7 20}%
      {3.7 10}{2.3 0}{2.3 350}{2.3 325}{2.3 300}%
      {3.7 290}{3.7 270}{3.7 250}{2.2 230}{2.2 220}{2.2 210}%
      {3.7 190}{3.7 180}{3.7 170}{3.7 160}{3.5 143}{3.2 140}{3 145}{3 150}
      {}
    }
    {\piclinewidth{15}
     \picstroke{
       \curvepath{2.8 170}{2.83 167}{3.1 166}{3.1 159}{3.1 153}
       {3.1 145}{3.5 145}{3.6 160}{3.6 174}{3.6 180}
       {3.3 190}{3.1 195}{3.1 175}{3.0 174}{2.83 173}{}
     }
     \picstroke{
       \curvepath{3.3 210}{3.27 211}{2.9 214}{2.9 220}
       {2.9 226}{3.27 229}{3.3 230}{3.27 231}{2.9 235}{2.9 240}
       {2.7 240}{2.5 232}
       {2.3 225}{2.3 220}{2.3 215}{2.4 210}{2.6 205}{2.7 203}
       {2.9 203}{2.9 206}{3.27 209}{}
     }
     \picstroke{
       \curvepath{2.9 127}{2.9 135}{2.7 137}{2.5 111}{2.7 96}{2.8 97}
       {3.3 98}{3.3 101}{3.0 101}{2.9 110}{2.9 117}{}
     }
     \picstroke{
       \curvepath{3.2 09}{3.4 11}{3.6 18}{3.6 24}{3.6 32}{3.4 38}{3.2 42}{3.1 45}
       {3.1 35}{3.1 22}{2.7 21}{2.7 19}{3.1 18}{3.1 14}{3.1 11}{}
     }
    }
  }
}
\]
\caption{\label{fig4}}
\end{figure}

We show that $D'$ is reduced. Assume that $p$ is a nugatory
crossing in $D'$. Then its two non-Seifert circle regions
coincide. If $p$ is an old crossing, the only way to make
them distinct by undoing the Hirasawa move is if they are
$A$ and $B$. However, the non-Seifert circle regions of any
crossing in $D$ lie either both inside or both outside of $s$.
If $p$ is new, then we see directly that the non-Seifert
circle regions lie one inside and one outside of $s$, and
hence are also distinct.

Now we examine how the $\sim$-equivalence relation has
been altered under the Hirasawa move.

Assume $p\sim q$ in $D'$. Then $p$ and $q$ have the same
pair of non-Seifert circle regions. (Since no crossing in $D'$
is nugatory, the two non-Seifert circle regions are always
distinct.)
{\nopagebreak
\def\labelenumi{\arabic{enumi})}\mbox{}\\[-18pt]
\def\theenumi{\arabic{enumi})}
\begin{enumerate}
\item\label{x_itema}
Assume first $p$ and $q$ are old crossings in $D'$ and $p\sim q$.
Then the non-Seifert circle regions at $p$ and $q$ have remained
the same when undoing the Hirasawa move (possibly the small
Seifert circle part are removed), except that the region $AB$
has been separated. Thus $p\sim q$ in $D$ also, except if
there is a region $C$ in $D$ such that the non-Seifert circle regions
at $p$ in $D$ are $A$ and $C$, and the non-Seifert circle regions
at $q$ are $B$ and $C$. However, the non-Seifert circle regions
of any crossing lie either both inside or both outside of $s$.
Since $A$ and $B$ lie on different sides of $s$, $C$ cannot exist.
Thus $p\sim q$ in $D$.
\item\label{x_itemb}
If $p$ and $q$ are both new, then the non-Seifert circle regions
can be explicitly given among the inner and outer adjacent
regions of $s$ (in our example the $A$s and $B$s). It is easy
to see that from the definition of the Hirasawa moves, all
pairs $(A_i,B_i)$ and $(A_i,B_{i+1})$ occurring this way are
distinct. (Here one needs that $D$ is prime; some of the $A_i$
or $B_i$ alone may coincide!) Thus $p\not\sim q$.

\item\label{x_itemc}
Finally consider the case $p$ new, $q$ old. The two
non-Seifert circle regions of $p$ lie one in- and one outside
of $s$ in $D$. By the argument with the region $C$ in case
\reference{x_itema}, we have $p\sim q$ only if $p$ and
$q$ are adjacent to the new region $AB$ in $D'$. From the
description of pairs of non-Seifert circle regions of the new crossings
and primeness of $D$, as mentioned in case \reference{x_itemb},
there are two such crossings $p$. They are the first and last
on the newly created strand. They are named $p_1$ and $p_2$
in figure \reference{fig4}. They are unique because of the
argument in case \reference{x_itemb}. 
\end{enumerate}
}

We proved that the only possible newly created pairs of
$\sim$-equivalent crossings are $(p_1,q_1)$ and $(p_2,q_2)$,
in which $q_1$ nd $q_2$ are taken from two fixed (distinct)
$\sim$-equivalence classes in $D$.
%

The addition of a crossing to a $\sim$-equivalence class
can reduce the crossing number of the underlying
generator by at most 2 (which occurs if the crossing
added is of odd number in its $\sim$-equivalence class).
Since $2\inx(s)-1\ge 5$ new crossings are added,
it follows that the crossing number of the
generator goes up under a Hirasawa move. 

\case Consider finally the slightly tricky case
$\inx(s)=2$. Now $2\inx(s)-5<0$ so that we cannot argue
as before. Let $n_i(s)$ and $n_o(s)$ be the number
of inner and outer crossings adjacent to $s$ resp.
Clearly $n_i(s),n_o(s)\ge \inx(s)=2$. If $n_i(s)=n_o(s)=2$, the two
internal and external crossings of $s$ are $\sim$-equivalent.
Then by assumption they are also neighbored equivalent,
and then we have the figure-8-knot diagram of genus
one, which we excluded. Thus assume w.l.o.g. (up to $S^2$-moves)
that $n_i(s)\ge 3$. Then we apply a \em{modified}
Hirasawa index-2-move:
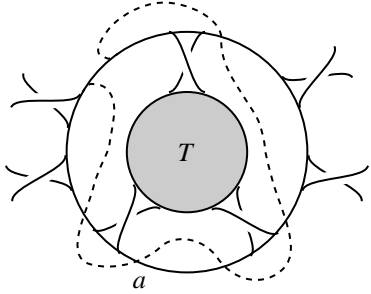
\begin{eqn}\label{mHm}
\diag{8mm}{6}{6}{
  \pictranslate{3 3}{
    \picfillgraycol{0.8}
    \picputtext{2.3 250 polar}{$a$}
    \picrotate{30}{\lbraid{2.5 0}{1 0.6}}
    \picrotate{90}{\lbraid{1.5 0}{1 0.6}}
    \picrotate{160}{\lbraid{2.5 0}{1 0.6}}
    \picrotate{190}{\lbraid{2.5 0}{1 0.6}}
    \picrotate{228}{\lbraid{1.5 0}{1 0.6}}
    \picrotate{310}{\lbraid{1.5 0}{1 0.6}}
    \picrotate{350}{\lbraid{2.5 0}{1 0.6}}
    \picfilledcircle{0 0}{1}{$T$}
    \piccircle{0 0}{2}{}
    {\piclinedash{0.1}{0.05}
     \picstroke{
       \opencurvepath{2 125}{2 130}{3.3 1 - 133}
            {2.4 125}{2.5 117}{3.6 1 - 108}{3.5 1 - 90}{3.3 1 - 75}
            {2.6 1 - 70}{2.4 1 - 55}{2.3 1 - 40}{2.3 1 - 30}{2.3 1 - 20}
            {2.3 1 - 10}{2.3 1 - 0}{2.3 1 - 359}{2.8 1 - 329}
	    {3.5 1 - 320}{3.6 1 - 300}{2.3 290}{2.4 1 - 280}
	    {2.4 1 - 260}{1.9 250}{3.3 1 - 240}{3.7 1 - 220}
	    {2.9 1 - 210}{2.5 1 - 190}{2.3 1 - 180}{2.3 1 - 170}
	    {2.3 1 - 160}{2.5 1 - 143}{2.7 1 - 140}{2 145}{2 150}
            {}
    }}
  }
}
\end{eqn}
This means, we create an additional trivial clasp in
one of the inner groups of at least 2 crossings.

It is easy to see that this move augments the
number of generator crossings (at least) by $1$.
\end{caselist}

The condition that $\sim$-equivalent crossings are neighbored
equivalent is not necessarily preserved under the
Hirasawa move. However, using flypes (which preserve
the number of crossings and $\sim$-equivalence classes),
one can again reestablish it. It is an easy observation
that the flypes needed reduce the number of separating
Seifert circles.

Then iterate Hirasawa moves (and possible flypes),
until $D$ has become special. Then we know, by this
argument, that the crossing number of the generator
of this special diagram is strictly greater
than this of $D$. This proves the second assertion
in the lemma.

The proof of the first part of the assertion goes
along similar lines, but is simpler. Since we showed
that the number of $\sim$-equivalence classes raises
by $2\inx(s)-3$ under each Hirasawa move, it is
sufficient to assume $\inx(s)\ge 2$. This makes
unnecessary the flyping argument, the need to care about
$\sim$-equivalent crossings being neighbored equivalent, or about
$\inx(s)=2$, and eliminates the exception of the 4-crossing diagram.
\qed

\begin{corr}
For given $n$ and $\chi$, except $(n,\chi)=(1,-1)$, 
a prime generator diagram $D$ of maximal crossing number is special,
and has maximal number of $\sim$-equivalence classes.
\end{corr}

\proof That $D$ is special follows from lemma
\reference{pq}. If $D$ is special and has not the maximal number
of $\sim$-equivalence classes, then by the work of \cite{SV}
we know that its unbisected Seifert graph $G'$ is either not
trivalent or not 3-connected. (The argument there was applied for
knots, but this condition was not used.) Latter case reduces
to former, since if $G'$ is not 3-connected, one can apply a flype
on $D$ to have a vertex of valence $\ge 4$ in $G'$. In \cite{SV}
we argued (under the exclusion of cut vertices), that one can apply a
decontraction on such a vertex, such that the created clasp is a new
$\sim$-equivalence class (and among the previous $\sim$-equivalence
classes no identifications occur). Then the crossing number of the
generator is also augmented. \qed

\subsection{Proof of the inequalities\label{413}}

The lemma does the main part of the proof of theorem
\reference{th1gen}, which will now be completed.
(Some of the arguments that follow appear also in very similar form
in \cite{SV}, so that the reader may consult there
for more details.)

\proof[of theorem \reference{th1gen}]
The lemma shows that for obtaining the stated estimates,
only special diagrams need to be considered.
(Since the 4-crossing diagram does not violate any
of the two assertions of the theorem, let us exclude
it in any further consideration.)

Assume that $D$ has a Seifert circle of valence $\ge 4$.
If we have a $\bm$-irreducible special diagram with a Seifert
circle of valence $\ge 4$, then one can always perform a Reidemeister
II move on any pair of non-neighbored edges in this Seifert circle
region. Then one obtains a special diagram of the same Euler
characteristic and two crossings more. By proper choice of
the pair of edges, the new crossings can be made to form a separate
new $\sim$-equivalence class, so that the diagram is still
$\bm$-irreducible (although not always any pair of such
edges will do).

This means that we need to consider only special diagrams
$D$ with Seifert circles of valence $2$ and $3$.
In other words, the Seifert graph (see \cite[\S 1]{Cromwell})
$G$ of $D$ is 2-3-valent (all its vertices have valence $2$ or $3$).
Removing vertices of valence $2$ in $G$ by unbisections means
identifying edges which correspond to $\sim$-equivalent
crossings in $D$. Let $G'$ be the 3-valent (cubic)
graph obtained this way. Since $D$ is reduced, $G'$ is
connected and contains no loop edges. Note that $G$ is a reduced
bisection of $G'$.

However, $G'$ may be a single loop, which occurs when $\chi(D)=0$.
Then one obtains the $(2,k)$-torus link diagrams with reverse
orientation (which are generated by the Hopf link and have one
$\sim$-equivalence class). Thus we can exclude the degenerate case
$\chi(D)=0$ in the rest of the proof.

Then, if $e(G')$ is the number
of edges of $G'$ and $v(G')$ the number of its vertices,
we have $v(G')=\myfrac{2}{3}e(G')$, and 
\[
-\frac{e(G')}{3}\,=\,v(G')-e(G')\,=:\,
\chi(G')\,=\,\chi(G)\,=\,\chi(D)\,.
\]
Since the number of $\sim$-equivalence classes in $D$ is at most
$e(G')$, we obtain the first result. 

%
%
%

Now consider the second statement in the theorem.
{}From Hirasawa's algorithm we know that the maximal
crossing number of a generator is achieved by a special
diagram (except for the case of knot diagrams of genus one),
and that among such diagrams, in a diagram whose Seifert graph
is 2-3-valent. Such a graph has $2-\chi(D)$ faces
(corresponding to regions of $D$ which do not contain a
Seifert circle). For a while forget about the orientation
of any component. If every $\sim$-equivalence class
of such a diagram $D$ had 2 crossings, the number
of components of $D$ would be equal to its number
$2-\chi(D)$ of non-Seifert circle regions. Since
the change of 2 crossings to 1 in each $\sim$-equivalence class
\begin{eqn}\label{yrep}
\diag{7mm}{2}{2}{
  \picPSgraphics{0 setlinecap}
  \pictranslate{1 1}{
    \picrotate{-90}{
      \rbraid{0 -0.5}{1 1}
      \rbraid{0 0.5}{1 1}
      \pictranslate{-0.5 0}{
    }
    }
    }
}
\quad\lra\quad
\diag{7mm}{1}{1}{
    \picmultiline{0.18 1 -1.0 0}{0 1}{1 0}
    \picmultiline{0.18 1 -1.0 0}{0 0}{1 1}
}
\end{eqn}
changes the number of components at most by $\pm 1$, we need at
least $2-\chi(D)-n(D)$ such replacements to obtain a diagram
of $n(D)$ components. Thus the maximal number of crossings
we can have is
\[
-6\chi(D)\,-\,\bigl(2-\chi(D)-n(D)\bigr)\,=\,-5\chi(D)+n(D)-2\,,
\]
and we showed
\begin{eqn}\label{eqn4}
c(D)\,\le\,\left \{
\begin{array}{c@{\quad\mbox{if}\es}c} -5\chi(D)+n(D)-2 & \chi(D)<0 \\ 
2 & \chi(D)=0 \end{array}\right.\,.
\end{eqn}
To obtain \eqref{sDi}, we must show that we can improve the estimate in
\eqref{eqn4} by one crossing if $\chi<0$ and
$n<2-\chi$ (we excluded knot diagrams of genus one).

Assume there is a diagram $D$ with $-5\chi+n-2$ crossings.
Let $G$ be its 2-3-valent Seifert graph, and $G'$ be the
3-valent graph obtained from $G$ by deleting valence-2 vertices. Call
the edges of $G'$ obtained under such un-bisections even, and the
others odd. $G'$ inherits a particular planar embedding from $G$,
and hence we can build its dual graph $G'^*$. Note that each edge
$e'$ of $G'$ corresponds bijectively to an edge $e'^*$ in $G'^*$.
Define for a subgraph $\Gm\subset G'$ the ``dual'' graph
\[
\Gm^*\,:=\,\{\,e'^*\,:e'\,\in \Gm\,\}\,.
\]
Let
\[
\Gm\,:=\,\{\,e'\,\in\,G'\,:\,e'\,\ \mbox{odd}\,\}\,.
\]

We claim that $\Gm^*$ is loopless. Assume that $\Gm^*$ contains
a loop made up of edges $l_1,\dots,l_k$. Take the diagram $D_0$
with the same graph $G'$, but such that all edges are even.
Then $D$ is obtained from $D_0$ by $2-n-\chi$ moves \eqref{yrep},
corresponding to the odd edges of $G'$. (The order of these moves is
irrelevant.) To have $n$ components in $D$ (and mot more), each
such move must reduce the number of components by one, so that
both strands on the left of \eqref{yrep} must belong to different
components. However, if one performs on $D_0$ the moves \eqref{yrep}
corresponding to $l_1,\dots,l_{k-1}$, and then one
applies the move for $l_k$, it is easy to see that both
strands on the left of \eqref{yrep} for this last move 
belong to the same component. This contradiction shows that
$\Gm^*$ is loopless.

Since $\Gm^*$ is loopless, $\Gm^*\subset \Gm'^*$ for some spanning tree
$\Gm'^*$ of $G'^*$. Then $G'\sm \Gm'$ is also a spanning tree (of $G'$).
Thus $G'$ has a spanning tree made up of even edges.

To show that this is impossible, consider the orientation of the
Seifert circles corresponding to the vertices of $G$. Sign the vertices
of $G$ positive or negative depending on this orientation. (This makes
$G$ bipartite.) Then this signing reduces also to a signing
of $G'$, with the property that even edges connect vertices
of the same sign, and odd edges connect vertices of opposite
sign. Since $n<2-\chi$, there must be vertices of $G'$ of both
signs. But $G'$ has a spanning tree made up of even edges,
and so all its vertices must have the same sign, a contradiction.
(Clearly adding just one edge to $\Gm$ may solve the problem,
because $\Gm^*$ may no longer be loopless, and $G'\sm \Gm$ may get
disconnected.)

This completes the proof of theorem \reference{th1gen}. \qed

\begin{rem}\label{rem23}
Note that we always have $2-\chi(D)\ge n(D)$, so that we can always
eliminate $n$ for $\chi$.
\end{rem}

\begin{rem}\label{remimp}
In practice, it will be often 
convenient to use \eqref{eqn4} rather than \eqref{sDi},
because we can absorb by this additional crossing two of the
exceptional cases (and avoid unpleasant case distinctions).
Nonetheless, the additional argument for \eqref{sDi} is useful
to settle completely the question on the maximal crossing
numbers of generators in many cases, in particular for knots.
\end{rem}

\begin{rem}
It is clear from the proof that there is no natural way to obtain an
upper restriction on the number of $\sim$-equivalence classes of
odd number of crossings. In fact, no such restriction exists: in
the case of knots and genus $g\ge 6$ we constructed in \cite{gen2}
generator diagrams with the maximal number of $-3\chi=6g-3$
$\sim$-equivalence classes, each one having a single crossing.
\end{rem}

\begin{rem}
If $G'$ is $3$-connected,
then it is easy to see that in fact different edges of $G'$
correspond to different $\sim$-equivalence classes in $D$,
so that one easily finds examples where the first bound $-3\chi$
in theorem \reference{th1gen} is realized sharply. We will study
such generators $D$ in more detail in a subsequent paper \cite{3ang}.
\end{rem}

\subsection{Applications and improvements}

We start with some first applications of the
crossing number estimates for generators. For this
we use the degree-2-Vassiliev invariant $v_2=\myfrac 12\Dl''(1)$.
Beside the independent interest of these inequalities, their proofs
also introduce an idea that, in modified form, will be important
for the later more general results.

\begin{theorem}\label{ppq}
If $K$ is a positive knot, then each reduced positive diagram $D$
of $K$ has $c(D)\le 9g(K)-8+2v_2(K)$ crossings, except if $K$ is
the trefoil (where $c(D)\le 4$).
\end{theorem}

\proof Let $D'$ be the positive generator in whose series $D$ lies.
We have $g(D')=g(D)=g(K)$, and we know that a positive $\bm$-twist
in a positive diagram augments $v_2$. If $g(D')>1$, then $c(D')\le
10g(D')-7$ by corollary \reference{cor4.1}. On the opposite side,
we have $c(D')\ge 2g(D')+1$, so that by theorem 6.1 of \cite{pos}
\[
v_2(D')\, \ge \, \frac{c(D')}{4}\,\ge\, \frac{g(D')+1}{2}\,
\]
(where in the second inequality we used the obvious integrality of
$v_2$). Now we can apply at most $v(D)-v_2(D')$ $\bm$-twists to $D'$ to
obtain $D$, so that
\[
c(D)\,\le\,c(D')+2\bigl(v_2(K)-v_2(D')\bigr)\,\le
\,10g(D')-7+2v_2(K)-g(D')-1\,=\,9g(K)-8+2v_2(K)\,.
\]
If $g(D)=g(D')=1$, then we must show $c(D)\le 2v_2(K)+1$. This can be
checked directly, since by theorems \reference{thge} and
\ref{thgen}, $D$ is either a $(p,q,r)$-pretzel
diagram for $p,q,r\ge 1$ odd, or a positive rational knot diagram
$C(p,-q)$, with $p,q\ge 2$ even. In former case $c(D)=p+q+r$ and
$v_2(D)=(pq+pr+qr+1)/4$, and the inequality easily follows. In
latter case
$c(D)=p+q$ and $v_2(D)=pq/4$, so the inequality holds unless $p=q=2$,
which is the positive 4-crossing trefoil diagram. \qed

\begin{theorem}\label{ppp}
If $K$ is an almost positive knot, then each reduced almost
positive diagram $D$ of $K$ has $c(D)\le 10g(K)+3+2v_2(K)$
crossings.
\end{theorem}

\proof $D$ has genus at most $g(K)+1$ by theorem \reference{Biq}.
Consider the almost positive generator $D'$, from which $D$ is
obtained by positive $\bm$-twists. Since almost positive genus 1
diagrams belong to positive knots, $g(D')>1$, so that by corollary
\reference{cor4.1},
$D'$ has at most $10g(D')-7\le 10g(K)+3$ crossings. Moreover, since
$g(D')>1$, $D'$ is not an unknotted twist knot diagram. Thus
$v_2(D')>0$ by theorem 4.1 of \cite{apos}.

Now we apply in a proper
order the $\bm$-twists taking $D'$ to $D$, keeping track of $v_2$.
We know from the Polyak-Viro formula for $v_2$ that a $\bm$-move
augments $v_2$ except if it is at a crossing $p$, linked with only
two other crossings, $q$ and $r$, one of which, say $q$, is the
negative one. (In this situation $v_2$ is preserved.) Call such a
crossing $p$ \em{thin}. In the proof of lemma 5.3 of \cite{apos} we
observed that then $q$ and $r$ form a trivial clasp in $D'$.
To avoid the elimination of $q$ in $D$, this clasp must be parallel,
and a twist at $r$ must be applied in the sequence of $\bm$-twists
taking $D'$ to $D$. Now we count how many different crossings occur as
$r$ in the above situation (for some $p$). Since $q\ssim r$, and
$\ssim$ is an equivalence relation, any two possible $r$ and $r'$ are
$\ssim$-equivalent. Then they intersect the same set of other
chords, in particular $p$. But $p$ intersects at most one other
chord different from $q$. So no two different $r$ exist, although
for the unique $r$ there may be different $p$. However, applying
first the $\bm$-twist at $r$, we eliminate all these crossings $p$
as thin crossings (and there are no other ones). So we ensure that
all following twists augment $v_2$. Since $v_2(D')>0$, at most
$v_2(K)=v_2(D)$ twists can be performed (the initial one included).
Thus
\begin{myeqn}{\qed}
c(D)\,\le\, c(D')+2v_2(K)\,\le\,10g(K)+3+2v_2(K)\,.
\end{myeqn}

Even albeit the inequalities in theorem \reference{th1gen}
are optimal in general, they can be improved under
additional conditions. One such condition is related to
the \em{signature} $\sg(L)$ of the underlying link $L$,
in comparison to the maximal \em{Euler characteristic} $\chi(L)$
of an orientable spanning surface of $L$. It is well-known
\cite{Murasugi3} that $|\sg(L)|\le 1-\chi(L)$, at least when
$L$ does not bound disconnected surfaces of small genus, so
for example for knots or (non-split) positive links. We can
say something about the situation when $\sg(L)\ll 1-\chi(L)$.

\begin{lemma}\label{lem7}
Assume we have a diagram $D$ of a link $L$ with $\chi(D)<0$ and
$a(D)$ negative crossings, and
\[
k(D)\,:=\,\frac{1-\chi(D)-2a(D)-\sg(D)}{2}\,\ge\,0\,.
\]
Then $D$ has at most $-3\chi(D)-\ds\frac{3}{2}k(D)$
$\sim$-equivalence classes, and its underlying generating
diagram has at most
\[
-5\chi(D)+n(D)-2-\frac{1}{2}k(D)
\]
crossings.
\end{lemma}

\proof If $D$ is a 4-crossing knot diagram, then $k(D)=0$,
so this case follows easily. Otherwise, let $\tl D$
be the special diagram obtained from $D$ by Hirasawa's
algorithm. Hereby, when applying a move at a separating
Seifert circle $s$, the new strand is placed
so as to create $\inx(s)-1$ (and not $\inx(s)$) new negative
(non-nugatory) crossings.

We assume first that $\inx(s)\ge 3$.
We have shown in the proof of theorem \reference{th1gen}
that the number of $\sim$-equivalence classes is augmented by 
\begin{eqn}\label{eqn1}
2\inx(s)-3\,\ge\,\frac{3}{2}\bigl(\inx(s)-1\bigr)\,,
\end{eqn}
and that of the crossing number of the underlying generating
diagram at least by
\begin{eqn}\label{eqn2}
2\inx(s)-5\,\ge\,\frac{1}{2}\bigl(\inx(s)-1\bigr)\,.
\end{eqn}
As for the positive special alternating diagram $\r D$
obtained from $\tl D$ by crossing changes, we have
from \cite{Murasugi3} $\sg(\r D)=1-\chi(\r D)=1-\chi(\tl D)$,
and $\sg$ changes at most by $2$ under a crossing change, we
have that $\tl D$ has at least
\begin{eqn}\label{7.5}
\frac{1-\chi(D)-\sg(D)}{2}
\end{eqn}
negative crossings. Thus $\tl D$ has at least $k(D)$
negative crossings added by the Hirasawa moves (and not
inherited from $D$). Then apply \eqref{eqn1},
\eqref{eqn2} and theorem \reference{th1gen}. 

This argument shows the assertion if no Hirasawa moves
at index-2 Seifert circles are applied. To deal with
$\inx(s)=2$ we consider again the modified Hirasawa move
\eqref{mHm}. Then we refine our argument as follows.
Mark in each modified Hirasawa move the negative crossing
in the additional clasp (this is crossing $a$ in \eqref{mHm}).
Then we claim that $\tl D$ has \eqref{7.5} many (not
only negative but) non-marked negative crossings. To see
this, it suffices to show that the switch of all
marked (positive) crossings in $\r D$ preserves $\sg$.
Now, switching all marked crossings in $\r D$ and
resolving the trivial clasps gives a diagram $\hat D$.
Clearly $\hat D$ is special and positive. But it has also the same
$\chi$ as $\r D$, since modulo crossing changes it is obtained
from $D$ by the ordinary (non-modified for $\inx(s)=2$)
Hirasawa moves. Thus $\sg(\r D)=\sg(\hat D)$ follows
from Murasugi's result $\sg=1-\chi$.

Therefore, \eqref{mHm} creates (with proper choice between putting
the arc above or below) only one non-marked negative crossing,
but at least two new $\sim$-equivalence classes, and augments the
generator crossing number at least by one. Then \eqref{eqn1} resp.
\eqref{eqn2} remain valid with the l.h.s.\ being the augmentation
of $\sim$-equivalence classes and generator crossing number resp.,
and the parenthetical term on the r.h.s.\ the number of
(non-marked) negative crossings contributing to \eqref{7.5}. \qed

\begin{rem}\label{ret}
Observe that, for knots, we can use instead of $\sg$ the
Ozsv\'ath-Szab\'o signature $\tau$, since it coincides (up
to a factor) with $\sg$ on (special) alternating knots,
and enjoys (up to that factor) the property \eqref{2a}.
However, Livingston in fact deduced the Rudolph-Bennequin
inequality \eqref{RBi} (in the original form of the estimate,
without the term $s_-(D)$ on the right) from $\tau$, and so
lemma \reference{lem7} becomes trivial for $\tau$. A similar
comment applies on the signature $s$ derived from Khovanov
homology (see Bar-Natan \cite{BarNatanKh} for example).
\end{rem}

Note that the signature improvement of theorem \reference{th1gen}
can be applied to theorems \reference{ppq} and \reference{ppp}.
However, we see in \cite{adeq} that, with further tools introduced,
the variety of possible modifications and improvements of such
inequalities grows considerably, so that we cannot discuss each
one in detail. We will thus in particular not elaborate much on the
use of $\tau$ and $s$, and leave the adaptation to an interested
reader.

} 
 
\section{Generators of genus 4\label{gen4}}

The methods developed in \S\reference{S33} allow us to push the compilation
of generators for genus 2 and 3 in \cite{gen2} one step further, and
to settle genus 4. Beside its algorithmical merit, this project
had the purpose of verifying (and in fact correcting minor mistakes
in) the theoretical thoughts, on which our approach bases. Several
applications presented below in this paper rely on this compilation.

Table \reference{tabg} gives the distribution of the prime
(knot) generators $K$ of genus $4$ by crossing number and number
of $\sim$-equivalence classes. For maximal generators, we know from
our previous work (theorem \reference{th1gen}) that the crossing
numbers in question lie between $21$ and $33$, but excluded $21$
in \cite{gen2} by a computational argument using the property
that the Seifert graph must be ($3$-valent and) bipartite.
For 22 crossings, still there are no maximal generators.

\begin{table}[ptb]
\captionwidth0.8\vsize\relax
\newpage
\vbox to \textheight{\vfil
\tabcolsep5.8pt
\rottab{%
\def\hh{\\[0.8mm]\hline[1.5mm]}%
\hbox to \textheight{\tiny\hss
\begin{mytab}{|c||r|r|r|r|r|r|r|r|r|r|r|r|r|r|r|r|r|r|r|r|r|r|r|r|r||r|}%
  { & \multicolumn{25}{|r||}{ } & }%
  \hline [1mm]%
  \rx{-3.6em}%
  \diag{4mm}{1.5}{1}{
    \piclinewidth{35}
    \picline{-0.15 1.05}{1.66 -0.25}
    \picputtext{1.2 0.7}{$c$}
    \picputtext{0.4 0.2}{$\#\sim$}
  }\rx{-3.6em}%
& 9 & 10 & 11 & 12 & 13 & 14 & 15 & 16 & 17 & 18 & 19 & 20 & 21 & 22 & 23 & 24 & 25 & 26 & 27 & 28 & 29 & 30 & 31 & 32 & 33 & \multicolumn{1}{|c|}{\rx{-0.1em}total\rx{-0.1em}} \\[2mm]%
\hline
\hline [2mm]%
8 &  &  &  &  &  &  &  & 29 &  &  &  &  &  &  &  &  &  &  &  &  &  &  &  &  &  & 29 \hh
9 & 1 & 2 & 10 & 28 & 71 & 104 & 147 &  &  & 145 &  &  &  &  &  &  &  &  &  &  &  &  &  &  &  & 508 \hh
10 &  & 21 & 72 & 210 & 356 & 557 & 660 & 819 & 1092 &  &  & 369 &  &  &  &  &  &  &  &  &  &  &  &  &  & 4156 \hh
11 &  &  & 48 & 257 & 766 & 1791 & 2942 & 3832 & 3080 & 2804 & 3188 &  &  & 447 &  &  &  &  &  &  &  &  &  &  &  & 19155 \hh
12 &  &  &  & 55 & 487 & 2033 & 4734 & 8585 & 12145 & 13523 & 8500 & 5168 & 4707 &  &  & 313 &  &  &  &  &  &  &  &  &  & 60250 \hh
13 &  &  &  &  & 56 & 548 & 3087 & 9661 & 19112 & 27552 & 31293 & 27717 & 14629 & 5427 & 3876 &  &  & 111 &  &  &  &  &  &  &  & 143069 \hh
14 &  &  &  &  &  & 46 & 590 & 3519 & 13251 & 32388 & 52870 & 61747 & 53398 & 35540 & 15787 & 3173 & 1827 &  &  & 20 &  &  &  &  &  & 274156 \hh
15 &  &  &  &  &  &  & 41 & 489 & 3584 & 14749 & 41049 & 78373 & 102880 & 95709 & 61646 & 28311 & 10626 & 965 & 465 &  &  & 1 &  &  &  & 438888 \hh
16 &  &  &  &  &  &  &  & 27 & 356 & 2814 & 13781 & 42566 & 90877 & 135278 & 138221 & 100392 & 48096 & 13094 & 4195 & 115 & 49 &  &  &  &  & 589861 \hh
17 &  &  &  &  &  &  &  &  & 14 & 231 & 1854 & 9704 & 34955 & 83859 & 141210 & 163710 & 125842 & 68515 & 24978 & 2942 & 837 &  &  &  &  & 658651 \hh
18 &  &  &  &  &  &  &  &  &  & 4 & 96 & 989 & 5258 & 20307 & 56939 & 110240 & 150023 & 136642 & 75688 & 27646 & 8127 & 172 & 46 &  &  & 592177 \hh
19 &  &  &  &  &  &  &  &  &  &  & 4 & 25 & 300 & 2109 & 8414 & 25220 & 57598 & 92985 & 105424 & 77316 & 29771 & 5059 & 1302 &  &  & 405527 \hh
20 &  &  &  &  &  &  &  &  &  &  &  &  & 6 & 52 & 401 & 2181 & 6905 & 17039 & 32977 & 45891 & 44939 & 27879 & 7828 &  &  & 186098 \hh
21 &  &  &  &  &  &  &  &  &  &  &  &  &  &  & 9 & 36 & 205 & 876 & 2328 & 4882 & 8272 & 10236 & 9024 & 5094 & 1332 & 42294 
\\[2mm]
\hline
\hline [2mm]%
total& 1 & 23 & 130 & 550 & 1736 & 5079 & 12201 & 26961 & 52634 & 94210 & 152635 & 226658 & 307010 & 378728 & 426503 & 433576 & 401122 & 330227 & 246055 & 158812 & 91995 & 43347 & 18200 & 5094 & 1332 & 3414819 \\[1.2mm]
\hline
\end{mytab}%
\hss}%
}{The number of $\bt$ irreducible prime genus 4 alternating knots
tabulated by crossing number $c$ and number of $\sim$-equivalence
classes ($\#\sim$). \protect\label{tabg}}
\vss}
\newpage
\end{table}

\begin{theo}
There are 3,414,819 prime knot generators of genus $4$. Among them
there are 1,480,238 special and 1,934,581 non-special ones.
\end{theo}

Note that the effort with every new
genus increases dramatically, and new theoretical insight was
needed to reduce calculations to a feasible magnitude.
In our proof we will explain the algorithm used now.
(Still it seems out of scope to attempt the case $g=5$,
already for reasons of storing and maintaining the result.)

\proof 
For the maximal generators we use the work in \cite{SV}. We have to
determine first the planar
$3$-connected graphs $G$ having knot markings. They can be obtained by
a (graphic version of the) $\gm$-construction in \cite{Vdovina}
from graphs $G''$ of genus 3:
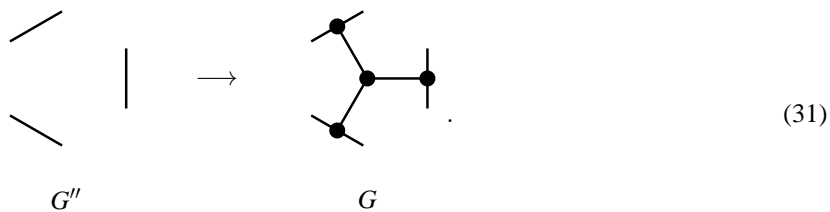
\begin{eqn}\label{gmG}
\diag{1cm}{2}{3}{
  \pictranslate{1 2}{
  \picputtext{0 -1.6}{$G''$}
  \picmultigraphics[rt]{3}{120}{
    \picline{0.8 0.4}{0.8 -0.4}
  }}
}\qquad\raisebox{0.5cm}{$\lra$}\qquad
\diag{1cm}{2}{3}{
  \pictranslate{1 2}{
  \picputtext{0 -1.6}{$G$}
  \picmultigraphics[rt]{3}{120}{
    \picline{0.8 0.4}{0.8 -0.4}
    \vrt{0.8 0}
    \picline{0.8 0}{0 0}
  }
  \vrt{0 0}
  }
}\,.
\end{eqn}
It is easy to see that if the graph $G$ on the right of \eqref{gmG} is
$3$-connected, then the graph $G''$ on the left must be at least $2$-%
connected. We determined the graphs $G''$ from the list of
maximal Wicks forms of genus $3$ compiled by A.\ Vdovina. (Note that
in the words an edge is a pair of inverse letters $a^{\pm 1}$,
while a vertex is a set of 3 letters $\{a,b,c\}$, such
that, for some choice of signs, $a^{\pm 1}b^{\pm 1}$,
$b^{\pm 1}c^{\pm 1}$ and $c^{\pm 1} a^{\pm 1}$ are subwords.)
From the graphs thus found, we selected only the planar
and $2$-connected ones, and obtained $34$ graphs $G''$. 

Then the graphs $G$ on the right of \eqref{gmG} are generated
by a $\gm$-construction. In applying the move \eqref{gmG},
some restrictions can be taken into account. In order
$G$ to become $3$-connected, on the left of \eqref{gmG}
not all $3$ segements can belong to the same edge in $G''$
(although two of them can). Moreover, whenever (a segement of)
one in a copy of a multiple (hence double) edge in $G''$ is
affected, the other copy cannot be. Otherwise, the graph $G$
is non-planar or still only $2$-connected, or $G''$ is not
$2$-connected.

All $G$ obtained were tested for planarity and $3$-connectedness
(and, of course, isomorphy). We obtained 50 graphs $G$. This
list was later verified by the program of Brinkmann and McKay
\cite{BMc} (which I discovered only afterwards and which made
it possible to deal with maximal generators of genus 5 and 6;
see \cite{3ang}).

Then for each graph $G$ on the right of \eqref{gmG}, we selected
all vertex orientations with the orientation of one vertex
fixed, and bisected edges between equally oriented vertices.
The new graph $G'$ corresponds to a knot iff the number of its
spanning trees is odd (see \cite{MS}).

All such graphs $G'$ were again tested for isomorphy.
We proved in \cite{SV} that the alternating diagrams with
these graphs as Seifert graphs do not admit flypes, so
that there is a one-to-one correspondence between alternating
knot, alternating diagram (up to mirroring and orientation) and
Seifert graph. Note that isomorphy between the $G'$ is preserved
under unbisecting edges, so that isomorphic $G'$ must
come from the same $G$, and the isomorphy check for $G'$
can be done for each $G$ separately. After duplications of
$G'$ were eliminated, the Dowker notations \cite{DT} were
generated, and so 42,294 maximal generators obtained (which
were later confirmed by the calculation in \cite{3ang}).

We also checked knot-theoretically that these generators are
distinct. The skein polynomial distinguished all except 70
pairs. Among these pairs, the Kauffman polynomial dealt
with 62. The remaining 8 pairs were distinguished by
the hyperbolic volume. 

We know from \cite{SV} that one can obtain the other,
not maximal, \em{special} generators from the maximal ones
by resolving reverse clasps. Then for each diagram of
a maximal generator we resolved any possible set of
(reverse) clasps, and identified among all these (alternating)
knots the generators of genus 4. They are 1,480,238.

Then we generated all possible special diagrams by flypes.
{}From these diagrams one can then obtain the other (non-special)
generating diagrams inductively by the number of separating
Seifert circles, using the (reverse to) the moves of Hirasawa
(figure \reference{fig4}), and the study of the effect of these moves on
generators, carried out in \S\reference{S33}.

We start with a generator diagram $D$ and want to determine
all generator diagrams $D'$ obtained from $D$ by
undoing a Hirasawa move. This means that we must identify the
segment $S$ in $D$ created by the move. We call $S$ below
the \em{Hirasawa segment}. There are several properties
of $S$ which follow from our treatise.
\begin{enumerate}
\item This is a segment of
odd length (i.e., odd number of crossings on it),
and no self-intersection.
\item Let $x$ and $y$ be the first and
last crossing on $S$, and call the other crossings \em{internal}.
By the work in \S\reference{S33}, we know that no internal crossings
should lie in a non-trivial $\sim$-equivalence class.
\item Moreover,
$x$ and $y$ have at most 2 $\sim$-equivalent crossings each.
So, prior to undoing the Hirasawa move, we have the option of
creating (by a $\bt$-move and possible flypes) two crossings
$\sim$-equivalent to $x$ and/or $y$, if $x$ and/or $y$ is in a
trivial $\sim$-equivalence class.
\item After the optional application
of such $\bt$-moves and flypes, we have the additional condition
that the length of the Hirasawa segment $S$ is strictly less
than the number of crossings outside $S$.
\item 
All possible such segments $S$ were generated, and their
crossings removed. Then the resulting Gau\ss{} codes of $D'$
were tested for realizability, $\bt$-reducedness,
primeness and for having the correct genus ($4$) and number of
separating Seifert circles (one more than in $D$).
\end{enumerate}

Note that we tested only necessary conditions for $S$. So some
of the diagrams $D'$ we obtain this way may actually not give
$D$ by a Hirasawa move. The difficulty is, though, to obtain
all diagrams $D'$ we need at least once, and we accomplish this,
since we reverse all correct Hirasawa moves. The decisive merit
of the Hirasawa move is that the way it provides to find all the
$D'$ is much more economical, in comparison to the previous
algorithm in \cite{STV,gen2}. Duplicated diagrams $D'$ can
be subsequently eliminated by bringing the
Dowker notation of any diagram to standard form. 

One can then iterate this procedure, undoing (potential)
Hirasawa moves (and performing the necessary consistency
checks), without the need to flype in between. The
process terminates, when for some number of separating
Seifert circles ($8$ here, and in general $2g$ by an
easy argument), no more (regular) diagrams arise. The
flypes are needed after the generation of the diagram
lists to identify the underlying knots, and the work is finished. \qed

In order to secure that this algorithm is practically
correct (theoretically it is justified by the proof of
theorem \reference{th1gen}), we verified that the subsets
of generators we obtained for $\le 18$ crossings coincide
with those that can be selected from the alternating knot
tables in \cite{KnotScape} and \cite{FR}. Prior to starting
genus $4$, we also applied this procedure for genus $3$. We
checked that the sets of generators we obtain coincide with
those generated previously in \cite{gen2}, by the algorithm
described there and in \cite{STV}. While the new implementation
dealt with genus $3$ about 50 times faster, it was still not
before several weeks of labour and computation (partly on
distributed machines) that table \ref{tabg} became complete.

A first small application is easy to obtain, and useful to mention.

\begin{corr}
There are 680 achiral generators of genus $4$.
The number of prime achiral alternating genus $4$ knots of crossing
number $n$ grows like $\ffrac{43}{2^8\,8!}n^8+O(n^7)$ for $n\to\infty$
even.
\end{corr}

\proof As in \cite{gen2}, one must find the achiral generators.
The use of Gau\ss{} sums (see proposition 13.3 of \cite{gen2})
pre-selects 798 generators, and 680 of them are indeed found
achiral by the hyperbolic symmetry program of J.\ Weeks. The check 
of the maximal number of $\sim$-equivalence classes of such generators 
shows it to be 18. This number is attained by 43 of the achiral knots, 
and none of them has non-trivial $\ssim$-equivalence classes, or 
admits (even trivial) type A flypes (as explained in figure
\ref{ffl2}). \qed

\section{Unknot diagrams, non-trivial polynomials and achiral knots
\label{S7}}

\subsection{Some preparations and special cases}

We apply now the classification of genus 4 generators to prove
theorem \reference{th1.6}. We begin with an easy proof
for 2-almost positive diagrams, that outlines the later methods. 

\begin{prop}
If $D$ is a $2$-almost special alternating diagram of the unknot,
then $D$ has a trivial clasp.
\end{prop}

\proof We use the Bennequin inequality stated in theorem \ref{Biq}.
It says that $g(D)\le 2$.

If $D$ is composite or of genus $1$, then the claim is easy to
obtain using the genus $1$ case \cite{apos}. So assume $g(D)=2$
and $D$ is prime.

If we have a $\sim$-equivalence class which
has crossings of oppsite (skein) sign, we either have a trivial 
clasp, or if not, we have a non-trivial $\sim$-equivalence class which
admits flypes. If after the flype and resolution of the clasp we still
have a diagram of genus two, we see that we did not start with a 
diagram of the unknot. By a direct check of the special genus two 
generators, we verify that this always happens.

We can thus assume that no $\sim$-equivalence class has crossings of 
oppsite (skein) sign. We call a $\sim$-equivalence class positive or 
negative after the sign of its crossings. 

Let $\hat D$ be the $2$-almost positive (generator)
diagram, obtained when removing in $D$ (by undoing $\bt$ twists)
crossings in (positive) $\sim$-equivalence classes, so that at most
2 crossings in each class remain. Then $\sg(\hat D)\le 0$.

Assume $\hD$ has two twist equivalent crossings $x$ and $y$ of
different sign. A reverse clasp in $\hD$ persists under $\bt$ twists,
even up to flypes, so $x$ and $y$ must be $\ssim$-equivalent.
Let us call their equivalence class \em{bad}.

One can flype the negative crossing $x\in\hD$ to form a parallel
clasp with $y$, and since $x$ is not twisted at in $\hD$ when
recovering $D$, this flype persists in $D$. Possible $\bt$ twists
in $D$ at $y$ create a $(2k+1,-1)$-rational tangle.
Then by a $(2k+1,-1)\to (2k,1)$ rational tangle move in $D$ (which
is a particular type of wave move),
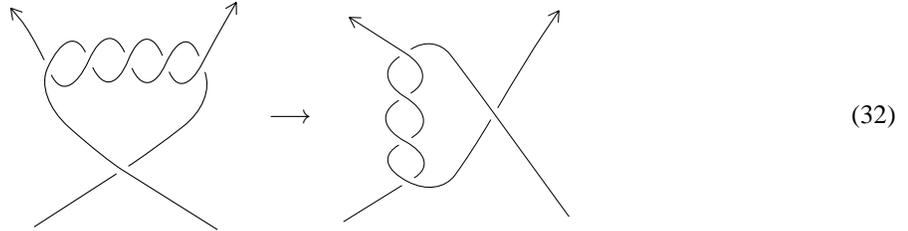
\begin{eqn}\label{rmv}
\epsfsv{3cm}{t1-twtg1}\quad\lra\quad
\epsfsv{3cm}{t1-twtg2}
\end{eqn}
we find a diagram $D_1$ with one crossing less and $\sg\le 0$, so
that $c(\hD_1)=c(\hD)+1$, and $\hD_1$ has one bad class less than
$\hD$. Still $D_1$ is $\le 2$-almost positive, though not necessarily
special. Repeat now such a move until $\hD_1$ has no bad class. 

Now by direct check all $2$-almost positive genus 2 generator diagrams
$\hD$ with $\sg\le 0$ and no bad class have $\le 7$ crossings and
are not special. (Actually $\sg=0$, since $\sg=4$ on positive genus
two diagrams.) They occur for the generators $6_3$, $7_6$ and $7_7$.


Since $D_1$ is not special, $D_1\ne D$, and $7\ge c(\hD_1)>c(\hD)$, so 
$c(\hD)\le 6$. The only special generator coming in question for $\hD$
is $5_1$. However, it is easy to check that $2$-almost positive
$5$-pretzel 
diagrams of the unknot have a clasp, which completes the proof. \qed

\begin{prop}\label{Za}
If $D$ is a $2$-almost positive diagram of the unknot, then $D$ can be 
made trivial by Reidemester moves not augmenting the crossing number. 
\end{prop}

\proof If $D$ is special, of genus $1$ or composite, the claim follows
from the almost positive case and the previous proposition. So we 
assume that $D$ is of genus two, prime and non-special. 

We take all non-special genus 2 generators, apply all possible flypes,
and switch crossings so that the diagrams are 2-almost positive.
We discard all diagrams of $\sg>0$, and from those of $\sg=0$
remove the duplicates. 

If now we have crossings of opposite sign forming a (trivial)
reverse clasp, this clasp will remain after flypes, and the
diagram is simplifiable using a Reidemeister II move. 
We discard thus diagrams with trivial reverse clasps.

If we have a trivial parallel clasp, we need to have a $\bt$-twist
at the positive crossing. We apply this twist in each such clasp, 
and discard all diagrams of $\sg>0$. There remain only 3 diagrams
in the series of $6_2$. They contain a tangle 
$\diag{8mm}{1.4}{1.6}{
  \pictranslate{-0.3 0}{
    \picveclength{0.13}\picvecwidth{0.11}
    \picline{1 1.1}{0.3 1.1}
    \picvecline{0.7 0.2}{1.3 0.8}
    \picmultivecline{0.17 1 -1.0 0}{0.7 0.8}{1.3 0.2}
    \picmulticurve{0.17 1 -1.0 0}{1.0 1.5}{0.7 1.5}{0.4 1.1}{0.7 0.8}
    \picmulticurve{0.17 1 -1.0 0}{1.3 0.8}{1.6 1.1}{1.3 1.5}{1.0 1.5}
    \picmultivecline{0.17 1 -1.0 0}{1 1.1}{1.7 1.1}
    \picline{0.9 1.1}{1.1 1.1}
    \piccurve{1.0 1.5}{0.7 1.5}{0.4 1.1}{0.7 0.8}
  }
}$
with a reverse clasp. This tangle can be simplified by Reidemeister III
unless there is a twist at the crossings in the reverse clasp.
Then again we check that $\sg>0$. \qed


We prepare some arguments that will be useful in the
proofs of theorems \reference{th1.6}, \reference{th1.7}
and \reference{th1.8}.

\begin{lemma}\label{lm1}
Assume $D$ is a $k$-almost positive diagram of genus $k$,
and $D$ depicts an achiral knot. Then $P(D)=1$.
\end{lemma}

\proof By Morton's inequality \eqref{mq}, we have $\md_l
P(D)\ge 0$. So achirality implies that $P(D)\in \bZ[l^2]$, i.e.
it has only a non-trivial $m$-coefficient in degree $0$.
The \cite[proposition 21]{LickMil} implies that $P(D)=1$.
\qed

\begin{lemma}\label{lm2}
Assume $D$ is a $k$-almost positive diagram of genus $k$,
and $P(D)=1$. Then for each pair $(a,b)$ of Seifert circles
in $D$, connected by a negative crossing, there are at least
two crossings connecting $(a,b)$.

Moreover, if $k\le 4$, then for at least one pair $(a,b)$
of Seifert circles connected by negative crossings,
there is a positive crossing connecting $(a,b)$.

\end{lemma}

Using lemma \reference{l45}, we have

\begin{corr}\label{cZ}
Any $k$-almost special alternating diagram with $P=1$ of genus
$k\le 4$ contains a trivial clasp up to flypes. \qed
\end{corr}

This immediately discards the special genus 4 generators.

\proof[of lemma \reference{lm2}]
We use \eqref{q2}. Since Morton's inequality is sharp, we see
that $\inx_-(D)=0$, so that each negative crossing has
a Seifert equivalent one. If no negative crossing has a
Seifert equivalent positive one, then for $k\le 4$, the
negative crossings lie in at most 2 Seifert equivalence classes.
Then an easy observation shows that the diagram $D$ is
$A$-adequate. So by theorem \reference{Vn1}, $V(D)\ne 1$, and
then the same is true for $P(D)$. \qed

\subsection{Reduction of unknot diagrams}


We prove here first a slightly weaker version of
theorem \reference{th1.6}, which allows for flypes.
We will explain in the next subsection how to remove the
need of flypes.

\begin{prop}\label{pz1}
For $k\le 4$, all $k$-almost positive unknot diagrams are
trivializable by crossing number reducing wave moves and flypes.
\end{prop}

\proof
We start with a $\le 4$-almost positive unknot diagram $D$, and
will show that one can perform, after possible flypes, a wave move
on $D$. In the following a wave move is always understood to be
one that \em{strictly reduces} the crossing number $c(D)$ of $D$.
Then we work by induction on $c(D)$. Note that a wave move does not
augment neither the number of positive, nor the number of negative
crossings (while it reduces one of both). So it does not lead out
of the class of $\le 4$-almost positive diagrams. It is possible
that the diagram genus becomes larger than $4$. However, in this
case we know that $D$ cannot be a diagram of the unknot. So we
need to be just interested in either showing that $D$ is not an
unknot diagram, or finding a wave move.

Next, genus and number of negative crossings is additive under
connected sum of diagrams. It is then sufficient to work with
prime diagrams $D$, with the following additional remark. If we
have a composite diagram $D=A\# B$, then flypes in $A$ and $B$
are not affected by taking the connected sum. However, we may have
a bridge $a$ in a factor diagram $A$, which is spoiled when taking
the connected sum with $B$. This can be remedied by displacing $B$
out of $a$. Such a move can formally be regarded as a flype, if in
figure \ref{fig_} one of $P$ or $Q$ has a diagram in which one
of its strands is an isolated trivial arc. (Let us call this move
\em{factor sliding flype}; it is as in figure \reference{figtan},
except that the tangle is additionally flipped.)

Let $D$ be the initial diagram. Clearly $g(D)\le 4$. Now let
$\hD$ be $D$ with \em{positive} $\sim$-equivalence classes reduced
(by $\bt$ moves) to 1 or 2 crossings. (That is, if a negative
$\sim$-equivalence class has $>2$ crossings, $\hD$ is not a
crossing-switched generator diagram.) Clearly $\sg(\hD)\le 0$. 

Let $k$ be the number of negative $\sim$-equivalence classes in $D$ (or
equivalently $\hD$), and $\ap=(a_1,\dots,a_k)$ be the vector of their 
sizes (number of crossings). We sort the integers $a_i$
non-increasingly, ignoring their order.

Let $D'$ be the positification of a diagram of a generator, 
whose series contains $D$ and $\hD$. That is, $D'$ is obtained
from $D$ or $\hD$ by replacing each $\sim$-equivalence class of
$n$ crossings with a class of one positive crossing if $n$ is odd
and 2 positive crossings if $n$ is even. Then $\sg(D')\le 2k$,
since twists in a $\sim$-equivalence class do not alter $\sg$
by more than $\pm 2$. In particular, we can exclude the case
$k=1$ (i.e. $\ap$ being one of $(4)$, $(3)$ or $(2)$), because 
for a positive diagram $D'$ the property $\sg=2$ implies (see
end of \S\reference{Sz}) that $g(D')=g(D)=1$, and then we are
easily done by \cite{gen1} identifying genus 1 unknot diagrams.

In the first step to generate $D$, we try to obtain $\hD$ from
$D'$, while $D'$ is obtained directly from the generator table by
flypes and positification. Flypes on the generator are necessary,
because they commute with the $\bt$-twists only up to mutations,
which we try to avoid (in our stated repertoire of moves simplifying
the unknot diagram).

For the possible $\ap$, the potential diagrams $\hD$ are determined by 
crossing switches (and in the case of $\ap=(3,1)$, one $\bt$ twist).

It is convenient to check the bound for $\sg$, which is a necessary
condition for $\sg(\hD)\le 0$, after each crossing switch. Moreover,
for given $\ap$ we select only the $D'$ with $\sg(D')\le 2k$. This
allows to discard 
a lot of possible continuations of the crossing switch procedure,
which are understood not to a lead to a proper $\hD$. Also, one should
order the $\sim$-equivalence classes of 1 and 2 elements and switch
w.r.t.
this order to avoid repetitions. Concretely, we proceeded as follows. 
\begin{itemize}
\item $\ap=(1,1,1,1)$. Change successively one crossing and check $\sg
 \le 6$, $\sg\le 4$,  $\sg\le 2$, $\sg\le 0$. 
\item $\ap=(2,1,1)$. Change a reverse clasp and two crossings and check 
 $\sg\le 4$,  $\sg\le 2$, $\sg\le 0$.
\item $\ap=(2,2)$. Change two reverse clasps and check $\sg\le 2$, 
 $\sg\le 0$.
\item $\ap=(3,1)$. Change twice one crossing and check $\sg\le 4$ and
 $\sg\le 2$. Then perform a $\bt$ at one of the negative crossings and
 test $\sg\le 0$.
\item The cases $\ap=(2,1)$ and $(1,1,1)$ are analogous.
\end{itemize}


The diagrams $\hD$ found are relatively few in comparison to the number 
of genertors $D'$. Also, $\sg$ is a mutation invariant. Therefore, to 
handle the flypes in practice, it is easier to start with one diagram 
$D'$ per generator, determine whether it gives a $\hD$, and only then to
process the other diagrams obtained by flypes from $D'$. Actually, one
can apply the flypes in $\hD$ directly, except for $\ap=(3,1)$. Latter
is easy to handle, since the $D'$ with $g(D')=2$ we need to consider
are few. Thus let us exclude this in the following.

If $\hD$ has a trivial clasp up to flypes, then all $D$ arising from 
this $\hD$ will admit a $(2k+1,-1)\to (2k,1)$ rational tangle move 
\eqref{rmv} up to flypes, and so we can discard such $\hD$. Note that 
only a negative crossing in a trivial $\sim$-equivalence class can 
form a parallel clasp up to flypes, and since in passing from $\hD$
to $D$ we do not twist at negative crossings, one can flype the 
negative crossing properly also in $D$ for \eqref{rmv} to apply.

Let $\hat\cD$ be the set of diagrams $\hD$ obtained/maintained so
far. In trying to reconstruct $D$ from $\hD$ we distinguish at
which positive $\sim$-equivalence classes we apply twists.

Let $X$ be a diagram and $\cX$ a set of diagrams. Fix for $X$ the set
$C_X$, consisting of one crossing in each positive
$\sim$-equivalence class in $X$. Let 
\begin{eqn}\label{*77}
X_*\,=\,\left\{\,X'\,:\,\begin{tabular}{c}
$X'$ is obtained from $X$ by 0 or 1\\ 
twist at each crossing in $C_X$
\end{tabular}\,
\right\}\,, \quad\mbox{and}\quad \cX_*\,=\,\bigcup_{X\in\cX}\,X_*\,.
\end{eqn}
(The choice to twist or not is made independently for each crossing,
so $|X_*|=2^{|C_X|}$.) The reason we will make use of this
construction shortly is the following 

\begin{lemma}\label{lm4}
If a diagram $D$ admits a wave move, and $D$ has a $\sim$-equivalence
class of at least $4$ crossings, then one can flype and reduce $D$ by
a $\bt$-move, so that the resulting diagram $D'$ still admits a wave
move. Contrarily, if a diagram $D'$ admits a wave move, and $D$ is
obtained from $D'$ by a $\bt$-move at a non-trivial $\sim$-equivalence
class, then $D'$ admits a wave move after flypes.
\end{lemma}

\proof We focus on the first claim. Let $D$ be presented w.r.t. four
of the crossings in the $\sim$-equivalence class $T$,
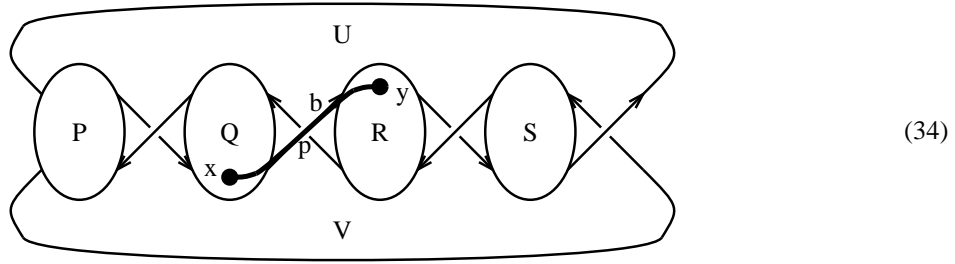
\begin{eqn}\label{fde}
\diag{1cm}{10}{5}{
  \picmultigraphics{2}{-4 0}{
    \picvecline{9 2}{8 3}
    \picmultivecline{-5 1 -1 0}{8 2}{9 3}
    \picvecline{6 3}{7 2}
    \picmultivecline{-5 1 -1 0}{7 3}{6 2}
  }
  \opencurvepath{9 3}{9.5 3.5}{9 4.3}{1 4.3}{0.5 3.5}{1 3}{}
  \opencurvepath{9 2}{9.5 1.5}{9 0.7}{1 0.7}{0.5 1.5}{1 2}{}
  \picmultigraphics{4}{-2 0}{
    \piccurve{8 3}{7.7 3.3}{7.3 3.3}{7 3}
    \piccurve{8 2}{7.7 1.7}{7.3 1.7}{7 2}
  }
  \picfilledellipse{7.5 2.5}{0.6 0.9}{S}
  \picfilledellipse{5.5 2.5}{0.6 0.9}{R}
  \picfilledellipse{3.5 2.5}{0.6 0.9}{Q}
  \picfilledellipse{1.5 2.5}{0.6 0.9}{P}
  \picputtext{5 3.8}{U}
  \picputtext{5 1.2}{V}
  \pictranslate{1 0}{
    \piclinewidth{15}
    \opencurvepath{2.5 1.9}{2.7 1.9}{3 2}{4 3}{4.3 3.1}{4.5 3.1}{}
    \picputtext{3.65 2.9}{b}
    \picputtext{3.5 2.3}{p}
    \picputtext{2.25 2.0}{x}
    \picputtext{4.8 3.0}{y}
    \point{2.5 1.9}
    \point{4.5 3.1}
  }
}
\end{eqn}
with the (possibly empty) tangles $P$, $Q$, $R$, $S$ being
enclosed by these crossings. Let $b$ be a bridge in $D$ that
can be shortened by a wave move, and $\len b$ the number of
overcrossings $b$ contains. $b$ is understood to start and end
on points $x$ and $y$ on the segment (edge) of the diagram
line after/before the beginning/ending undercrossing. Now the
existence of a wave move means that, if we delete $b$ from the
plane curve of $D$, the distance between the regions $X$ and
$Y$ that contain the endpoints $x,y$ of the resulting arc is
smaller than $\len b$. Distance between regions $X$ and $Y$
is meant as the shortest length of a sequence of successively
neighbored regions, starting with $X$ and ending with $Y$,
counting \em{one} of $X$ and $Y$, and neighbored regions means
regions sharing an edge.

It is an easy observation that $b$ contains at most one crossing
in $T$ (assuming these crossings have all the same sign). Let
this be the crossing $p$ between $Q$ and $R$. (The case that $b$
does not contain any crossing in $T$ is dealt with a similar,
though slightly simpler argument.) Now let $\gm=(X,\dots,Y)$
be a shortest length sequence of neighbored regions connecting
$X$ and $Y$. If $\gm$ contains a region $W$ \em{completely}
inside $P$ or $S$, then it must contain also $U$ and $V$. However,
the distance of $U$ and $V$ is two, and moving from $U$ to $V$,
it is not necessary to pass via $W$. So $\gm$ can be chosen not to
contain any regions completely inside $P$ or $S$ (though it may 
have to contain a region whose corner is one of the crossings in $T$).

This means now that the flype that makes $P$ and $S$ empty (moving
their interior into $Q$ or $R$), does not change the distance
between $X$ and $Y$. Likewise, it does not affect the bridge $b$.
So $b$ remains too long after such flypes. Now when $P$ and $S$ are
empty, one can reduce $D$ by a $\bt$ move, and so we are done
with the first statement.

The proof of the second statement is similar. Again a wave
move in $D'$ can be chosen so that it does not affect the
diagram outside $Q$, $R$ and $p$. Then in $D$ one can adjust
by flypes the new crossings created by the $\bt$-move at a crossing
$q$ to lie outside that portion of the diagram if $p=q$, and
entirely within $R$ or $Q$ if $p\ne q$.
\qed

Now from our set $\hat\cD$ we generate the set $\cD=(\hat\cD)_*$.
Our understanding is that one can obtain all $D$ from an $X\in
\cD$ by $\bt$ moves at a non-reduced $\sim$-equivalence class.
The point in this restriction comes from the previous lemma,
which says that the applicability of a wave move up to flypes
persists under such $\bt$ moves. So if some $X\in \cD$ admits a
wave move, one can discrard all $D$ obtained from $X$ by twising
at a non-reduced $\sim$-equivalence class. (In fact, the lemma
shows that it is necessary to twist only at a crossing in a
trivial positive $\sim$-equivalence class, but we will need
\eqref{*77} it is given form shortly below.) Similarly if $\sg(X)>0$,
the same holds for all $D$ obtained from $X$ by positive $\bt$
twists. Such $X$ can be likewise discarded.

Therefore, discard all diagrams $X$ in $\cD$ admitting a wave
move or having a positive signature, obtaining a subset $\cE$ of 
$\cD$. We deal with $\cE$ using the maximal degree of the Jones
polynomial. This procedure was described in \cite{canon}. What
we need in our case is the following property.

\begin{lemma}(regularization; \cite{canon})\label{lrez}
Assume that there is an integer $n$ such that 
\begin{eqn}\label{rez}
n=c(D'')-\Md V(D'')\ \mbox{for all}\ D''\in D_*\,.
\end{eqn}
Then \eqref{rez} holds for all diagrams $D'$ obtained by positive
$\bt$-twists and flypes from $D$.
\end{lemma}

The proof is a simple inductive application of the skein relation
of $V$. We called in \cite{canon} property \eqref{rez} \em{regularity},
and its verification \em{regularization}. We will will use below the
same terms. (There is a similar property $P$, which we will encounter
in some form in \S\reference{S6}.)

We check first that $\Md V(E)>0$ for all diagrams $E\in \cE$. Then
split $\cE$ into sets $\cE_n$ of diagrams $E$, accodring to $n=c(E)
-\Md V(E)$. (In practice the occurring values of $n$ are between 3
and 6.) Now $D$ is supposed to be obtained from some diagram in
some $\cE_n$ by some number of positive $\bt$-twists. To rule
this out, build $(\cE_n)_*$ and check for each $E\in (\cE_n)_*$ that
$n=c(E)-\Md V(E)$. Using the skein relation of $V$ one can then
show that this property is preserved under further positive $\bt$-%
twists. So if $D$ is a diagram not admitting a wave move, then 
\[
\Md V(D)= c(D)-n\ge c(E)-n=\Md V(E)>0\,.
\]
Thus no unknot diagram can occur as $D$.
With this proposition \ref{pz1} is established.
\qed

\subsection{Simplifications}

There are a few shortcuts in the above procedure which we
mention. In building $\cE$ from $\cD$, one can also discard
all semiadequate diagrams $X$. Semiadequacy persists under
twising, and by \cite{Thistle} (or by theorem \reference{Vn1})
the unknot has no non-trivial
semiadequate diagrams. Since we will come to talk about it
in the sequel, we call this step shortly the \em{semiadequacy
shortcut}. For knots with $P=1$, it can always be applied.

When switching crossings one by one to obtain $\hD$ from $D'$,
for $c_-(D)=g(D)\in \{3,4\}$, we can use lemma \reference{lm2}.
We discard special generators immediately and for non-special
generators switch only crossings that have a Seifert equivalent one.
For $g=4$, we are left, after flyping and discarding $\hD$ with
trivial clasps, with 31 possible diagrams. Using the semiadequacy
shortcut, we find $\#\cE=4$, and regularization applies with $n=5$.

Another shortcut comes from the Rudolph-Bennequin inequality
\eqref{RBi}. It implies the following:

\begin{lemma}\label{lm7.}
If $D$ is a diagram of a knot $K$, and $k$ the number of
$\sim$-equivalence classes of $D$ which contain a negative
crossing, then $g_s(K)\ge g(D)-k$.
\end{lemma}

\proof We flype $D$ so that between the $n$ negative crossings in each
such $\sim$-equivalence class, we have $n-1$ negative Seifert circles
(of valence 2). The estimate then follows from \eqref{RBi}. \qed

This means that for given $\ap$ (index of distribution of negative
crossings into $\sim$-equivalence classes) of length $k$, one
needs to consider only $D'$ with $g(D')\le k$. We
call this step the \em{Rudolph-Bennequin shortcut}.

The reason for not including these shortcuts into the proof
of proposition \reference{pz1} is that in the modifications of
the proof for the related theorems their legitimity will vary
from case to case.

We now set out to prove theorem \reference{th1.6} in its full form.
The basic point of removing flypes from proposition \reference{pz1}
lies in the construction \eqref{*77} and lemma \reference{lm4}.
We first improve lemma \reference{lm4}.

\begin{lemma}\label{lm4'}
Assume a diagram $D$ is given in the presentation \eqref{fde}, with
at least 2 tangles $R$, $Q$. Assume $D$ admits a reducing wave move
shrinking a bridge $b\subset R\cup Q\cup \{p\}$ to a bridge $b'$.
Then this move can be chosen so that $b$ lies entitrely within
$R$ or $Q$, or in $R\cup\{p\}$ or $Q\cup\{p\}$.
\end{lemma}

\proof If $b$ does not contain $p$, then $b$ lies entitrely within
$R$ or $Q$. So assume again that $b$ passes $p$. It is obvious
that $b'$ can be chosen not to pass any crossing outside of
$R\cup Q\cup \{p\}$. If $b'$ passes $p$, then the wave move
between $b$ and $b'$ would split into two parts, one which
modifies the part of $b$ between $x$ and $p$, and another
one which acts on the remaining part between $p$ and $y$.
At least one of these would be a reducing move, and so we
are done.

Thus we may assume that $b'$ does not pass $p$. Then there
are 2 options how $b'$ looks outside $R$ and $Q$:
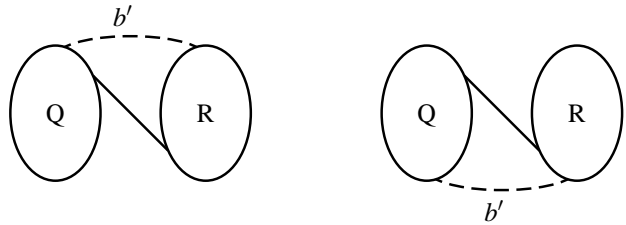
\[
\diag{1cm}{4.2}{3}{
  \picline{2.1 2}{3.1 1}
  {\piclinedash{0.2 0.1}{0.25}
   \piccurve{3.6 2.3}{3.3 2.6}{1.9 2.6}{1.6 2.3}
  }
  \picputtext{2.5 2.8}{$b'$}
  \picfilledellipse{3.6 1.5}{0.6 0.9}{R}
  \picfilledellipse{1.6 1.5}{0.6 0.9}{Q}
}
\qquad
\diag{1cm}{4.2}{3}{
  \picline{2.1 2}{3.1 1}
  {\piclinedash{0.2 0.1}{0.25}
   \piccurve{3.6 0.7}{3.3 0.4}{1.9 0.4}{1.6 0.7}
  }
  \picputtext{2.5 0.2}{$b'$}
  \picfilledellipse{3.6 1.5}{0.6 0.9}{R}
  \picfilledellipse{1.6 1.5}{0.6 0.9}{Q}
}
\]
They are symmetric, so look just at the first one.
By drawing a circle $\gm$ like
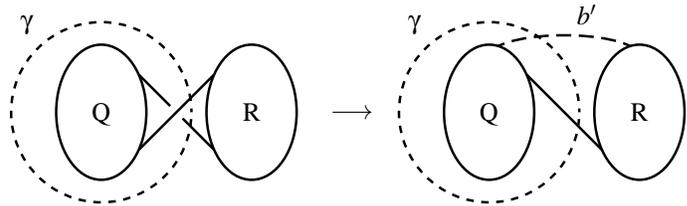
\[
\diag{1cm}{3.7}{3}{
 \pictranslate{-0.5 0}{
  \picline{2.1 2}{3.1 1}
  \picmultiline{-7 1 -1 0}{2.1 1}{3.1 2}
  {\piclinedash{0.1 0.1}{0.05}
   \piccircle{1.6 1.5}{1.2}{}
  }
  \picputtext{0.6 2.7}{$\gm$}
  \picfilledellipse{3.6 1.5}{0.6 0.9}{R}
  \picfilledellipse{1.6 1.5}{0.6 0.9}{Q}
 }
}
\quad\lra\quad
\diag{1cm}{3.7}{3}{
 \pictranslate{-0.5 0}{
  \picline{2.1 2}{3.1 1}
  {\piclinedash{0.1 0.1}{0.05}
   \piccircle{1.6 1.5}{1.2}{}
  }
  {\piclinedash{0.2 0.1}{0.25}
   \piccurve{3.6 2.3}{3.3 2.6}{1.9 2.6}{1.6 2.3}
  }
  \picputtext{2.9 2.8}{$b'$}
  \picfilledellipse{3.6 1.5}{0.6 0.9}{R}
  \picfilledellipse{1.6 1.5}{0.6 0.9}{Q}
  \picputtext{0.6 2.7}{$\gm$}
 }
}\,,
\]
we see that again $\gm$ splits the wave move into two parts,
for which the bridge lies within $Q\cup\{p\}$ and $R$ resp., and
at least one of which is reducing. \qed

Now we need a more delicate treatment of the freedom
to apply flypes at a crossing of a $\sim$-equivalence class.
It will be convenient to slightly change the nomenclature of
tangles from \eqref{fde}.

\begin{defi}\label{ddg}
Let the \rm{(flyping) degree} $\deg C$ of a $\sim$-equivalence class
$C$ in $D$ be the maximal number of tangles $T_i$ in a presentation
of $D$ like \eqref{fde'} (where it is 4), such that no $T_i$ contains
crossings in $C$, and each crossing not contained in any tangle
belongs to $C$. We call such $T_i$ \em{essential} (w.r.t. $C$).
\begin{eqn}\label{fde'}
\diag{1cm}{10}{5}{
  \picmultigraphics{2}{-4 0}{
    \picvecline{9 2}{8 3}
    \picmultivecline{-5 1 -1 0}{8 2}{9 3}
    \picvecline{6 3}{7 2}
    \picmultivecline{-5 1 -1 0}{7 3}{6 2}
  }
  \opencurvepath{9 3}{9.5 3.5}{9 4.3}{1 4.3}{0.5 3.5}{1 3}{}
  \opencurvepath{9 2}{9.5 1.5}{9 0.7}{1 0.7}{0.5 1.5}{1 2}{}
  \picmultigraphics{4}{-2 0}{
    \piccurve{8 3}{7.7 3.3}{7.3 3.3}{7 3}
    \piccurve{8 2}{7.7 1.7}{7.3 1.7}{7 2}
  }
  \picfilledellipse{7.5 2.5}{0.6 0.9}{$T_4$}
  \picfilledellipse{5.5 2.5}{0.6 0.9}{$T_3$}
  \picfilledellipse{3.5 2.5}{0.6 0.9}{$T_2$}
  \picfilledellipse{1.5 2.5}{0.6 0.9}{$T_1$}
}
\end{eqn}
We call $C$ \em{(flyping) inactive} if it has degree $1$. Hereby
we insist on the strand orientations given in \eqref{fde'}. Note
that in a trivial $\sim$-equivalence class, a crossing may admit
type A flypes. In that case we still regard its class as inactive.
\end{defi}

The following lemma shows how much the degree of a
$\sim$-equivalence class can grow. It is useful (and was used)
as a consistence test during the clacultion.
An \em{$n$-twist tangle} is a tangle (with a diagram) of $n$ crossings,
containing $n-1$ clasps. Again we can distinguish between parallel
and reverse twists if $n>1$. A \em{$(n,m)$-pretzel tangle} is the
sum of an $n$- and $m$-twist tangle, performed so that we do not
obtain a $n+m$-twist tangle.

\begin{lemma}\label{lmda}
In a knot diagram of genus $g$, each $\sim$-equivalence class
$C$ has degree at most $g$. Moreover, if $\deg C=g$ and $C$ is
non-trivial, then one of the $T_i$ is a twist tangle, or a
$(n,m)$-pretzel tangle for $n,m$ odd (and reverse twists in
either case).
\end{lemma}

\proof Note that the Seifert circles in \eqref{fde'} (thickened
below) are of the shape
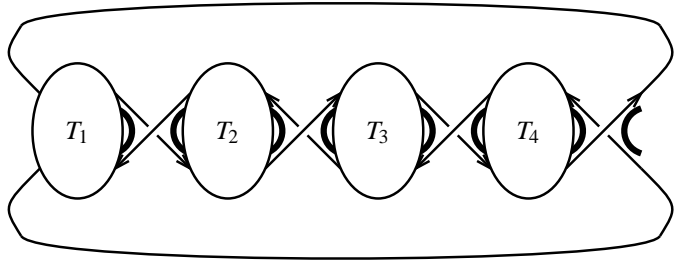
\[
\diag{1cm}{10}{5}{
  \picmultigraphics{2}{-4 0}{
    \picvecline{9 2}{8 3}
    \picmultivecline{-5 1 -1 0}{8 2}{9 3}
    \picvecline{6 3}{7 2}
    \picmultivecline{-5 1 -1 0}{7 3}{6 2}
  }
  \opencurvepath{9 3}{9.5 3.5}{9 4.3}{1 4.3}{0.5 3.5}{1 3}{}
  \opencurvepath{9 2}{9.5 1.5}{9 0.7}{1 0.7}{0.5 1.5}{1 2}{}
  \picmultigraphics{4}{-2 0}{
    {\piclinewidth{12}
     \piccurve{8.0 2.8}{8.3 2.7}{8.3 2.3}{8 2.2}
     \piccurve{9.0 2.8}{8.7 2.7}{8.7 2.3}{9 2.2}
    }
  } 
  \picfilledellipse{7.5 2.5}{0.6 0.9}{$T_4$}
  \picfilledellipse{5.5 2.5}{0.6 0.9}{$T_3$}
  \picfilledellipse{3.5 2.5}{0.6 0.9}{$T_2$}
  \picfilledellipse{1.5 2.5}{0.6 0.9}{$T_1$}
}\,.
\]
We look at the Euler characteristic. First $\chi(D)=1-2g$. We 
apply next the replacement 
\begin{eqn}\label{*&}
\diag{1cm}{2}{2}{
  \pictranslate{0 -0.5}
  {\piclinewidth{12}
   \piccurve{1.5 1.8}{1.8 1.7}{1.8 1.3}{1.5 1.2}
   \piccurve{0.5 1.8}{0.2 1.7}{0.2 1.3}{0.5 1.2}
  }
  \picfilledellipse{1 1}{0.6 0.9}{$T_i$}
}\quad\lra\quad
\diag{1cm}{2}{2}{
  {\piclinewidth{12}
   \picellipse{1 1}{0.8 0.5}{}
  }
}
\end{eqn}
(in other words, ``emptying out'' $T_i$ and making it into a single
Seifert circle). One observes that, because $T_i$ is essential,
\eqref{*&} augments $\chi$ by at least 1, and by exactly 1 only if
$T_i$ is a $(n,m)$-pretzel tangle or a $2n$-twist tangle.
Connectivity of $D$ implies that in former case $n,m$ are odd, and
in either case that at most one such $T_i$ can occur. Thus
\eqref{*&} for all
but at most one $i$ augments $\chi$ by at least 2. Finally, after all
instances of \eqref{*&} are performed, the resulting diagram $D_0$
is still connected, so $\chi(D_0)\le 1$.

The first claim of the lemma follows then by a simple estimate. 
The second claim follows similarly, because if $C$ is non-trivial,
then $\chi(D_0)=0$. Thus at least once in \eqref{*&}, $\chi$ must
go up only by $1$. \qed


\begin{defi}\label{*77a}
Let $X$ be a diagram. Fix for $X$ a set $C_X$, consisting of
one positive crossing in each non-negative $\sim$-equivalence
class in $X$. (Non-negative means here that the class
contains at least one positive crossing.) Let $p\,:\,C_X\,\to
\,\bN$ be a function, defined on $c\in C_X$ depending on its
$\sim$-equivalence class $C$ as follows:
\[
p(c)\,=\,\left\{\es
\begin{array}{c@{\es\mbox{if}\ }l}
1 & \mbox{$C$ is trivial and inactive,} \\
0 & \mbox{$C$ is non-trivial and inactive,} \\
0 & \mbox{$C$ is active and $|C|\ge \deg C$,} \\
\Br{\ffrac{\deg C-|C|+1}{2}} & \mbox{$C$ is active and $|C|<\deg C$.}
\end{array}\,\right\}\,.
\]
Then let
\[
X_+\,=\,\left\{\,X'\,:\,\begin{tabular}{c}
$X'$ is obtained from $X$ by $0,\dots,p(c)$ \\ 
twists at $c$ for each $c\in C_X$ \em{and type} \\
\em{B flypes}, so that $X'$ has no trivial clasp
\end{tabular}\,
\right\}\,.
\]
(The choice of twists is again understood to be independent for each
crossing, so that there are $\prod_c (p(c)+1)$ choices. Now, however,
$|X_+|$ will in general be larger, due to the flypes.)
Again if $\cX$
is a set of diagrams, we let $\cX_+\,=\,\bigcup_{X\in\cX}\,X_+$\,.
\end{defi}

The next lemma now shows up to how many crossings we must twist
in each $\sim$-equivalence class, in order to be sure that
reducing wave moves will presist (without need of flypes)
if more crossings are added.

\begin{lemma}
Assume that for a generator $D$, every diagram in $D_+$ admits
a reducing wave move. Then so does every diagram in its series
$\br{D}$ obtained by type B flypes and positive $\bt$-twists from $D$.
\end{lemma}

\proof Let $D'\in \br{D}$. We fix a $\sim$-equivalence class $C$
in $D'$ that contains at least one positive crossing. Then
we will obtain below a diagram $D''\in \br{D}$ from $D'$
by flypes and removing pairs of crossings in $C$ by undoing twists.
In this diagram $D''$ the class $C$ will have at most $\deg C+1$
crossings if $C$ is active, and at most $3$ crossings if $C$ is
inactive. We will argue that if $D''$ admits a reducing wave move,
then so does $D'$. The lemma follows then by induction on the
$\sim$-equivalence classes of $D$.

In order to find this diagram $D''$, first assume that $C$ is
inactive. Let $T=T_1$ be the complementary tangle to $C$.
We may then assume that $C$ is positive, since otherwise
we have a trivial clasp. (In the absence of flypes, one
cannot put tangles ``between'' clasps.)
If $|C|\le 3$, we may set $D''=D'$, so let $|C|\ge 4$.
Now remove by undoing twists pairs of crossings in $C$ until
$2$ or $3$ remain. Let this be $D''$. 

If $D''$ admits a reducing wave move, then lemma \ref{lm4'},
applied on $T_1=T$ and $T_2$ consisting of the $|C|-2$ internal
crossings in $C$, shows that this move avoids $T_2$. (Note
that, unlike definition \ref{ddg}, lemma \ref{lm4'} makes no
assumption whatsoever on $T_i$.) Then $T_2$ can be replaced by
arbitrary many crossings in $C$, and our claim for $D''$ is
justified.

Next let $C$ be active.
We assume that $|C|\ge \deg C+2$, and that between the
$T_i$ only crossings of the same sign occur in $D'$.
(Otherwise we have again a trivial clasp.)
We call such a collection of crossings a \em{group}.
We attach a non-zero integer to a group by saying
that it is a group of $-t_i$ if it has $t_i$
negative crossings, and a group of $t_i$ if it has $t_i$
positive crossings.

If there is a triple of crossings, i.e.\ a group of $t_i$ for
$|t_i|\ge 3$, then we can undo a $\bt$-twist. Each $T_i$
will have a crossing of the same sign on either side, so
by lemma \reference{lm4'} a wave move in the simplified diagram
would give a wave move in the original diagram.

If there are two pairs of crossings, then we remove one crossing
in each pair, and flip the portion of the $T_i$ on one side
around. (This is a flype followed either by undoing a twist, or
removing a trivial clasp.) The same argument applies, only that
bridges might have turned into tunnels. (Keep in mind that
we regard tunnels as equivalent to bridges.)

With this we managed to simplify $C$ by two crossings, so
to find $D''$, repeat this until $|C|-\deg C\in\{0,1\}$.

Let us finally briefly argue why we can assume that $D''$
has no trivial clasps. First, we assumed that in $D'$
every group of a $\sim$-equivalence class is signed,
and constructed $D''$ so that it has the same property.
This is why $D''$ has no trivial reverse clasp.
If $D''$ had a trivial parallel clasp, then the two
crossings in it would be in a trivial and inactive
$\sim$-equivalence class. Thus the only difference
between $D''$ and $D'$ in that classes could be that at the
positive crossing, some $\bt$-twists are applied.
We would then have a variant of the move \eqref{rmv} in $D'$.
\qed

The following modification of the previous lemma is 
the analogue of the test given in the paragraph after
the proof of lemma \reference{lm4}, which we apply now
to avoid flypes.

\begin{lemma}
Assume that $D$ is a (not necessarily alternating) generator
diagram, or obtained from such a diagram by one negative $\bt$-twist,
and possible flypes. Assume that each diagram $D''$ in $D_+$ admits
a reducing wave move, or $\sg(D'')>0$. Then some of these two
alternatves applies to every diagram $D'\in \br{D}$ 
obtained by type B flypes and positive $\bt$-twists from $D$.
\end{lemma}

\proof Remark that in the previous proof undoing 
a $\bt$-twist at a positive crossing, as well as a flype,
does not augment $\sg$.

Then we consider only the situation where we need
to undo a negative $\bt$-twist. Let $C$ again be
the $\sim$-equivalence class modified in passing from
$D'$ to $D''$. Following the argument in the previous
proof, we see that a negative $\bt$-twist could occur
(after a possible flype)
only if in $C$ we have a group of $t_i\le -3$,
or two groups of $t_i=-2$. Then we have $\ge 3$ negative
crossings (at at least one positive crossing) in $C$.

So assume that $D$ has a $\sim$-equivalence of $\ge 3$
negative crossings. This means, by the Rudolph-Bennequin
shortcut, that we need to consider only generators of genus
at most 2. The description of unknot diagrams of genus
1 is easy.

We thus assume $g(D')=2$. By lemma \ref{lmda}, we have
$\deg C\le 2$. We assumed that $C$ has a positive crossing,
and thus if $D'$ has no trivial clasp, at least one of the two
groups in $C$ has $t_1>0$. Because of the negative crossings in
$C$, the other one must have $t_2\le -3$, and in particular
$\deg C=g(=2)$. The second part of the lemma then shows that one
of the $T_i$ is a twist tangle resp. pretzel tangle. (One can
use alternatively theorem \reference{thgen} directly.) This
tangle, say $T_1$, is connected by a crossing in $C$ on one side
in a non-alternating way. Thus for a twist tangle $T_1$, we have
a version of the move \eqref{rmv} (in which the twist tangle
would have an even number of crossings and orientation would
be slightly different).  

For a pretzel tangle $T_1$, we can apply a similar move if one
of $m$ or $n$ is equal to $1$. If $m,n\ge 3$, we will argue that
the diagram is knotted. We have in $D'$ a tangle like
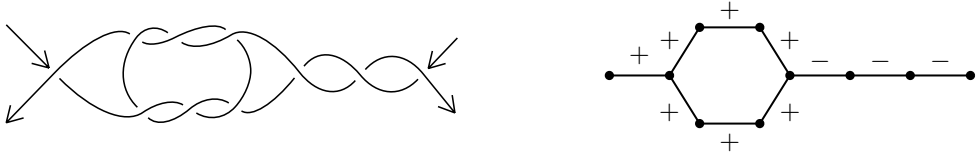
\[
\epsfsv{3cm}{t1-4apunkn}\kern2cm
\diag{8mm}{6}{3}{
  \pictranslate{0 0.5}{
    \ppt{6 1} \ppt{5 1} \ppt{4 1} \ppt{3 1} \ppt{1 1} \ppt{0 1}
    \ppt{1.5 1.8} \ppt{2.5 1.8} \ppt{1.5 0.2} \ppt{2.5 0.2}
    \picline{6 1}{3 1}
    \picline{0 1}{1 1}
    \picline{2.5 1.8}{1.5 1.8}
    \picline{2.5 0.2}{1.5 0.2}
    \picline{3 1}{2.5 1.8}
    \picline{3 1}{2.5 0.2}
    \picline{1 1}{1.5 1.8}
    \picline{1 1}{1.5 0.2}
    \picputtext{5.5 1.25}{$-$}
    \picputtext{4.5 1.25}{$-$}
    \picputtext{3.5 1.25}{$-$}
    \picputtext{0.5 1.35}{$+$}
    \picputtext{1.0 0.4}{$+$}
    \picputtext{1.0 1.6}{$+$}
    \picputtext{3.0 0.4}{$+$}
    \picputtext{3.0 1.6}{$+$}
    \picputtext{2 2.1}{$+$}
    \picputtext{2 -.1}{$+$}
  }
}
\]
Its Seifert graph is shown on the right.
We use now the inequality \eqref{q2} of Murasugi-Przytycki.
(The problematic cases discussed below in \S\reference{S61}
do not occur here.) By successively contracting the leftmost
vertex of the negative edges, we see that $\inx_-(D')\ge 3$.
Then \eqref{q2}, slightly rewritten using \S\ref{kkd}, yields
\[
\md_l P(D')\,\ge\,2g(D')-2c_-(D')+2\inx_-(D')\,\ge\,2\,.
\]
The skein polynomial thus shows that the unknot cannot occur
in such a diagram. \qed

\proof[of theorem \reference{th1.6}]
The case $k=1$ was known, and $k=2$ can be easily concluded
from the proof of proposition \reference{Za}. Let us thus
assume $k=3,4$. Moreover, by theorem \reference{thge}, we may
assume $g(D)\ge 2$.

The work we did allows us to proceed now similarly to
proposition \ref{pz1}. We consider generators of
genus $2$ to $k$. The factor slide allows us to
discard composite generators, and corollary \reference{cZ}
special ones of genus $k$. We take the positfication of
one diagram per generator. We switch
crossings (successively testing $\sg$ to save work) to
find all $k$-almost positive diagrams of $\sg\le 0$. We must
here, though, allow for crossings of opposite sign in the
same $\sim$-equivalence class, incl. for trivial clasps.

Our understanding is that each diagram we have to treat
is obtained from the grenerator diagram first by type A
flypes, then $\bt$-twists and finally type B flypes.
Thus first we apply type A flypes. 

Additionally, we must consider for $k=4$ the diagrams, where the
negative crossings split into $\sim$-equivalence classes $(3,1)$.
We take all genus 2 generators, apply type A flypes, make positive,
switch two crossings, and test $\sg\le 2$. Then we apply one negative
$\bt$-twist, and test $\sg\le 0$.

Let $\cD$ be the set of diagrams that we obtained so far.
Our understanding is now that every diagram we need to
consider is obtained from those in $\cD$ by positive $\bt$-twists
and type B flypes. For each diagram $D$ in $\cD$, we then
replace $D_*$ in the context of \eqref{*77} by $D_+$. 
(This then removes the diagrams with trivial clasps.)
We discard all diagrams in $\cD_+$ admitting a reducing wave
move or having $\sg>0$.

Finally, we apply Jones polynomial regularization (lemma \ref{lrez})
on the remaining set $\cE$. Here still $(\cE_n)_*$ must be
used (an \em{not} $(\cE_n)_+$). However, now crossings of either
sign may occur in a $\sim$-equivalence class, so the set $C_X$ in
\eqref{*77} must be understood to contain one positive crossing
in each $\sim$-equivalence class that has a positive crossing.
\qed

\subsection{Examples}


\begin{exam}\label{xzq}
It is to be expected that the simplification of unknot diagrams
grows considerably more complicated when the number of negative
crossings increases. We give a few examples illustrating this.

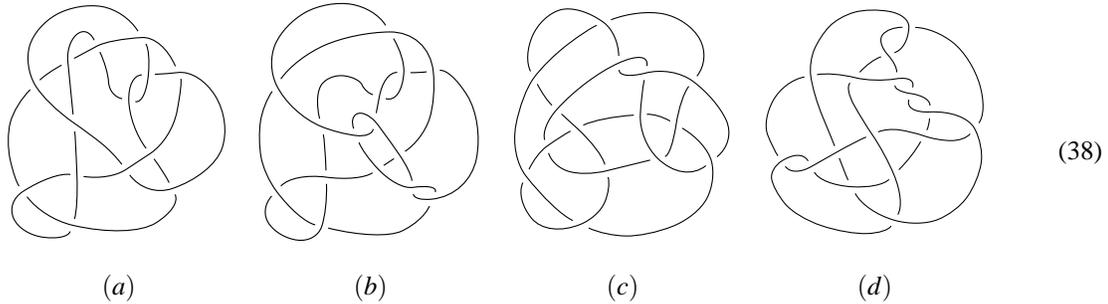
\begin{eqn}\label{Zb}
\begin{array}{*4c}
\epsfsv{3cm}{t1-unkndg} & \epsfsv{3cm}{t1-unkndg2} &
\epsfsv{3cm}{t1-unkndg4} & \epsfsv{3cm}{t1-unkndg3} \\[1.8cm]
(a) & (b) & (c) & (d)
\end{array}
\end{eqn}

The 5-almost positive 14 crossing diagram (a) does not allow
a reducing wave move, and so theorem \reference{th1.6} fails
in this form for $k=5$. Our example allows a reducing wave
move after a flype. We checked then that up to 18 crossings,
all 5-almost positive unknot diagrams still were found to
admit a reducing wave move after flypes, so that at least
proposition \reference{pz1} might be true for $k=5$. For
6-almost positive diagrams, even flypes are not enough.
The 15 crossing diagram (b) admits a flype, but no wave
move applies in any of the two diagrams obtained by flypes.
The diagram (c) is similar, though now no flypes apply either.
The 16 crossing diagram (d) is like (c) but is additionally
special.
\end{exam}


{}From the point of view of alternation, we have from theorem
\ref{th1.6}:

\begin{corr}\label{cr!}
For $k\le 4$, a special $k$-almost alternating unknot diagram is
trivializable by crossing number reducing wave moves and factor
slices. \qed
\end{corr}

\hbox to \textwidth{\parbox{12.4cm}{
\begin{exam}\label{z1}
We made again some experiments to gain evidence as to how (in)essential
is the speciality assumption is. Tsukamoto's theorem \ref{TTs}
for $k=1$ motivates that flypes must be included. We found that
all $2$-almost alternating unknot diagrams up to 18 crossings
admit, after possible flypes, a reducing wave move, and in fact
such a move that leads \em{to a $2$-almost alternating diagram}.
In contrast, there is a 3-almost alternating 17 crossing
diagram (on the right) which does not admit neither a flype
nor a (whatsoever) reducing wave move.
\end{exam}
}\hss\hbox{
\epsfsv{3cm}{t1-17_3aal} 
}\hss}

As already highlighted, both corollary \ref{cr!} and the above
computation deserve to be accompained by the warning that wave
moves do not (in general) preserve the property $\le k$-almost
alternating. (In particular, we did not claim in the corollary
that diagrams obtained intermediately after moves during the
process of reduction are either special or $\le 4$-almost
alternating.) Thus there seems some meaning in asking what moves
in $\le k$-almost alternating diagrams should naturally play the
role of wave moves in $\le k$-almost positive diagrams.

A common feature of these examples is that still one can simplify them
by (reducing) wave moves, after first performing some crossing number
preserving ones. There seems some folklore belief that this might be
true in general.

\begin{conj}\label{cu}
Crossing number preserving and crossing number reducing wave
moves in combination suffice to trivialize unknot diagrams.
\end{conj}

I had only an imprecise source for this conjecture, which also
reports that it was checked up to relatively large crossing
number (above 20). Within my computing capacity, I was able to
confirm it up to 17 crossings.
We explained already that (since factor sliding is a
preserving wave move), theorem \reference{th1.6} in particular
settles the conjecture for $\le 4$-almost positive diagrams.

\begin{rem}\label{rDn}
We do not know if factor sliding is essential in theorem
\ref{th1.6}. The calculation we did was based on the Dowker
notation \cite{DT}. It would make a treatment of composite
diagrams indefinitely more complicated, too much for the
minor (potential) benefit of removing the factor sliding.
\end{rem}

\subsection{Non-triviality of skein and Jones polynomial}

\proof[of theorem \ref{th1.7}]
The non-triviality result for the skein polynomial is almost
immadiate from the proof of proposition \reference{pz1}.
Assume $P=1$. Then in particular $V=1$, so $K$ is not 
semiadequate by theorem \ref{Vn1}, and $\Dl=1$, so $\sg=0$. Also
by $\md_lP=0$ and Morton's inequalities a $\le 4$-almost 
positive diagram $D$ must have $g(D)\le 4$. The rest of the
argument and calculation is as in the proof of proposition
\reference{pz1}. Here the semiadequacy shortcut is allowed,
but not the Rudolph-Bennequin shortcut. This deals with the
statement in theorem \ref{th1.7} about the skein polynomial.

We make now more effort for the more interesting problem
regarding $V$.
%
%
The cases $k\le 1$ are well-known. First assume
$k=2$. If $V=1$, then by \cite{restr} a $2$-almost positive
diagram $D$ has $g(D)\le 3$. Also $V=1$ implies $\dt(L)=1$,
so that $8\mid \sg(L)$, and using $g(L)\le g(D)\le 3$, we
have again $\sg=0$. 

If $D$ is composite, then one is easily out using the
result for almost positive knots in \cite{restr}.
(Note that when $\Md V>0$, then $V$ is not a unit in $\bZ[t^{\pm 1}]$;
see \cite[\S 12]{Jones2}.)
If $D$ has genus $1$, then $K$ is a rational or pretzel
knot, so a Montesinos knot. By \cite{LickThis}, it is
semiadequate, and thus by theorem \reference{Vn1}, $V\ne 1$.

So start with positive generators $D'$ of genus 
2 or 3 and $\sg\le 4$, switch two crossings to find 
$\hD$, and check $\sg(\hD)\le 0$ and that $\hD$
is not semiadequate. We obtain again a set $\hat\cD$.
We construct again $\cD=(\hat\cD)_*$ and discard
semiadequate, positive signature or wave move admitting
diagrams $X$ therein. Then we obtain again $\cE$ and 
apply a $\Md V$ regularization. 

For the rest of the proof we consider $k=3$.
By \cite{restr} we have $g(D)\le 5$. Again the composite
case can be easily reduced to the prime one, so assume $D$ is prime.

\begin{figure}[htb]
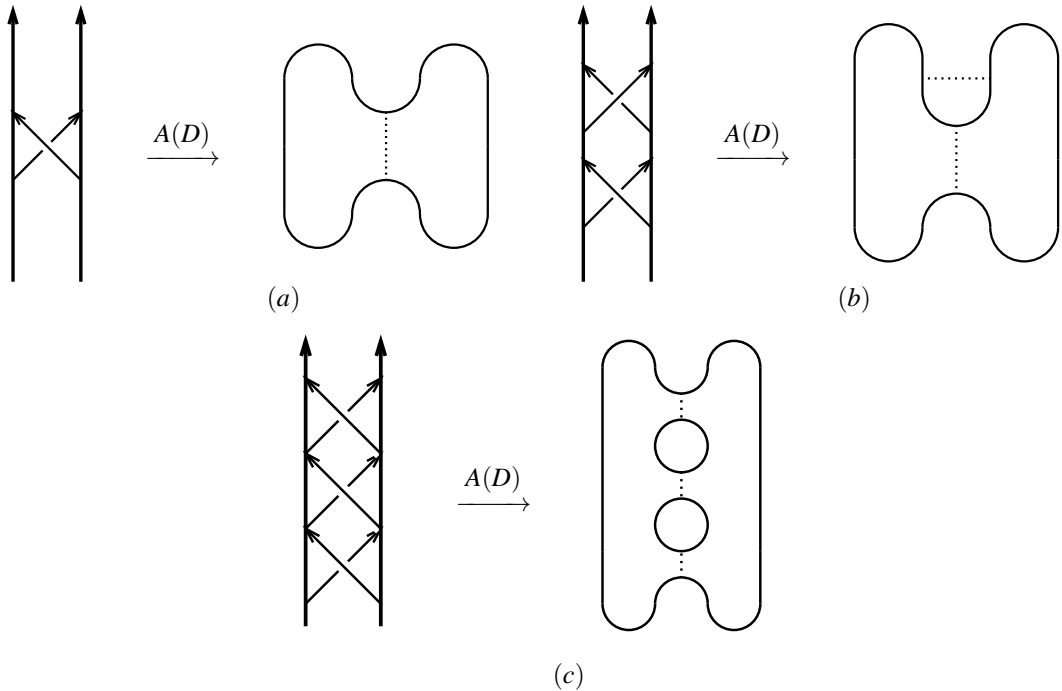

\[
\begin{array}{cc}
\diag{9mm}{8}{4}{
  \pictranslate{0 1.5}{
    \picmultivecline{0.18 1 -1.0 0}{0 0}{1 1}
    \picmultivecline{0.18 1 -1.0 0}{1 0}{0 1}
  }
  {\piclinewidth{20}
   \picvecline{0 0}{0 4}
   \picvecline{1 0}{1 4}
  }
  \picputtext{2.5 1.8}{\rato{1cm}}
  \picputtext{2.5 2.2}{$A(D)$}
  \pictranslate{5.5 2}{
    \picmultigraphics[S]{2}{1 -1}{
      \picmultigraphics[S]{2}{-1 1}{
	\picline{1.5 0}{1.5 1}
	\piccirclearc{1.0 1}{0.5}{0 180}
	\piccirclearc{0 1}{0.5}{270 0}
      }
    }
    {\piclinedash{0.03 0.06}{0.01}
     \picline{0 0.5}{0 -0.5}
    }
  }
}
&
\diag{9mm}{8}{4}{
  \pictranslate{0 0.8}{
    \picmultivecline{0.18 1 -1.0 0}{0 0}{1 1}
    \picmultivecline{0.18 1 -1.0 0}{1 0}{0 1}
  }
  \pictranslate{0 2.2}{
    \picmultivecline{0.18 1 -1.0 0}{1 0}{0 1}
    \picmultivecline{0.18 1 -1.0 0}{0 0}{1 1}
  }
  {\piclinewidth{20}
   \picvecline{0 0}{0 4}
   \picvecline{1 0}{1 4}
  }
  \picputtext{2.5 1.8}{\rato{1cm}}
  \picputtext{2.5 2.2}{$A(D)$}
  \pictranslate{5.5 1.8}{
    \picmultigraphics[S]{2}{1 -1 x}{
      \picscale{-1 1 x}{
	\picline{1.5 0}{1.5 1.0}
	\piccirclearc{1.0 1}{0.5}{0 180}
	\piccirclearc{0 1}{0.5}{270 0}
      }
      \picline{1.5 0}{1.5 1.5}
      \piccirclearc{1.0 1.5}{0.5}{0 180}
      \picline{0.5 1.5}{0.5 1.0}
      \piccirclearc{0 1}{0.5}{270 0}
    }
    {\piclinedash{0.03 0.07}{0.01}
     \picline{0 0.5}{0 -0.5}
     \picline{1.2 0.5 x}{1.2 -0.5 x}
    }
  }
} \\
(a) & (b)
\\[3mm]
\multicolumn{2}{c}{
\diag{1cm}{7}{4}{
  \picmultigraphics{3}{0 1}{
    \pictranslate{0 0.5}{
      \picvecline{0 0}{1 1}
      \picmultivecline{0.18 1 -1.0 0}{1 0}{0 1}
    }
  }
  {\piclinewidth{20}
   \picvecline{0 0.2}{0 4}
   \picvecline{1 0.2}{1 4}
  }
  \picputtext{2.5 1.8}{\rato{1cm}}
  \picputtext{2.5 2.2}{$A(D)$}
  \pictranslate{5.0 1.2}{
   \picscale{0.7 0.7}{
    \picmultigraphics[S]{2}{1 -1 x}{
      \picscale{-1 1 x}{
	\picline{1.5 0}{1.5 1.0}
	\piccirclearc{1.0 1}{0.5}{0 180}
	\piccirclearc{0 1}{0.5}{270 0}
      }
      \picline{1.5 0}{1.5 3.5}
      \piccirclearc{1.0 3.5}{0.5}{0 180}
      \piccirclearc{0 3.5}{0.5}{270 0}
    }
    \piccircle{0 2.0}{0.5}{}
    \piccircle{0 0.5}{0.5}{}
    \pictranslate{0 0.5}{
      \piclinedash{0.03 0.07}{0.01}
      \picline{0 2.5}{0 2.0}
      \picline{0 1.0}{0 0.5}
      \picline{0 -1.0}{0 -0.5}
    }
   }
  }
}} \\
\multicolumn{2}{c}{\ry{1.5em}(c)}
\end{array}
\]
\caption{\label{figAD}Some possible loops in the $A$-state of $D$.
The left sides display the location of crossings for a pair of
Seifert circles (thickened), and the right part how the loops
are connected in $A(D)$.}
\end{figure}

\begin{caselist}
\case 
Clearly we must get disposed of the (computationally
impossible) case $g(D)=5$ first.

If the 3 negative
crossings are divided into Seifert-equivalence
classes 1-1-1, then the situation in part (a) of
figure \reference{figAD} occurs 3 times; possibly the
loops on the right of (a) are connected to one.
(Here we must take into account that the Seifert graph
of $D$ is bipartite, so any cycle of non-parallel
edges has length at least 4.)
So the $A$-state of $D$ has $|A(D)|=s(D)-3$ loops. 
Converting this to $V$ using \eqref{XA},
we see that it contributes in degree $m(D)=g(D)-3=2$.
So $\md V\ge 2$, and $V\ne 1$.

If the 3 negative crossings are divided 2-1 in Seifert-equivalence
classes, with the same argument and $|A(D)|=s(D)-1$, we have
$\md V\ge m(D)=1$.

Then all 3 negative crossings are in a common Seifert-equivalence
class. So we have a loop in $A(D)$ like in part (c) of figure
\reference{figAD}. Now if there are no more (positive) crossings in
that Seifert-equivalence class, $D$ is $A$-adequate. If there are
such crossings, then in the $A$-state of $D$ there are isolated
traces connecting the same loops. By the analysis of \cite{Bmo} in
lemma \reference{lBM}, we see that the coefficient vanishes in the
degree in which the $A$-state contributes. But since $|A(D)|=s(D)+1$,
this degree is $m(D)=g(D)-5=0$. So $\md V>0$. This finishes the
case $g(D)=5$.

\case 
Let $g(D)=4$.
Again one cannot have the 3 negative crossings divided into
Seifert-equivalence classes 1-1-1, because then $\md V\ge m(D)=g(D)-3=1$.

\begin{caselist}
\case 
$D$ is special.

If all 3 negative crossings are Seifert-equivalent, then either
the diagram is $A$-adequate (if no further crossings are
Seifert-equivalent), or simplifies after a possible
flype (if there are such crossings). 

So assume the 3 negative crossings are divided into
Seifert-equivalence classes 2-1. Now $D$ would simplify
if there is a positive Seifert-equivalent crossing to
a negative one. Otherwise the class that has the single negative
crossing $p$ must be trivial. Then in the $A$-state of $D$ the
crossing $p$ leaves an isolated self-trace (as in part (a)
of figure \reference{figAD}).
By lemma \reference{lBM}, we see again that the coefficient of
the $A$-state vanishes. But since $|A(D)|=s(D)-1$, this degree is
$m(D)=g(D)-4=0$. So $\md V>0$. This finishes the case $g(D)=4$, $D$
special.

\case
$D$ is non-special.

\begin{caselist}
\case 
The 3 negative crossings are Seifert-equivalent. Since we can assume
$D$ is not $A$-adequate, we must have a positive Seifert-equivalent
crossing. Now the canonical surface of $D$ has then a compressible
Murasugi summand. By Gabai's work \cite{Gabai2,Gabai3}, this surface
is not of minimal genus, so $g(L)\le 3$. Now again, since $8\mid
\sg$, the condition $\sg=0$ emerges, and we can use a $\sg$ test.

We need to start only with positivized generator diagrams $D'$ for
which $\sg(D')\le 6$. (This drops drastically the work~-- see the
proof of proposition \reference{45657}.) To obtain $\hD$, we switch
3 crossings in a Seifert-equivalence class of at least 4, and check
if $\sg(\hD)\le 0$. (Since $\sg$ is unaltered under mutation, one
diagram per generator knot will do.) 10 diagrams $\hD$ remain; after
discarding semiadequate diagrams 4 remain. Then we build $\cX_*$ as
in \eqref{*77} for the set $\cX$ of these 4 diagrams. After discarding
diagrams in $\cX_*$ admitting a wave
move and such of $\sg>0$, the diagram set becomes empty.

\case 
The 3 negative crossings are in Seifert-equivalence classes 2-1. The
degree $m(D)$ in \eqref{XA} of the $A$ state is again $0$, and by
lemma \ref{lBM}, we see that, for $\md V=0$, there must be positive
crossings Seifert-equivalent to the single negative crossing in the one
class (as in part (b) of figure \reference{figAD}), and no positive
crossings Seifert-equivalent to the pair
of negative crossings in the other class.

Again Gabai and $8\mid \sg$ imply that $\sg=0$. So we need only
the positivized generators with $\sg(D')\le 6$. It turned out that
after switch of the pair of crossings in a Seifert-equivalence class
of $2$, we never had $\sg\le 2$. In particular we cannot build a
$\hD$ by one more crossing switch so that $\sg(\hD)\le 0$. So we
are done.
\end{caselist}
\end{caselist}

\case $g(D)\le 3$. Now $8\mid \sg$ again implies $\sg=0$.
So we could proceed as for theorem \reference{th1.6},
just avoiding the use of lemma \reference{lm2}. \qed
\end{caselist}

\begin{rem}
For $k=2$ and the Alexander polynomial a similar result
follows from \cite{2apos}, because we showed (as likewise
announced, but not given account on by Przytycki) that
2-almost positive knots of zero signature are only the twist
knots (of even crossing number with a negative clasp).
Later we will derive this result in a far more elegant
way from our setting, using some of the work in \cite{gsigex}.
\end{rem}

\begin{exam}
The $(-3,5,7)$-pretzel knot is an example of a knot with 
trivial Alexander polynomial, which is (by the preceding remark)
3-almost positive. We do not know whether it is the only one.
(For this particular case one can, again, use the methods
of this paper to easily verify 3-almost positivity
directly, rather than appealing to \cite{2apos}.)
\end{exam}

\begin{corr}
If $L$ is a $\le 4$-almost positive smoothly slice knot, then
$V(L)\ne 1$.
\end{corr}

\proof Let $D$ be a $\le 4$-almost positive diagram of $K$. Theorem
\ref{th1.7} easily deals with the case that $D$ is composite
(because, as noted it its proof, if $V\ne 1$, then $V$ is not a
unit). Considering $D$ to be prime, we have $\sg=0$, and
lemma \ref{lm7.} implies that $g(D)\le 4$. Then a calculation
like in the proof of theorem \reference{th1.7} applies. \qed

In special diagrams, one can extend theorem \ref{th1.6} to $k=5$,
if one allows for flypes, and similarly the part of theorem
\ref{th1.7} for the skein polynomial.

\begin{theo}\label{t5p}
Let $D$ be a $k$-almost special alternating knot diagram, for
$k\le 5$. Then the following conditions are equivalent:
\def\labelenumi{(\arabic{enumi})}
\def\theenumi{(\arabic{enumi})}
\begin{enumerate}
\item $D$ is unknotted,
\item $P(D)=1$, \label{z2}
\item $D$ can be trivialized by reducing wave moves and flypes.
  \label{z3}
\end{enumerate}
\end{theo}

Remembering the diagram (d) in \eqref{Zb}, we see that the theorem
is not true for $k=6$. On the other hand, for $k=5$ it might
we true for $k$-almost positive diagrams (without the assumption
the diagram to be special), as explained in example \ref{xzq}.
The warning formulated below example \ref{z1} is valid here, too.

\proof Clearly we need to prove only \ref{z2}$\ \So\ $\ref{z3}.
We may further assume that $k=5$, since the cases $k\le 4$ are
contained in theorems \ref{th1.6} and \ref{th1.7}. Next we may
asume thet $D$ is prime, because (again, for example using
\eqref{PtoV} and properties of $V$), if $P(D)$ is a unit in
$\bZ[l^{\pm 1},m^{\pm 1}]$, then $P(D)=1$.

Since $k=c_-(D)=5$, we have $g(D)\le 5$ by \eqref{mq}.
Next, $P=1$ implies $\Dl=1$,
which in turn implies $\sg=0$, and the signature test is allowed.
Likewise is the semiadequact test, since $P=1$ implies $V=1$.

Assume $g(D)=5$. Since $D$ is special, by lemma \reference{l45},
each $\ssim$-equivalence class is signed, or we can simplify
$D$ after a flype. By \eqref{q2}, since $\md_l P(D)=1$ and
$k=c_-(D)=5$, we must have $\inx_-(D)=0$. This means that there
is no negative $\ssim$-equivalence class of a single crossing.
Thus the 5 negative crossings in $D$ must lie either in a single
negative class, or split into two such classes 2-3. In either
case $D$ is $A$-adequate, and so $V(D)\ne 1$ by theorem \ref{Vn1}.

Then we have $g(D)\le 4$, and can use the generator table. If
$g(D)=1$, then by theorem \ref{thge}, the claim can be checked
easily. (All such $D$ are semiadequate, except if they have
a trivial clasp.)

We consider thus $g(D)=2,3,4$. The rest of the proof is as for
proposition \ref{pz1}, though with a fifth negative crossing.
It may be necessary (and thus we do) consider non-special
diagrams. Again, at most one negative $\bt$-twist must we
dealt with. In the case of such a twist, here we must test
all the generators of genus 2 and 3, and those of genus 4
whose positification has $\sg\le 6$. These cases are easily
finished.

We consider the case of no negative $\bt$-twist. We proceed
as before, until we build again $\cD=(\hat\cD)_*$ (and not
$(\hat\cD)_+$). Next, to obtain $\cE$, we discard now not only
diagrams $D\in \cD$ admitting a wave move, but also such that
do so after a flype. (Our understanding is that we apply $\bt$
moves only on non-trivial $\sim$-equivalence classes, so that
we may use also type A flypes.)

Regularization of $V$ applies with $n=6,7,8$, but 
we obtain now for $n=7,8$ a total of 330 diagrams $E$ in $\cE$
on which regularization fails. We can fix this, though, by
restricting the condition \eqref{rez} to diagrams $E\in \cE$
(in our current notation) with $\sg(E)\le 0$. Then \eqref{rez}
is always satisfied.
\qed

The calculation is now longer than for proposition \ref{pz1}, due
to the extra negative crossing. By \eqref{q2} and semiadequacy, we
have restrictions also when $g(D)\le 4$. We left many of these tests
out, though, adopting a simplistic attitude. As soon as the calculation
became manageable (say, within a couple of days), we made very
little effort to speed it further up by additional tests. We were
aware that every single new step put into our procedure augments
the risk of error. And when all diagrams can be verified, then
certainly a subset can be, too, without that we get aware of
overlooking the others. Indeed, a number of iterations were
necessary in order to have the single steps working well and in
the correct order.

\subsection{On the number of unknotting Reidemeister moves\label{sNR}}

There has been recently some interest in the literature to estimate
the minimal number of Reidemeister moves needed for turning an unknot
diagram $D$ of $n$ crossings into the trivial (0-crossing) diagram
(see e.g. \cite{HL,HN,Hayashi}). Let us call this quantity here $r(D)$.

It was long known from Haken theory that $r(D)$ should be estimable
from above, and in \cite{HL} this was made explicit, by showing that
\begin{eqn}\label{cHL}
r(D)\,\le\,O(C^n)\,,
\end{eqn}
where $n=c(D)$, and $O(\dots)$ means `at most \dots\,times a constant
independent on $n$'. Unfortunately, the exponential base $C$ is a number
with billions of digits, which renders the estimate quite impracticable.

Our above results imply the following statement, which can be
regarded also as a continuation to proposition \reference{Za}:

\begin{prop}\label{pRei}
If $D$ is an $n$ crossing $k$-almost positive unknot diagram for
$k\le 4$ or $k$-almost special alternating unknot diagram for
$k\le 5$, then $r(D)=O(n^p)$, where $p$ depends only on $k$.
\end{prop}

This is in sharp contrast to the Hass-Lagarias bound in the general
case. Since in proposition \ref{pRei} we consider only small $k$,
we could take $p$ also to be a global constant. Underscoring the
dependence on $k$ is made in order to indicate that our approach
could work also for higher $k$, provided proposition \ref{pz1} or
theorem \ref{t5p} can be extended. The values of $p$ we will obtain
for given $k$ (the largest of which is $144$) can likely be lowered
by further direct calculation (as in proposition \ref{Za} when $k\le
2$). However, we opted for a more conceptual argument, formulating
three lemmas below that have also some independent meaning.

For many diagrams $D$ we have $r(D)=O(n)$. Only recently some
diagrams $D$ were found where $r(D)$ is quadratic in $n$ \cite{HN}.
Even although this
still leaves a large gap between lower and upper bounds, we do not
believe a polynomial upper bound is possible in general. We should
note, though, that any set of (not-too-fancy) crossing number
non-augmenting unknotting moves will yield much more tangible
estimates in \eqref{cHL}. For example, would conjecture \ref{cu}
be true, the below lemmas \ref{lcw} and \ref{lfw}, together with
the result in \cite{nlpol}, would justify the constant $C=10.399$
in \eqref{cHL}. See also remark \reference{yuuu}.

Note also that, by lemma \reference{lcw}, the difference between
studying Reidemeister moves in $S^2$ or $\bR^2$ is immaterial from
the point of view of proposition \reference{pRei} (it affects $p$
at most by $\pm 1$). 

Proposition \ref{pRei} follows from the below three lemmas.

\begin{lemma}\label{lnd}
For fixed $k$, the number of $k$-almost positive unknot diagrams of $n$
crossings is polynomial in $n$.
\end{lemma}

\proof 
There are ${n\choose k}=O(n^k)$ choices of negative crossings in an
$k$-almost positive diagram, once the underlying alternating diagram 
$\hD$ is given.

Moreover $g=g(D)=g(\hD)\le k$, if $D$ is unknotted. Thus by theorem
\ref{th1gen}, $D$ has at most $6g-3$ $\sim$-equivalence classes.
Next, by lemma \ref{lmda}, each $\sim$-equivalence class has at most
$g$ twist equivalence classes, so that $\hD$ has at most $g(6g-3)$
twist equivalence classes. 

Thus the number of $\hD$ (regarded in $S^2$) of $n$ crossings is at
most $O(n^{g(6g-3)-1})$, with $g\le k$, and so the number of $D$ is at
most $O(n^{k(6k-3)-1+k})$.
\qed

\begin{lemma}\label{lcw}
A wave move in an $n$-crossing diagram can be realized by $O(n)$
Reidemeister moves. \qed
\end{lemma}

\begin{lemma}\label{lfw}
A flype in an $n$-crossing diagram can be realized by $O(n^4)$
Reidemeister moves. 
\end{lemma}

\proof It is enough to prove the claim for the tangle isotopy
in (a) below:
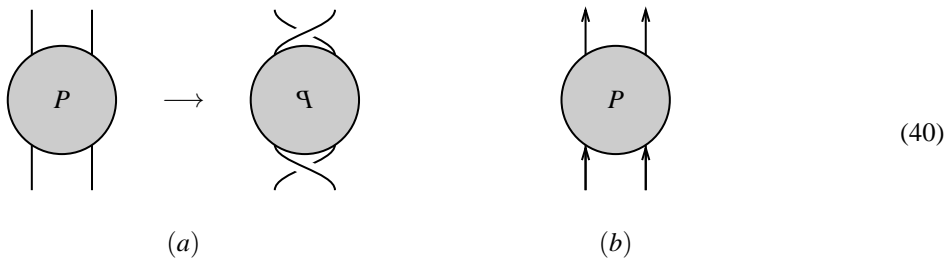
\begin{eqn}\label{tis}
\begin{array}{c@{\kern1.7cm}c}
  \diag{8mm}{3}{3}{
  \picfillgraycol{0.8}
  \picline{1 0}{1 3}
  \picline{2 0}{2 3}
  \picfilledcircle{1.5 1.5}{0.9}{$P$}
}\,\lra\,
\diag{8mm}{3}{3}{
  \picfillgraycol{0.8}
  \rbraid{1.5 0.38}{1 0.76}
  \lbraid{1.5 2.62}{1 0.76}
  \picscale{-1 1}{
    \picfilledcircle{-1.5 1.5}{0.9}{$P$}
  }
} & 
\diag{8mm}{3}{3}{
  \picfillgraycol{0.8}
  \picvecline{1 0}{1 0.7}
  \picvecline{1 0}{1 3}
  \picvecline{2 0}{2 0.7}
  \picvecline{2 0}{2 3}
  \picfilledcircle{1.5 1.5}{0.9}{$P$}
} \\[15mm]
(a) & (b)
\end{array}
\end{eqn}
Here the second tangle is obtained by a flip (rotation by
$180^\circ$) on the vertical axis.

For any tangle $P$, there are two ways to group its 4 ends into two
neighbored pairs. The results of the flips with respect to these
two choices differ by two wave moves. Thus, by lemma \ref{lcw},
we may pick any of these two choices, and with this freedom we
can then orient the tangle strands so that the orientation of the
ends on the left of \eqref[a]{tis} is like in \eqref[b]{tis}.

Now we would like to apply the procedure of Vogel \cite{Vogel}
outside a trivial tangle (a) below,
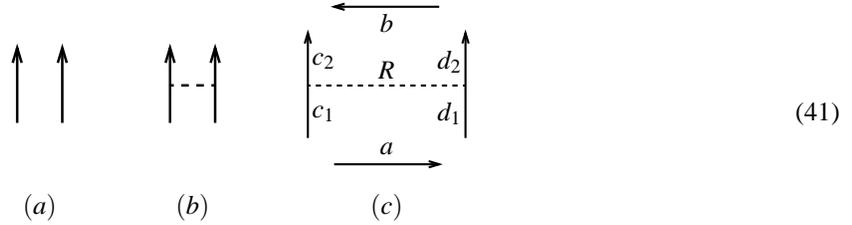
\begin{eqn}\label{hz0}
\begin{array}{c@{\kern1cm}c@{\kern1cm}c}
\diag{1.0cm}{1}{1}{
  \picvecline{0.2 0}{0.2 1}
  \picvecline{0.8 0}{0.8 1}
} &
\diag{1.0cm}{1}{1}{
  \picvecline{0.2 0}{0.2 1}
  \picvecline{0.8 0}{0.8 1}
  {\piclinedash{0.1 0.1}{0.05}
   \picline{0.2 0.5}{0.8 0.5}
  }
} &
\diag{7mm}{3}{3}{
  \piclinewidth{25}
  \picvecline{0 0.5}{0 2.5}
  \picvecline{3 0.5}{3 2.5}
  {\piclinedash{0.1 0.1}{0.05}
   \picline{0 1.5}{3 1.5}
  }
  \picvecline{0.5 0}{2.5 0}
  \picvecline{2.5 3}{0.5 3}
  \picputtext{1.5 0.3}{$a$}
  \picputtext{1.5 2.7}{$b$}
  \picputtext{0.3 2.0}{$c_2$}
  \picputtext{0.3 1.0}{$c_1$}
  \picputtext{2.7 2.0}{$d_2$}
  \picputtext{2.7 1.0}{$d_1$}
  \picputtext{1.5 1.8}{$R$}
} \\[12mm]
(a) & (b) & (c)
\end{array}
\end{eqn}
complementary to \eqref[b]{tis}, to obtain a diagram of $P$
as a partially closed braid $\be$. (We will use below the
terminology of \S\reference{SBR} for braids.)
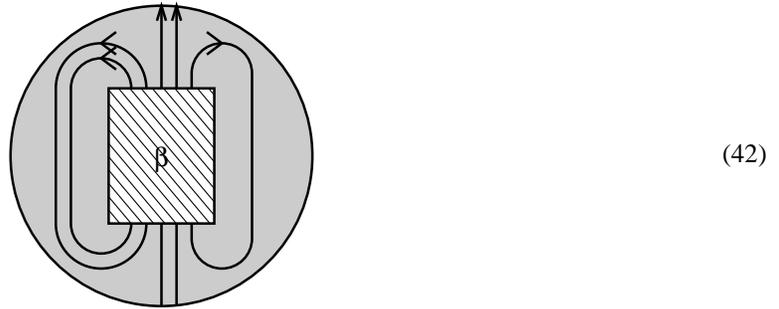
\begin{eqn}\label{hz}
\diag{1cm}{5}{5}{
  \pictranslate{2.5 2.5}{
    {\picfillgraycol{0.8}
     \picfilledcircle{0 0}{2}{}
    }
    \picovalbox{-0.8 0}{1.2 3.0}{0.6}{}
    \picovalbox{-0.8 0}{0.8 2.6}{0.4}{}
    \picvecline{0 -2}{0 2}
    \picvecline{0.2 -2}{0.2 2}
    \picvecwidth{0.15}
    \picovalbox{0.8 0}{0.8 3}{0.4}{}
    {\piclinewidth{70}
     \picshade{
       \picfilledbox{0 0}{1.4 1.8}{}
     }{40}{-8}
    }
    \picbox{0 0}{1.4 1.8}{$\be$}
    \picvecline{-0.8 1.5}{-0.801 1.5}
    \picvecline{-0.8 1.3}{-0.801 1.3}
    \picvecline{0.8 1.5}{0.801 1.5}
  }
}
\end{eqn}

The Vogel move is a Reidemeister II move of a special type.
It affects two edges in the boundary of a region $R$, which
have the \em{same orientation} as seen from inside $R$, and
belong to \em{distinct Seifert circles}, and creates a
reverse trivial clasp inside $R$.

In order to realize these moves in the tangle complementary to
\eqref[a]{hz0}, we draw a dashed line in \eqref[a]{hz0}, as in 
\eqref[b]{hz0}, and need to be concerned only with the Vogel
moves that affect this line.

Let $R$ be the region containing the line. The two arrows in \eqref[b]
{hz0} are oriented oppositely as seen from inside $R$, and it is easy to
see that for orientation reasons, these arrows must belong to distinct
Seifert circles. 
If now the pair $(a,b)$ of segments in \eqref[c]{hz0} allows for
a Vogel move, then
at least one of the four pairs $(a,c_1)$, $(a,d_1)$, $(b,c_2)$ and
$(b,d_2)$ does also. So one can find a Vogel move that avoids the
dashed line.

Vogel proves that, no matter what choice of move is taken at each
step, after $O(n^2)$ applications no further moves are possible, and
then the diagram is in braid form.

So we achieved the form \eqref{hz} for $P$ by $O(n^2)$ Reidemeister
moves. Since the number of Seifert circles is not changed, $\be$ has
$l=O(n)$ strands and $m=O(n^2)$ crossings. (All strands of $\be$ are
closed except $2$. There is no evidence how many strands are closed
on the left and on the right.)

The reverse procedure can be applied after the flip with the same cost,
so it is enough to look from now on at the case that $P$ is in the
form \eqref{hz}.

The flip of $P$ now transcribes into the conjugation of $\be$ with
the half-twist $\dl\in B_l$ (whose square is generating the
center of $B_l$), plus the interchange of the strands
closed on the left and on the right in \eqref{hz}. Latter step can
be accomplished by $O(n^3)$ moves, since by lemma \ref{lcw},
moving one strand between either side takes $O(m)=O(n^2)$ moves.

It remains to look at the flip of the braid $\be$, i.e.\ the
conjugation with $\dl$. With 
\[
\gm_k=\sg_1\dots\sg_{k}\,,
\]
we can write
\begin{eqn}\label{sii0}
  \dl\,=\,\gm_{l-1}\dots\gm_1\,.
\end{eqn}

To simplify things slightly, let us introduce the \em{bands}
\begin{eqn}\label{s_}
\sg_{i,j}^{\pm 1}\,=\,\sg_i\dots\sg_{j-1}\sg_{j}^{\pm 1}
  \sg_{j-1}^{-1}\dots\sg_i^{-1}\,,
\end{eqn}
for $1\le i<j\le l$, so that
\begin{eqn}\label{siij}
\sg_i=\sg_{i,i+1}\,.
\end{eqn}
Let 
\[
\pi\,:\,B_l\to S_l\,,\qquad\,\pi(\sg_i)=\tau_i\,=\,(i\ i+1)\,
\]
be the permutation homomorphism. It now is easy to check that,
\em{unless} $i<k+1<j$,
\begin{eqn}\label{sii1}
\gm_k\sg_{i,j}\,=\,\sg_{\tau(i),\tau(j)}\gm_k\,,
\end{eqn}
with $\tau=\pi(\gm_k)^{-1}$, and that this identity is
accomplished by $O(l)$ YB relations (recall \S\ref{SBR}).

Now, write $\be$ as a word in $\sg_{i,j}$ using \eqref{siij}.
Next, replace this word by $\dl^{-1}\dl\be$ by $O(l^2)=O(n^2)$
Reidemeister II moves, using \eqref{sii0} to express $\dl$.
Then, by using \eqref{sii1}, one can ``pull'' each $\gm_k$
through $\be$ (and the situation $i<k+1<j$ never occurs!).
Thus by $O(m\cdot l)=O(n^3)$ applications of \eqref{sii1}
one changes $\dl\be$ to $\bar\be\dl$ (with $\bar\be$ being
$\be$ with $\sg_i$ changed to $\sg_{l-i}$), and this requires
then $O(n^3\cdot l)=O(n^4)$ Reidemeister moves.

After that, there are two half-twists of opposite sign on the
left and right closed strands, which cancel by $O(n^2)$ moves,
and the interchange of the left and right closed strands, which,
as argued, costs $O(n^3)$ moves.
\qed

\begin{rem}
It appears that the braid form is not essential, and in the statement
of the lemma possibly $O(n^4)$ could be made into $O(n^3)$. It seems,
though, very awkward to write a concrete procedure down without using
braids.
\end{rem}

\proof[of proposition \ref{pRei}]
With proposition \ref{pz1}, theorem \ref{t5p}, and lemmas \ref{lnd}
to \ref{lfw}, after switching between $O(n^{6k^2-2k-1})$ different
$k$-almost positive unknot diagrams, each time using $O(n^4)$, so
totally $O(n^{6k^2-2k+3})$, Reidemeister moves, we can reduce the
crossing number. Then this procedure can be iterated. So we need at
most $O(n^{6k^2-2k+4})$ moves, that is, we can set $p=6k^2-2k+4$
(for moves in $S^2$). \qed

\begin{rem}\label{yuuu}
An analogon of proposition \ref{pRei} holds for almost alternating
diagrams by theorem \ref{TTs} and using lemmas \ref{lcw} and \ref{lfw}.
Lemma \ref{lnd} is not true for almost alternating diagrems, but
its use can be avoided. While by flypes one can in general obtain
exponentially many diagrams, only a polynomial number of flypes is
needed to pass between any two diagrams. This can be seen from the
structure of flyping circuits discussed in \cite{SunThis} (see also
the proof of lemma \reference{pq}).
\end{rem}

\subsection{Achiral knot classification}

We turn, concluding this section,
to achiral knots. The proof of theorem \ref{th1.8}
divides into two main cases, for prime and composite knots.
For the prime case we have done the most work with the previous
explanation. However, for composite knots some further arguments 
are necassary, the problem being that prime factors of achiral knots
are not necessarily achiral. We will need some of the results in
\cite{gen2} and Thistlethwaite's work \cite{Thistle}.

\proof[of theorem \ref{th1.8} for prime knots]
Let $D$ be a $\le 4$-almost positive diagram and $K$ the achiral
knot it represents. We will work by induction on the crossing number
of $D$, showing that each time we can find a wave move to simplify 
$D$, unless $D$ has small crossing number. (For at most 16 crossings
the claim of the theorem was verified in \cite{restr}.) Now assuming
that $K$ 
is prime, we know that if $D$ is composite, all but one of the prime
factors of $D$ represent the unknot. Using the result for the unknot
we already proved, we can trivialize these factors, and then go over
to deal with the non-trivial one. With this argument we can w.l.o.g.
assume that $D$ is prime.

The proof now follows the one for theorem \ref{th1.6} with the following
modifications. To determine $\hD$, the Rudolph-Bennequin shortcut
is allowed. The inequality \eqref{RBi} must hold for an achiral knot
$K$ by replacing $g_s(K)$ by $0$, even if the knot $K$ is not slice.
This is because one of $K\# K$ or $K\#-K$ is slice, and the estimate
on the right is not sensitive w.r.t. orientation and additive under
(proper) connected diagram sum. The $\sg$ tests in determining $\hD$
are allowed because $\sg=0$ by achirality. 

Lemma \reference{lm1} shows that one can
apply the semiadequacy shortcut if $c_-(D)=g(D)$ also here.
So for such $D$ we are immediately done with the
verification in the proof of theorem \reference{th1.7}.
Even if $c_-(D)>g(D)$, the semiadequacy test is legitimized
partly, namely if $X\in \cD$ satisfies $\md V(X)\ge -1$.
This can be explained as follows from the work in \cite{Thistle}.

If $X$ is $B$-adequate, then it has the minimal number of
positive crossings, which is at most $4$ since $K$ is achiral and
$\le 4$-almost positive. But $X$ has also at most 4 negative crossings,
so $c(X)\le 8$, and we are easily out. If $X$ is $A$-adequate, then
positive $\bt$-twists preserve $A$-adequacy and $\md V$. So
$\md V(D)\ge -1$, but $V(D)$ is reciprocal; if $\md V\ge 0$,
then $V=1$ and contradiction to 
the semiadequacy. Then $\md V(D)= -1$; moreover, $[V(D)]_{t^{-1}}=
[V(D)]_{t^{1}}=\pm 1$ by semiadequacy. Then we use the property 
$V(1)=1$ (see \cite[\S 12]{Jones2}), which shows that only two
polynomials are possible, $-1/t+3-t$ and $1/t-1+t$. These are ruled out
using $V(i)=\pm 1$ and $V(e^{\pi i/3})\ne 0$. 

Regularization (as described at the end of the proof of
theorem \reference{th1.6}) applies similarly, and again we have
$3\le n\le 6$. Now $\Md V\ne 0$, but $\Md V\le 6$, using that
$\md V\ge -6$ by \cite{restr} and reciprocity of $V$. So regularity
shows that we can deal with all diagrams $D$, unless $c(D)\le n+\Md
V\le 12$. In that case the check of \cite{restr} for $\le 16$
crossings finishes the proof. \qed

We turn now to the two remaining knots. While $4_1\#4_1$ is easy to
deal with inductively, we will need some work for $3_1\#!3_1$.

\begin{lemma}\label{f1}
Let $K'$ be a chiral connected sum factor of a 
$\le 4$-almost positive achiral knot $K$. Then $\Md_z F(K')\le 3$.
\end{lemma}

\proof We apply \eqref{QFc}, so in our case $\le 4$-almost
positivity and achirality show $\Md_z F(K)\le 7$. Now since
$K$ is achiral but $K'$ is not, $K$ must have $!K'$ as a factor too,
i.e. $K=K'\#!K'\#L$, where $L$ is achiral. As $\Md_zF(L)\ge 0$,
we have $7\ge \Md_z F(K)\ge 2\Md_z F(K')$, and so $\Md_z F(K')\le 3$.
\qed

\begin{lemma}\label{f2}
Let $D'$ be a genus one digram occurring as a prime factor of $D$.
Assume $D'$ represents a chiral knot $K'$. Then this knot $K'$
is one of the trefoils.
\end{lemma}

\proof By \cite{gen1}, $K'$ is a rational or pretzel knot, so a
Montesinos knot. By \cite{LickThis} then we have $\Md_zF(K')\ge
c(K')-2$. By the previous lemma $c(K')-2\,\le\, \Md_z F(K')\le 3$,
so $c(K')\le 5$. Now $4_1$ is achiral, and the 5-crossing knots are
excluded because $\Md_zF=4$. So $K'$ must be one of the trefoils. \qed

\proof[of theorem \ref{th1.8} for composite knots]
Let $D$ be a $\le 4$-almost positive diagram and $K$ the 
composite achiral knot it represents. We have $g(D)\le 4$.
We assume w.l.o.g. that $D$ cannot be further simplified by
flypes and wave moves.

Let first $D$ be prime. As $D$ admits no wave move after flypes,
the proof in the prime knot case shows that $c(D)\le 12$. By
lemma \reference{f1}, we have $\Md_zF\le 7$. Then we see by 
check of the maximal degree of the $F$ polynomial of achiral
prime knots of at most 12 crossings that if $K$ contains
two achiral connected sum factors, then $K=4_1\# 4_1$. If $K$ 
has a chiral factor $K'$, then $K=K'\#!K'\#L$, where $L$ is 
achiral. Then we use lemma \ref{f1} and verify that only the 
trefoils have $\Md_z F(K')\le 3$ among connected sum factors 
of knots of $\le 12$ crossings, and that $L$ must be trivial. 

So for the rest of the proof we assume that $D$ is composite. If all 
prime factors of $D$ are achiral, then again we can deal with these
factors separately. With the argument in the preceding paragraph,
$D$ depicts some connected sum of trefoils and figure-8-knots.
Because of $\Md_zF(D)\le 7$, there remain only the options that we 
want to show. So assume below that at least one prime factor of $D$ 
depicts a chiral knot. 

\begin{caselist} %
\case
Let $D$ have first at least 3 prime factors $D=D_1\#D_2\#D_3$.
Since $4\ge g(D)=\sum g(D_i)$, we may assume w.l.o.g. that
$g(D_1)=g(D_2)=1$. Then the knots $K_{1,2}$ these diagrams
represent are trefoils (by lemma \ref{f2}). If they are same sign
trefoils, then $K$ must factor as $3_1\#3_1\#!3_1\#!3_1\#L$ for
an achiral knot $L$, and $\Md_zF\ge 8$, a contradiction to
$\Md_z F(K)\le 7$. So $K_{1,2}$ are trefoils of opposite sign.
Since $a_-(!3_1)=3$ (see \eqref{a+-}), one of $D_{1,2}$ has at
least three negative crossings. So $D_3$ is has at most one
negative crossing, but represents an achiral knot. We know,
though, that non-trivial positive or almost positive knots are
not achiral. So $K=3_1\#!3_1$, and we are done.

\case Now we consider the situation that $D=D_1\#D_2$ has two
prime factors $D_{1,2}$. Let $K_{1,2}$ be the knots of $D_{1,2}$.
We can group the mirrored chiral
prime factors of $K$ which occur on different sides of the decomposition
$D=D_1\#D_2$ into a (possibly composite) factor knot $K'$ of $K$, and
assume that $K_1=K'\#L$ and $K_2=!K'\#M$ with $L$ and
$M$ achiral (and possibly trivial) and $K'$ chiral (and non-trivial by
the assumption preceding the case distinction). 
We have by lemma \reference{f1} that $\Md_zF(K')\le 3$.

\begin{caselist} %
\case 
$D_1$ is positive, $D_2$ is $\le 4$-almost positive. 

\begin{caselist} %
\case 
$g(D_1)=1$.
Then $g(K'\#L)=1$ so $L$ is trivial, and $K'$ is a positive knot of
genus $1$. By lemma \reference{f2} it is the positive (right-hand)
trefoil. Now consider $D_2$, with $g(D_2)\le 3$.

\begin{caselist} %
\case 
$g(D_2)=1$. Then by lemma \reference{f2}, $K_2$ is a trefoil, and by 
achirality of $K$ it must be negative, so we are done.

\case \label{c112}
$g(D_2)=3$. Since $D_2$ is $\le 4$-almost positive, the left
Morton inequality in \eqref{mq}
shows $\md_lP(D_2)\ge -2$. Moreover $D_1$ is positive of genus 1, so
$\md_lP(D_1)=2$. It follows that $\md_lP(D)\ge 0$. Now $K$ is achiral,
and then \cite[proposition 21]{LickMil} shows $P(D)=1$. But the trefoil 
polynomial must divide $P(D)$, a contradiction.

\case \label{c113}
$g(D_2)=2$. We repeat the argument in the previous case. Now $\md_lP(D_2)\ge 
-4$, so $\md_lP(D)\ge -2$. The case $\md_lP(D)\ge 0$ is ruled out as before, 
so assume that $\md_lP(D)= -2$, and then $\Md_lP(D)= 2$. Again $K=K'\#!K'\#
M$, with $M$ achiral. Comparison of the $P$-degrees and \cite[proposition 21]
{LickMil} shows then that $P_M=1$. So we seek an achiral knot $M$ with trivial 
skein polynomial, whose connected sum with a negative trefoil $!K'$ is a knot 
$K_2$ in a $\le 4$-almost positive genus 2 diagram $D_2$. However, in that case 
$V(M)=1$, so $\Md V(K_2)=-1$, and also $\sg(K_2)<0$. Then, assuming w.l.o.g. 
that we reduced $D_2$ by flypes and wave moves, we checked in the proof of 
theorem \ref{th1.6} that $c(D_2)-\Md V(D_2)=n\le 6$, so $c(D_2)\le 5$. Then 
clearly $M$ is trivial, and we are done.
\end{caselist}

\case 
$D_1$ is positive and $g(D_1)=2$. In the proof, given in \cite{gen2}, of the fact 
that positive knots of genus 2 have minimal positive diagrams, we verified that
all positive diagrams $D_1$ of genus 2 can be reduced to diagrams $D_1'$
with $\Md Q(D_1')\ge c(D_1')-2$. Now $\Md Q(D_1')\le \Md _zF(D_1')\le \Md _zF(D)\le 7$.
So $c(K_1)\le c(D_1')\le 9$, and
$K'$ is a chiral knot with $\Md_zF(K')\le 3$ occurring as a
factor of a knot $K_1$ with at most 9 crossings. Then again $K'$ is a trefoil. But
then, looking at the signature $\sg$, and using that $L$ is achiral, we have 
$\sg(D_1)=\sg(K')+\sg(L)=\pm 2$, while for a positive diagram $D_1$ of genus 2 
we showed in \cite{gen2} that $\sg=4$. This contradiction finishes 
the case.

\case 
$D_1$ is positive and $g(D_1)=3$. So $g(D_2)=1$ and by lemma \reference{f2}, $K'$ is a trefoil 
(and $M$ is trivial). Now again $\sg(D_1)\ge 4$ by \cite{apos}, while $\sg(D_2)=
\sg(K')=\pm 2$, and $\sg(D)=\sg(D_1)+\sg(D_2)>0$, a contradiction.

\end{caselist} 

\case 
$D_1$ is almost positive, $D_2$ is $\le 3$-almost positive.

\begin{caselist} %
\case 
$g(D_2)=3$. Then $g(D_1)=1$ and $K'$ is a trefoil, which must be
right-hand, since $D_1$ is almost positive. Using \eqref{mq}, we have
$\md_l P(D_2)\ge 0$, so $\md_l P(D)> 0$, and a contradiction to
achirality.

\case 
$g(D_2)=2$. By Morton's inequality
$\md_l P(D_2)\ge -2$. Also we know from \cite{gen2}
that there is no almost positive knot of genus 1, so if
$g(D_1)=1$, then $K_1$ is positive. Then by Morton $\md_l
P(D_1)\ge 2$, and $\md_l P(D)\ge 0$. So $P(D)=1$. However,
if $\md_l P(D_1)\ge 2$, then $P(D_1)$, which divides $P(D)$,
cannot be a unit in $\bZ[l^{\pm 1},m^{\pm 1}]$ by
\cite[proposition 21]{LickMil}, and we have a contradiction.

\case 
$g(D_2)=1$. Then by lemma \reference{f2}, $!K'$ is a trefoil (and $M$ is
trivial). It must be a negative trefoil, because otherwise $\sg(D_2)$ and 
$\sg(D_1)$ are both positive, and $\sg(D)=\sg(D_1)+\sg(D_2)>0$.
So $K'$ is a positive trefoil.
\begin{caselist}

\case If $g(D_1)=1$, we are easily done by lemma \reference{f2}. 

\case If $g(D_1)=3$, then by Morton's inequality $\md_l P(D)\ge 0$,
and we have
a contradiction along the above lines in case \reference{c112}. 

\case So assume $g(D_1)=2$. Now $D_1$ contains a positive trefoil
factor. By Morton, $\md_l P(D)\ge -2$, so by achirality and \cite
[proposition 21]{LickMil}, and then along the lines in case \ref{c113},
we see that $L$ must have trivial skein polynomial. So $K_1$ has the
trefoil polynomial, and an almost positive diagram $D_1$. Then by
\cite{posqp}, we have $2=\Md_m P(K_1)=2g(K_1)=2g(K')+2g(L)$. (One
can use alternatively $\Dl$ and the theorem that $\Md\Dl=1-\chi$ in 
\cite{canon}.) This implies that $L$ must be trivial, and we are done.
\end{caselist}

\end{caselist} 

\case 
$D_1$ and $D_2$ are both $2$-almost positive.

If $g(D_1)=1$, then ($L$ is trivial and) $K'$ is a chiral knot with
$\Md_z F\le 3$ and a $2$-almost positive genus 1 diagram. Up to
reducible diagrams or diagrams with trivial clasps (which we can
exclude) there are no such knots. Similarly we argue if $g(D_2)=1$.
So assume $g(D_1)=g(D_2)=2$. Then by Morton's inequality $\md_l
P(D_{1,2})\ge 0$, and the meanwhile well-known way how to obtain
a contradiction.

\end{caselist} 
\end{caselist} 

The case distinction and the proof are complete. \qed\\

\hbox to \textwidth{%
\parbox[b]{12.4cm}{

The proof of theorem \reference{th1.8} suggests that, beyond
the description of the occurring knots, we can also obtain a
partial simplification statement for their diagrams. We make
this explicit in the following way.

\begin{corr}\label{Cz}
Let $K$ be a $k$-almost positive achiral knot for $k\le 4$, and
$K\ne 3_1\# !3_1$. Then one can simplify a $\le 4$-almost positive
diagram $D$ of $K$ by reducing wave moves and factor slides either
to a minimal crossing diagram, or to the 8 crossing diagram of
$K=6_3$ shown on the right.
\end{corr}
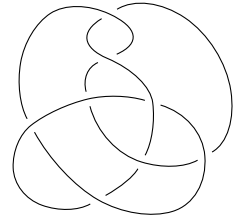
\vspace{-2em}%

\mbox{}%
}\hss
\epsfsb{3cm}{t1-6_3_8cr}%
\kern-1mm\hss}\vspace{\parskip}

\proof The exclusion of $K=3_1\# !3_1$ was made in order to spare
us the technical complications of dealing with non-amphicheiral
prime factors. Thus we can again consider w.l.o.g. prime diagrams.

Now the proof of theorem \reference{th1.6} can be carried
through, with the only difference coming in the Jones polynomial
regularization of $\cE_n$ at the very end. (As before, the cases
$c_-(D)\le 2$ and $g(D)=1$ can be checked directly.)

We can discard all alternating diagrams $E$ in $\cE$, since
$\bt$-twists preserve alternation, and these are then minimal
crossing diagrams, as we claimed. From the rest we check that 
\begin{eqn}\label{Vv1}
\Md V(E)+\md V(E)\ge -1\,.
\end{eqn}
(It is enough to consider only $E\in \cE$ with $\Md V(E)\le 6$.)

The property \eqref{Vv1}, together with the growth of $\Md V(E)$
under twists (which was the subject of regularization), implies
that $V$ will not be self-conjugate on diagrams that are obtained
by positive $\bt$-twists from diagrams $E\in \cE$.

It is thus enough to look at $E$ themselves. We identify all diagrams
that depict achiral knots, and this leads to the shown 8 crossing
diagram of $6_3$ (which admits no reducing wave move, and a few
diagrams of $3_1\# !3_1$, which we chose not to dwell upon). \qed

{}From the work in \S\reference{sNR} we have then also:

\begin{corr}
Let $K$ be a $k$-almost positive achiral knot for $k\le 4$, and $K\ne
3_1\# !3_1$. Then one can simplify an $n$-crossing $\le 4$-almost
positive diagram $D$ of $K$ to a minimal crossing diagram by
$O(n^p)$ Reidemeister moves, where $p$ depends only on $k$. \qed
\end{corr}

Apart from the generator compilation, the most decisive
reason for the success of such proofs is the easy comutability
of the invariants we use, primarily the signature $\sg$ and
Jones polynomial $V$. It still took some days of work on a
computer to complete the proofs. However, rather than thinking
of such an approach as tedious, we feel that is should be considered
paying high tribute to this efficiency, which makes it not too
hard to evaluate the invariants on thousands of diagrams. (The
most complicated ones we encountered have 35 crossings.) Recently,
a lot of attention is given to ``relatives'' of $\sg$ and $V$ defined
in terms of homology theory. The comutability of these new invariants
is still sufficiently difficult, so that performing with them
a procedure similar to ours seems, unfortunately, out of question
for quite a while.

\section{The signature\label{Ss}}

As a last application related to positivity, we can finally settle
a problem, initiated in \cite{2apos}, and whose solution has been
suspected for a while. First we have

\begin{prop}\label{45657}
The knot $14_{45657}$ (depicted in figure 4 of \cite{2apos}) is the
only positive knot of genus $4$ with $\sg=4$. It has only one positive
(reduced) diagram, its unique 14 crossing diagram.
\end{prop}

This is the first \em{non-alternating} positive knot, which
now is known to have only one positive diagram.

\proof Let us seek a positive diagram of such a knot. Series of
composite and special generators have $\sg\ge 6$ and are clearly
ruled out. Checking the $\sg$ of the positively crossing-switched
non-special generators reveals that $\sg=8$ occurs 1,927,918
times, $\sg=6$ occurs 6662 times, and $\sg=4$ only once, for
the positive diagram of $14_{45657}$. Applying a $\bt$ twist
at whatever crossing of this diagram gives $\sg=6$. \qed

In general, the problem what is the minimal signature of positive
knots of some genus is very difficult to study. We devoted to
it a separate paper \cite{gsigex}. Some of the theoretical
thoughts there can be applied in practice here to prove

\begin{theo}\label{tsg4}
The positive knots of $\sg=4$ are:
\begin{mylist}{\arabic}
\myitem all genus 2 knots,
\myitem an infinite family of genus 3 knots, which is scarce, though,
in the sense that asymptotically for $n\to\infty$ we have
\begin{eqn}
\label{op}
\frac{\#\,\{\,K\,:\,\mbox{$K$ positive, $g(K)=3$, $\sg(K)=4$, 
$c(K)=n$}\,\}}{\#\,\{\,K\,:\,\mbox{$K$ positive, $g(K)=3$,
$c(K)=n$}\,\}}\,=\,O\left(\frac{1}{n^{10}}\right)\,,
\end{eqn}
and
\myitem the knot $14_{45657}$ (of genus 4).
\end{mylist}
\end{theo}

\proof
The genus 2 case is obvious, and let us first argue briefly
about the estimate \eqref{op} for genus 3. First, by the result
of \cite{SV} we know that the number of special alternating
knots of genus 3 and $n$ crossings behaves asymptotically like
a constant times $n^{14}$. That these are asymptotically dense
in the set of positive genus 3 knots follows from the result
of \cite{SV} that maximal generators are special alternating,
and an estimate in \cite{adeq}, which shows that
\begin{eqn}\label{oq}
c(D)-c(K)\,\le\,2g(K)-1\,,
\end{eqn}
for a positive diagram $D$ of a (positive) knot $K$. This
argument determines the behaviour of the denominator on the
left of \eqref{op} for $n\to\infty$.

To estimate (asymptotically,
from above) the numerator, we use again \eqref{oq}. This
basically (up to a constant) allows one to go over from counting
positive knots to counting positive diagrams. Now the positive
diagrams of genus $3$ with $\sg=4$ were described in \S 4 of
\cite{2apos} (see the remarks following proposition 4.1 there).
We know that we have infinite degree of freedom to apply
$\bt$-twists in at most 5 $\sim$-equivalence classes. Thus
the number of relevant diagrams of $n$ crossings is $O(n^4)$,
and the behaviour of the numerator of \eqref{op} is also clarified. 

With this argument, and proposition \ref{45657}, for the rest of the
proof we can assume we settled $g\le 4$, so we want to show that
there is no positive knot of genus $g\ge 5$ with $\sg=4$.

Let $D$ be a positive knot diagram. Then we consider a
sequence of diagrams $D=D_0,\dots,D_n=\bigcirc$, which is
created as follows:
\begin{mylist}{\arabic}
\myitem (``generalized clasp resolution'')
If $D_i$ has two equivalent crossings (i.e. two crossings 
that form a clasp, parallel or reverse, after flypes; see
definition \reference{dqa}) then we 
change one of these two crossings, apply a possible flype that 
turns the crossings into a trivial clasp, and resolve the clasp.
(This is the move (17) of \cite{gsigex}.) If there is no flype
necessary, nugatory crossings may occur after the clasp resolution. 
Then possibly remove these nugatory crossings to obtain $D_{i+1}$.
\myitem (``shrinking a bridge'')
If $D_i$ has no two equivalent crossings, then we choose a 
piece $\gm$ of the line of $D_i$ (with no self-intersections),
such that the endpoints of $\gm$ lie in neighbored regions 
of $D_i\sm\gm$. This means that after proper crossing changes 
on $\gm$, so that it becomes a bridge/tunnel, we can apply 
a wave move that shrinks it to a bridge/tunnel $\gm_1$ to 
length $1$ in $D_{i+1}$. We choose between bridge/tunnel 
so that the resulting single crossing of $\gm_1$ is positive. 
Herein we assume that the length of $\gm$, which is the number 
of crossings it passes, is bigger than $1$. We call such
a curve \em{admissible}. Among all such possible admissible 
$\gm$ (we showed in \cite{gsigex} that some always exist) 
we choose one of minimal length. If several minimal length 
bridge/tunnels $\gm$ are available, we choose $\gm$ so that 
$g(D_i)-g(D_{i+1})$ is minimal.
\end{mylist}

Now assume there is a $D$ with $\sg(D)\le 4$ and $g(D)\ge 5$. 
We will derive a contradiction. By considering the proper 
diagrams of the sequence found for $D$, we may assume 
w.l.o.g. that $D_1=D'$ has genus $\le 4$ and $\sg(D')\le 4$.
(Note that always $\sg(D_{i+1})\le \sg(D_{i})$.)

If $D'$ differs from $D$ by a clasp resolution, then $g(D)\le 
g(D')+1$. So the only option is that $g(D)=5$ and $g(D')=4$. Since 
$\sg(D')\le 4$, we checked that $D'$ is the diagram of $14_{45657}$. 
(We write below the knot for its positive diagram, since latter is 
unique.) Then $D$ is obtained from $D'$ by creating possibly first 
a certain number of nugatory crossings, and then a clasp (so that 
nugatory crossings become non-nugatory), and an optional flype (if 
there were no nugatory crossings). It is easy to see that the 
previously nugatory crossings can be switched in $D$ so that the 
clasp to become one that does not require nugatory crossings
to be added in $D'$. Thus if we check that all positive diagrams
obtained from $14_{45657}$ by adding a clasp \em{without} nugatory 
crossings have $\sg=6$, we are done. Now, if the genus remains $4$, 
we already checked it. If $g=5$, we have a prime positive diagram of 
16 crossings. Thus we seek a prime (\cite{Ozawa}) positive knot of
$g=5$ and $\le 16$ crossings. If the knot is alternating, then by
\cite{Murasugi5} we have $\sg=10$, so consider 
only non-alternating knots. A pre-selection from the table of 
\cite{KnotScape} using the (necessary) skein polynomial condition 
$\md_lP=\Md_mP=10$, shows that all these knots have $\sg\ge 6$.

For the rest of the argument assume that $D'$ is obtained from $D$
by shrinking a bridge. We would like to show that this second move 
is needed only in very exceptional cases, and they do not occur here.

We proved in \cite{gsigex} that $D'$ has at most one $\sim$-equivalence
class of $\ge 3$ elements, and they are $3$. Now we use the argument
in that proof to show even stronger restrictions.

\begin{lemma}
(1) Assume $D'$ has $\ge 2$ disjoint clasps (meaning that the pairs
of crossings involved are disjoint; the clasps may be parallel or 
reverse, but \em{not} up to flypes). Then $g(D)\le g(D')+1$.\\
(2) Assume $D'$ has $\ge 3$ disjoint clasps. Then $D'$ cannot be 
obtained from a diagram $D$ by shrinking a minimal bridge.
\end{lemma}

\proof 
(1) There is at least one clasp $a$ that does not contain the
crossing of the shrinked bridge/tunnel $\gm_1$. Since $D$ has
no equivalent crossings by assumption, the curve $\gm$ must pass
through the clasp $a$. We argued in \cite{gsigex} that it intersects
the two edges of the clasp only once. Then $D$ has an admissible
curve of length $3$, and shrinking this curve to a one crossing
curve reduces the genus by at most $1$. By choice of admissible
curve, we have the first claim.

(2) There are at least two clasps that do not contain the
crossing of the shrinked bridge/tunnel $\gm_1$. Since $D$ has
no equivalent crossings by assumption, the curve $\gm$ must
pass through both clasps. Then it intersects at least 4 edges, 
but we saw that at each clasp we have an admissible curve of 
length $3$ in $D$. This is a contradiction to the minimality 
of the admissible curve. \qed

The rest is a simple electronic check. We determined the prime
generators of genus 2 (they are 24) and 3 whose positification has 
$\sg=4$ (they are 13, and described in \cite{2apos}). Apply flypes to 
them, an optional $\bt$ twist, and again flypes. Then check that all 
the resulting diagrams have at least 2 disjoint clasps, and we see 
that $g(D)\le 4$ by part (1) of the lemma, so are done. The composite 
generators are only of genus 2 and directly ruled out the same way.
Likewise, the diagram of $14_{45657}$ has 4 clasps, and we apply 
part (2) of the lemma. This completes the proof of theorem \ref{tsg4}.
\qed

\begin{rem}
We remarked in \cite{gsig} the relation of the sequence of moves
we described to Taniyama's partial order. Indeed, one can apply
the above sort of proof in that context. For example, one can 
show that $4_1$ dominates all knots except connected sums of 
$(2,n)$-torus knots.
\end{rem}

\begin{prop}
If $K$ is an almost positive knot of genus $g\ge 3$, then $\sg\ge 4$.
\end{prop}

\proof If the almost postive diagram is composite, we can conclude 
the claim from the fact that $\sg>0$ when a knot is almost positive. 
So assume $D$ is prime. Since $g(K)\ge 3$, also $g(D)\ge 3$. We make 
almost postive the 13 generators with $\sg=4$ and $14_{45657}$.
In latter case always $\sg=4$, and in former case this is also 
true, up to 8 diagrams of $10_{145}$, which has $g=\sg=2$.

Now we build for these 8 almost positive diagrams $E$ the set $E_*$ as
in \eqref{*77}. Again if the negative crossing is equivalent to a 
positive one, we can discard the case by assuming we work 
with a least crossing almost postive diagram, or because the
diagram becomes positive. We can also work with one diagram 
per generator, because the signature is mutation invariant, 
and so is the genus of the knots if their mutated diagrams are 
almost positive, as can be concluded easily from the work in 
\cite{posqp}. Now we check that for all diagrams in all $E_*$ 
we have $\sg\ge 4$ or that $\Md_mP=2$. Latter option shows by 
\cite{posqp} that $g=2$, and the work there shows also that $g$ 
will not change under $\bt$ twists at a (positive) $\sim$-equivalence 
class of more than $1$ crossing (in an almost positive diagram). 
Since all diagrams we need to check are obtained from some diagram 
in $E_*$ by twists at a $\sim$-equivalence class of more than $2$
crossings, we see that if $\sg=2$, then $g=2$. \qed

\begin{corr}
If $K$ is a 2-almost positive knot, then $\sg>0$ except if $K$ is
a (non-positive) twist knot.
\end{corr}

\proof If $D$ is a 2-almost postive prime diagram of genus at least 
3, then we can switch one more crossing in the above set of 8 diagams
of $10_{145}$ and check that $\sg=2$. It remains to deal with
$g(D)\le 2$, and this was done in \S 6 of \cite{gen2} (see remark
6.1 therein). The composite case $D$ follows again easily from the 
prime one.
\qed

Even though the check in \S 6 of \cite{gen2} is somewhat tedious,
we must consider this proof as a considerable simplification of, and
far more elegant than the one in \cite{2apos}. It was the lack of
such proof that was so bothering while writing \cite{2apos}.


% ...

\section{Braid index of alternating knots\label{S6}}

\subsection{Motivation and history\label{mhs}}

We tried to apply a variant of the regularization of lemma \ref{lrez}
to the skein polynomial. Our motivation was to address the problem
studied by Murasugi \cite{Murasugi4} on
the sharpness of the MWF inequality \eqref{MWF} on alternating links.


If we have equality in \eqref{MWF}, then
one can determine the braid index from the skein
polynomial. Observably this often occurs. However, the examples
of unsharpness, albeit sporadic, are diversely distributed,
and it seems very difficult to make meaningful statements as to
nice classes on which the inequality would be sharp. It was a
reasonable conjecture that it would be so on alternating links.
Murasugi worked on this conjecture, and proved that it holds if
the link is rational or fibered alternating. His results remain
the most noteworthy
ones. Then he and Przytycki proved in \cite{MP} the cases when
$|\Mc\Dl\,|<4$, but found a counterexample of a 15 crossing
4-component link and an 18 crossing knot. 

The extreme paucity of
such examples (see proposition \ref{p72}) suggests still that
sharpness statements are possible for large classes of alternating
links. In that spirit, we show theorem \reference{tBi}.
%
%
%
Note that the Murasugi-Przytycki link counterexample is of genus
3, so that the theorem does not hold for links in this form. Their
knot counterexample has genus 6. 

To prove theorem \ref{tBi}, we need first to have a good control on
the MWF bound. We would like to have a growth in degree under
$\bt$-twists, and we can formulate analogous initial conditions
to the case of $V$ that ensure inductively under the skein relation
such degree growth. Originally this enabled us to deal with genus 2
and, with considerable effort, genus 3. However, we always needed
the initial two twist vectors for each $\sim$-equivalence class.
The resulting tremendous growth of the number of generators and
crossings in the initial diagrams made the case of genus 4 intractable,
except in special cases.

Then Ohyama's paper \cite{Ohyama} turned
our attention to \cite{MP}. There an efficient graph theoretic
machinery is developed to study the MWF bound. With the help
of this machinery we were able to considerably improve our work
and deal with genus 4 completely, and this is what we describe in
\S\reference{S6}. (We will, however, unlikely be able to do so for
genus 5.) Later, in \S\reference{SMS}, we will show how to modify
the arguments to obtain also a Bennequin surface.

The following clarification is to be put in advance. During our study
of \cite{MP}, we found a gap, which is explained in \S\ref{S61}. It
occurred when we wanted to understand the diagram move
of Figure 8.2 of \cite{MP}. Murasugi-Przytycki seem to assume
that Figure 8.2 is the general case, but we will explain that
it is not. And taking care of the missing case leads to a modified
definition of index, which we call $\inx_0$ (see definition
\reference{di0_}). Roughly speaking, the correction needed is
that in certain cases some edges in $\str v$ are not contracted
(cf. definition \reference{d2.3}). So Murasugi-Przytycki's diagram
move just proves instead of \eqref{q3} that 
\begin{eqn}\label{iis-1}
b(L)\le s(D)-\inx_0(D)\,.
\end{eqn}
Then the question is how do $\inx(D)$ and $\inx_0(D)$ relate to
each other. We will argue that 
\begin{eqn}\label{iis0}
\inx(D)\le \inx_0(D)\,,
\end{eqn}
which justifies \eqref{q3} (and its applications in Murasugi-%
Przytycki's Memoir). After we found this argument, we speculated,
based on our computational evidence, whether in fact always
\begin{eqn}\label{iis1}
\inx(D)=\inx_0(D)\,.
\end{eqn}
Later this was indeed established by Traczyk, who proved
in \cite{Traczyk2} the reverse inequality to \eqref{iis0}.
Still our (much more awkward) definition of index needs
(at least temporary) treatment, 
in order to prove \eqref{iis0} or \eqref{iis1} and fix the gap
in \cite{MP}. (A minor modification of $\inx_0$ will also be
used for the Bennequin surfaces.) Also, one must realize that
Murasugi-Przytycki's definition of index loses its geometric
meaning \em{per s\'e}. It simplifies the true transformation of
the Seifert graph under their diagram move, in a way which is
\em{a priori} incorrect but (fortunately) \em{a posteriori} turns
out to still give the right quantity. If one likes to keep the
correspondence between (Seifert) graph and diagram, one must live
with the circumstance that (in general) not all of $\str v$ is
to be contracted.

\subsection{Hidden Seifert circle problem\label{S61}}

\begin{figure}
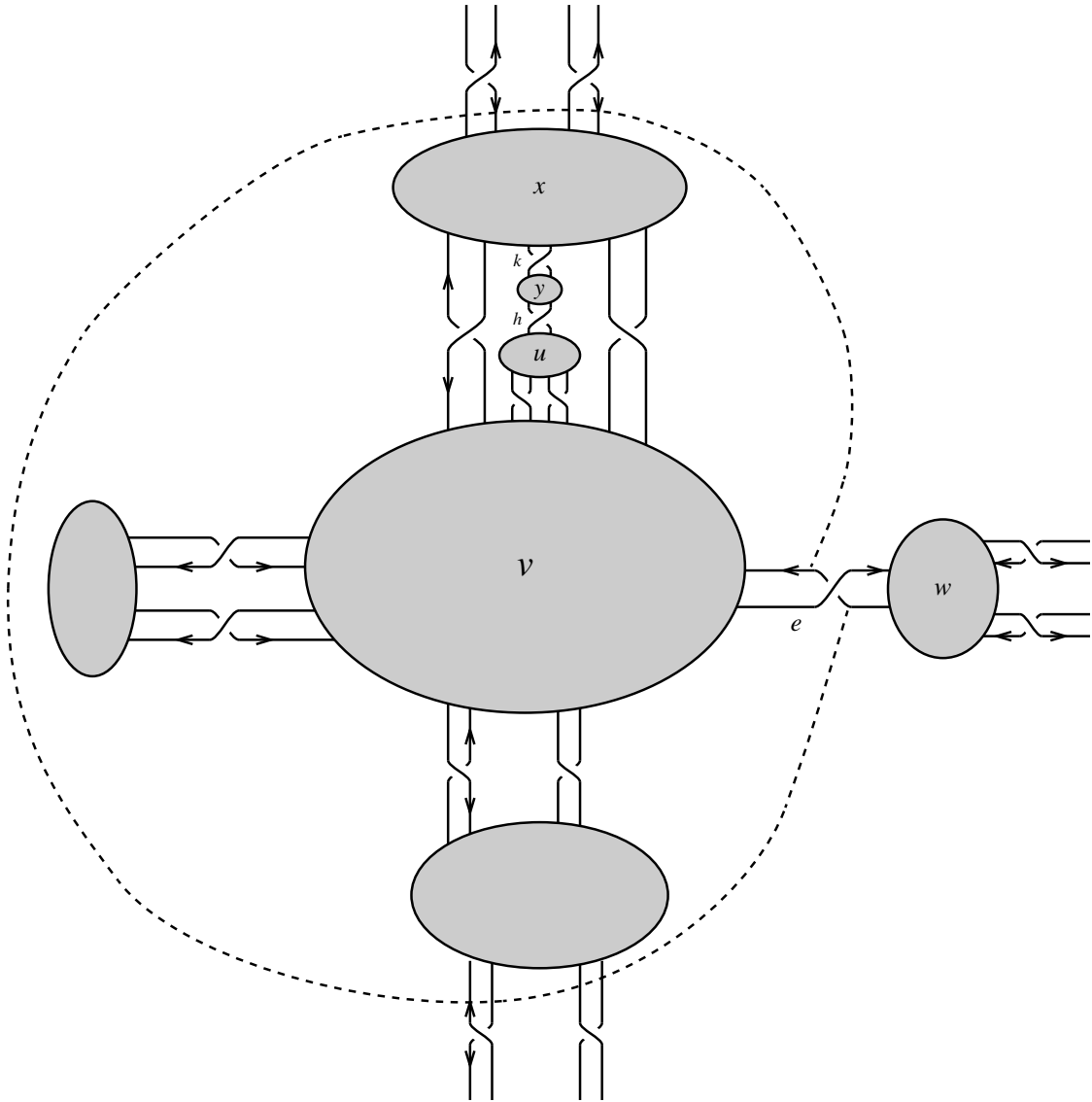

\[
\diag{1cm}{15}{15}{
  \picPSgraphics{0 setlinecap}
  \picclip{\picPSgraphics{0 0 m 0 15 l 15 15 l 15 0 l cp}}{
     \rvbr{11.5 7}{0}{0.5}{3}
     \lvbr{14.2 6.5}{180}{0.3}{2.2}
     \lvbr{14.2 7.5}{180}{0.3}{2.2}
     \lvbr{6.7 14}{-90}{0.4}{2}
     \lvbr{8.1 14}{-90}{0.4}{2}
     \lvbr{6.5 10.5}{90}{0.5}{3.8}
     \rbr{8.7 10.5}{90}{0.5}{3.8}
     \rbr{7.25 9.6}{90}{0.25}{1.2}
     \rbr{7.75 9.6}{90}{0.25}{1.2}
     \lbr{7.5 11.5}{-90}{0.3}{0.6}
     \lbr{7.5 10.7}{-90}{0.3}{0.6}
     \rvbr{6.4 4.5}{-90}{0.3}{2.8}
     \rbr{7.9 4.5}{-90}{0.3}{2}
     \rvbr{6.7 0.9}{90}{0.3}{2}
     \rbr{8.2 0.9}{90}{0.3}{2}
     \rvbr{3.2 7.5}{180}{0.4}{3}
     \rvbr{3.2 6.5}{180}{0.4}{3}
     \picfillgraycol{0.8}
     \picfilledellipse{7.5 12.5}{4 2 : 1.6 2 :}{$x$}
     \picfilledellipse{7.5 11.1}{0.6 2 : 0.4 2 :}{\small$y$}
     \picfilledellipse{7.5 10.2}{1.1 2 : 0.6 2 :}{$u$}
     \picfilledellipse{7.3 7.3}{4 2 : 6 2 : x}{\Large$v$}
     \picfilledellipse{7.5 2.8}{3.5 2 : 2.0 2 :}{}
     \picfilledellipse{1.4 7}{1.2 2 : 2.4 2 :}{}
     \picfilledellipse{13 7.0}{1.5 2 : 1.9 2 :}{$w$}
  }
  \picputtext{7.2 11.5}{\scriptsize$k$}
  \picputtext{7.2 10.7}{\scriptsize$h$}
  \picputtext{11 6.5}{$e$}
  {\piclinedash{0.1 0.1}{0.05}
   \opencurvepath{11.2 7.3}{11.5 8}{11.8 9}{11.8 10}{11.5 11}
                 {11 12}{10.1 13}{9 13.5}{7.5 13.6}
                 {6 13.5}{3.6 12.9}{1.5 10.9}
                 {1.0 10}{0.3 9}{0.2 5}{1 4}{2.5 2}
                 {5 1.2}{9 1.5}{10.5 3}{11.2 5}{11.7 6.7}
                 {}
  }
}
\]
\caption{A move of Murasugi-Przytycki, where the relayed strand 
(dotted line) does not go along a Seifert circle (denoted 
as $u$) adjacent to $v$. The Seifert circles are depicted
in gray to indicate that their interior may not be empty.}
\label{f10}
\end{figure}

\begin{figure}
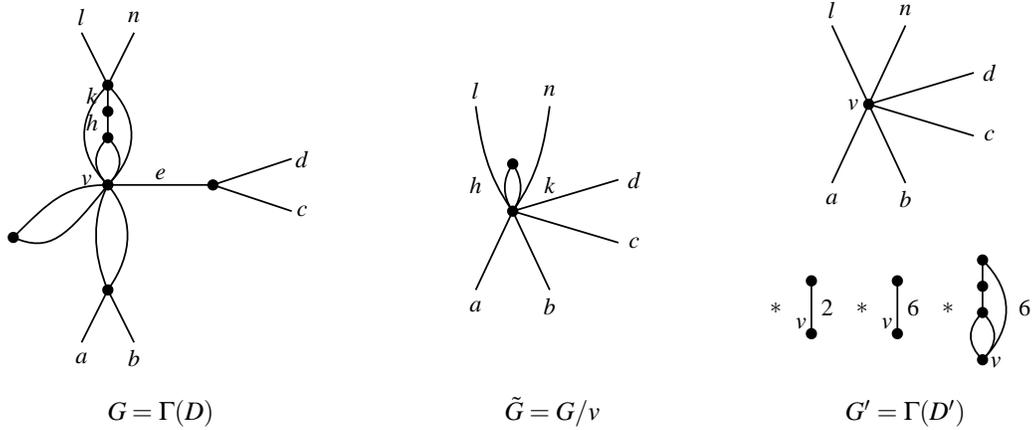

{\small
\[
\begin{array}{c@{\kern15mm}c@{\kern15mm}c}
\diag{7mm}{7}{7}{
  \vrt{2.5 3.5}
  \picputtext{2.1 3.6}{$v$}
  \piccurve{2.5 3.5}{2.8 3.8}{2.8 4.1}{2.5 4.4}
  \piccurve{2.5 3.5}{2.2 3.8}{2.2 4.1}{2.5 4.4}
  \vrt{2.5 4.4}
  \picline{2.5 4.4}{2.5 5.4}
  \picputtext{2.2 4.7}{$h$}
  \picputtext{2.2 5.2}{$k$}  
  \vrt{2.5 4.9}
  \vrt{2.5 5.4}
  \piccurve{2.5 3.5}{3.1 4.1}{3.1 4.8}{2.5 5.4}
  \piccurve{2.5 3.5}{1.9 4.1}{1.9 4.8}{2.5 5.4}
  \picline{2.5 5.4}{2.0 6.4}
  \picputtext{2.0 6.7}{$l$}
  \picline{2.5 5.4}{3.0 6.4}
  \picputtext{3.0 6.7}{$n$}
  \piccurve{2.5 3.5}{1.7 3.5}{1.4 3.2}{0.7 2.5}
  \piccurve{2.5 3.5}{1.7 2.5}{1.4 2.2}{0.7 2.5}
  \piccurve{2.5 3.5}{2.2 2.9}{2.2 2.2}{2.5 1.5}
  \piccurve{2.5 3.5}{3.0 2.9}{3.0 2.2}{2.5 1.5}
  \vrt{2.5 1.5}
  \vrt{0.7 2.5}
  \vrt{4.5 3.5}
  \picline{2.5 1.5}{2.0 0.5}
  \picputtext{2.0 0.2}{$a$}
  \picline{2.5 1.5}{3.0 0.5}
  \picputtext{3.0 0.2}{$b$}
  \picline{2.5 3.5}{4.5 3.5}
  \picputtext{3.5 3.7}{$e$}
  \picline{4.5 3.5}{6 4}
  \picputtext{6.2 4}{$d$}
  \picline{4.5 3.5}{6 3}
  \picputtext{6.2 3}{$c$}
}
&
\diag{7mm}{3.5}{5}{
  \vrt{1 2}
  \pictranslate{-1.5 d}{
     \piccurve{2.5 3.5}{2.7 3.8}{2.7 4.1}{2.5 4.4}
     \piccurve{2.5 3.5}{2.3 3.8}{2.3 4.1}{2.5 4.4}
     \vrt{2.5 4.4}
     \picputtext{1.8 4.0}{$h$}
     \picputtext{3.2 4.0}{$k$}  
  }
  \piccurve{1 2}{0.5 2.7}{0.4 3.3}{0.3 4}
  \picputtext{0.3 4.3}{$l$}
  \piccurve{1 2}{1.5 2.7}{1.6 3.3}{1.7 4}
  \picputtext{1.7 4.3}{$n$}
  \picline{1 2}{0.3 0.5}\picputtext{0.3 0.2}{$a$}
  \picline{1 2}{1.7 0.5}\picputtext{1.7 0.2}{$b$}
  \picline{1 2}{3 2.6}\picputtext{3.3 2.6}{$d$}
  \picline{1 2}{3 1.4}\picputtext{3.3 1.4}{$c$}
}
&
\
\vcbox{
\shortstack{
\diag{7mm}{3.5}{4}{
  \vrt{1 2}
  \picline{1 2}{0.3 3.5}\picputtext{0.3 3.8}{$l$}
  \picline{1 2}{1.7 3.5}\picputtext{1.7 3.8}{$n$}
  \picline{1 2}{0.3 0.5}\picputtext{0.3 0.2}{$a$}
  \picline{1 2}{1.7 0.5}\picputtext{1.7 0.2}{$b$}
  \picline{1 2}{3 2.6}\picputtext{3.3 2.6}{$d$}
  \picline{1 2}{3 1.4}\picputtext{3.3 1.4}{$c$}
  \picputtext{0.7 2}{$v$}
}\\[3mm]
$*$
\es
\diag{7mm}{0.8}{1}{
  \vrt{0.2 0}
  \vrt{0.2 1}
  \picline{0.2 0}{0.2 1}
  \picputtext{0.5 0.5}{$2$}
  \picputtext{0.0 0.2}{$v$}
}\es
$*$
\es
\diag{7mm}{0.8}{1}{
  \vrt{0.2 0}
  \vrt{0.2 1}
  \picline{0.2 0}{0.2 1}
  \picputtext{0.5 0.5}{$6$}
  \picputtext{0.0 0.2}{$v$}
}\es
$*$\es
\diag{7mm}{1.5}{2.4}{
  \pictranslate{-2.2 -3.3}{
    \vrt{2.5 3.5}
    \piccurve{2.5 3.5}{2.8 3.8}{2.8 4.1}{2.5 4.4}
    \piccurve{2.5 3.5}{2.2 3.8}{2.2 4.1}{2.5 4.4}
    \vrt{2.5 4.4}
    \picline{2.5 4.4}{2.5 5.4}
    \vrt{2.5 4.9}
    \vrt{2.5 5.4}
    \piccurve{2.5 3.5}{3.1 4.1}{3.1 4.8}{2.5 5.4}
    \picputtext{3.3 4.5}{$6$}
  }
  \picputtext{0.55 0.16}{$v$}
}}
}
\\[2.4cm]
\ry{1.7em}\mbox{\normalsize $G=\Gm(D)$} &
\mbox{\normalsize $\tG=G/v$} & \mbox{\normalsize $G'=\Gm(D')$}
\\[2mm]
\end{array}
\]
}
\caption{The various Seifert graphs of the diagrams related to
the move of Murasugi-Przytycki in figure \reference{f10},
in the case when the relayed
strand does not go along all Seifert circles adjacent to $v$.
The graph of $D'$ is given in its block decomposition, which
corresponds to the Murasugi sum decomposition along the newly
created Seifert circle. For simplicity, we display a multiple
edge by attaching the multiplicity to the edge drawn as simple
(otherwise, a letter attached just indicates the name).}
\label{f11}
\end{figure}

Recall definition \reference{d2.3}.
Now we must understand the move of Murasugi-Przytycki that
corresponds to the choice of a simple edge $e$ and the
contraction of the star of $v$ in $G$. (To set the record
straight, we should say that this move was considered,
apparently simultaneously and independently, also by Chalcraft
\cite{Chalcraft}, although merited there only with secondary
attention. With this understanding, we will refer to it
below still as the Murasugi-Przytycki move.) This move is
shown in figure 8.2 of \cite{MP}. Let $D$ be the diagram 
before the move and $D'$ the diagram resulting from it.
Let us for simplicity identify an edge with its crossing
and a vertex with its Seifert circle (see the remark above
proposition \reference{p71}). In this language, the move
of Murasugi-Przytycki eliminates one crossing, corresponding
to $e$. The crossings of the other edges $e'\ne e$,
incident to $v$, do not disappear under the Murasugi-Przytycki
move. Instead, they become in $D'$ parts of join factors of
$\Gm(D')$ that correspond to a Murasugi summand on the opposite
side of the modified Seifert circle. See the proof of lemma
8.6 in \cite{MP} and figures \reference{f10} and \reference{f11}.

The subtlety, which seems to have been overlooked
in the proof of \cite{MP}, is illustrated in figures \ref{f10} and
\ref{f11}. The Seifert circles adjacent to $v$ may be nested in $D$
in such a way that relaying the arc of $v$ by the move, one does (and
can) \em{not} go along \em{all} Seifert circles adjacent to $v$. In
the Seifert graph $G'=\Gm(D')$ of $D'$ some of the edges incident to
$v$ in $G=\Gm(D)$ may not enter, as written in the proof of lemma
8.6 in \cite{MP}, into block components that are
2-vertex graphs (with a multiple edge).

Still we see that contracting the star of $v$ in $G=\Gm(D)$, we obtain
a graph $\tG=G/v$, which is a contraction of $G'=\Gm(D')$. (We will
later describe \em{exactly} how $G'$ is constructed from $G$, but let 
us for the time being use the easier to obtain $G/v$ instead.) Here
contraction of a graph means that $\tG$ is obtained from $G'$ by
contracting some (possibly several or no) edges, and we allow
multiple edges in $G'$ to be contracted (by doing so
simultaneously with all simple edges they consist of).

More precisely, the difference between the block component
of $\tG$ and $G'$ is that in the last block component $X$ of
$G'$ in figure \ref{f11} the star of $v$ is contracted  to
obtain the block of $\tG$.
So for the proof of lemma 8.6 in \cite{MP} and \eqref{q3},
we actually need the following lemma.

\begin{lemma}\label{lmQ}
If a graph $H'$ is a contraction of $H$, then $\inx(H')\le \inx(H)$.
\end{lemma}

This lemma can be proved easily by induction on the number of
vertices, using the definition of the index. (The main point is that
decontraction does not increase edge multiplicity.) 
\em{Still is should be understood that the contraction of a vertex is
not fully correct as modelling the Murasugi-Przytycki diagram move.}

\subsection{Modifying the index\label{mht}}

In the following we will need to define several
``relatives'' of Murasugi-Przytycki's index. This results
from our desire first to fix the aforementioned error,
and then to keep track of the braided surfaces. We summarize
the various indices in the scheme \eqref{xs} and the remarks
after it. The reader may consult this in advance to avoid
confusion.

It becomes necessary to understand exactly the 
transformation of the Seifert graph $G$ under the move 
of Murasugi-Przytycki. We describe it now, also filling
in the detail overlooked by them.

Now we use graphs, in which we mark edges (each edge has a
$\bZ_2$-grading: it is either
marked or not). In the initial graph all edges are unmarked.
A marked edge is to be understood as one that cannot be
chosen as an edge $e$. It corresponds to a multiple edge.
\em{We assume} for the rest of the note
that \em{$G$ is bipartite}, thus $G$ has no cycles of length 3, 
which avoids some technical difficulties.

We choose a non-marked edge $e$ and a vertex $v$ of $e$. Let $w$ 
be the other vertex of $e$ (see figure \ref{f10}). We define the
notion \em{on the opposite side to $e$} as follows.

\begin{defi}\label{defopp}
A vertex $y\ne v,w$ is on 
the oppisite side to $e$ if there is a vertex $x\ne v,w,y$
adjacent to $v$ such that $y$ and $w$ are in different
conected components of $(G\sm v)\sm x$.
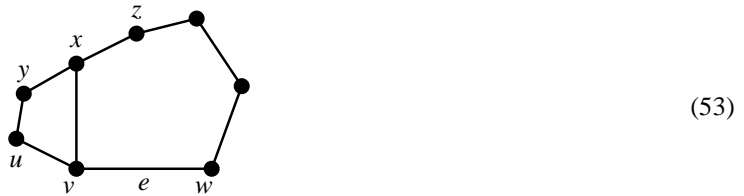
\begin{eqn}\label{yyy}
\diag{1cm}{3}{2.6}{
  \lvrt{0.8 0.5}{-0.1 -0.25}{$v$}
  \picputtext{1.7 0.3}{$e$}
  \lvrt{2.6 0.5}{-0.1 -0.25}{$w$}
  \vrt{3 1.6}
  \vrt{2.4 2.5}
  \lvrt{1.6 2.3}{0 0.3}{$z$}
  \lvrt{0.8 1.9}{0 0.3}{$x$}
  \lvrt{0.1 1.5}{0 0.3}{$y$}
  \lvrt{0.0 0.9}{0 -0.3}{$u$}
  \picline{0.8 0.5}{2.6 0.5}
  \piclineto{3 1.6}
  \piclineto{2.4 2.5}
  \piclineto{1.6 2.3}
  \piclineto{0.8 1.9}
  \piclineto{0.1 1.5}
  \piclineto{0.0 0.9}
  \piclineto{0.8 0.5}
  \piclineto{0.8 1.9}
}
\end{eqn}
(Here `$\sm$' stands for the deletion of a vertex together with
all its incident edges~-- but not its adjacent vertices.) Otherwise 
we say $y$ is \em{on the same side as $e$}.
\end{defi}

The meaning of 
this distinction is that the Murasugi-Przytycki move lays 
the arc along a Seifert circle $x$ adjacent to (the Seifert 
circle of) $v$, if $x$ is on the same side as $e$. This move 
affects the crossings that connect $x$ 
to $v$, or to a Seifert circle $z$ on the same side as $e$.

\begin{defi}\label{Gve}
We define now the marked graph $G\_ev$. The vertices of
$G\_ev$ are those of $G$ except $w$. The edges and markings
on them are chosen by copying those in $G$ as follows. Let an edge
$e'$ in $G$ connect vertices $v_{1,2}$.

\begin{caselist}
\case $v$ is among $v_{1,2}$, say $v=v_1$.

\begin{caselist}
\case If the other vertex $v_2$ of $e'$ is $w$ (i.e. $e=e'$), 
then $e'$ is deleted.

\case If $v_2$ is on the opposite side to $e$, then $e'$ 
is retained in $G\_ev$ with the same marking.\label{CG}

\case If $v_2$ is on the same side as $e$, then $e'$ 
is retained in $G\_ev$, but marked. 
\end{caselist}

\case $v$ is not among $v_{1,2}$.

\begin{caselist}
\case If none of $v_{1,2}$ is adjacent to $v$, then
$e'$ retains in $G\_ev$ the same vertices and marking.

\case One of $v_{1,2}$, say $v_1$, is adjacent to $v$.
(Then $v_2$ is not adjacent to $v$ by bipartacy.)

\begin{caselist}
\case If $v_1=w$, then change $v_1$ to $v$ in $G\_ev$, and retain
the marking. 

\case So assume next $v_1\ne w$. If $v_2$ is on the opposite
side to $e$, then retain $v_{1,2}$ and the marking.

\case If $v_2$ is on the same side as $e$, then 
we change $v_1$ to $v$, and retain the marking. 
(Note that by bipartacy, if $v_2$ is on the same side as $e$,
then so must be $v_1$.)
\end{caselist}

\end{caselist}

\end{caselist}
\end{defi}

Since a marking will indicate for us only that the edge cannot
be chosen as $e$, the resulting graph $G\_ev$ may be reduced
by turning a multiple edge into a single marked one. (This also 
makes it irrelevant to create multiple edges in case \ref{c1.3}.) 

\begin{theorem}
If $D'$ arises from $D$ by a Murasugi-Przytycki diagram move at a
crossing $e$ and Seifert circle $v$, then $\Gm(D')=\Gm(D)\_ev$. \qed
\end{theorem}

\begin{defi}\label{di0_}
Replace in definition \ref{d2.3} the two occurrences of $G_{v_i}=G/v_i$
by $G\_ev_i$, as given in definition \ref{Gve}. Then we define the
corresponding notions of \em{$0$-independent} edges and index
$\inx_0(D)$.
\end{defi}

With this definition, we obtain \eqref{iis-1}.
The property \eqref{iis0} can be proved by induction over
the edge number, looking at a maximal set of independent
(not $0$-independent) edges, using lemma \ref{lmQ} and
the observation that $G/v$ is a contraction of $G\_ev$.
This fixes Murasugi-Przytycki's proof of lemma 8.6 in \cite{MP}.


We speculated whether in fact \eqref{iis-1} can be stronger than
\eqref{q3}. If \eqref{iis1} is false, then conjecture \ref{C7} is
also false, so we were wondering whether this is a way to find
counterexamples to that conjecture. We explained, though, that
indeed \eqref{iis1} is true (and proved by Traczyk \cite{Traczyk2}
as a followup to our discovery of Murasugi-Przytycki's gap). One
could then say that this makes our (much more awkward) definition
of index obsolete. Still its treatment is necessary in order to
prove \eqref{iis0} or \eqref{iis1}, and fix the gap in \cite{MP}.

Also, if one likes to keep the correspondence between (Seifert)
graph and diagram, one must live with the circumstance that (in
general) not all of $\str v$ is to be contracted. The idea of
using vertex contraction (straightforwardly, following \cite{MP}) 
appeared in at least one further paper, \cite{MT}. We feel this
misunderstanding may cause a problem at some point, which was
additional motivation for us to provide the present correction.

We should stress that, while our proof of \eqref{iis0} might
imply that independent edges (in Murasugi-Przytycki's sense
of definition \reference{d2.3}) are $0$-independent (i.e.
corresponding to a set of single crossings admitting Seifert
circle reducing operations), we \em{do not know} if the converse
is true, despite \eqref{iis1}. Thus proofs (like those we give
below) relying on $\inx$ (rather than $\inx_0$) are to some
extent non-diagrammatic.

The important difference of $\inx_0$ to $\inx$ lies in not 
affecting edges in case \reference{CG}. The treatment of 
vertices on the opposite side to $e$, the technical detail missed 
by Murasugi-Przytycki, does not affect the result by \eqref{iis1}, 
yet it creates a lot of calculation overhead (which we experienced
in attempts to use the possibly better estimate \eqref{iis-1}
prior to Traczyk's proof of \eqref{iis1}). Note, however, that it 
implies the additivity of $\inx_0$ under block sum in an easier
(and much more natural) way 
than Murasugi-Przytycki's corresponding statement for $\inx$.

\begin{defi}
A marked graph is \em{not $2$-connected} if it has an \em{unmarked}
edge whose deletion disconnects it. If $G$ is not $2$-connected, there
is a plane curve intersecting $G$ in a single, and unmarked, edge.
We call such a curve a \em{separating curve}.
\end{defi}

Note that the initial (unmarked Seifert) graph of $D$ is 2-connected
because $D$ has no nugatory crossings. 

\begin{lemma}
If $G$ is 2-connected, so is $G\_ev$.
\end{lemma}

\proof We assume to the contrary that $G\_ev$ is not 2-connected. Let
$e'$ be a disconnecting edge. So there is a separating curve $\gm$ that
intersects $G\_ev$ only in $e'$. The only edges in $G\_ev$ which
do not exist in $G$ are of the type $vz$ in \eqref{yyy}
($z$ is a vertex on the same side as $e$, adjacent to a vertex
$x$ adjacent to $v$ in $G$). By definition \ref{defopp},
$vz$ belongs to a cycle, and so cannot disconnect $G\_ev$.

Therefore, $e'$ persists in $G$. It must be unmarked in $G$,
since the move from $G$ to $G\_ev$ never deletes markings.
Thus the curve $\gm$ must intersect $G$ in some other
edge. The only edges added in $G$ when recovering it from
$G\_ev$ (except that $e$ is decontracted) are of the form
$xz$ in \eqref{yyy} ($x$ is a vertex adjacent to $v$, and
$z$ is a vertex adjacent to $x$ on the same side as $e$). 
Then $\gm$ passes in $G$ through a cycle as the right one in
\eqref{yyy} (the one containing $z,x,v,w$ in consecutive order;
note that $z\ne w$ by bipartacy).
In $G\_ev$ this cycle is affected only by replacing $zx$,$xv$
by $zv$ (and contracting $e$). So $\gm$ must pass through
$e'=zv$ in $G\_ev$. But by
construction $zv$ is marked in $G\_ev$, and $\gm$ is not a
separating curve, a contradiction. \qed

It is easy to see that $G_1*G_2$ is 2-connected iff both
$G_1$ and $G_2$ are so.

\begin{lemma}\label{LAD0}
If $G_{1,2}$ are $2$-connected, then $\inx_0(G_1*G_2)=
\inx_0(G_1)+\inx_0(G_2)$.
\end{lemma}

\proof It is enough to see that the contraction procedure of an
edge $e$ in $G_1$ does not affect edges or markings in $G_2$,
except possibly the change of vertex at which the block sum
$G_1*G_2$ is performed. 

Let $v,w$ be the ends of $e$, and we consider the building of $G\_ev$
for $G=G_1* G_2$. Let $z$ be the vertex at which the block sum
$G_1* G_2$ is performed.

If $z\ne v$ is not adjacent to $v$, then nothing is changed in $G_2$
when building $G\_ev$.

Next assume $z=v$. The vertex $v$ must be adjacent to at least
one more vertex $x\ne w$ in $G_1$ (else $G_1$ is not $2$-connected
or $e$ is multiple). Then we see with this
choice of $x$ in definition \reference{defopp}
that the vertices in $G_2$ except $v$ lie on
the opposite side to $e$. Thus building $G\_ev$ does not
affect $G_2$. 

Finally assume $z\ne v$, but $z$ is adjacent to $v$.
If $z=w$ is the other end of $e$, then in $G\_ev$ all
edges incident in $G_2$ to $w$ are redirected to $v$
with the same marking, and so $G_2$ is not affected.
If $z\ne w$, then choosing $z$ for $x$ in definition
\reference{defopp}, we see that all vertices of $G_2$
except $z$ are on the opposite side to $e$. Thus none
of these edges is affected by building $G\_ev$. \qed


\subsection{Simplified regularization\label{SMT}}

To speed up the test for genus 4 generators, we make heavy use of
Murasugi-Przytycki's work. We recall the inequalities \eqref{q1}~--
\eqref{q3}, and the exactness of inequality \eqref{q4} in alternating
diagrams.

The proof of theorem \reference{tBi} will demonstrate the efficiency
of proposition \reference{p71} as a tool in determining the braid
index. To give a first glimpse of that capacity, we mention that
we have also obtained the following by computer verification
of the tables in KnotScape \cite{KnotScape} and Knotilus \cite{FR}:

\begin{prop}\label{p72}
If $K$ is an alternating knot of $\le 18$ crossings, then MWF is
sharp, except if $K$ is the (18 crossing) Murasugi-Przytycki knot
or its mutant. \qed
\end{prop}
 
This fact also shows the depth of the insight Murasugi-Przytycki
must have had in picking up exactly these two knots as counterexample
candidates!

The first easy lemma gives an upper control on the growth
of the braid index under $\bt$ twists.

\begin{lemma}\label{l72}
If $\tD$ is obtained from $D$ by a $\bt$-twist, then
$\inx(\tD)\ge \inx(D)+1$.
\end{lemma}

\proof Let $e$ be the egde in the Seifert graph $\Gm(D)$ of $D$ 
which is bisected twice to obtain the Seifert graph $\Gm(\tD)$ 
of $\tD$. Double bisection means to put two valence-2-vertices 
$v_{1,2}$ on $e$, dividing it into three edges $e_{1,2,3}$.
Let $e_2$ be the (middle) edge connecting $v_{1,2}$ in $\Gm(\tD)$. 

For an independent edge set $S$ in $\Gm(D)$ we construct an
independent edge set $S'$ in $\Gm(\tD)$ by keeping $S\sm\{e\}$
and including two of $e_{1,2,3}$ into $S'$ if $e\in S$, and one
of these three edges otherwise. This shows $\inx(\tD)\ge \inx(D)+1$.

(If $e$ connects a vertex of valence $2$, then one easily sees
that $\inx(\tD)=\inx(D)+1$. However, otherwise we cannot exclude
the possibility that $\inx(\tD)= \inx(D)+2$.)
\qed

\begin{rem}\label{XXX}
Let us stress again that any choice of two resp. one edge(s) among
$e_{1,2,3}$ to include into $S'$ will do when $e\in S$ resp. $e\not
\in S$, provided we specify
the to-be-contracted vertices correctly. This observation
will be important both when $e\in S$ (for the proof of corollary
\reference{tcr2}) and $e\not\in S$ (for the proof of lemma
\reference{L52}).
\end{rem}

The next step is to control the degrees of the skein polynomial
to estimate the braid index from below.

\begin{lemma}\label{L52}
Let $D$ be an alternating diagram with (see \eqref{q3})
\begin{eqn}\label{_str}
\mwf(D)\,=\,\mpb(D)\,.
\end{eqn}
Let $e$ be an edge in the Seifert graph $\Gm(D)$, which is
not contained in at least one maximal independent set. (That
is, there is a maximal independent set $C$ with $e\not\in C$.)
Assume further that $e$ is either a simple edge (i.e. its
crossing has no Seifert equivalent one) or that $D$ is special.

Let $\tD$ be obtained from $D$ by a $\bt$ twist at (the crossing
corresponding to) $e$. Then
\begin{enumerate}
\item \eqref{_str} holds for $\tD$\,,
\item $\inx(\tD)=\inx(D)+1$ and $\mwf(\tD)\,=\,\mwf(D)+1$, and
\item $e$ is not contained in some maximal independent set
  of $\Gm(\tD)$.
(Here $e$ is to be considered as a crossing in $\tD$ and identified
with some of the three crossings after the twist.)
\end{enumerate}
\end{lemma}

\proof The equality \eqref{_str}, together with \eqref{q1}~--
\eqref{q3} and the exactness of inequality \eqref{q4}, imply that
\eqref{q1}, \eqref{q2} are exact for $D$. Now let $D'$ be $D$ with
$e$ switched. Assume w.l.o.g. $e$ is positive in $D$ (the negative
case is handled analogously).

We want to show first
\begin{eqn}\label{&&}
\md_l P(\tD)\,=\,\md_l P(D)\,.
\end{eqn}
If $D$ is special (alternating), \eqref{&&} holds, because
both hand sides equal $1-\chi(D)$. So assume $e$ is simple.
Then $\inx_-(D')=\inx_-(D)+1$, because in $\Gm(D')$ we have a
single negative edge in a block. Denote by $q_1(D)$ and $q_2(D)$ the
right hand-sides of the inequalities \eqref{q1} and \eqref{q2}.
Since $w(D')=w(D)-2$, it follows that $q_2(D')=q_2(D)$, and since
\eqref{q2} is exact for $D$, we have 
\begin{eqn}\label{cf}
\md_l P(D')\,\ge\,\md_l P(D)\,.
\end{eqn}
Since (by the skein relation)
\begin{eqn}\label{xx}
P(\tD)\,=\,(l^2+1)P(D)\,-\,l^2P(D')\,,
\end{eqn}
we see that \eqref{&&} holds.

Now let us look at \eqref{q1}. For an alternating diagram $D$,
as remarked, all blocks of $G=\Gm(D)$ have edges of the same
sign. Thus theorem \ref{T25} implies that the indersection
of each maximal independent set of $G$ with the positive
blocks also realizes $\inx_+(G)$. The assumption that $e$ is
not contained in some maximal independent set of $G$ (and the
exactness of \eqref{q4} in alternating diagrams) means then
that $\inx_+(D')\ge \inx_+(D)$. Thus $q_1(D')< q_1(D)$. Since
\eqref{q1} is exact for $D$, we see from \eqref{xx} that 
\begin{eqn}\label{xyz}
\Md_l P(\tD)\,=\,\Md_l P(D)+2\,.
\end{eqn}
Now, \eqref{xyz} and \eqref{&&} imply that $\mwf(\tD)=\mwf(D)+1$,
and because $s(\tD)=s(D)+2$ and \eqref{_str}, we have that
$\inx(\tD)\le \inx(D)+1$. The opposite inequality was shown in
lemma \reference{l72}. So $\inx(\tD)=\inx(D)+1$, and the second
claim in our lemma follows. The first claim is then also clear.
The third claim is a consequence of remark \reference{XXX}.
%
\qed

Note that condition $3$ implies that one can iterate the twisting.
So we have a test for sharpness of MWF on a given series.

%

\begin{corr}\label{tcr2}
Assume $D$ is a special
alternating generator diagram such that \eqref{_str} 
holds, and let $S$ be the intersection of all maximal independent
sets of $D$. Let for $S'\subset S$ the diagram $D_{S'}$ be obtained
from $D$ by applying one $\bt$ twist at each crossing in $S'$.
Assume that for all $S'$
\begin{eqn}\label{tst}
\mwf(D_{S'})=\mwf(D)+|S'|\,.
\end{eqn}
Then the MWF inequality is exact on, and conjecture \reference{C7}
is true for the series of $D$.
\end{corr}

\proof As before, we split the series of $D$ into $2^{|S|}$
subseries, depending on whether we twist or not at each crossing 
in $S$. We know that each twist augments the index by at least one.
The equalities \eqref{_str} and \eqref{tst} then ensure that
it goes up by exactly one under a $\bt$ twist at some choice
of crossings $s$ from $S$. So for each of the 3 $\sim$-equivalent
crossings $s'$ obtained from $s$ under this twist, one can choose a
maximal independent set not to contain $s'$ by remark \ref{XXX}.
Then one can apply
further $\bt$ twists at $s$, and the sharpness of MWF is preserved.
\qed

For non-special generators we must take into account Seifert
equivalence classes.

\begin{corr}\label{tcr2'}
Assume for a non-special alternating generator diagram $D$
that \eqref{_str} holds, and let $S$, $S'$ and $D_{S'}$ be as
in corollary \ref{tcr2}. Let $S''$ be a set of crossings
of $D$, such that the intersection of $S''$ with each
Seifert equivalence class of $n$ crossings contains at most
$n-1$ elements. (In particular, $S''$ is disjoint from trivial
Seifert equivalence classes, and so $S''\cap S=\vn$.)
Assume that for all such sets $S'$ and $S''$ we have
\begin{eqn}\label{tst'}
\mwf(D_{S'\cup S''})=\mwf(D)+|S'|+|S''|\,.
\end{eqn}
Then the MWF inequality is exact on, and conjecture \reference{C7}
is true for the series of $D$.
\end{corr}

\proof The proof is completely analogous. When twisting further
at the classes twisted in $S''$ (resp.\ the one remaining
element in a Seifert equivalence class not in $S''$), all edges are
simple. \qed

Note the following slight simplification. Each of the Seifert
equivalence classes $X_i$ of $D$ decomposes by lemma \ref{l45}
completely into $\ssim$-equivalence classes $Y_{i,j}$. If for two sets
$S''$ and $\tS''$ we have $|S''\cap Y_{i,j}|=|\tS''\cap Y_{i,j}|$
for all applicable $i,j$, then $D_{S'\cup S''}$ and $D_{S'\cup \tS''}$
differ by a mutation. So we need to choose just one $S''$ for
given tuple of sizes $\{|S''\cap Y_{i,j}|\}_{i,j}$.

\proof[of theorem \reference{tBi}]
This proof consists now in a verification that bases heavily on
corollaries \reference{tcr2} and \reference{tcr2'}. First consider
one diagram $D$ of each generator knot $K$.

The decisive merit of these corollaries is that they provide a condition
to test sharpness of MWF by calculation of $P$ only on very few
diagrams in the series, and without adding a high number of crossings
by $\bt$ twists. The previous verification needs $2^{t(D)}$
diagrams, with $t(D)$ being the number of $\sim$-equivalece classes.
Note that by theorem \ref{th1gen} for $g(D)=4$ we have $c(D)\le 33$
and $t(D)\le 21$, so with that naive method, we would have up
to $2^{21}$ diagrams per generator, with up to $75$ crossings.

The graph theoretic determination of $S$ is not computationally
trivial, but takes only a fraction of the time that would have been
required to calculate the large number of extra polynomials. It
can be done by essentially the same recursion (and thus takes about
as long) as the calculation of the index:
\[
S\,=\,S_G\,=\,\bigcap\,\left\{\es\{e\}\cup S_{G/v} \, :\,
\begin{array}{c}
\mbox{$e$ simple edge of $G$, $v$ vertex}\\
\mbox{of $e$,\ } \inx (G/v)\,=\,\inx (G)-1
\end{array}\es\right\}\,.
\]
(And of course $S=\vn$ if $G$ has no simple edge.)
This calculation reveals (see \eqref{zsz}) that $S$ is small; always 
$|S|\le 5$. For most generators $S$ is actually empty.

For special generators the calculation of the index from the
definition is more time consuming, because they have more crossings
than the non-special ones, and admit no Murasugi-sum decomposition.
This suggested to seek some further simplifications. We found an easy
sufficient condition, which we describe in proposition \ref{psim},
to assure that $S$ is empty. In case of a special generator $D$ with
empty $S$ only the polynomial of $D$ needs to be calculated.

Proposition \ref{psim} is stated and proved in \S\ref{SMS}, since
it relates to some further notions treated there. The generator
$12_{1202}$ of genus two easily shows that the inequality \eqref{**}
in the proposition is false for non-special generators. On the other
hand, \eqref{oop} holds, and \eqref{**} becomes an equality for
all special genus two generators. For genus four, it is still exact
for slightly more than one half (748,193) of the special generators.

Even when proposition \ref{psim} does not apply directly, its proof
shows some priority that should be given to the contraction of
vertices of valence two. Note that the only difficulty for
non-special generators is that $F$ may not be a forest,
and we need to take care of not contracting vertices of
valence 2 incident to a double edge. With this additional
restriction, one can deal with non-special generators, too,
and this speeds up considerably their calculation.

The following table shows the sizes of intersections $S$ of maximal
independent sets for special and non-special genus 4 generators;
the last column gives the total number of extra, non-generator,
knots whose polynomials needed to be calculated. The most complicated
such knots have 35 crossings. For a special generator, their number
is $2^{\#S}-1$.

\begin{eqn}\label{zsz}
\begin{array}{c|*{6}{|c}||c}
\ry{1.2em}%
\#S & 0 & 1 & 2 & 3 & 4 & 5 & \mbox{total test knots} \\[1.4mm]
\hline\ry{1.2em}%
\mbox{special} & 1,082,270 & 269,406 & 106,204 & 20,676 & 1649 & 33 & 758,508
\\[1mm]
\hline\ry{1.2em}
\mbox{non-spec} & 1,573,426 & 290,668 & 63,821 & 6,488 & 178 & \mbox{---} &
 15,331,751
\\[1mm]
\end{array}
\end{eqn}

Even though the combination of $S'$ with the sets $S''$ creates
extra work for the non-special generators, this overhead was
manageable. We used the simplification for the choice of $S''$
explained after the proof of corollary \ref{tcr2'}. (There are
theoretic ways to further reduce the calculation, but they result
in too technical conditions, whose implementation augments the
risk of an error.) In total, non-special generators required
less time to deal with, because the calculation of the sets $S$
for special generators took about 15 times as long as for the
non-special ones.

With this, however, the work is not yet finished. The lack of
(at least confirmed) flype invariance of $\inx(D)$ requires
us to deal with diagrams $D_0$ obtained by flypes from the
generator diagrams $D$ we considered. To simplify the occurring
problem, first note that type B flypes (see figure \ref{ffl2})
commute with $\bt$ twists, so it is enough to generate
$D_0$ from $D$ by type A flypes. (This leads already
only to relatively few new diagrams $D_0$, about 58000 for
special generators $K$ and about 150000 for non-special ones.)
Also, the condition \eqref{oop} is invariant under flypes, so
no $D_0$ coming from diagrams $D$ we discarded using proposition
\reference{psim} need to be considered. Moreover, the mutation
invariance of $P$ ensures that the MWF bound will grow (by 1)
for each $\bt$ twist we apply on $D_0$. Since we argued that the
diagram index will grow also (at least) by 1, the verification
we need to perform on the remaining $D_0$ reduces to the check
of \eqref{_str}. (We determined in fact also the interesctions
$S$ of maximal independent sets for all the diagrams $D_0$, and
found that they have the same size as for the corresponding
diagrams $D$ that differ by type A flypes.)

We must at last admit that, even after these various arguments
and simplifications, still several weeks were necessary to
complete the work. \qed

\begin{rem}
It is clear from our proofs that the conjecture \reference{C7} of
Murasugi-Przytycki is also confirmed for the knots in theorem
\ref{tBi} and proposition \ref{p72}.
\end{rem}

Using \cite{STV}, we have for example:

\begin{corr}
The number of alternating genus 2 resp. 3 knots of braid index
$n$ grows like $O(n^8)$ resp. $O(n^{14})$. The number of achiral
such knots grows like $O(n^2)$ resp. $O(n^{5})$. 
\end{corr}

\proof Regularity of all generators implies that, up to an additive
constant, the braid index behaves like half of the crossing number on
such knots. Then the results on enumeration by crossing number 
extend directly by replacing crossing number by braid index. \qed

\subsection{A conjecture\label{OhS}}

The above work suggests that for fixed genus, $b$ should behave 
similarly to $c/2$. A calculation of the range of the difference $c-2b$ 
for alternating knots $K$ of genus $g\le 4$ reveals that we have
\begin{eqnarray}
g-2\,\stackrel{(*)}{\le}\,c-2b\,\le\,2g-3 & & \mbox{if $K$ is
special alternating,} \label{XX} \\
-2\,\le\,c-2b\,\le\,2g-4 & & \mbox{if $K$ is
not special alternating,} \label{XY} 
\end{eqnarray}
with all inequalities realized sharply. In fact, both estimates 
from above are shown in \cite{mwf} for all alternating \em{knots}
$K$ (of all $g$), and we see that they are the best possible (at 
least for $g\le 4$). The lower estimate in \eqref{XY} follows by 
the result of Ohyama \cite{Ohyama}. The improved lower bound (*) 
in \eqref{XX} seems unclear in contrast. It would be implied by 
the inequality
\begin{eqn}\label{ci}
4\inx(G)+e(G)\,\ge\,3(v(G)-1)\,,
\end{eqn}
for a (planar bipartite) graph $G$ with $e(G)$ edges and 
$v(G)$ vertices, and with \em{an odd number of spanning trees}. 
The spanning tree condition is necessary, as show simple 
examples, and its use makes an approach to prove \eqref{ci} 
difficult. Still we verified (*) and \eqref{ci} for (the 
Seifert graphs of alternating diagrams of) all special alternating
knots of up to 18 crossings and found no counterexample. If it 
holds, \eqref{ci} would in fact then imply (*) more generally
for alternating knots all whose Murasugi atoms (of the alternating
diagram; see above theorem \reference{T25}) are knots (and
none are links). 

One could conjecture (*) also for links:

\begin{conj}
For a special alternating link $L$ we have $c(L)-g(L)\ge 2(b(L)-1)$.
\end{conj}

In the case of links \eqref{ci} will take its extended form 
for general (planar bipartite) graphs $G$ 
\begin{eqn}\label{ci'}
4\inx(G)+e(G)\,\ge\,3v(G)-2-n(D)\,,
\end{eqn}
where $D$ is the (special alternating) diagram with $\Gm(D)=G$ and
$n(D)$ the number of its components. (Note that $G$ has odd number 
of spanning trees iff $n(D)=1$; see \cite{MS}.) The other 3
inequalities in \eqref{XX} and \eqref{XY} hold also for links, but
the upper estimates need some change depending on $n(D)$.


There is one more important special case in which one can
show our conjecture.

\begin{prop}
If $L$ is a special alternating arborescent link (\em{without} 
hidden Conway spheres) then $c(L)-g(L)\ge 2(b(L)-1)$.
\end{prop}

\begin{rem}
The hidden Conway spheres of Menasco \cite{Menasco} refer to
the links which are arborescent and alternating (i.e. have 
an arborescent diagram and an alternating diagram), but not
alternating\em{ly} arborescent (i.e. do not have a diagram
which is \em{simultaneously} arborescent and alternating).
Even though this family of exceptions is relatively small,
it exists, and is often carelessly overlooked by authors.
See \cite{Thistle2} for an explanation.
\end{rem}

\proof We show \eqref{ci'}
for a special alternating arborescent diagram $D$ and $G=\Gm(D)$.
For a special diagram $D$, the Seifert graph $G$ coincides with 
one of the checkerboard graphs. If $D$ is arborescent, 
latter are series-parallel. This means (see e.g. \cite{achir}) 
that they can be obtained from\es $\diag{8mm}{1}{0.6}{
\svrt{0 0.3}\svrt{1 0.3}\picline{0 0.3}{1 0.3}}$\es by edge 
bisections (putting a vertex of valence two on some edge) and 
doublings (adding a new edge between the same two vertices). 
So we assume now $G$ is (bipartite and) series-parallel. 
We use induction over the number $e_s(G)$ of simple edges of 
$G$. (Note that multiple edges enter into the count $e(G)$ 
with their multiplicity.)

First assume $e_s(G)=0$. Since decreasing by one the multiplicity 
$>2$ of a multiple edge changes $n(D)$ by $\pm 1$, it is enough to 
deal with the case that all edges of $G$ are double. Next one can 
successively delete double edges (which preserves $n(D)$ and 
$\inx(G)=0$), until one has a tree of double edges. Then the 
inequality \eqref{ci'} we claimed is easy to see (as an equality;
with $D$ being a connected sum of Hopf link diagrams).

Now assume $e_s(G)>0$. We consider a simple edge $e$ of $G$ which is
the first in a maximal independent sequence. 

This edge $e$ can, as a first case, be adjacent to
a vertex $v$ connecting only two vertices
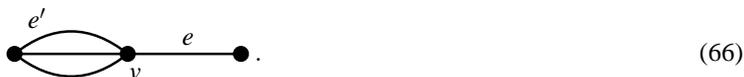
\begin{eqn}\label{*a}
\diag{1cm}{3}{1.2}{
  \lvrt{0 0.6}{0.3 0.48}{$e'$}
  \lvrt{3 0.6}{-0.7 0.2}{$e$}
  \lvrt{1.5 0.6}{0.1 -0.26}{$v$}
  \picline{0 0.6}{3 0.6}
  \piccurve{1.5 0.6}{1 0.2}{0.5 0.2}{0 0.6}
  \piccurve{1.5 0.6}{1 1.0}{0.5 1.0}{0 0.6}
}\es.
\end{eqn}
We want to use the definition of index and contract $v$. 

If $e'$ is simple, then the contraction preserves $n(D)$ (it 
trivializes a reverse clasp), reduces $v(G)$ by 2 and $\inx(G)$
at least by 1. The resulting graph is either series-parallel or
a block sum of such. So again we can argue by induction and the
stability of \eqref{ci'} under connected/block sum. If $e'$ is 
multiple, then this can be easily reduced to the argument when 
$e'$ is simple. This deals with the case \eqref{*a}.

\def\dblwedge#1#2#3{
  \pictranslate{#1}{
    \picrotate{#1 #2 2{1 i 4 i -}R 
      2 c d * x d * + sqrt /rr x D x atan 5 1 r 4{p}R}{
      \piccurve{0 0}{rr #3 * d N}{rr 1 #3 - * rr #3 N *}{rr 0}
      \piccurve{0 0}{rr #3 * d}{rr 1 #3 - * rr #3 *}{rr 0}        
    }
  }
}%
\def\dbledge#1#2{\dblwedge{#1}{#2}{0.2}}
\def\dbledgex#1#2{\dblwedge{#1}{#2}{0.3}}

If we have no fragment \eqref{*a}, then we have a simple edge
$e$ of the following form:
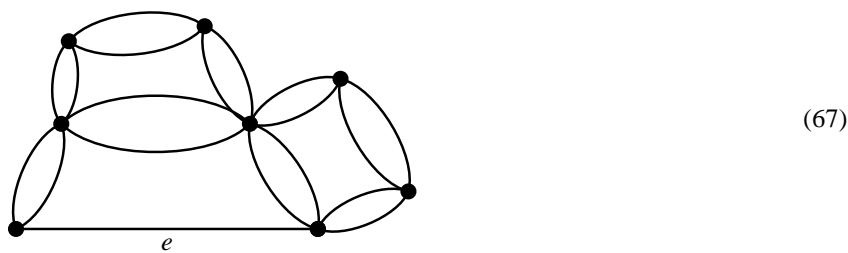
\begin{eqn}\label{**a}
\diag{1cm}{6}{3.5}{
  \vrt{0 0.3}
  \vrt{4 0.3}
  \vrt{0.6 1.7}
  \vrt{0.7 2.8}
  \vrt{2.5 3}
  \vrt{3.1 1.7}
  \vrt{4 0.3}
  \vrt{4.3 2.3}
  \vrt{5.2 0.8}
  \picline{0 0.3}{4 0.3}
  \picputtext{2 0.1}{$e$}
  \dbledge{0 0.3}{0.6 1.7}
  \dbledge{0.6 1.7}{0.7 2.8}
  \dbledge{0.7 2.8}{2.5 3}
  \dbledge{0.6 1.7}{3.1 1.7}
  \dbledge{2.5 3}{3.1 1.7}
  \dbledge{4.3 2.3}{3.1 1.7}
  \dbledge{5.2 0.8}{4.3 2.3}
  \dbledge{4 0.3}{5.2 0.8}
  \dbledge{4 0.3}{3.1 1.7}
}
\end{eqn}
The edges except $e$ are all multiple, but they can be reduced to
double ones using the previoius argument. Now from \eqref{**a}
one sees that there is always a cycle in which only two vertices 
are connected (not necessarily by simple edges, as shown below)
to vertices outside the cycle. 
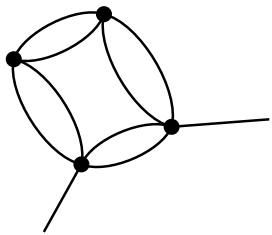
\[
\diag{1cm}{3.5}{3.2}{
  \pictranslate{-3 0.6}{
    \vrt{3.1 1.7}
    \vrt{4 0.3}
    \vrt{4.3 2.3}
    \vrt{5.2 0.8}
    \dbledge{4.3 2.3}{3.1 1.7}
    \dbledge{5.2 0.8}{4.3 2.3}
    \dbledge{4 0.3}{5.2 0.8}
    \dbledge{4 0.3}{3.1 1.7}
    \picline{5.2 0.8}{6.5 0.9}
    \picline{4 0.3}{3.5 -0.6}
  }
}
\]
By contractions of two double edges like
\[
\diag{1cm}{5}{1}{
  \picline{0 0.5}{1 0.5}
  \vrt{1 0.5}
  \vrt{2 0.5}
  \vrt{3 0.5}
  \vrt{4 0.5}
  \picline{5 0.5}{4 0.5}
  \dbledgex{1 0.5}{2 0.5}
  \dbledgex{3 0.5}{2 0.5}
  \dbledgex{3 0.5}{4 0.5}
}
\kern1cm\lra\kern1cm
\diag{1cm}{3}{1}{
  \picline{0 0.5}{1 0.5}
  \picline{3 0.5}{2 0.5}
  \vrt{1 0.5}
  \vrt{2 0.5}
  \dbledgex{1 0.5}{2 0.5}
}
\]
we can lose $2$ vertices, $4$ edges, and $n(D)$ goes down 
by 2. The index goes not go up by lemma \ref{lmQ}. It is 
clear that by such contractions and reducing a multiplicity 
$>3$ of edges by $2$, one can turn \eqref{**a} into\es
$\diag{1cm}{1}{0.8}{
  \svrt{0 0.4}\svrt{1 0.4}
  \picline{0 0.4}{1 0.4}
  \pictranslate{0 0.4}{
    \piccurve{0 0}{0.3 -0.4}{0.7 -0.4}{1 0}
    \piccurve{0 0}{0.3 0.4}{0.7 0.4}{1 0}
  }
  \picputtext{0.7 0}{$e$}
}$~.
(Keep in mind that the graph is bipartite.) Then $e$ is no 
longer simple, and again we can use induction. \qed

\begin{corr}
The inequality $c(L)-g(L)\ge 2(b(L)-1)$ holds for special alternating
Montesinos links $L$.
\end{corr}

\proof Such links are arborescent and have no hidden Conway spheres.
\qed

\section{Minimal string Bennequin surfaces\label{SMS}}

\subsection{Statement of result}

We apply the work in the previous section to the problem, what knots
have a braided surface of minimal genus on the minimal number of
strands. With the notation of \S\reference{SBR}, let $\sg_{i,k}\in
B_n$ for $1\le i<k\le n$ be the bands we introduced in \eqref{s_}:
\[
\sg_{i,k}\,:=\,
\sg_i\dots\sg_{k-2}\sg_{k-1}\sg_{k-2}^{-1}\dots\sg_{i}^{-1}\,.
\]
If one represents a link $L$ as the closure of a braid $\be\in B_n$
which is written as a product of $l$ of the $\sg_{i,k}$ and their
inverses, one obtains a Seifert surface $S$ of $L$ consisting of
$n$ disks and $l$ bands. $S$ is called \em{braided surface} of $L$.
If $S$ has minimal genus (i.e. equal to the genus of $L$), then $S$
is called a \em{Bennequin surface}. This terminology was coined by
Birman and Menasco (see for example \cite{BirMen}), and relates to
the work of Bennequin. He showed in \cite{Bennequin} that such
a surface exists for 3-braid links on a 3-string braid. Rudolph
\cite{Rudolph2} showed that \em{every} (not necessarily minimal
genus) Seifert surface can be made into a braided surface.
In particular a Bennequin surface always exists on a braid of
some (possibly very large) number of strings. A natural question
was whether the minimal number of strings (i.e. the braid
index of $L$) are always enough to span a Bennequin surface.
With M.\ Hirasawa \cite{HS} we showed later that there are knots
of genus 3 and braid index 4, which have no minimal (i.e. 4-)string
Bennequin surface. Hirasawa has also shown, contrarily, that for
2-bridge links \em{any} minimal genus surface is a braided surface on
the minimal number of strands. (The same holds for 3-braid links,
because in \cite{mwf} it was shown that such links have a single
minimal genus, in fact even just incompressible, surface.)
Here we want to show

\begin{theo}\label{thbs}
Any alternating knot of genus up to 4 or at most 18
crossings has a minimal string Bennequin surface.
\end{theo}

\subsection{The restricted index}

We define now a third variant of graph index, this time one which
keeps track of surfaces. As for $\inx_0$ in \S\reference{mht},
we use marked edges. In the initial graph $G=\Gm(D)$ all
edges are unmarked. A marked edge is to be understood as one that
cannot be chosen as an edge $e$. It corresponds to a crossing that
is grouped with other crossings to form a band, or which connects
the same 2 Seifert circles as some other crossing. \em{We assume}
again for the rest of the treatment that \em{$G$ is bipartite}, and
thus has no cycles of length 3, which would create some problems.

First we \em{reduce} the graph $G$ by turning a multiple edge 
into a simple marked one.

Next we choose a non-marked edge $e$ and a vertex $v$. Let $w$ 
be the other vertex of $e$ (see figure \ref{f10}). Recall the
notion on the opposite side to $e$ from definition \ref{defopp}.
The meaning of the distinction between the side of $e$ and
the opposite side was that the Murasugi-Przytycki move lays 
the arc along a Seifert circle $x$ adjacent to (the Seifert 
circle of) $v$, if $x$ is on the same side as $e$. This move 
affects the crossings (possibly in bands) that connect $x$ 
to $v$, or to a Seifert circle $z$ on the same side as $e$.

\begin{defi}\label{Gve'}
We define now the marked graph $G\^ev$. The vertices of
$G\^ev$ are those of $G$ except $w$. The edges and markings
on them are chosen by copying those in $G$ as follows.
Let an edge $e'$ in $G$ connect vertices $v_{1,2}$.

\begin{caselist}
\case $v$ is among $v_{1,2}$, say $v=v_1$.

\begin{caselist}
\case If the other vertex $v_2$ of $e'$ is $w$ (i.e. $e=e'$), 
then $e'$ is deleted.

\case If $v_2$ is on the opposite side to $e$, then $e'$ 
is retained in $G\^ev$ with the same marking.

\case If $v_2$ is on the same side as $e$, then $e'$ 
is retained in $G\^ev$, but marked. \label{c1.3}
\end{caselist}

\case $v$ is not among $v_{1,2}$.

\begin{caselist}
\case If none of $v_{1,2}$ is adjacent to $v$, then
$e'$ retains in $G\^ev$ the same vertices and marking.

\case One of $v_{1,2}$, say $v_1$, is adjacent to $v$.
(Then $v_2$ is not adjacent to $v$ by bipartacy.)

\begin{caselist}
\case If $v_1=w$, then change $v_1$ to $v$ in $G\^ev$, and retain
the marking. 

\case So assume next $v_1\ne w$. If $v_2$ is on the opposite
side to $e$, then retain $v_{1,2}$ and the marking.

\case If $v_2$ is on the same side as $e$, then 
we change $v_1$ to $v$, and put a marking. \label{hj}
(Note that by bipartacy, if $v_2$ is on the same side as $e$,
then so must be $v_1$.)
\end{caselist}

\end{caselist}

\end{caselist}
\end{defi}

Since a marking will indicate for us only that the edge cannot
be chosen as $e$, the resulting graph $G\^ev$ may be again reduced
by turning a multiple edge into a single marked one. (This also
makes it irrelevant to create multiple edges in case \ref{c1.3}.) 

\begin{defi}
We define the \em{restricted index} $\inx_b(G)$ of a graph $G$
like the (Murasugi-Przytycki) index $\inx(G)$ in definition
\ref{d2.3}, replacing $G/v$ by $G\^ev$. Again set $\inx_b(D)=
\inx_b(\Gm(D))$ for a diagram $D$.
\end{defi}

\begin{rem}\label{r6.3'}
If in case \reference{hj} of the definition of $G\^ev$,
we retain the old marking (and do not necessarily put one),
then we ignore the restriction coming from keeping bands.
So marked edges become the equivalent of multiple ones.
Thus the corresponding index is exactly the previously
defined $\inx_0$, correcting Murasugi-Przytycki's
definition of index to reflect their diagram move.
\end{rem}

Again we have

\begin{lemma}\label{LAD}
If $G_{1,2}$ are $2$-connected, then $\inx_b(G_1*G_2)=
\inx_b(G_1)+\inx_b(G_2)$.
\end{lemma}

\proof Similar to lemma \reference{LAD0}. \qed
%
%
%
%

\begin{prop}\label{oqa}
Let $D$ be a minimal genus diagram of a link $L$.
(I.e. the canonical surface of $D$ is a minimal genus surface
of $L$.) Then $L$ has a Bennequin surface on a braid of
$s(D)-\inx_b(D)$ strings.
\end{prop}

\proof We call a \em{band} a set of crossings in a diagram
which looks locally like a braid of the form 
\[
\sg_1\dots\sg_{k-1}\, \sg_n^{-1}\dots\sg_{k+1}^{-1}\,
\sg_{k}^{\pm 1}\sg_{k+1}\dots\,\sg_n\sg_{k-1}^{-1}\dots\sg_1^{-1}
\]
for some $1\le k\le n$. (This is a fragment of a diagram
that contains parts of $n+1$ Seifert circles.)

A marked edge corresponds to a crossing which either belongs
to a band (of $>1$ crossing), or which connects the same
Seifert circles as at least one other crossing. With this
understanding, one has to verify that the definition
of $\inx_b$ using $G\^ev$ reflects precisely the effect
on the diagram by the Murasugi-Przytycki move, keeping
track of markings. (Note, however, that a contraction at
$v$ can delete marked edges incident to $v$.)

Also, each Murasugi-Przytycki move acts within a single 
block. So we see that, starting with $D$, we can apply 
\em{at least} $\inx_b(D)$ Murasugi-Przytycki moves, wich 
each move reducing the number of Seifert circles and bands 
by 1. Finally, the moves of Yamada \cite{Yamada} can likewise 
be applied so that the number of (Seifert circles and)
bands is preserved. Then we have a braid diagram of $L$, 
which gives the desired surface. With this proposition 
\reference{oqa} is proved. \qed

\subsection{Finding a minimal string Bennequin surface}

It is clear now that with proposition \ref{oqa} we have a tool for
a computational verification of theorem \ref{thbs}. Unfortunately,
the description of $\inx_b$ makes its calculation very time
consuming. In practice, we tried to further speed up the recursion.
We feel details too technical and too long to present here, but
roughly speaking, we tried to take advantage
of the observation that handling of vertices `far apart' in the
graph does not depend on their order. We still could not avoid
that some sequences of edges were left unverified, so that we may
obtain smaller values than $\inx_b$. But this computation was much
faster, and when its bound already coincides with the one of MWF,
we know that the one of $\inx_b$ must coincide too. This dealt with
most generators much more quickly than $\inx_b$ would have done,
and we had to calculate $\inx_b$ (following the definition) only for
the handful of remaining ones.


We explain first the following partial case.

\begin{prop}\label{prth}
Alternating knots of genus $\le 3$ have a minimal
string Bennequin surface.
\end{prop}

\proof First one checks (by computer) for all generators $K$
and all their diagrams $D$ the property
\begin{eqn}\label{DD*}
\inx(D)=\inx_b(D)\,.
\end{eqn}
It is helpful to distinguish again flypes of type A and type B
as in figure \ref{ffl2}.
Flypes of type B commute with $\bt$ twists and type A flypes. So
it is enough to generate the diagrams obtainable from a given diagram
of a generator $K$ by type A flypes. For each such diagram $D$,
it is enough to find a diagram differing from $D$ by type B flypes,
which satisfies \eqref{DD*}. The following argument shows that
\eqref{DD*} is preserved under $\bt$ twists.

When \eqref{DD*} holds, then it is easy to see, either by a graph
theoretic
argument, or by deforiming the canonical surface geometrically, that
a $\bt$ move at a crossing not deleted by a Murasugi-Przytycki move
adds one disk and one band. Consider the case the crossing we
twist at is deleted by a Murasugi-Przytycki move. So the edge $e$ is
in a maximal independent set, and we chose a vertex $v$ to contract.
Then we apply first a move 
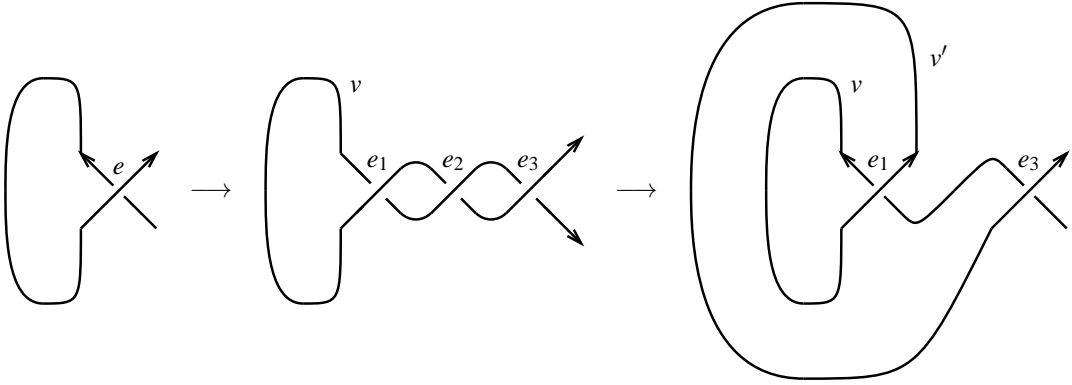
\[
\diag{1cm}{2}{3}{
  \opencurvepath{1 2}{1 3}{0 3}{0 2}{0 1}{0 0}{1 0}{1 1}{}
  \picvecline{2 1}{1 2}
  \picmultivecline{-6 1 -1 0}{1 1}{2 2}
  \picputtext{1.5 1.8}{$e$}
}\quad\lra\quad
\diag{1cm}{4.2}{3}{
  \opencurvepath{1 2}{1 3}{0 3}{0 2}{0 1}{0 0}{1 0}{1 1}{}
  \picputtext{1.2 2.9}{$v$}
  \pictranslate{1 1}{
    \picmultigraphics{3}{1 0}{
      \picline{0.2 0.8}{0.8 0.2}
      \picmultiline{-8.5 1 -1.0 0}{0.2 0.2}{0.8 0.8}
    }
    \picmultigraphics{2}{1 0}{
      \piccirclearc{1 0.6}{0.28}{45 135}
      \piccirclearc{1 0.4}{0.28}{-135 -45}
    }
    \picline{0 0}{0.3 0.3}
    \picline{0 1}{0.3 0.7}
    \picvecline{2.7 0.3}{3.2 -0.2}
    \picvecline{2.7 0.7}{3.2 1.2}
  }
  \picputtext{1.5 1.9}{$e_1$}
  \picputtext{2.5 1.9}{$e_2$}
  \picputtext{3.5 1.9}{$e_3$}
}\quad\lra\quad
\diag{1cm}{5.2}{5}{
  \pictranslate{1 1}{
    \picvecline{1.8 1.2}{1 2}
    \picmultivecline{-6 1 -1 0}{1 1}{2 2}
    \picline{4 1}{3.2 1.8}
    \picmultivecline{-6 1 -1 0}{3 1}{4 2}
    \picputtext{1.5 1.9}{$e_1$}
    \picputtext{3.5 1.9}{$e_3$}
    \opencurvepath{1.8 1.2}{2.0 1.0}{3 2}{3.2 1.8}{}
    \picputtext{1.2 2.9}{$v$}
    \opencurvepath{1 2}{1 3}{0 3}{0 2}{0 1}{0 0}{1 0}{1 1}{}
    \picputtext{2.3 3.3}{$v'$}
    \opencurvepath{2 2}{2 4}{-1 4}{-1 2}{-1 1}{-1 -1}{2 -1}{3 1}{}
  }
}
\,.
\]
(In the notation we identified edges with crossings and vertices with
Seifert circles.) So we alter
the middle of the 3 crossings after the $\bt$ move. (We lay the
rest of the new Seifert circle $v'$ below all other crossings attached
to $v$ from outside.) We obtain now at the site of the $\bt$ twist
locally a braid $\sg_2\sg_1$. Then delete by
a Murasugi-Przytycki move the edge $e_3$, contracting (the
vertex of) $v'$. So the result follows. 
\qed

The extension to genus $4$ requires both more calculation and more
delicate implementation techniques of the index. Let us state, however,
a simple and useful test that was quoted already in \S\ref{SMT}.

\begin{prop}\label{psim}
Assume $D$ is a special knot generator diagram, and $t(D)$
is the number of $\sim$-equivalence classes of $D$. Then
\begin{eqn}\label{**}
\mpb(D)\,=\,s(D)-\inx(D)\,\le\,t(D)+\chi(D)\,.
\end{eqn}
Also, $\inx_b(D)\ge c(D)-t(D)$. Moreover, the intersection of all
independent sets of size $c(D)-t(D)$ is empty. In particular, if
\begin{eqn}\label{oop}
t(D)+\chi(D)\,=\,\mwf(D)\,,
\end{eqn}
then MWF is exact and conjecture \reference{C7} is true
on the series of $D$.
\end{prop}

\proof Up to flypes assume that each non-trivial $\sim$-equivalence
class of $D$ forms a clasp. $D$ has $c(D)$ crossings and $c(D)+
\chi(D)$ Seifert circles. There are $c(D)-t(D)$ Seifert circles
of valence 2. We claim that for each such Seifert circle we can
take one of the edges adjacent to it into an independent set.

Now, because $D$ is a knot diagram, and is special, the edges
incident to valence 2 vertices in the Sefiert graph $\Gm$ of $D$
form a graph $F$, which is a forest.
Moreover, for each region $R$ of $\bR^2\sm \Gm$ there
are at least two edges in $\Gm\cap\partial R$ not incident to a
valence 2 vertex. (There is at least one such edge because $D$
is a knot diagram, and then at least one other, because $\Gm$ is
bipartite.) So one can contract all valence 2 vertices, obtaining
$c(D)+\chi(D)-(c(D)-t(D))=t(D)+\chi(D)$ vertices.

To find two disjoint independent sets, just choose the opposite
edge at each valence 2 vertex. To see that $\inx_b(D)\ge 
c(D)-t(D)$, choose a root $r$ in each tree of $F$, which is
not of valence $2$. Include in the independent set an edge
which is incident to a valence $2$ vertex $v$ adjacent to
$r$, but not connecting $r$ and $v$. Then contract $v$ and
proceed inductively. \qed

\proof [of theorem \ref{thbs}]
For genus 4, again we need to test all generator diagrams and 
all those obtained from them by type A flypes, with the 
freedom to apply type B flypes on each diagam before 
calculating $\inx_b$.

Some of the special generators to test with the most crossings required 
up to about half a day of calculation, and this made the work painful.
In practice it thus turned out necessary to speed up the procedure.


Some heuristical ways to choose the sequence of vertices $v$ and
edges $e$ more efficiently were used. Briefly speaking, contracting
edges far away from each other commutes, so we need to try only
one order of contractions. Accounting for such redundancies made
the calculation considerably faster at some places. It still took
about 2 weeks to verify \eqref{DD*} for all generators,
and thus to complete the work for the genus 4 knots.

In order to give an impression why this simplification is necessary,
we mention that much later we tried to recompute the result, for
verification purposes, using the definition of $\inx_b$ as given here.
This confirmed the outcome, but took $3\myfrac{1}{2}$ months, after
splitting the generator list into 100 equal parts and processing them
simultaneously. For some generators the calculation took several weeks.

The knots up to 18 crossings served again as a test case and were
checked several times with the different index implementations for
runtime performance purposes.
This concludes the proof of theorem \ref{thbs}. \qed

\begin{rem}
There is a question of Rudolph whether a strongly quasipositive
knot (or link) always has a strongly quasipositive (Bennequin)
surface on the minmial number of strings. While we expect this
not to be the case in general, we know of no counterexamples.
The theorem shows that the answer is positive for (special)
alternating knots of genus up to 4, or up to 18 crossings.
\end{rem}

In an attempt at clarification, we summarize the 3 types of
index (of \em{bipartite} graphs or diagrams) that occurred above:
\begin{eqn}\label{xs}
\def\se{\scbox{!}}
\def\gee{\stackrel{\se}{\ge}}
\def\rrr#1{\rbox{\rbox{\rbox{#1}}}}
\def\rrrr{\rrr{$\stackrel{\rbox{\scbox{!}}}{\ge}$}}
\begin{array}{*7c}
             &      & \inx_b(D)   &    &            &     &          \\
[3mm]        &      & \rrr{$\ge$}&    &       &          \\
s(D)-\mwf(D) & \gee & \inx_0(D)& =& \inx(D)          &          \\
\end{array}
\end{eqn}
The indices of the first, resp.\ second, line take, resp do not
take, into account keeping the bands. Thus these indices can,
resp.\ can not, be used to estimate Bennequin surface string numbers.
The indices of the left, resp. right, side take, resp.\ do not take,
into account the distinction of vertices on the opposite side to an
edge. The three indices are arguedly additive under join (theorem
\reference{T25}, lemmas \ref{LAD0} and \ref{LAD}). The inequality
marked with a `$!$' is known to become strict (i.e. not an equality)
for some diagrams/graphs. The equality on the right uses Traczyk's
aforementioned argument \cite{Traczyk2}, which we did not present here.

In \cite{mwf} we proved the existence of the minimal string
Bennequin surface for alternating links of braid index 4.
Note that our knot examples in \cite{HS} have braid index 4,
crossing number 16 and genus 3. So, if we exclude alternation,
they would fall into \em{any} of the three categories we
confirmed the Bennequin surface on. Thus alternation is a
crucial assumption. However, the proof of theorem \ref{thbs} 
and the Murasugi-Przytycki examples \cite[\S 19]{MP} of unsharp 
MWF inequality ($5=\mwf<b=6$) show the difficulties with 
the extension of the Bennequin surface result for alternating 
knots both of higher genus and braid index ($>4$). Still, 
even if our construction of the Bennequin surface would
fail giving minimum number of strings on alternating diagrams
(of which we have no example yet),
this would by far not mean that there is no minimal string
surface. So, with some courage, we could ask:

\begin{question}
Does every alternating knot have a minimal string Bennequin surface?
\end{question}

Note that our tests are more general than the cases we applied 
them on. For example, the sharp MWF and minimal string Bennequin 
surface criteria apply also for an alternating pretzel link $(x_1,
\dots,x_n)$, with $x_i$ odd and the twists therein reverse.

Dealing with canonical surfaces leads to some related question:
Does even every knot whose genus equals the canonical genus have a
minimal string Bennequin surface? It is interesting to remark that for
none of the examples reported in \cite{HS} these genera are equal.

%
%
%

\section{\label{SHo_}The Alexander polynomial of alternating knots}

While it is very well known what Alexander polynomials occur for
an arbitrary knot, the question about those occurring for an
alternating knot is much harder. We will apply our work to this
problem in the last section.

\subsection{\label{SHo}Hoste's conjecture}

There is a conjecture I learned from personal communication with
Murasugi, who attributes it to Hoste.

\begin{conj}(Hoste)
If $z\in\bC$ is a root of the Alexander polynomial of an
alternating knot, then $\R z > -1$.
\end{conj}

Here we will use again our generator classification and an
appropriate calculation to show

\begin{theo}\label{thHo}
Hoste's conjecture is true for knots up to genus 4 (or 
equivalently for Alexander polynomials of maximal degree up to 4).
\end{theo}

The parenthetic rephrasing refers to the normalization of
$\Dl$ with $\Dl(t)=\Dl(1/t)$ and $\Dl(1)=1$. (Beware that 
we will change this normalization during the proof.)

Some properties of the Alexander polynomials of alternating
knots are known; see \cite{Crowell,Murasugi2,MS}. In particular
the following holds:

\begin{theo}\label{thCM}(\cite{Crowell,Murasugi2})
The coefficients of $\Dl_K$ for an alternating knot $K$ alternate
in sign, i.e. $[\Dl_K]_{i}[\Dl_K]_{i+1}\le 0$ for all $i\in \bZ$.
\end{theo}

The zeros seem less easy to control than the coefficients, 
and Hoste's conjecture appears open even for example for 
2-bridge knots. Still theorem \ref{thCM} will be a useful
ingredient in the proof of theorem \ref{thHo}. Also, there is 
one important tool to control some zeros of $\Dl$, the signature.

\begin{theo}\label{thsg}
For any knot $K$ the number of zeros of $\Dl_K$ in $\{\,z\in\bC\,:\,|z|=1,
\,\I z>0\,\}$, counted with multiplicity, is at least $\ffrac{1}{2}
|\sg(K)|$.
\end{theo}

This is a folklore fact; see \cite{Gar,spec} for some explanation.
In \cite{spec} it was proved alternatively for special alternating 
knots, and we deduced:

\begin{corr}(\cite{spec})
If $K$ is special alternating, then all zeros of $\Dl$ lie on $S^1=
\{\,z\in\bC\,:\,|z|=1\,\}$, and so Hoste's conjecture holds for $K$.
\end{corr}

(The value $\Dl(-1)$ is the determinant, and non-zero for example 
by theorem \ref{thCM}; but in fact it is never zero for any knot.)

One further important tool, Rouch\'e's Theorem, will be introduced 
during the proof.

\proof[of theorem \ref{thHo}]
Assume $\Dl(z)=0$ for $\R z\le -1$. Since by theorem \ref{thCM},
$\Dl$ has no zeros on the negative real line, and $\Dl$ is real
and reciprocal, we see that $\Dl$ must have the $4$ distinct
zeros $z^{\pm 1},\bar z^{\,\pm 1}$. Then by theorem \ref{thsg},
we see that Hoste's conjecture is true whenever $|\sg|\ge 2g-2$,
and in particular the case $g=1$ is finished.

Consider next $g=2$, $\sg=0$. Let $\Dl_{[i]}=[\Dl]_{\md \Dl+i}$.
Then we have
\begin{eqn}\label{star+}
\sum_{\Dl(z)=0}\,z\,=\,-\frac{\Dl_{[1]}}{\Dl_{[0]}}\,\ge\,0\,.
\end{eqn}
However, if $\R z\le -1$, then the same holds for $\bar z$, and also
$\R z^{-1}, \R \bar z^{-1}$ are negative, so that the left sum
in \eqref{star+} is in fact $<-2$, a contradiction. This
finishes $g=2$.

The argument using \eqref{star+} works also for $g=3$, $\sg=2$, 
since the extra pair of zeros on $S^1$ augments the sum on the
left of \eqref{star+} by at most two, and so it is still negative. 

For $g=3$ it remains to consider $\sg=0$, the first non-trivial case. 

Note that $\sg$ is constant on each series of alternating diagrams 
(as can be inferred from the combinatorial formula of $\sg$,
given for example in \cite{Kauffman}). Thus we need to consider
the series of $\sg=0$ generators. Composite knots (inductively)
and mutations are immaterial from the point of view of Hoste's
conjecture, so again it suffices to consider prime generators, 
and one diagram per generator.

Among the 4017 generators of genus 3, only 210 have $\sg=0$.
Now note (e.g. using \eqref{star+}) that if $\Dl(z)=0$, $\R z
\le -1$, then the remaining two zeros of $\Dl$, different from 
$z^{\pm 1}$ and their conjugates, must be real positive
(and mutually inverse). Let $z'$ be the (real) zero $>1$.

Now for any diagram $D$ in the series $\br{D'}$ of $D'$
(resuming the notation of definition \reference{dbr}), the
Alexander polynomial has the form
\begin{eqn}\label{star++}
\Dl(D)\,=\,\sum_{i=1}^n\, a_i\,\Dl(D_i)\cdot (1-t)^{n_i}\,,
\end{eqn}
where $a_i\in\bN$, and $D_i$ are obtained by smoothing out 
one crossing in $n_i$ different $\sim$-equivalence classes of 
$D'$. Hereby the convention for $\Dl$ was \em{changed}, so that 
$\md \Dl=0$ and $\Dl_{[0]}=\Dl(0)>0$. An explanation of formula
\eqref{star++} is given \cite{STW}. It was clarified that all
terms have leading coefficients of equal sign (and then the same 
applies to all other coefficients using theorem \ref{thCM}).

The polynomials
\begin{eqn}\label{DDL}
\Dl_i\,:=\,\Dl(D_i)\cdot (1-t)^{n_i}\mbox{\quad and\quad}
\tl\Dl_i\,:=\,\Dl(D_i)
\end{eqn}
will be the central object of attention (and calculation) 
from now on. Note also that 
\begin{eqn}\label{star+++}
n_i=\Md \Dl(D')-\Md \Dl(D_i)\,=\,2g-\Md \Dl(D_i)\,=\,6-\Md \tl\Dl_i\,.
\end{eqn}
In particular, $\Dl_i=0$ if $n_i>2g(=6)$, 
and these $D_i$ and $\Dl_i$ can be discarded in the sum of 
\eqref{star++}. (We will soon see that it is very helpful to 
choose the number $n$ of terms in this sum as small as 
possible.) 

With this reduction, it is easy to see that 
\begin{eqn}\label{star4+}
\begin{array}{c}
m_1\\M_1
\end{array}
\,:=\,\begin{array}{c}
\inf \\
\sup
\end{array}\,\left\{ \,\left|\,\ffrac{\Dl_{[1]}(D)}{\Dl_{[0]}(D)}\,
  \right|\ :\ D\in\br{D'}\ \right\}\,=\,
\begin{array}{c}
\min\\
\setbox\@tempboxa=\hbox{$\max$}
\vbox{\copy\@tempboxa\vbox to \z@{\hbox to 
  \wd\@tempboxa{\scriptsize \hss $i$\ry{1.0em}\hss}\vss}}
\end{array}\,\gm(D_i)\,,
\end{eqn}
where
\begin{eqn}\label{star4+.}
\gm(D_i):=\frac{(\Dl_i)_{[1]}}{(\Dl_i)_{[0]}}=
\frac{(\tl\Dl_i)_{[1]}}{(\tl\Dl_i)_{[0]}}+6-\Md \tl\Dl_i\,.
\end{eqn}
Then \eqref{star+} implies that if there is a $z$ of $\R z\le -1$ with 
$\Dl_{D}(z)=0$, then the real zero $z'>1$ of $\Dl_{D}$ satisfies
\begin{eqn}\label{star5+}
z'+\frac{1}{z'}>2+m_1\,.
\end{eqn}
Since $z'$ is the unique zero $>1$ (and because $2+m_1\ge 2$),
\eqref{star5+} is equivalent to
\begin{eqn}\label{star6+}
\Dl(D)(z_0)<0\,,\quad\mbox{where}\quad z_0+\frac{1}{z_0}=2+m_1
\mbox{\ and\ }z_0>1\,.
\end{eqn}
(Keep in mind that we normalized $\Dl$ so that $\Dl(0)>0$, so
the leading coefficient is also positive.) Now if we check that 
contrarily
\begin{eqn}\label{star7+}
\Dl(D_i)(z_0)>0
\end{eqn}
for all $i=1,\dots,n$, then \eqref{star6+} clearly cannot hold.
Let us call \eqref{star7+} the \em{positive zero test}.

For each generator $D'$ the $D_i$ can be generated, $m_1$ and 
$z_0$ calculated and then \eqref{star7+} checked. (We excluded 
$n_i\ge 7$, and if $n_i=6$ the only polynomial coming in question 
is $\tl\Dl_i=1$, so we can assume $n_i\le 5$.)

The positive zero test took a few minutes, and succeeded on the 
210 generators. With this $g=3$ is also finished.

For $g=4$ more work is needed. As before $\sg\ge 6$ is trivial,
but for $\sg=4$ we have something to check. Now with two pairs
of (conjugate) zeros of $\Dl$ on $S^1$, \eqref{star+} and 
\eqref{star4+} give the condition
\begin{eqn}\label{c2}
\,\left|\,\ffrac{\Dl_{[1]}(D)}{\Dl_{[0]}(D)}\,\right|\,<\,2\,.
\end{eqn}
To see that this is false, by theorem \ref{thCM} and \eqref{star+++} 
it is enough to consider for each generator diagram $D'$ the 
diagrams $D_i$ with $n_i\le 1$, and calculate that $\gm(D_i)\ge 2$
for $\gm(D_i)$ in \eqref{star4+.}. There are about 500,000 generators
of $g=4$ and $\sg=4$, and this test took a few hours, but succeeded.

For $g=4$, $\sg=2$ we have a modification of the positive zero
test used for $g=3$, $\sg=0$. Most thoughts apply still, except
that in \eqref{star+++} and \eqref{star4+.} the `6' (standing
for $2g$) becomes `8', and in \eqref{star5+} and \eqref{star6+}
the additive $2$ on the right disappears, since again we have
a new pair of zeros $z$ on $S^1$, giving an additional contribution
$2\R z$ up to 2 in the sum of \eqref{star+}. (Now we must also
take care that $m_1\ge 2$, so that $z_0$ in \eqref{star6+} is
still real.)

The positive zero test thus now weakens (i.e. the condition tested is
stronger and may more likely not hold), and there are more generators
(about 190,000). There turn up to be 1157 cases of failure. In an 
attempt to handle these generators $D'$, we enhanced the positive 
zero test.

Apply the test on the diagrams obtained from $D'$ by $\bt$ twisting
once (separately) at each $\sim$-equivalence class of $D'$. If the 
test succeeds on $D'_i$, obtained from $D'$ by a $\bt$ move at the 
$i$-th class, then we can exclude twists in $D'$ at this class. With 
restrictions thus obtained we return again to the positive zero test 
of $D'$. The fewer classes to twist at mean that the set of diagrams 
$D_i$ in \eqref{star++} becomes smaller. Then $m_1$ in \eqref{star4+}, 
and with also $z_0$ in \eqref{star7+}, goes up, while there are 
fewer positivity tests to perform. 

With this trick from the 1157 cases the number of difficult ones was 
reduced to 304. These, as well as the $\sg=0$ knots, require a new 
method, which we sought for a while, and found in Rouch\'e's Theorem.

\begin{theo}(Rouch\'e's Theorem; see e.g. \cite[\S 1]{RahSch})\label{RT}
If two holomorphic functions $f,g$ inside and on some piecewise
smooth closed contour $C\subset \bC$ satisfy 
\begin{eqn}\label{0fg}
0<|g(z)|<|f(z)|\mbox{\ \ \ for all $z\in C$}\,,
\end{eqn}
then $f$ and $f+g$ have the same number of zeros (with
multiplicity) inside $C$.
\end{theo}

For us $f,g$ will be always sums of the type \eqref{star++}. We 
also would like to choose 
\begin{eqn}\label{C}
C=\{\,\R z=-1\,\}\,.
\end{eqn}
Now, the contour is not closed. So let us argue why this choice
is admissible.

In our case $f,g$ will be polynomials of the same degree with
real leading coefficients $\Mc f$ and $\Mc g$, which do not
cancel out. Moreover, all ratios between coefficients are
bounded, in a way depending only on the generator $D'$. Thus
it is not hard to a find a constant $R$ (depending only on 
$D'$), such that the norms of $f,g$ compare uniformly on
$\{\,\R z\le -1, |z|=R\,\}$ in the same way as their absolute 
leading coefficients $|\Mc f|$ and $|\Mc g|$ compare. If $\Mc f \ne
\Mc g$, we can achieve \eqref{0fg} after possibly swapping $f,g$.
If $\Mc f=\Mc g$, we take $g$ in \eqref{0fg} to be our present $g/2$,
and apply the theorem twice.

With this understanding we assume that $C$ is chosen as in \eqref{C}.

We will consider the property that
\begin{eqn}\label{C_}
\mbox{$\R(f\bar g)$ is positive on all of $C$}.
\end{eqn}
This positivity condition
is insensitive for the tricky case we decided to divide $g$
by $2$. If \eqref{C_} is satisflied, Rouch\'e's Theorem implies
that $f+g$ has the same number of roots $z$ with $\R z\le -1$ as
\em{either} of $f$ or $g$ do. For $f$ this is precisely theorem
\ref{RT}, assuming with the preceding argument that \eqref{0fg}
holds. For $g$ one applies that theorem on $f+g$ and $-g$, and
uses \eqref{C_} to ensure \eqref{0fg}. (Clearly, $z\in C$ is not
a problem in either case, because of \eqref{0fg}.)

Therefore, to confirm Hoste's conjecture on $\br{D'}$, it suffices
that (1) the constant polynomial $1$ is among the $\tl\Dl_i$, and
that (2) all conditions
\[
\R \bigl(\Dl_i(z)\Dl_j(\bar z)\bigr)>0\quad\mbox{%
for $z\in C$ and all $1\le i<j\le n$}\,,
\]
are sastisfied, with $\Dl_i$ as in \eqref{DDL}. Then we can conclude
the property that the number of roots $z$ left of $C$ is $0$ by
induction on, say, $\Dl_{[0]}$. We call the above two conditions
in the following the \em{Rouch\'e test}.

Now the number of checks grows quadratically in $n$, 
which makes the reduction of the number of $\Dl_i$ 
much more urgent than before. We already discarded 
zero polynomials. Clearly, duplicate $\Dl_i$, even up
to scalars, can be discarded too, and one easily sees that 
in fact this is true also up to multiplication by powers
of $1-t$. So we divided all factors $1-t$ out of the
$\Dl_i$ prior to discarding duplicates. Since the origin 
of the values $n_i$ in \eqref{DDL} will be irrelevant for 
the rest of the calculation, assume w.l.o.g. that $\tl
\Dl_i$ are the so reduced polynomials (indivisible by 
$t-1$). By symmetry one sees then that $\Md\tl\Dl_i$ is 
an even number (between 0 and 8).

For the purpose of testing the 304 remaining $\sg=2$ generators, 
we used MATHEMATICA\TM{} \cite{Wolfram} and the Rouch\'e test with 
the so far reduced sets of $\tl\Dl_i$ for each $D'$. In determining 
the $\tl\Dl_i$, we also used the restrictions, obtained in the
enhanced positive zero test, on which $\sim$-equivalence classes
of the generator twists are necessary. After an hour and 15 minutes 
MATHEMATICA reported success on all 304 generators $D'$, thus 
completing the $\sg=2$ case.

For $\sg=0$, however, there are about 60,000 generators (with
no restrictions on $\sim$-equivalence classes to twist at), and
MATHEMATICA\TM{}, albeit very intelligent, is known to pay for
its intelligence with its speed. This made further reductions
of the set of $\tl\Dl_i$ necessary. 

In this vein, note that if some $\tl\Dl_i$ lies in the convex
hull of others, it is also redundant. (It is enough that 
$a_i$ remain positive in \eqref{star++}; integrality is 
inessential.) Unfortunately, the effort of identifying exactly 
convex linear combinations as a pre-reduction to the Rouch\'e 
test would take not less effort than the test itself. 

In practice we did the following for polynomials $\tl\Dl_i$ of 
degrees $d=2,4$. We determined the minimal and maximal ratios
of $(\tl\Dl_i)_{[j]}/(\tl\Dl_i)_{[0]}$ for $1\le j\le d/2$.
Then we replaced the set of all $\tl\Dl_i$ of degree $d$ by
a collection of $2^{d/2}$ polynomials given by taking for each 
coefficient $\Dl_{[j]}$ with $1\le j\le d/2$ once the minimal, 
and once the maximal ratio, setting $\Dl_{[0]}=1$ (and completing 
the other coefficients by symmetry $\Dl_{[d-j]}=\Dl_{[j]}$).

For $d=2$ this replacement is equivalent (under convex linear
combinations), so it is a genuine simplification, but this did
not speed up the check sufficiently. Contrarily, for $d=4$ the
new polynomials have a larger convex hull, which augments the
risk of failure of the test. However, using only the 4 polynomials
instead of the dozens of others made the test about 10 times
faster. It still took a couple of days to carry out. 

MATHEMATICA reported success of the last, fastest, form of the test 
on all $\sg=0$ generators, except only one, the knot $12_{1039}$
of \cite{KnotScape}. There the test worked out at least using the
full list of polynomials $\tl\Dl_i$ of degree 4, rather than their
4 substitutes.

With this, after a few weeks of work, the check of
Hoste's conjecture in genus $\le 4$ was completed. \qed

\begin{rem}
After we finished our proof, we learned of Ozsv\'ath and Szab\'o's
inequality \eqref{ppOS__}. It would replace the r.h.s. of
\eqref{star+} by 2, thereby settling $g=3$, $\sg=0$ directly,
and sharpening the positive zero test for $g=4$. Also,
\eqref{c2} would be excluded straightforwardly.
Still both the nature of the tests and the extent of
calculation would not be reduced significantly.
\end{rem}

Note that \eqref{star4+} can be extended easily to
$\big|\Dl_{[i]}/\Dl_{[0]}\big|$ for $i>1$. In particular we
have

\begin{prop}\label{ppbd}
For given $i\ge 1$ and $g>0$, the ratios
$\bigl\{\,\big|\Dl_{[i]}/\Dl_{[0]}\big|\,:\,g(K)=g\,\bigr\}$
are bounded. \qed
\end{prop}

Such a property is not at all evident without using generators.
Of course, it is possible, using e.g. corollary \reference{cor4.1},
to give explicit estimates in terms of $g$. Moreover, we gain a
practical method to calculate (even sharp) upper and lower bounds
on this ratio for given genus when the generators are available.

\begin{exam}
For example, for $g=4$, $\sg=0$ prime alternating knots we found
\[
3\le -\frac{\Dl_{[1]}}{\Dl_{[0]}}\le 20\,,\quad
5\le  \frac{\Dl_{[2]}}{\Dl_{[0]}}\le 122\,,\quad
7\le -\frac{\Dl_{[3]}}{\Dl_{[0]}}\le 333\,,\quad\mbox{and}\quad
8\le  \frac{\Dl_{[4]}}{\Dl_{[0]}}\le 461\,,
\]
and these bounds are the best possible (except I did not check 
if equalities are attained, or the inequalities are always strict).
They were found by unifying the analogous bounds that can be
obtained for each generator separately. Using latter, we also
estimated for each generator the Jensen integral for the Euclidean
Mahler measure $M(\Dl)$, the product of norms of zeros of $\Dl$
outside the unit circle (see \cite{STW}). In combination, we have 
\[
M(\Dl)\le 638.21
\]
for the Euclidean Mahler measure of $g=4$, $\sg=0$ prime alternating
knots. (This estimate may not be very sharp.) 
\end{exam}

Such data were obtained in an earlier attempt at 
this case of Hoste's conjecture, prior to using Rouch\'e's 
theorem. They can be extended with the proper (not small, 
but manageable) amount of calculation to all $g=4$ alternating 
knots. Note again that the special alternating knots (here 
$\sg=8$) are mostly clear from this point of view. Even for 
arbitrary genus $g$ it is easy to deduce from the preceding 
discussion that
\[
0\le (-1)^j\frac{\Dl_{[j]}}{\Dl_{[0]}}<{2g\choose j}\,,
\]
with the right inequality the best possible and the left bound 
not improvable beyond $1$. (The bound is exactly $1$ for $j=1$ 
by the work in \cite{MS}; apparently this is true also for higher
$j$, and would follow from Fox's Trapezoidal Conjecture, discussed
in \cite{HS,spec} and right below.) Also $M(\Dl)\equiv 1$.

By modifying the contour $C$ in Rouch\'e's test one could
obtain further more specific information about location of
$\Dl$-zeros on certain families of knots.

\subsection{The $\log$-concavity conjecture}

The last open problem we consider is the $\log$-concavity conjecture.
This conjecture, made in \cite{spec}, states:

\begin{conj}\label{lcc}
Call a polynomial $X\in \bZ[t^{\pm 1}]$ to be \em{$\log$-concave},
if $X_{[k]}:=[X]_k$ are $\log$-concave, i.e.
\begin{eqn}\label{_lcc}
X_{[k]}^2\ge X_{[k+1]}X_{[k-1]}
\end{eqn}
for all $k\in \bZ$. Let
$n(K)$ be the number of components of a link $K$.
\def\labelenumi{(\arabic{enumi})}
\begin{enumerate}
\item\label{lc1} If $K$ is an alternating link,
then $t^{(1-n(K))/2}\Dl_K(t)$ is $\log$-concave.
\item\label{lc2} If $K$ is a positive link, then
$t^{(1-n(K))/2}\nb_K(\sqrt{t})$ is $\log$-concave.
\end{enumerate}
\end{conj}
We will refer to both properties as `$\Dl$-$\log$-concavity'
resp. `$\nb$-$\log$-concavity'.

We remind that by Crowell-Murasugi \cite{Crowell,Murasugi2}
for part \ref{lc1}, when $\Dl_K(t)$ is normalized so that
$\md\Dl_K=0$ and $[\Dl_K]_0>0$, the polynomial $\Dl_K(-t)$
is positive, i.e. all its coefficients are non-negative.
The same property holds by Cromwell \cite{Cromwell} for part
\reference{lc2} and $t^{(1-n(K))/2}\nb_K(\sqrt{t})$.

Part \reference{lc1} is a natural strengthening of \em{Fox's
Trapezoidal conjecture}.

\begin{conj}\label{Fxc}(Fox)
If $K$ is an alternating knot and $\Dl_K$ is normalized so
that $\md\Dl_K=0$, then there is a number $0\le n\le g(K)$
such that for $\Dl_{[k]}:=\big|[\Dl_K]_k\big|$ we have
\begin{eqn}\label{Fxq}
\begin{array}{cl}
\Dl_{[k]}>\Dl_{[k-1]} & \mbox{for $k=1,\dots,n$}, \\
\Dl_{[k]}=\Dl_{[k-1]} & \mbox{for $k=n+1,\dots,g(K)$}\,.
\end{array}
\end{eqn}
We call polynomials of this form \em{trapezoidal}. A similar
property can be conjectured for (non-split) alternating links,
replacing $g(K)$ by $\Br{\spn\Dl_K/2}$.
\end{conj}

The Trapezoidal conjecture has received some treatment in the
literature, being verified for 2-bridge knots \cite{Hartley}
(see also \cite{Burde}) and later for a larger class of
alternating algebraic knots \cite{Murasugi6}. More recently,
Ozsv\'ath and Szab\'o used their knot Floer homology \cite{OS}
to derive a family of linear inequalities on the coefficients of
$\Dl$ for an alternating knot.

\begin{prop}(Ozsv\'ath and Szab\'o)\label{ppOS}
Let $K$ be an alternating knot of signature $\sg=\sg(K)$, and genus
$g=g(K)$, and let $\Dl=\Dl(K)$ be normalized so that $\md\Dl=0$
and $[\Dl]_0>0$. Then for each integer $s\ge 0$,
\begin{eqn}\label{ppOS_}
(-1)^{s+g}\,\sum_{j=1}^{g-s}\,j\cdot [\Dl]_{s+g+j}\,\le\,
(-1)^{s+\sg/2}\,\max \left(0,\BR{\frac{|\sg|-2s}{4}}\right)\,.
\end{eqn}
\end{prop}

(Note that when $\Dl$ is normalized so that $\Dl(1)=1$ and $\Dl(t)=
\Dl(1/t)$, then the sign of $[\Dl]_{\pm g}$ is $(-1)^{g+\sg/2}$.)

For genus $g=2$ the inequalities \eqref{ppOS_} are very similar to
(and slightly stronger than) the Trapezoidal conjecture, and for
$s=g-2$ yield
\begin{eqn}\label{ppOS__}
-[\Dl]_1\,\ge\,2[\Dl]_0+\,\left\{\,\begin{array}{cl}
-1 & \mbox{if $|\sg|=2g$} \\
+1 & \mbox{if $|\sg|=2g-2$} \\
 0 & \mbox{otherwise}
\end{array}\,\right.,
\end{eqn}
which settles in \eqref{Fxq} the case $k=1$ (for knots). For $g>2$ and
$s<g-2$, the inequalities \eqref{ppOS_} do not relate directly to the
Trapezoidal conjecture.

In-Dae Jong \cite{Jong} has proved independently the Trapezoidal
conjecture up to genus 2 using the generator description in theorem
\reference{thgen}, and observed that for genus 2 the $\log$-concavity
of $\Dl$ easily follows from trapezoidality.

In \cite{spec} we showed part \reference{lc2} of conjecture \ref{lcc}
for special alternating knots, by giving a new proof of theorem
\reference{thsg} in this case.

We have now

\begin{theorem}\label{tcc}
Both part \ref{lc1} and \ref{lc2} of conjecture \ref{lcc} (and
therefore also Fox's conjecture) hold for knots of genus at most 4.
\end{theorem}

The preceding work easily settles the positive case.

\proof[of part \ref{lc2} of theorem \ref{tcc}]
{}From theorem \ref{thsg} one easily sees that all zeros of
$\nb(\sqrt{t})$ are real if $\sg(K)\ge 2g(K)-2$. With the explanation
in \cite{spec} we have then $\log$-concavity of $\nb$. Thus we
need to check just the polynomial of $14_{45657}$ directly. \qed

Part \reference{lc1} is much harder to check, since we have no
easy sufficient conditions at hand. We explain how we proceeded.

\proof[of part \reference{lc1} of theorem \reference{tcc}]
We had to test all generators up to genus 4. A small
defect of $\log$-concavity is that it is not straightforwardly
preserved under product (as is trapezoidality), so composite
generators must be handled either.

Due to its volume, the check had to be optimized strongly.
Consider a particular generator diagram $D$. We use the notation
of \eqref{DDL} in the proof of theorem \ref{thHo}. Since one can
recover $\Dl_i$ by
\begin{eqn}\label{DDD}
\Dl_i\,=\,(1-t)^{g(D)-\Md\tl\Dl_i/2}\tl\Dl_i
\end{eqn}
even when $\tl\Dl_i$ is taken up to powers of $1-t$, we can assume
w.l.o.g. that we divide out all factors $1-t$ in $\tl\Dl_i$ so that
$\tl\Dl_i(1)\ne 0$, and $[\tl\Dl_i]_0>0$. First we designed the
calculation of $\tl\Dl_i$ so that diagrams $D_i$ of genus 0 or
split ones (where $\tl\Dl_i$ is zero or a scalar) are detected
and discarded in advance. For the calculated $\tl\Dl_i$
again those of small degree are treated extra. For degree 0 it
is enough to keep the polynomial $1$. For degree 2 we need only
two polynomials, those whose ratio $[\tl\Dl_i]_1/[\tl\Dl_i]_0$ is
minimal and maximal.

The $\tl\Dl_i$ of degree 4 lead to the main
reduction. In this case we have a number of points
\begin{eqn}\label{plll}
\bigl(\,[\tl\Dl_i]_1/[\tl\Dl_i]_0,\,
  [\tl\Dl_i]_2/[\tl\Dl_i]_0\,\bigr)\,\in\,\bR^2
\end{eqn}
in the plane. We need to keep only the \em{extremal} ones, i.e.
those which span the convex hull, and discard the others (which
lie in the convex hull). In practice it turned out that this
reduced the number of $\tl\Dl_i$ of degree 4 up to a factor of
150, and the number of all $\tl\Dl_i$ up to a factor of 11.5.

With the so reduced set of $\tl\Dl_i$, we recover $\Dl_i$ by
\eqref{DDD} and must test for $1\le k\le g(D)$
\begin{eqn}\label{(st)}
[\Dl_i]_{k}^2\ge [\Dl_i]_{k-1}[\Dl_i]_{k+1}\,,
\end{eqn}
and the set of such conditions for convex linear combinations of $\Dl
_i$ and $k<g(D)$. (The case $k=g(D)$ for linear combinations follows
directly from \eqref{(st)} and symmetry of $\Dl_i$.) Latter again turns
into deciding whether the convex polygon $\Sg\subset \bR^2$ spanned by
points 
\[
\bigl(\,-[\Dl_i]_{k+1}/[\Dl_i]_k,\,[\Dl_i]_{k+2}/[\Dl_i]_k\,\bigr)
\,\in\,\bR^2
\]
is contained in the region $R=\{\,(x,y)\,:\,x\ge 0,\,y\le x^2
\,\}$. Clearly it is enough to test the inclusion $\partial\Sg
\subset R$ for the boundary $\pa\Sg$ of $\Sg$. This boundary
is a certain collection of edges between extremal points, and
can be determined similarly to the convex hull reduction of the
$\tl\Dl_i$ of degree 4. Finally, we have to check that the function
$\bR^2\ni (x,y)\mapsto x^2-y$ is non-negative on these edges.
In practice, we made the tests by an $\eps$ stricter, to avoid
floating point number deviations.

For a handful of pairs $(D,k)$ we have $\Sg\not\subset R$. Let us call
such pairs and the corresponding generators $D$ \em{irregular}. For
genus 3 the irregular generators are $7_1$ and $8_5$. They still do
not mean that we have a counterexample, since the points coming from
$\Dl$ polynomials of diagrams in the series of $D$ form only a subset
of $\Sg$. In such cases again we check $\log$-concavity of the
generator
polynomial $\Dl(D)$, and then apply the described test to diagrams
$D'$ obtained by one $\bt$ twist from $D$. Iterating this a few times
shows that only one explicit one-parameter family, of $(3,3,2m)$
pretzel diagrams ($m>0$), remains, which can be verified directly. 

Among prime generators of genus 4, only 3 are irregular, and they are
$9_1$, $10_2$ and $10_{46}$. Using the twisted diagram check, they
leave out only one one-parameter family each, consisting of the 
rational knots $2m-1,8$ and the pretzel knots $P(2m,1,7)$ and $P(2m,
3,5)$ (with $m>0$) resp. Again latter are easily settled directly.
Among composite generators up to genus 4, no irregular ones occur.

Despite its theoretical simplicity, the check involved much skillful
implementational effort. Among others, we decided to write an own
(optimized) procedure for calculating the Alexander polynomial, instead
of using the functionality of \cite{KnotScape}, which applies a
substitution to the skein polynomial obtained by the Millett-Ewing
algorithm.

With all optimization, we had a speed-up by a factor of almost
5 in comparison to the initial implementation. Finally, while
genus 3 could be done in about 8 minutes, we could reduce the effort
for genus 4 only to about a second for the simplest generators,
and up to 30 seconds for the most complicated ones. This meant that
the verification of the full list of generators had to be done in
a few dozen parts, a couple of them running over several weeks. \qed

\begin{rem}
Note that for Fox's conjecture the test could be simplified,
since dealing with linear combinations of the $\Dl_i$ becomes
obsolete. (The difficulty of $\log$-concavity is that \eqref{_lcc}
is not a linear condition.) Still this would not have reduced the
magnitude of calculation significantly, since the most time
is required for the calculation of the $\Dl_i$.
\end{rem}

With some effort to adapt the computation to links, we obtained
the following outcome.

\begin{theorem}\label{lcL}
The $\Dl$-$\log$-concavity (and hence Fox's conjecture \ref{Fxc})
holds for $n$-component alternating links $L$ with $\spn \Dl(L)+n
\le 9$ (or equivalently $g(L)+n\le 5$). 
\end{theorem}

\proof We recurred the calculation from generators of these links to
knot generators by using claspings \eqref{y}. The check turned out
far more laborious than for knots, but there is no real difference
in the method, so we skip details to save space. \qed

\begin{rem}
It is suspectable that in conjecture \reference{Fxc} we always have
for a knot $K$
\begin{eqn}\label{ineqn}
n\,\ge\,g(K)-|\sg(K)/2|\,.
\end{eqn}
It is easily observed that it suffices to check \eqref{ineqn} for the
generators, and that prime ones are enough. Moreover, the inequality is
trivial for special alternating generators, so consider only the non-%
special ones. We tested that all such generators up to genus 4 satisfy
\eqref{ineqn}. (For $g=2$, and also $n=0$, this property is implied
by Ozsv\'ath--Szab\'o's inequality \eqref{ppOS__}.) Actually, very few
generators have $n<g(K)$: only 4 of the non-special prime generators
of genus 3 and 30 of those of genus 4 have coefficients showing a true
``trapez'', rather than a ``triangle''. We may mention that we did not
test any analogon of \eqref{ineqn} for links.
\end{rem}

\begin{rem}
A minor addition to part \reference{lc1} of the $\log$-concavity
conjecture \ref{lcc} is that if we have equality in \eqref{_lcc},
then all 3 of $X_{[k]}$, $X_{[k\pm 1]}$ are equal. Our test verified
this extra property, too, for knots up to genus 4 and the links in
theorem \ref{lcL}.
\end{rem}

\subsection{\label{rcs}Complete linear relations by degree}


Let us remark that our proof of the Fox conjecture can be conceptually
extended. While it is clear that one can verify a particular sort of
linear inequality between coefficients of $\Dl$ for given genus, we
can in fact determine the \em{complete list} of such inequalities.

\begin{defi}
A \em{(convex) polytope} $\Sg$ is the convex hull of a finite number
of points in $\bR^n$ and \em{rays} $\bR_+\cdot \vec b$ for $\vec b
\in\bR^n$. Alternatively, one may describe a polytope by a set of
linear inequalities (or as intersection of half-spaces of $\bR^n$).
If no rays occur, we call $\Sg$ \em{bounded}. If at most one point
is there (with no point meaning that the vertex is the origin), then
$\Sg$ is a \em{cone}. A \em{facet} of $\Sg$ is a codimension-1 subset
of the boundary $\bd\Sg$ of $\Sg$, given by intersection of $\Sg$ 
with a hyperplane.
\end{defi}

There seems no uniform rule for the usage of `polytope' and
`polyhedron' in the literature. Our attitude here is to avoid latter
term throughout. The word `convex' will be usually omitted since all
polytopes we will deal with are such. Also note that, in contrast
to us, many authors consider a polytope to be bounded per definition.
Our terminology needs one more important clarification.

\noindent{\bf Convention.} In the following the term \em{convex hull}
$\conv(D)$ of a set $D\subset \bR^n$ is meant as \em{the closure} of
the set 
\[
\Bigl\{\,\sum_{i=1}^m\,\lm_ix_i\,:\,\sum_{i=1}^m\,\lm_i=1,\, \lm_i\ge 0
  ,\,x_i\in D\,\Bigr\}\,.
\]
The reason for taking closures is easily seen from the
example of a line and a point
(outside the line) in the plane. This phenomenon wildens in higher
dimension: for a non-closed convex polytope only parts of some
facets would be there. Formulated in terms of the linear inequalities
describing the polytope, certain points for which some inequality
is exact (i.e. an equality) are allowed, while others are not. The
discussion of these cases puts no reasonable use-over-effort ratio
in prospect. So we assume \em{all convex sets are closed} (and the 
inequalities describing the polytope are always non-strict).

We shall treat here genus 3 exemplarily. For technical reasons
we formulate the result in terms of $\nb$.

\begin{prop}
The signed coefficients $\nb_k=[\nb]_k\cdot \sgn([\nb]_6)$ of
$\nb(K)$ for an alternating knot $K$ of genus 3 satisfy the
following 6 inequalities:
\begin{eqn}
\label{=} \begin{array}{rcl@{\qquad}rcl@{\qquad}rcl}
\ds \nb_4 -\frac{\nb_2+5\nb_6}{2} & \le &\ds \frac{1}{2}\,, &
\ds \nb_4 -\nb_2-3\nb_6 & \le &\ds 2\,,&
\ds \nb_4 +3\nb_2+9\nb_6 & \ge &\ds -4\,,  \\[3mm]
\ds \nb_4 +\frac{\nb_2}{2}+4\nb_6 & \ge &\ds -\frac{1}{2}\,, &
\ds \nb_4 -\nb_2+7\nb_6 & \ge &\ds -7\,, &
\ds \nb_2 -3\nb_6 & \le &\ds 4\,. \\[1mm]
\end{array}
\end{eqn}
This set of linear inequalities is complete up to constants, in the
sense that any other linear inequality valid for all (or all but
finitely many) alternating genus 3 knots is, up to a worse absolute
term, a consequence of those above. 
\end{prop}

We will explain later how to remove the defect of the constants, but
it requires a bit more sophisticated calculation.

{\def\pt#1{{\picfilledcircle{#1}{0.1}{}}}
\def\mydraw#1{\pt{#1 d 2{3 1 r :}R x}}

\proof If, for each geneartor $D$, we take all $\tl\Dl_i\ne \Dl(D)$,
then all have degree at most 4. (We will deal with the $\Dl(D)$
later.) Let $\nb_i(z)$ be given by
\[
\nb_i(t^{1/2}-t^{-1/2})=\Dl_i(t)\cdot t^{-\Md\Dl_i/2}\,.
\]
The reason for using $\nb_i(z)$ is that the factor $1-t$ in
\eqref{DDD} turns into $z$, and so we have just a degree shift,
which is faster on the computer. (Thus, in fact, we used $\nb_i$
instead of $\Dl_i$ almost throughout the calculation for theorem
\ref{tcc}.) Now $\Md\tl\Dl_i=\spn \nb_i$ (recall \S\ref{SLP}).

Now we can do the previous sort of determination of extremal
polynomials $\nb_i$ of span 4,
but this time with all generators taken together. We collect
again the corresponding (extremal) points in the plane
\begin{eqn}\label{pl4}
(x,y)=\bigl(\,[\nb_i]_2/[\nb_i]_6,\,
  [\nb_i]_4/[\nb_i]_6\,\bigr)\,\in\,\bR^2\,,
\end{eqn}
as in \eqref{plll}, together with the two points (with
$x=0$) describing the convex hull (which now is an interval)
for the $\nb_i$ of span 0 and 2. Here is the result:
\begin{eqn}\label{pdg}
\diag{7mm}{10}{10}{
  \pictranslate{5 5}{
\picmultigraphics[rt]{2}{90}{
  {\piclinewidth{90}
   \picmultigraphics{11}{0 1}{\picline{-5 -5}{5 -5}}
  }
  \picvecline{-5 0}{5 0}
}
\picvecline{0 -5}{0 5}
\picvecline{0 -5 x}{0 5 x}
\picputtext{4.7 0.3}{$x$}
\picputtext{4.7 -0.3 x}{$y$}
\picfillgraycol{0.9 0.0 0.1}
\mydraw{0 -4 1 }
\mydraw{0 5 2 }
\mydraw{-3 0 1 }
\mydraw{-2 -7 2 }
\mydraw{-2 -3 1 }
\mydraw{-2 1 1 }
\mydraw{-1 2 1 }
\mydraw{1 3 1 }
\mydraw{2 -9 2 }
\mydraw{2 -5 1 }
\mydraw{3 -4 1 }
\mydraw{3 -2 1 }
\mydraw{3 4 1 }
\mydraw{4 7 2 }
}}
\end{eqn}
The convex hull $\tl\Sg$ is a hexagon, desribed by equations (one for
each boundary edge) in which the r.h.s. of an equation in \eqref{=}
is set to 0, $\nb_2$ is replaced by $x$, $\nb_4$ by $y$, and $\nb_6$
by $1$. To obtain \eqref{=}, we need then to calculate the minimal
or maximal value of the l.h.s. over all generator polynomials
$\Dl(D)$ (which we previously omitted from the calculation).

To see completeness, note that for each $\Dl_i\ne \Dl(D)$ one has
a sequence of diagrams $\bar D_j\in\br{D}$ with 
\[
\lim_{j\to\infty}\,\frac{\Dl(\bar D_j)}{[\Dl(\bar D_j)]_0} =
\frac{\Dl_i}{[\Dl_i]_0}\,,
\]
so the set of points \eqref{pl4} is not reducible. \qed

\begin{figure}[htb]
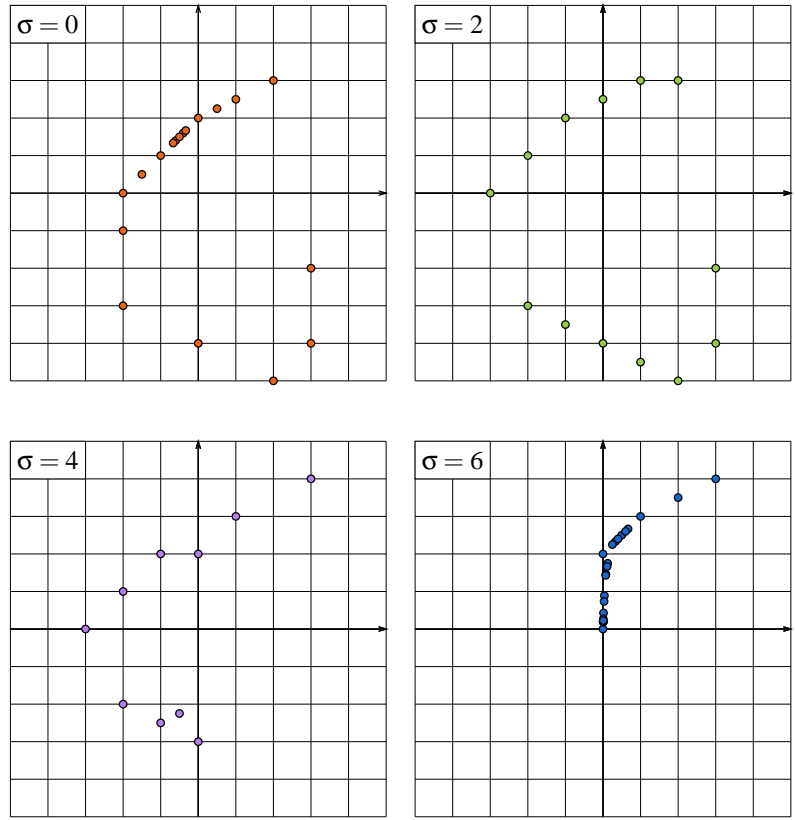

\[
\diag{5mm}{10}{10}{
\pictranslate{5 5}{
\picmultigraphics[rt]{2}{90}{
  {\piclinewidth{90}
   \picmultigraphics{11}{0 1}{\picline{-5 -5}{5 -5}}
  }
  \picvecline{-5 0}{5 0}
}
{\piclinewidth{90}\picfilledbox{-4 4.5}{2 1}{$\sg=0$}}
\picfillgraycol{0.9 0.4 0.1}
\mydraw{0 -4 1 }
\mydraw{0 2 1 }
\mydraw{-3 1 2 }
\mydraw{-3 7 5 }
\mydraw{-2 -3 1 }
\mydraw{-2 -1 1 }
\mydraw{-2 0 1 }
\mydraw{-2 4 3 }
\mydraw{-2 8 5 }
\mydraw{-1 1 1 }
\mydraw{-1 3 2 }
\mydraw{-1 5 3 }
\mydraw{2 -5 1 }
\mydraw{2 3 1 }
\mydraw{2 5 2 }
\mydraw{2 9 4 }
\mydraw{3 -4 1 }
\mydraw{3 -2 1}
}}
\quad
\diag{5mm}{10}{10}{
\pictranslate{5 5}{
\picmultigraphics[rt]{2}{90}{
  {\piclinewidth{90}
   \picmultigraphics{11}{0 1}{\picline{-5 -5}{5 -5}}
  }
  \picvecline{-5 0}{5 0}
}
\picvecline{0 -5}{0 5}
\picvecline{0 -5 x}{0 5 x}
{\piclinewidth{90}\picfilledbox{-4 4.5}{2 1}{$\sg=2$}}
\picfillgraycol{0.6 0.8 0.3}
\mydraw{0 -4 1}
\mydraw{0 5 2}
\mydraw{-3 0 1}
\mydraw{-2 -7 2}
\mydraw{-2 -3 1}
\mydraw{-2 1 1}
\mydraw{-1 2 1}
\mydraw{1 3 1}
\mydraw{2 -9 2}
\mydraw{2 -5 1}
\mydraw{2 3 1}
\mydraw{3 -4 1}
\mydraw{3 -2 1}
}}
\]

\[
\diag{5mm}{10}{10}{
\pictranslate{5 5}{
\picmultigraphics[rt]{2}{90}{
  {\piclinewidth{90}
   \picmultigraphics{11}{0 1}{\picline{-5 -5}{5 -5}}
  }
  \picvecline{-5 0}{5 0}
}
\picvecline{0 -5}{0 5}
\picvecline{0 -5 x}{0 5 x}
{\piclinewidth{90}\picfilledbox{-4 4.5}{2 1}{$\sg=4$}}
\picfillgraycol{0.7 0.5 0.9}
\mydraw{0 -3 1}
\mydraw{0 2 1}
\mydraw{-3 0 1}
\mydraw{-2 -9 4}
\mydraw{-2 -5 2}
\mydraw{-2 -2 1}
\mydraw{-2 1 1}
\mydraw{-1 2 1}
\mydraw{1 3 1}
\mydraw{3 4 1}
}}
\quad
\diag{5mm}{10}{10}{
\pictranslate{5 5}{
\picmultigraphics[rt]{2}{90}{
  {\piclinewidth{90}
   \picmultigraphics{11}{0 1}{\picline{-5 -5}{5 -5}}
  }
  \picvecline{-5 0}{5 0}
}
\picvecline{0 -5}{0 5}
\picvecline{0 -5 x}{0 5 x}
{\piclinewidth{90}\picfilledbox{-4 4.5}{2 1}{$\sg=6$}}
\picfillgraycol{0.1 0.4 0.8}
\mydraw{0 0 1}
\mydraw{0 2 1}
\mydraw{1 3 1}
\mydraw{1 5 2}
\mydraw{1 7 3}
\mydraw{1 9 4}
\mydraw{1 14 8}
\mydraw{1 15 9}
\mydraw{1 19 13}
\mydraw{1 20 14}
\mydraw{1 24 27}
\mydraw{1 24 128}
\mydraw{1 25 34}
\mydraw{1 26 134}
\mydraw{1 28 140}
\mydraw{1 30 70}
\mydraw{1 30 145}
\mydraw{1 31 115}
\mydraw{1 31 137}
\mydraw{2 8 3}
\mydraw{2 12 5}
\mydraw{3 4 1}
\mydraw{3 13 5}
\mydraw{4 7 2}
}}
\]
\caption{\label{FG5}Extremal points of the polygons $\tl\Sg$
for genus $g=3$ and given signature.}
\end{figure}
}


Let us write more generally $\tl\Sg_g\subset \bR_+^{g-1}$ for the
convex sets obtained as in the preceding proof for generators of
genus $g$. For $g>3$ 
the procedure remains largely the same. Note that these polytopes
will always be bounded, due to proposition \ref{ppbd}.

It is natural to restrict oneself also to knots of given signature,
and figure \ref{FG5} shows the pictures corresponding to \eqref{pdg}
for $g=3$. For example, for $\sg=6$ (the special alternating knots),
one derives the (complete up to constants linear) inequalities:
\begin{eqn}\label{==}
\nb_4-\frac{4}{3}\nb_2\,\ge\,-3\,,\qquad
\nb_2\,\ge\,3\,,\qquad
\nb_4-\nb_2-2\nb_6\,\le\,-2\,,\qquad
\nb_4-\frac{\nb_2+5\nb_6}{2}\,\le\,-\frac{1}{2}\,.
\end{eqn}
Let us denote the corresponding (bounded) polytopes as $\tl\Sg_{g,\sg}$.

The graphics of figure \ref{FG5} were displayed to show the difference
to the properties considered unrelatedly to our approach. Roughly
(i.e.\ ignoring the additive terms independent on $\nb_k$), the
conditions \eqref{ppOS_} of Ozsv\'ath--Szab\'o (only $s=0,1$ are
relevant here) would lead to the restrictions
\begin{eqn}\label{OSS}
y\le 4\mbox{\quad and\quad}y\le \frac{x+6}{2}\,,
\end{eqn}
while the trapezoidal inequalities \eqref{Fxq} (with the options
`$>$' and `$=$' merged into `$\ge$') transcribe into 
\begin{eqn}\label{trS}
y\,\le\,\min\,\left(\,5,\,\frac{x+5}{2},\,\frac{x}{3}+3\,\right)\,.
\end{eqn}
These inequalities are (as should be) easily seen to follow from
the ones describing our polygons $\tl\Sg$, but as much weaker
consequences. Note in particular that, even taken together, \eqref{OSS}
and \eqref{trS} restrict the possible $(x,y)$ only to an
infinite region of the plane. Also, the dependence of the shape of
$\nb$ on the signature becomes much more apparent in our polygons
than in \eqref{ppOS_}. This shows how much further information
one can obtain with our method, if one considers a fixed genus.
Still \eqref{Fxq} or \eqref{ppOS_} may define (without abusing
constants) a facet of $\Sg$. We will test this explicitly in
below for $g=3$.


If one wants to remove the inaccuracy up to constants, we have the following
statement, in general dimension, which is worth taking record of.

\begin{theorem}\label{thSgz}
The (closed) convex hull $\Sg_g$ of the set of Conway polynomials of
alternating knots of given genus $g$, or $\Sg_{g,\sg}$ given genus and
signature, is a convex polytope in $\bR^g$.
\end{theorem}

\proof To determine $\Sg_g$, one must find the convex hull $\tl\Sg_D$ of
the points \eqref{pl4} for each generator $D$ separately, take the cone
$\bR_+\cdot (\tl\Sg_D\times\{1\})$ of $\tl\Sg_D$ in $\bR^g$, and
translate this cone to have the summit $\nb(D)$. Then $\Sg_g$ is the
convex hull of the union of all these translated cones taken
over all generators $D$. This can be rewritten as
\begin{eqn}\label{chr}
\Sg_g=\bR_+\cdot (\tl\Sg_g\times\{1\})+\conv\left(\,\bigl\{\,
\nb(D)\,:\,\mbox{$D$ is a generator of genus $g$}\,\bigr\}\,\right)
\end{eqn}
(with the '$+$' being the Minkowski set-theoretic sum
$A+B\,=\,\{\,a+b\,:\,a\in A,\,b\in B\,\}$).
The result $\Sg_g$ is a(n unbounded) convex polytope, as follows
from the Farkas-Minkowski-Weyl Theorem for convex polytopes (see e.g.
\cite{Sch86}, Corollary 7.1a). \qed

\begin{rem}
Note that one can replace `alternating' by `positive' in theorem
\reference{thSgz}. (For knots of given genus and signature a bit
extra argument is needed.) We chose not to elaborate on the positive
case here, though.
\end{rem}

Let us call a facet $F\subset\pa\Sg$ of $\Sg\subset \bR^g$ to
be \em{$n$-open} if it contains an affine-translated copy of
$\bR_{+}^{n}\times \{0\}^{\times g-n}$, but not one of $\bR_
{+}^{n+1}\times \{0\}^{\times g-n-1}$. Say $F$ is \em{largest}
if it is $(g-1)$-open. Then inequalities like the above \eqref{=}
or \eqref{==}, said to be complete up to constants, are in
fact those for the largest facets of $\Sg_g$.

Theorem \ref{thSgz} clearly also gives a practical way to determine
$\Sg_g$, though it requires to work in one dimension up in comparison
to $\tl\Sg_g$. Latter would likely be quite more complicated
to describe precisely for high genus. (In calculating $\tl\Sg_3$
above, we took advantage of the 2-dimensional convex hull algorithm
we had implemented for theorem \reference{tcc}.) Nevertheless this
is feasible (even by hand) for genus 2, and carried out by In-Dae
Jong \cite{Jong2}. The result can be stated as follows.

\begin{theorem}
The complete set of linear inequalities satisfied by the coefficients
$\nb_i=[\nb]_i\cdot \sgn([\nb]_4)$ of Conway polynomials of alternating
genus 2 knots of given signature is:
\[
\nb_4\ge 1 \mbox{\ \ (for all $\sg$),\es\  and\ \ \ }
\left\{\ \begin{array}{r@{\,\,}c@{\,\,\,}ll}
-2\nb_4-1 & \le\, \nb_2\,\le & \nb_4+1   & \mbox{for\ }\sg=0\,, \\
-2\nb_4+1 & \le\, \nb_2\,\le & 2\nb_4-1 & \mbox{for\ }|\sg|=2\,, \\
        2 & \le\, \nb_2\,\le & 2\nb_4+1  & \mbox{for\ }|\sg|=4\,
\end{array}\ \right\}\,.
\]
\end{theorem}

For genus $3$ one must use a computer. There is software for converting
a convex-hull (vertex-ray) representation of a polytope like
\eqref{chr} into a half-space intersection (linear inequalities)
representation. We used the programs {\tt cdd/cdd+} of Komei Fukuda
\cite{Fukuda} and {\tt lrs} of David Avis \cite{Avis} to find a
linear inequality representation for $\Sg_{3,\sg}$ from \eqref{chr}
and the calculation of figure \ref{FG5}. (The programs do essentially
the same, but use different algorithms, and we applied them both for
consistency security.)

The description of $\Sg_{3,\sg}$ is shown in the following table in
columns 2 and 3. An entry of the form `$m$v$n$r' in the third column
means that $\Sg_{3,\sg}$
is the convex hull of $m$ vertices and $n$ rays. (For example,
`1v2r' is a planar angle segment, and `2v1r' is a plane half-strip.)
The second column gives the number of facets of (or minimal number
of linear inequalities describing) $\Sg_{3,\sg}$.

The following columns show the efficiency of the conditions
\eqref{ppOS_} (with 'OS$n$' standing for the case $s=n-1$)
and \eqref{Fxq} (with 'T$m$' standing for the inequality obtained
by putting $k=m$ and joining the two alternatives into a non-strict
inequality, without regard to the other $k$) for $g=3$ and given
signature $\sg$. An entry of the form `$m$v$n$r' is as in column 2
and describes here the intersection of $\Sg_{3,\sg}$ with the
hyperplane defined by the inequality being exact (i.e. an equality).
In case the intersection is empty
(the inequality is never exact), a bracketed number
`[$n$]' means that the smallest defect of the absolute term
is $n$. (That is, in an inequality $x\le y$ with $x\in\cL\{\nb_i\}$
and $y\in \bR$ the largest value of the l.h.s. on $\Sg_{3,\sg}$
is $y-n$.) One sees that both the trapezoidal and Ozsv\'ath--%
Szab\'o inequalities are (not more and not less than) moderately
good as conditions on the Alexander polynomial for genus 3, and
their sharpness increases with increasing (absolute) signature.

\[
\def\hh{\\[0.8mm]\hline[1.5mm]}%
\begin{mytab}{|c||c|c||c|c|c|c|c|}%
  { & & & \multicolumn{5}{|r|}{ } }%
  \hline [1mm]%
  \rx{-0.1em}%
  \diag{8mm}{1.5}{1}{
    \piclinewidth{42}
    \picline{-0.24 1.05}{1.73 -.04}
    \picputtext{1.2 0.8}{ineq.}
    \picputtext{0.4 0.3}{$\sg$}
  }\rx{-0.1em}%
& f & vr & OS1 & OS2 & T1 & T2 & T3 \\
\hline
\hline [2mm]%
0 & 10 & 5v7r &  [2]  &  [1] &  [2] &  [2] &  [2] \\[1mm]
2 & 12 & 6v8r &  [1]  &  1v  &  [1] &  [1] & 2v2r \\[1mm]
4 &  8 & 4v5r &  [1]  & 2v1r &  [2] & 1v1r & 2v2r \\[1mm]
6 &  6 & 3v4r &  2v   & 1v1r &   1v & 1v1r & 2v2r \\[1mm]
\hline
\hline[2mm]%
\mbox{total} & 12 & 6v6r & -- & -- & 1v & 2v1r & -- \\[1mm]
\hline
\end{mytab}%
\]

The determination of $\Sg_{3,\sg}$ also leads to a description
of $\Sg_3$, stated in the below theorem. Some data is given also
in the bottom line of the table. (Note that \eqref{ppOS_} depend,
at least in the constant terms, slightly on $\sg$, and so does T3
on $\sg\bmod 4$ when we sign $\nb_i$ so that $\nb_6>0$. Thus we
cannot compare $\Sg_3$ to these inequalities directly.) 

\begin{theorem}
The complete linear inequalities satisfied by the signed
Conway coefficients $\nb_k=[\nb]_k\cdot \sgn([\nb]_6)$ of
an alternating knot of genus 3 are those in \eqref{=}, together
with the following 6 other inequalities:
\[
\hbox to \textwidth{\hfill $\ds
\begin{array}{rcl@{\qquad}rcl@{\qquad}rcl}
40\nb_6 -\nb_4 -12\nb_2 & \ge & -37\,, &
34\nb_6 +8\nb_4+3\nb_2  & \ge & -1\,, &
5\nb_6 +\nb_4+\nb_2 & \ge & 0\,, \\[3mm]
3\nb_6 +\nb_2 & \ge &-1\,, &     9\nb_6 -3\nb_4+\nb_2 & \ge &0\,, &
\nb_6 & \ge & 1\,.\\[1mm]
\end{array}
$ \hfill
$\begin{array}{@{}c@{}}\mbox{\phantom{\rlap{$40\nb_6 -\nb_4 -12\nb_2$}}%
\strut}\\[3mm]%
\strut\phantom{\rlap{$40\nb_6 -\nb_4 -12\nb_2$}}\Box \\[1mm]\end{array}$}%
\]
\end{theorem}

We conclude with the table for $g=4$.

\[
\def\hh{\\[0.8mm]\hline[1.5mm]}%
\begin{mytab}{|c||c|c||c|c|c|c|c|c|c|}%
  { & & & \multicolumn{7}{|r|}{ } }%
  \hline [1mm]%
  \rx{-0.1em}%
  \diag{8mm}{1.5}{1}{
    \piclinewidth{42}
    \picline{-0.24 1.05}{1.73 -.04}
    \picputtext{1.2 0.8}{ineq.}
    \picputtext{0.4 0.3}{$\sg$}
  }\rx{-0.1em}%
& f & vr & OS1 & OS2 & OS3 & T1 & T2 & T3 & T4 \\
\hline
\hline [2mm]%
0 & 105 & 26v26r & 1v    & [2] & 1v   & [1] & [1]  & [1]  & [1]  \\[1mm]
2 & 115 & 29v31r & [1]   & [2] & [1]  & [2] & [2]  & [2]  & 4v4r \\[1mm]
4 &  81 & 21v24r & 1v    & [1] & 1v   & [1] & [1]  & 2v2r & 5v5r \\[1mm]
6 &  41 & 11v14r & 2v1r  & [1] & 2v1r & [2] & 1v1r & 2v2r & 5v5r \\[1mm]
8 &  18 &  5v8r  & 1v1r  & 3v  & 1v1r & 1v  & 1v1r & 2v2r & 2v4r \\[1mm]
\hline
\hline[2mm]%
\mbox{total}
  & 100 & 27v26r &  --   &  -- &  --  & 1v  & 2v1r & 4v2r & --   \\[1mm]
\hline
\end{mytab}%
\]

\noindent{\bf Acknowledgement.}
I wish to thank L.~Kauffman for inviting me to Waterloo in fall 2003
and letting me there present some the work of this paper in a talk. 
The material in \S\ref{S7} was motivated by some discussion with 
T.~Tsu\-ka\-mo\-to, and the one in \S\ref{sNR} by discussion with
A.~Holmsen. Furthermore,
S.~Rankin, M.~Thistlethwaite and K.~Fukuda provided some further
technical help for using their programs.

Part of this research was
supported by Postdoc grant P04300 of the Japan Society for the
Promotion of Science (JSPS) at the University of Tokyo (I wish to
thank my host T.~Kohno for his encouragement), the 21st Century
COE program at RIMS (Kyoto, Japan) and the BK21 Project at KAIST
(Daejeon, Korea).

\end{document}

%% file: myeqn.tex
\newenvironment{eqn}{\begin{equation}}{\end{equation}\@ignoretrue}

\newenvironment{myeqn*}[1]{\begingroup\def\@eqnnum{\reset@font\rm#1}%
\xdef\@tempk{\arabic{equation}}\begin{equation}\edef\@currentlabel{#1}}
{\end{equation}\endgroup\setcounter{equation}{\@tempk}\ignorespaces}

\newenvironment{myeqn}[1]{\begingroup\let\eq@num\@eqnnum
\def\@eqnnum{\bgroup\let\r@fn\normalcolor 
\def\normalcolor####1(####2){\r@fn####1#1}%
\eq@num\egroup}%
\xdef\@tempk{\arabic{equation}}\begin{equation}\edef\@currentlabel{#1}}
{\end{equation}\endgroup\setcounter{equation}{\@tempk}\ignorespaces}

\newenvironment{myeqn**}{\begin{myeqn}{(
\theequation)\es\es\mbox{\qed}}\edef\@currentlabel{\theequation}}
{\end{myeqn}\stepcounter{equation}}